\documentclass[12pt]{amsart}
\usepackage{amsmath,amsfonts}
\usepackage{amssymb}
\usepackage{amscd}
\usepackage{amsthm}
\usepackage{yhmath}
\usepackage{hyperref}
\usepackage{epsfig,subfigure}
\usepackage{import}
\usepackage[utf8]{inputenc}

\usepackage{todonotes}
\usepackage{tikz}
\usetikzlibrary{cd}

\usepackage[all]{xy}
\usepackage{color} 
\usepackage{soul}
\usepackage{tikz}
\usetikzlibrary{decorations.markings}

\numberwithin{equation}{section}

\newtheorem{proposition}[equation]{Proposition}

\theoremstyle{remark}

\theoremstyle{definition}

\newcommand{\salto}{$\\ $ \phantom{z_5=}}
\newcommand{\saltox}{$\\ $ \phantom{X_5=}}
\newcommand{\saltot}{$\\ $ \phantom{\theta_5=}}

\tikzset{->-/.style={decoration={
  markings,
  mark=at position #1 with {\arrow{>}}},postaction={decorate}}}

\def\XXint#1#2#3{{\setbox0=\hbox{$#1{#2#3}{\int}$}
	\vcenter{\hbox{$#2#3$}}\kern-.5\wd0}}

\newcommand{\R}{\mathbb R}

\newcommand{\g}{\mathfrak{g}}

\renewcommand{\span}{\operatorname{span}}

\begin{document}

\tikzset{->-/.style={decoration={
  markings,
  mark=at position #1 with {\arrow{>}}},postaction={decorate}}}
  \tikzset{-<-/.style={decoration={
  markings,
  mark=at position #1 with {\arrow{<}}},postaction={decorate}}}
\title{A cornucopia of Carnot groups in low dimensions}
	\author[E. Le Donne and F. Tripaldi]{Enrico Le Donne and Francesca Tripaldi}

\address{\textsc{Francesca Tripaldi}: 
University of Bologna, Department of Mathematics, Piazza di Porta S. Donato 5, 40126 Bologna, Italy}
\email{francesca.tripaldi2@unibo.it}
\address{\textsc{Enrico Le Donne}: 
Dipartimento di Matematica, Universit\`a di Pisa, Largo B. Pontecorvo 5, 56127 Pisa, Italy \\
\& \\
University of Jyv\"askyl\"a, Department of Mathematics and Statistics, P.O. Box (MaD), FI-40014, Finland}
\email{enrico.ledonne@unipi.it}
 \thanks{
 E.L.D. and F.T were partially supported by the Academy of Finland (grant
288501
`\emph{Geometry of subRiemannian groups}' and by grant
322898
`\emph{Sub-Riemannian Geometry via Metric-geometry and Lie-group Theory}')
and by the European Research Council
 (ERC Starting Grant 713998 GeoMeG `\emph{Geometry of Metric Groups}'). F.T. was also partially supported by the University of Bologna, funds for selected research topics, and by the European Union's Horizon 2020 research and innovation programme under the Marie Sk\l{}odowska-Curie grant agreement No 777822 GHAIA (`\emph{Geometric and Harmonic Analysis with Interdisciplinary Applications}').
}

\begin{abstract}
 Stratified groups  are those  simply connected Lie groups whose Lie algebras admit a derivation for which the eigenspace with eigenvalue 1 is Lie generating.
 When a stratified  group is equipped with a left-invariant path distance that is homogeneous with respect to the automorphisms induced by the derivation, this metric space  is known as  Carnot group. 
 Carnot groups appear in several mathematical contexts.
 To understand their algebraic structure, it is useful to study some examples explicitly.
 In this work, we provide a list of low-dimensional stratified groups, express their Lie product, and present a basis of left-invariant vector fields, together with their respective left-invariant 1-forms, a basis of right-invariant vector fields, and some other properties.  
 We exhibit all stratified  groups in dimension up to 7 and also study some free-nilpotent groups in dimension up to 14.
\end{abstract}

\renewcommand{\subjclassname} {\textup{2010} Mathematics Subject Classification}

\date{\today}
 
  \maketitle
  \tableofcontents
  
  
\section{Introduction}


Stratified groups, equipped with their homogeneous metrics, appear in several mathematical contexts.
Such metric groups arise in harmonic analysis, in the study of hypoelliptic differential operators, and as boundaries of strictly pseudo-convex complex domains, see the books \cite{Stein:book, Capogna-et-al} as initial references.
When  
 equipped with Carnot-Carath\'eodory metrics,
stratified groups are also  known as Carnot groups 
and they
appear in geometric group theory as asymptotic cones of nilpotent finitely generated groups, see
\cite{Gromov1, Pansu}.
Sub-Riemannian stratified groups are limits of Riemannian manifolds and are metric tangents of  sub-Riemannian manifolds. 
Sub-Riemannian geometries arise  in many areas of pure  and
applied  mathematics (such as algebra, geometry, analysis, mechanics, control theory, mathematical physics), as well as
 in applications (e.g., robotics), for references see the book \cite{Montgomery}. 
The literature on geometry and analysis on stratified groups is plentiful. In addition to the previous references, we also cite few more: \cite{
Roth:Stein,
Folland-Stein,
nagelstwe,
Koranyi-Reimann85,
varsalcou,
Heinonen-calculus,
magntesi,
Vittone-tesi,
Bonfiglioli:et:al,
jeancontrol,Rifford:book,
agrachev2019comprehensive}.

Stratified groups
 are simply connected Lie groups for which the Lie algebra admits
 a  special grading, called  a stratification,
  and is equipped with one such a stratification.
Namely, a grading is  a stratification if the degree-one layer of the  grading is  Lie generating. The presence of a positive grading implies that the Lie algebra is nilpotent. Not all nilpotent Lie algebras admit a stratification; however, free-nilpotent Lie algebras do.
Each positive grading induces a one-parameter 
family of automorphisms, giving rise to the consideration of homogeneous distances, which are unique up to biLipschitz equivalence, see \cite{ledonne_primer} for an introduction.
 Hence, stratified groups have many analogies with (finite-dimensional) vector spaces. However, their possible non-commutativity provides some crucial differences.

In this paper, we provide the list of all stratifiable Lie algebras of dimension up to 7, as well as the free-nilpotent Lie algebras of dimension up to 14. Moreover, up to dimension 6, we also consider those nilpotent Lie algebras that are not stratifiable, and since they all happen to be positively gradable, we provide one such grading for them.

Let us recall the terminology for stratifiability, stratifications, and gradings.
We first stress that all Lie algebras considered here are over $\R$ and finite-dimensional. 
Also, given two subspaces $V, W$ of a Lie algebra, we set
$
[V,W] := {\rm Span}\{[X,Y] ;\; X\in V,\ Y\in W\} .
$
 A \textit{stratification of step $s$} of a Lie algebra $\g$ is a direct-sum decomposition
$$\g = V_1 \oplus V_2 \oplus \cdots \oplus V_s,
  \label{def:stratifiable-algebras}
  $$
for some integer $s\geq 1$, where $V_s \not= \{0\}$ and $[V_1,V_j] = V_{j+1}$ for all integers $j\in \{1,\dots,s\}$ and where we set $V_{s+1} = \{0\}$. 
We say that a Lie algebra is \textit{stratifiable} if there exists a stratification of it.
Equivalently, as pointed out in \cite{Cornulier_SBE}, stratifiable algebras 
 are those nilpotent Lie algebras that possess a derivation inducing the identity map modulo the derived subalgebra.
  Stratifiable algebras are also called {\em Carnot algebras}.
We say that a Lie algebra is
\textit{stratified} when it is stratifiable and endowed with a stratification.
 We should stress that in the case of a stratifiable algebra, the choice of a stratification is essentially unique: every two stratifications of $\mathfrak{g}$ differ by a Lie algebra automorphism of $\mathfrak{g}$, see \cite[Proposition~2.17]{ledonne_primer} for a reference.

A  stratification is a particular example of grading. Indeed, it is a grading where the layer of degree one is Lie generating.
  A \emph{positive grading} of a Lie algebra $\g$ is a
family $(V_t)_{t\in (0,+\infty)}$ of linear subspaces of $\g$, where all but finitely many of the $V_t$'s are $\{0\}$, such that $\g$ is their direct sum
$$\g=  \bigoplus_{t\in (0,+\infty)} V_t $$
  and where 
  $$[V_t, V_u]\subset V_{t+u}, \qquad  \text{ for all } t,u >0.$$
We say that a Lie algebra is \emph{positively gradable} if there exists a positive grading of it. 
   We say that a Lie algebra is \textit{graded} (or {\em positively graded}, to be more precise) when it is positively gradable and endowed with a positive grading.
 More considerations on this subject can be found in \cite{LeDonne_Rigot_BCP_graded_groups}.

As usual, we only consider those Lie algebras that are 
{\em indecomposable}, i.e., those that are not the direct sum of two nontrivial Lie algebras. 
The classification of stratified algebras that we provide in this paper will give rise to the following consequence:

\begin{proposition}
There are 4 indecomposable stratified algebras in dimension 5, and 13 in dimension 6. All nilpotent Lie algebras of dimension less than or equal to 6 are positively gradable; but 2 in dimension 5 and 11 in dimension 6 are not stratifiable.
In dimension 7, there are two one-parameter families of indecomposable stratified algebras, plus 45 more single examples.
\end{proposition}

For the list of nilpotent Lie algebras of dimension less than or equal to 7, we follow Gong's classification from his thesis \cite{Gong_Thesis}. However, in our paper we shall also point out for each Lie algebra 
what the corresponding name/number is in the classifications present in de Graaf \cite{deGraaf_classification}, Magnin \cite{magnin}, and Del Barco \cite{delBarco}, respectively.
This will provide the reader with a database to navigate between the different notations. 
One should stress that the list in Magnin's paper \cite{magnin} consists of indecomposable nilpotent Lie algebras with complex structural constants (\ref{structural coeff}), and is therefore shorter than the other ones. This is due to the fact that some Lie algebras are isomorphic over algebraically closed fields, but are not isomorphic over $\mathbb{R}$ (we refer to Chapter 7 in \cite{Gong_Thesis} for a more thorough explanation of this fact).
 Regarding the free nilpotent algebras of dimension greater than 7, we shall use Hall's construction from \cite{Hall_basis}.

For each group in our list, we exhibit a basis of its Lie algebra following
the presentation in 
\cite{Gong_Thesis}
and we study the differential structure in exponential coordinates. More precisely, we calculate the group law, the vector fields that are left-invariant extensions of the given basis, and the respective left-invariant 1-forms. We also provide the expressions for the right-invariant vector fields.
One should stress that nilpotent groups that are not stratifiable, do not have a canonical sub-Riemannian structure. For these groups, we shall also {add a subsection} to present a possible grading and discuss which polarizations give rise to a maximal Hausdorff dimension for their respective sub-Riemannian distance, and calculate the tangent space. We also compute the asymptotic cones of all non-stratifiable nilpotent Lie groups of dimension less than or equal to 6.

\subsection{Notations and differential objects considered}

For the groups discussed in this paper, we shall use the following notation for  
describing
their differential structure 
in exponential coordinates.

We present each Lie algebra with a choice of a basis denoted by $X_1,\ldots, X_n$, and provide the list of non-zero bracket relations in the form
\begin{equation}\label{structural coeff}
    [X_i,X_j]=\sum_{k=1}^nc_{ij}^kX_k\,,
\end{equation}
where $c_{ij}^k\in\mathbb{R}$ are called \textit{structural constants}, and of course we only present those for which $i<j$.

Since for such low-dimensional Lie algebras the list of equations (\ref{structural coeff}) is rather short, we aim to provide a visualization of the grading of the Lie algebra through a diagram as follows. For Carnot groups up to dimension 7, the sum in (\ref{structural coeff}) is of at most one addend and most of the time the only non-zero coefficient is 1 or $-1$. Hence, if for  given $i,j$ there exists $ k$ such that $[X_i,X_j]=  X_k$,  we will then visualize this relation as
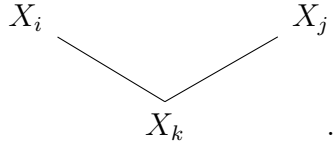
\begin{center}
    \begin{tikzcd}[end anchor=north]
    X_i \ar[dr, no head]& &X_j\ar[dl, no head]\\
    &X_k &\quad\;.
    \end{tikzcd}
\end{center}
In other words, the bracket relation expressed in the diagram should always be read \textit{from left to right}.
If the diagram should be read differently, we shall use the following notation. 
Namely, if  the bracket relation is   $[X_i,X_j]=-X_k$, we will draw the diagram as 
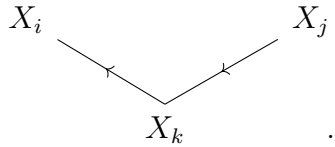
\begin{center}
    \begin{tikzcd}[end anchor=north]
    X_i \ar[dr, no head,-<-=.5]& &X_j\ar[dl, no head,->-=.5]\\
    &X_k &\quad\;.
    \end{tikzcd}
\end{center}
If instead the bracket relation is   $[X_i,X_j]=c X_k$ for some $c\in \R$, we will then write
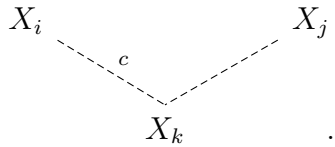
\begin{center}
    \begin{tikzcd}[end anchor=north]
    X_i \ar[dr, no head, dashed,"c"]& &X_j\ar[dl, no head, dashed]\\
    &X_k &\quad\;.
    \end{tikzcd}
\end{center}

In the case of Carnot algebras where the given basis can be adapted to a stratification, we will also draw the diagram by rows according to the different strata. For example,   a diagram of the form
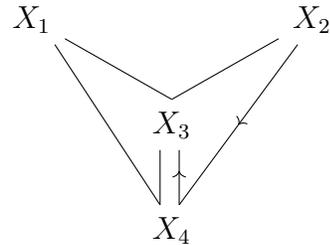
\begin{center}
    \begin{tikzcd}[end anchor=north]
    X_1 \ar[dr, no head]\ar[ddr, no head,end anchor={[xshift=-3.3ex]north east}]& &X_2\ar[dl, no head]\ar[ddl, no head, end anchor={[xshift=-1.9ex]north east}, ->-=.5]\\
    &X_3\ar[d, no head, start anchor={[xshift=-3.3ex]south east}, end anchor={[xshift=-3.3ex]north east}]\ar[d, no head, start anchor={[xshift=-1.9ex]south east}, end anchor={[xshift=-1.9ex]north east}, -<-=.5] &\\
    & X_4 &
    \end{tikzcd}
    
\end{center}
means that $X_1,X_2$ is a basis of the first stratum $V_1$, $X_3$ is a basis of the second stratum $V_2$, $X_4$ is a basis of the third stratum $V_3$, and the bracket relations are $[X_1,X_2]=X_3$, $[X_1,X_3]=X_4$, and $[X_2,X_3]=X_4$.

Moreover, other information is also readily available from simply looking at the diagram, such as the rank of the Lie algebra, which is equal to the number of vectors in the first row, and the nilpotency step, which is given by the number of rows in the diagram. 

Given a basis $X_1,\ldots, X_n$ of a nilpotent Lie algebra $\mathfrak{g}$, there exists a unique {(up to isomorphism)} simply connected Lie group $\mathbb{G}$ with $\mathfrak{g}$ as Lie algebra. Moreover, the exponential map $\exp\colon\mathfrak{g}\to\mathbb{G}$ is a diffeomorphism, see \cite{Corwin-Greenleaf} for a reference. We shall then  parametrize $\mathbb{G}$ via the exponential map and our choice of basis. Namely, we will use the identification
\begin{eqnarray}\label{identification}
\mathbb{R}^n&\longleftrightarrow&\mathbb{G}\nonumber\\
    (x_1,\ldots,x_n)&\longmapsto&\exp\Big(\sum_{i=1}^nx_iX_i\Big)\,.
\end{eqnarray}

Since in nilpotent groups the Baker-Campbell-Hausdorff formula converges globally, the identification above allows us to write the group product. In fact, for all $\mathbf{x},\mathbf{y}\in\mathbb{R}^n$, there exists a unique  $\mathbf{z}\in\mathbb{R}^n$ such that 
\begin{equation}\label{group law in G}
    \exp\Big(\sum_{i=1}^nx_iX_i\Big)\exp\Big(\sum_{i=1}^ny_iX_i\Big)=\exp\Big(\sum_{i=1}^nz_iX_i\Big)\,.
\end{equation}
Via the identification (\ref{identification}), one can write the group law in (\ref{group law in G}) as
\begin{equation}
    \mathbf{x}\ast\mathbf{y}=\mathbf{z}\,.
\end{equation}

Hence, we have a group law $\ast$ on $\mathbb{R}^n$ that makes $(\mathbb{R}^n,\ast)$ a simply connected Lie group with Lie algebra $\mathfrak{g}$, whose identity element is $\mathbf{0}$.

The basis $X_1,\ldots,X_n$ of $\mathfrak{g}$ induces left-invariant vector fields on $\mathbb{R}^n$ via the formula
\begin{equation}\label{leftinvariant vf}
    \mathbf{x}\mapsto \mathbf{d}\big(L_\mathbf{x}\big)_\mathbf{0}\mathbf{e}_i\;,\;i=1,\ldots,n\,,
\end{equation}
where $L_\mathbf{x}(\mathbf{y}):=\mathbf{x}\ast \mathbf{y}$ is the \emph{left translation}, and by $\mathbf{e}_i$ we denote the $i$-th vector of the standard basis of $\mathbb{R}^n$. We will still denote by $X_i$ the left-invariant vector fields given by equation (\ref{leftinvariant vf}).
One should stress that each vector field $X_i$ is represented by the $i$-th column of the matrix $\mathbf{d}\big(L_\mathbf{x}\big)_\mathbf{0}$, for $i=1,\ldots,n$. For better readability, however, in our paper we will provide both $\mathbf{d}\big(L_\mathbf{x}\big)_\mathbf{0}$ as a matrix, as well as the vector fields $X_1,\ldots,X_n$ written as derivations.

We also provide the explicit expression in coordinates of the basis $\theta_1,\ldots,\theta_n$ of left-invariant 1-forms that is dual to the basis of left-invariant vector fields $X_1,\ldots,X_n$, that is
\begin{equation}\label{leftinvariant form}
    \theta_i(X_j)=\delta_{ij}\;,\; \text{ for } i,j=1,\ldots,n\,.
\end{equation}

Likewise, one can repeat the same procedure for right translations. For shortness, we will only provide the differential at $\mathbf{0}$ of right translations $R_\mathbf{x}(\mathbf{y}):=\mathbf{y}\ast\mathbf{x}$. The reader can then deduce the right-invariant vector fields (and subsequently the right-invariant 1-forms) from the columns of the matrix $\mathbf{d}\big(R_\mathbf{x}\big)_\mathbf{0}$.

 In the case of non-stratifiable nilpotent groups, we will insert an extra subsection to discuss the possible Carnot-Carath\'eodory structures that  maximize the Hausdorff dimension. In order to do so, we now recall the notion of polarization and the method to calculate the dimension of the metric space that it defines (up to biLipschitz equivalence).  
Given $\mathfrak{g}$ a nilpotent Lie algebra, we denote by $\mathbb{G}$ the simply connected Lie group that has $\mathfrak{g}$ as  Lie algebra. By applying left-translations, we have that any subspace $V$ of $\mathfrak{g}$ induces a left-invariant subbundle $\Delta$ of $T\mathbb{G}$.  Following Gromov's terminology, we call the pair $(\mathbb{G},\Delta)$ a \emph{polarization} of $\mathbb{G}$. We only focus on the case where $V$ is Lie generating, which is equivalent to saying that $\Delta$ is a  bracket-generating distribution. One can check that, since $\mathfrak{g}$ is nilpotent, a subspace $V\subset \mathfrak{g}$ is Lie generating if and only if \begin{equation}
    \label{complementary_Delta}
V+[\mathfrak{g},\mathfrak{g}]=\mathfrak{g}.
\end{equation}
Once we fix a left-invariant norm $\Vert\cdot\Vert$ on $\Delta$, we get 
that the triple $(\mathbb{G},\Delta,\Vert\cdot\Vert)$ is an example of a
Carnot-Carath\'eodory space.
We refer to Gromov seminal paper \cite{Gromov1} for the theory of Carnot-Carath\'eodory spaces, also called {\em CC-spaces}.
 Every CC-space is a metric space when equipped with the control distance. The Hausdorff dimension with respect to this distance can be expressed as
\begin{equation}\label{Hausdorff dim}
    \sum_i\, i\,\big(\mathrm{rank\,}\Delta^i-\mathrm{rank\,}{\Delta}^{i-1}\big),
\end{equation}
where $\mathrm{rank\,}\Delta^0=0$, $\mathrm{rank\,}\Delta^1=\mathrm{dim\,}V$, and 
\begin{equation*}
    \mathrm{rank\,}\Delta^i=\mathrm{dim\,}\big(V+[V,V]+\cdots+\underbrace{[V,[V,\cdots,[V,V]]]}_{i-1}\big)\,.
\end{equation*} 
Let us point out that if $\mathfrak{g}$ is stratifiable, then a Lie generating subspace $V\subset\mathfrak{g}$ maximizes the Hausdorff dimension if and only if $V$ is the first layer of a stratification. On the other hand, this characterization is not present for non-stratifiable Lie algebras. For this reason in our paper, when presenting a non-stratifiable nilpotent Lie algebra $\mathfrak{g}$, we will add an extra subsection to discuss for which choice of polarization $(\mathbb{G},\Delta)$ we obtain maximal Hausdorff dimension. 
In low dimension, except for a few cases, such polarizations are  unique up to automorphism. Namely,
in  dimension up to 5,
all possible $\Delta$ with \eqref{complementary_Delta}   differ  by an automorphism. 
In dimension 6, except for $N_{6,1,2}$ and $N_{6,1,4}$, 
polarizations of maximal  dimension 
are unique up to automorphism.
For the considered polarizations  we calculate the tangent cone, which is also known as symbol, and by a theorem of Mitchell is a very easily computable Carnot group, see \cite{Mitchell}.

By a theorem of Pansu, the asymptotic cone of every nilpotent Lie group equipped with a CC-metric is a Carnot group, which does not depend on the choice of metric.
Thus, for every non-stratifiable nilpotent Lie algebra of dimension 5 or 6 we shall also determine its asymptotic cone.
The calculation is done via the {\em associated Carnot-graded Lie algebra}, for which we refer to \cite{Pansu-croissance}.


Finally, after the list of Carnot groups of dimension 7, we analyze free-Carnot groups of low dimension. Regarding the step-2 case, one can easily write down the product law and the left-invariant vector fields in arbitrary dimension (this calculation is not original and can also be found in \cite{LeDonne_Rigot_BCP_graded_groups, MR3549917, Bonfiglioli:et:al}).
Regarding the step-3 case, the free-Carnot groups of rank 2 is 5-dimensional, so it is already included in Section 4 (see Lie algebra ${N_{5,2,3}}$ on page 18).

In addition, we will present the rank-3 step-3 free-Carnot group, which has dimension 14, and the rank-2 
free-Carnot groups o step at most 5, which have dimensions 5, 8, and 14, respectively. We will not discuss the rank-4 step-3 case, which has dimension 30, nor the one with rank 3 and step 4, which has dimension 32.

Of the free-Carnot groups above, we will also provide exponential coordinates of the second type, together with the change of variables with respect to the ones of first type. We shall add an $s$, as an exponent, to  those differential objects that are expressed in exponential coordinates of the second type.


\section{1D--4D nilpotent Lie algebras}
  
The nilpotent Lie groups of dimension up to 4 are well known, and they are all stratifiable.  The list of the stratifiable Lie algebras of dimension less than or equal to 4 is: $\mathbb{R}$, $\mathbb{R}^2$, $\mathbb{R}^3$, $N_{3,2}$, $\mathbb{R}^4$, $N_{3,2}\times\mathbb{R}$, and $N_{4,2}$, where here $\mathbb{R}^n$ denotes the $n$-dimensional  abelian Lie algebra, $N_{3,2}$ is the first Heisenberg Lie algebra, and $N_{4,2}$ is the Engel Lie algebra.

We shall now recall the non-zero brackets of these last two Lie algebras in order to help the reader get familiar with our diagram notation.

The algebra $N_{3,2}$ and its stratification are represented as
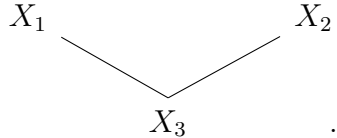
\begin{center}
    \begin{tikzcd}[end anchor=north]
    X_1 \ar[dr, no head]& &X_2\ar[dl, no head]\\
    &X_3 &\quad\; .
    \end{tikzcd}
\end{center}
Whereas the algebra $N_{4,2}$ and its stratification are represented as
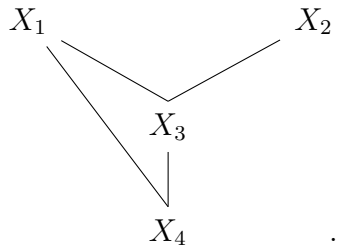
\begin{center}
    \begin{tikzcd}[end anchor=north]
    X_1\ar[dr, no head]\ar[ddr, no head] & &X_2\ar[dl, no head]\\
    &X_3\ar[d, no head] & \\
    & X_4 &\quad\;.
    \end{tikzcd}
\end{center}

In the case of the first Heisenberg group and the Engel group, it is sometimes convenient to work in exponential coordinates of the first kind, and some other times in exponential coordinates of the second kind. 

  \subsection*{$\bf N_{3,2}$ } 
  The following  Lie algebra is denoted as $N_{3,2}$ by Gong in \cite{Gong_Thesis}, as $L_{3,2}$ by de Graaf in \cite{deGraaf_classification}, as $\mathfrak{h}_3$ by Del Barco in \cite{delBarco}, and as $\mathcal{G}_3$ by Magnin in \cite{magnin}.
  
  The only non-trivial bracket is the following:
\begin{eqnarray*}
    [X_1,X_2]=X_3\,.
\end{eqnarray*}
This is a nilpotent Lie algebra of rank 2 and step 2 that is stratifiable, also known as the first Heisenberg algebra. The Lie brackets can be pictured with the diagram:
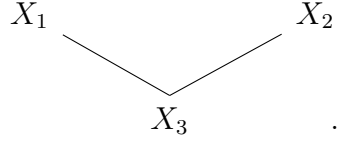
\begin{center}
 
\begin{tikzcd}[end anchor=north]
    X_1 \ar[dr, no head]& &X_2\ar[dl, no head]\\
    &X_3 &\quad\; .
    \end{tikzcd}
   
\end{center}



The composition law (\ref{group law in G}) of $N_{3,2}$ is given by:
\begin{itemize}
    \item $z_1=x_1+y_1$;
    \item $z_2=x_2+y_2$;
    \item $z_3=x_3+y_3+\frac{1}{2}(x_1y_2-x_2y_1)$.
\end{itemize}

Since
\begin{eqnarray*}
   \mathrm{d} (L_{\bf x})_{\bf 0}= \left[\begin{matrix} 
   1 & 0 & 0\\
   0 & 1 & 0\\
   -\frac{x_2}{2} & \frac{x_1}{2} &1
   \end{matrix}\right]\,,
\end{eqnarray*}
 the induced left-invariant vector fields \eqref{leftinvariant vf} are: 
\begin{itemize}
\item $X_1=\partial_{x_1}-\frac{x_2}{2}\, \partial_{ x_3}\,;$
\item $X_2=\partial_{x_2}+\frac{x_1}{2}\, \partial_{ x_3}\,;$
\item $X_3=\, \partial_{ x_3}\,,$
\end{itemize}
and the respective left-invariant 1-forms \eqref{leftinvariant form} are:
\begin{itemize}
\item $\theta_1=dx_1$;
\item $\theta_2=dx_2$;
\item $\theta_3=dx_3-\frac{x_1}{2} dx_2+\frac{x_2}{2}dx_1\,.$
\end{itemize}
 
 Finally, we have
 \begin{eqnarray*}
   \mathrm{d} (R_{\bf x})_{\bf 0}= \left[\begin{matrix} 
   1 & 0 & 0\\
   0 & 1 & 0\\
   \frac{x_2}{2} & -\frac{x_1}{2} &1
   \end{matrix}\right]\,.
\end{eqnarray*}

One can also consider the exponential coordinates of the second kind. In this case, we obtain the following expression for the left-invariant vector fields:

\begin{itemize}
\item $X_1^s=\partial_{x_1}\,;$
\item $X_2^s=\partial_{x_2}+x_1\, \partial_{ x_3}\,;$
\item $X_3^s= \partial_{ x_3}\,.$
\end{itemize}

  \subsection*{$\bf N_{4,2}$ } 
  The following  Lie algebra is denoted as $N_{4,2}$ by Gong in \cite{Gong_Thesis}, as $L_{4,3}$ by de Graaf in \cite{deGraaf_classification}, as (2) by Del Barco in \cite{delBarco}, and as $\mathcal{G}_4$ by Magnin in \cite{magnin}.
  
  The non-trivial brackets are the following:
\begin{eqnarray*}
    [X_1,X_2]=X_3\,,\,[X_1,X_3]=X_4\,.
\end{eqnarray*}
This is a nilpotent Lie algebra of rank 2 and step 3 that is stratifiable, also known as the filiform Lie algebra of dimension 4. The Lie brackets can be pictured with the diagram:
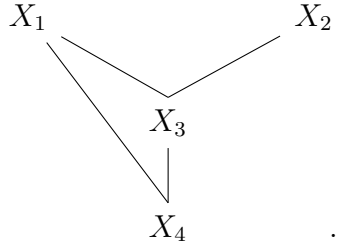
\begin{center}
    \begin{tikzcd}[end anchor=north]
    X_1\ar[dr, no head]\ar[ddr, no head] & &X_2\ar[dl, no head]\\
    &X_3\ar[d, no head] & \\
    & X_4 & \quad\;.
    \end{tikzcd}
\end{center}



The composition law (\ref{group law in G}) of $N_{4,2}$ is given by:
\begin{itemize}
    \item $z_1=x_1+y_1$;
    \item $z_2=x_2+y_2$;
    \item $z_3=x_3+y_3+\frac{1}{2}(x_1y_2-x_2y_1)$;
    \item $z_4=x_4+y_4+\frac{1}{2}(x_1y_3-x_3y_1)+\frac{1}{12}(x_1-y_1)(x_1y_2-x_2y_1)$.
\end{itemize}

Since
\begin{eqnarray*}
   \mathrm{d} (L_{\bf x})_{\bf 0}= \left[\begin{matrix} 
   1 & 0 & 0& 0 \\
   0 & 1 & 0& 0\\
   -\frac{x_2}{2} & \frac{x_1}{2} &1& 0\\
  -\frac{x_1x_2}{12} -\frac{x_3}{2} & \frac{x_1^2}{12} &\frac{x_1}{2}& 1 
   \end{matrix}\right]\,,
\end{eqnarray*}
the induced left-invariant vector fields $\eqref{leftinvariant vf}$ are: 
\begin{itemize}
\item $X_1=\partial_{x_1}-\frac{x_2}{2}\,\partial_{x_3}-\big(\frac{x_3}{2}+\frac{x_1x_2}{12}\big)\partial_{x_4}\,;$
\item $X_2=\partial_{x_2}+\frac{x_1}{2}\,\partial_{x_3}+\frac{x_1^2}{12}\,\partial_{ x_4}\,;$
\item $X_3=\partial_{x_3}+\frac{x_1}{2}\,\partial_{ x_4}\,;$
\item $X_4={\partial}_{ x_4}\,,$
\end{itemize}
and the respective left-invariant 1-forms \eqref{leftinvariant form} are:
\begin{itemize}
\item $\theta_1=dx_1$;
\item $\theta_2=dx_2$;
\item $\theta_3=dx_3-\frac{x_1}{2} dx_2+\frac{x_2}{2}dx_1$;
\item $\theta_4=dx_4-\frac{x_1}{2}dx_3+\frac{x_1^2}{6}dx_2+\big(\frac{x_3}{2}-\frac{x_1x_2}{6}\big)dx_1\,.$
\end{itemize}
 
 Finally, we have
 \begin{eqnarray*}
   \mathrm{d} (R_{\bf x})_{\bf 0}= \left[\begin{matrix} 
   1 & 0 & 0& 0 \\
   0 & 1 & 0& 0\\
   \frac{x_2}{2} & -\frac{x_1}{2} &1& 0\\
 \frac{x_3}{2} -\frac{x_1x_2}{12}  & \frac{x_1^2}{12} &-\frac{x_1}{2}& 1 
   \end{matrix}\right]\,.
\end{eqnarray*}

One can also consider the exponential coordinates of the second kind. In this case, we obtain the following expression for the left-invariant vector fields:

\begin{itemize}
\item $X_1^s=\partial_{x_1}\,;$
\item $X_2^s=\partial_{x_2}+x_1\, \partial_{ x_3}+\frac{x_1^2}{2}\partial_{x_4}\,;$
\item $X_3^s= \partial_{ x_3}+x_1\partial_{x_4}\,;$
\item $X_4^s=\partial_{x_4}\,.$
\end{itemize} 
\newpage
  \section{5D indecomposable nilpotent Lie algebras}
  
  Among all the indecomposable nilpotent Lie algebras of dimension 5, there are 4 Carnot algebras and 2 more nilpotent Lie algebras, which are gradable. Moreover, there are other 3 decomposable nilpotent algebras, which are stratifiable: the abelian $\mathbb{R}^5$, the direct products of the first Heisenberg group times $\mathbb{R}^2$, and the Engel group times $\mathbb{R}$.
  
  \subsection*{$\bf  N_{5,1}$ non-stratifiable}
The following Lie algebra is denoted as $N_{5,1}$ by Gong in \cite{Gong_Thesis}, as $L_{5,6}$ by de Graaf in \cite{deGraaf_classification}, as (1) by Del Barco in \cite{delBarco}, and as $\mathcal{G}_{5,6}$ by Magnin in \cite{magnin}.

The non-trivial brackets are the following:
\begin{equation*}
    [X_1,X_2]=X_3\,,\,[X_1,X_3]=X_4\,,\,[X_1,X_4]=[X_2,X_3]=X_5\,.
\end{equation*}
This is a nilpotent Lie algebra of rank 2 and step 4 that is positively gradable, yet not stratifiable. The Lie brackets can be pictured with the  diagram:
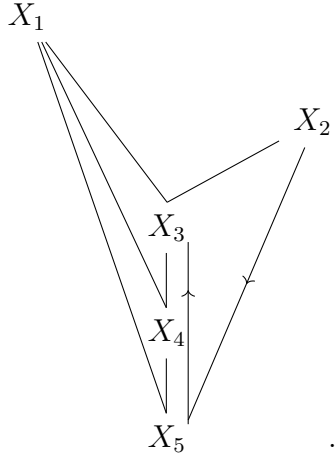
\begin{center}
 
	\begin{tikzcd}[end anchor=north]
		X_1\ar[ddr,no head] \ar[dddr,no head, end anchor={[xshift=-2.5ex]north east}]\ar[ddddr,no head, end anchor={[xshift=-2.5ex]north east}]& &\\
		&& X_2\dlar[no head]\ar[->-=.5,no head,dddl, end anchor={[xshift=-.9ex]north east},start anchor={[xshift=-3.ex]south east }, end anchor={[yshift=-.5ex]north east}]\\
		& X_{3}\dar[no head, end anchor={[xshift=-2.5ex]north east},start anchor={[xshift=-2.5ex]south east }]\ar[-<-=.3,no head,dd, end anchor={[xshift=-0.9ex]north east},start anchor={[xshift=-0.9ex]south east }, end anchor={[yshift=-0.8ex]north east},start anchor={[yshift=0.8ex]south east }]&\\
		&X_4\dar[no head, end anchor={[xshift=-2.5ex]north east},start anchor={[xshift=-2.5ex]south east }]&\\
		&X_5&  \quad\;.
	\end{tikzcd}

\end{center}


The composition law (\ref{group law in G}) of $N_{5,1}$ is given by:


 
\begin{itemize}
    \item $z_1=x_1+y_1$;
    \item $z_2=x_2+y_2$;
    \item $z_3=x_3+y_3+\frac{1}{2}(x_1y_2-x_2y_1)$;
    \item $z_4=x_4+y_4+\frac{1}{2}(x_1y_3-x_3y_1)+\frac{1}{12}(x_1-y_1)(x_1y_2-x_2y_1)$;
    \item $z_5=x_5+y_5+\frac{1}{2}(x_1y_4-x_4y_1+x_2y_3-x_3y_2)+\frac{1}{12}(x_1-y_1)(x_1y_3-x_3y_1)\salto +\frac{1}{12}(x_2-y_2)(x_1y_2-x_2y_1)-\frac{1}{24}x_1y_1(x_1y_2-x_2y_1)$.
\end{itemize}

Since
\begin{eqnarray*}
   \mathrm{d} (L_{\bf x})_{\bf 0}= \left[\begin{matrix} 
   1 & 0 & 0& 0 & 0\\
   0 & 1 & 0& 0& 0\\
   -\frac{x_2}{2} & \frac{x_1}{2} &1& 0& 0\\
  -\frac{x_1x_2}{12} -\frac{x_3}{2} & \frac{x_1^2}{12} &\frac{x_1}{2}& 1 & 0\\
  -\frac{x_1x_3+x_2^2}{12}-\frac{x_4}{2}& \frac{x_1x_2}{12}-\frac{x_3}{2} &\frac{x_1^2}{12}+\frac{x_2}{2} & \frac{x_1}{2}& 1
   \end{matrix}\right]\,,
\end{eqnarray*}
the induced left-invariant vector fields \eqref{leftinvariant vf} are:
\begin{itemize}
\item $X_1={\partial}_{x_1}-\frac{x_2}{2}{\partial}_{x_3}-\big(\frac{x_3}{2}+\frac{x_1x_2}{12}\big){\partial}_{x_4}-\big(\frac{x_4}{2}+\frac{x_1x_3+x_2^2}{12}\big){\partial}_{x_5}\,;$
\item $X_2={\partial}_{x_2}+\frac{x_1}{2}{\partial}_{x_3}+\frac{x_1^2}{12}{\partial}_{x_4}+\big(\frac{x_1x_2}{12}-\frac{x_3}{2}\big){\partial}_{x_5}\,;$
\item $X_3={\partial}_{x_3}+\frac{x_1}{2}{\partial}_{x_4}+\big(\frac{x_1^2}{12}+\frac{x_2}{2}\big){\partial}_{x_5}\,;$
\item $X_4={\partial}_{x_4}+\frac{x_1}{2}{\partial}_{x_5}\,;$
\item $X_5={\partial}_{x_5}\,,$
\end{itemize}
and the respective left-invariant 1-forms \eqref{leftinvariant form} are: 
\begin{itemize}
\item $\theta_1=dx_1$;
\item $\theta_2=dx_2$;
\item $\theta_3=dx_3-\frac{x_1}{2} dx_2+\frac{x_2}{2}dx_1$;
\item $\theta_4=dx_4-\frac{x_1}{2}dx_3+\frac{x_1^2}{6}dx_2+\big(\frac{x_3}{2}-\frac{x_1x_2}{6}\big)dx_1$;
\item $\theta_5=dx_5-\frac{x_1}{2}dx_4+\big(\frac{x_1^2}{6}-\frac{x_2}{2}\big)dx_3+\big(\frac{x_1x_2}{6}+\frac{x_3}{2}-\frac{x_1^3}{24}\big)dx_2+\big(\frac{x_4}{2}-\frac{x_1x_3+x^2_2}{6}+\frac{x_1^2x_2}{24}\big)dx_1\,.$
\end{itemize}

Finally, we have
\begin{eqnarray*}
   \mathrm{d} (R_{\bf x})_{\bf 0}= \left[\begin{matrix} 
   1 & 0 & 0& 0 & 0\\
   0 & 1 & 0& 0& 0\\
   \frac{x_2}{2} & -\frac{x_1}{2} &1& 0& 0\\
  \frac{x_3}{2}-\frac{x_1x_2}{12} & \frac{x_1^2}{12} &-\frac{x_1}{2}& 1 & 0\\
  \frac{x_4}{2}-\frac{x_1x_3+x_2^2}{12}& \frac{x_1x_2}{12}+\frac{x_3}{2} &\frac{x_1^2}{12}-\frac{x_2}{2} & -\frac{x_1}{2}& 1
   \end{matrix}\right]\,.
\end{eqnarray*}

\subsubsection{Grading, polarizations of maximal dimension, and asymptotic cone}
The Lie algebra $N_{5,1}$
   is not stratifiable, but it is gradable as
   \begin{eqnarray*}
   V_i &= & \span \{  X_i \} \,,\,i=1,\ldots,5\,.
   \end{eqnarray*}

      We claim that in this Lie algebra every two complementary subspaces to the  derived subalgebra, as in \eqref{complementary_Delta}, differ by an automorphism. Indeed, when $u_i^j$ varies, 
      the matrix
$$    \left[  \begin{array}{cccccc}
     1 &  0 &  0 &  0 &  0   \\
 0 &  1 &  0 &  0 &  0 \\
 u_1^3 &  u_2^3 &  1 &  0 &  0 \\
 u_1^4 &  u_2^4 &  u_2^3 &  1 &  0 \\
 u_1^5 &  u_2^5 & u_2^4-u_1^3 & u_2^3&  1 
      \end{array}\right]$$
         is a Lie algebra automorphism and sends 
         the subspace $\span\{X_1, X_2\}$ to $\span\{X_1+u_1^3X_3+u_1^4X_4+u_1^5X_5, X_2+u_2^3X_3+u_2^4X_4+u_2^5X_5\}$, which is an arbitrary subspace as in \eqref{complementary_Delta}.
       In particular, every $\Delta$ as in \eqref{complementary_Delta} gives maximal dimension.
The tangent cone of each of such polarizations has Lie algebra isomorphic to $N_{5,2,3}$, see page \pageref{N523}.

 The asymptotic cone 
 of the Lie group with Lie algebra
 $N_{5,1}$ has Lie algebra isomorphic to $N_{5,2,1}$, which is the filiform algebra of step 4.

 

   \subsection*{$\bf N_{5,2,1}$ } 
  The following  Lie algebra is denoted as $N_{5,2,1}$ by Gong in \cite{Gong_Thesis}, as $L_{5,7}$ by de Graaf in \cite{deGraaf_classification}, as (2) by Del Barco in \cite{delBarco}, and as $\mathcal{G}_{5,5}$ by Magnin in \cite{magnin}.
  
  The non-trivial brackets are the following:
\begin{eqnarray*}
    [X_1,X_2]=X_3\,,\,[X_1,X_3]=X_4\,,\,[X_1,X_4]=X_5\,.
\end{eqnarray*}
This is a nilpotent Lie algebra of rank 2 and step 4 that is stratifiable, also known as the filiform Lie algebra of dimension 5. The Lie brackets can be pictured with the  diagram:
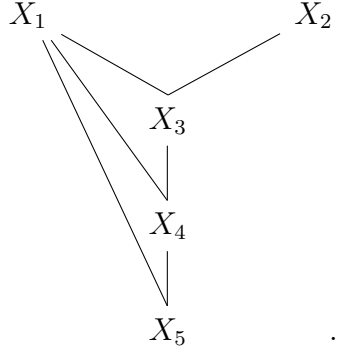
\begin{center}
 
	\begin{tikzcd}[end anchor=north]
		X_1\drar[no head] \ar[ddr,no head, end anchor={[xshift=-2.5ex]north east}]\ar[dddr,no head, end anchor={[xshift=-2.5ex]north east}]& & X_2\dlar[no head]\\
		& X_{3}\dar[no head, end anchor={[xshift=-2.5ex]north east},start anchor={[xshift=-2.5ex]south east }]&\\
		&X_4\dar[no head, end anchor={[xshift=-2.5ex]north east},start anchor={[xshift=-2.5ex]south east }]&\\
		&X_5& \quad\; .
	\end{tikzcd}
   
\end{center}



The composition law (\ref{group law in G}) of $N_{5,2,1}$ is given by:
\begin{itemize}
    \item $z_1=x_1+y_1$;
    \item $z_2=x_2+y_2$;
    \item $z_3=x_3+y_3+\frac{1}{2}(x_1y_2-x_2y_1)$;
    \item $z_4=x_4+y_4+\frac{1}{2}(x_1y_3-x_3y_1)+\frac{1}{12}(x_1-y_1)(x_1y_2-x_2y_1)$;
    \item $z_5=x_5+y_5+\frac{1}{2}(x_1y_4-x_4y_1)+\frac{1}{12}(x_1-y_1)(x_1y_3-x_3y_1) -\frac{1}{24}x_1y_1(x_1y_2-x_2y_1)$.
\end{itemize}

Since
\begin{eqnarray*}
   \mathrm{d} (L_{\bf x})_{\bf 0}= \left[\begin{matrix} 
   1 & 0 & 0& 0 & 0\\
   0 & 1 & 0& 0& 0\\
   -\frac{x_2}{2} & \frac{x_1}{2} &1& 0& 0\\
  -\frac{x_1x_2}{12} -\frac{x_3}{2} & \frac{x_1^2}{12} &\frac{x_1}{2}& 1 & 0\\
  -\frac{x_1x_3}{12}-\frac{x_4}{2}& 0 &\frac{x_1^2}{12} & \frac{x_1}{2}& 1
   \end{matrix}\right]\,,
\end{eqnarray*}
the induced left-invariant vector fields \eqref{leftinvariant vf} are:
\begin{itemize}
\item $X_1={\partial}_{x_1}-\frac{x_2}{2}{\partial}_{x_3}-\big(\frac{x_3}{2}+\frac{x_1x_2}{12}\big){\partial}_{x_4}-\big(\frac{x_4}{2}+\frac{x_1x_3}{12}\big){\partial}_{x_5}\,;$
\item $X_2={\partial}_{x_2}+\frac{x_1}{2}{\partial}_{x_3}+\frac{x_1^2}{12}{\partial}_{x_4}\,;$
\item $X_3={\partial}_{x_3}+\frac{x_1}{2}{\partial}_{x_4}+\frac{x_1^2}{12}{\partial}_{x_5}\,;$
\item $X_4={\partial}_{x_4}+\frac{x_1}{2}{\partial}_{x_5}\,;$
\item $X_5={\partial}_{x_5}\,,$
\end{itemize}
and the respective left-invariant 1-forms \eqref{leftinvariant form} are:
\begin{itemize}
\item $\theta_1=dx_1$;
\item $\theta_2=dx_2$;
\item $\theta_3=dx_3-\frac{x_1}{2} dx_2+\frac{x_2}{2}dx_1$;
\item $\theta_4=dx_4-\frac{x_1}{2}dx_3+\frac{x_1^2}{6}dx_2+\big(\frac{x_3}{2}-\frac{x_1x_2}{6}\big)dx_1$;
\item $\theta_5=dx_5-\frac{x_1}{2}dx_4+\frac{x_1^2}{6}dx_3-\frac{x_1^3}{24}dx_2+\big(\frac{x_4}{2}-\frac{x_1x_3}{6}+\frac{x_1^2x_2}{24}\big)dx_1\,.$
\end{itemize}
 
 Finally, we have
 \begin{eqnarray*}
   \mathrm{d} (R_{\bf x})_{\bf 0}= \left[\begin{matrix} 
   1 & 0 & 0& 0 & 0\\
   0 & 1 & 0& 0& 0\\
   \frac{x_2}{2} & -\frac{x_1}{2} &1& 0& 0\\
 \frac{x_3}{2} -\frac{x_1x_2}{12}  & \frac{x_1^2}{12} &-\frac{x_1}{2}& 1 & 0\\
  \frac{x_4}{2}-\frac{x_1x_3}{12}& 0 &\frac{x_1^2}{12} & -\frac{x_1}{2}& 1
   \end{matrix}\right]\,.
\end{eqnarray*}
  \subsection*{$\bf N_{5,2,2}$ non-stratifiable}  The following Lie algebra is denoted as $N_{5,2,2}$ by Gong in \cite{Gong_Thesis}, as $L_{5,5}$ by de Graaf in \cite{deGraaf_classification}, as (4) by Del Barco in \cite{delBarco}, and as $\mathcal{G}_{5,3}$ by Magnin in \cite{magnin}.
  
  The non-trivial brackets are the following:
  \begin{eqnarray*}
[X_1,X_2]=X_4\,,\,[X_1,X_4]=[X_2,X_3]=X_5\,.
\end{eqnarray*}
This is a nilpotent Lie algebra of rank 3 and step 3 that is positively gradable, yet not stratifiable. The Lie brackets can be pictured with the diagram:

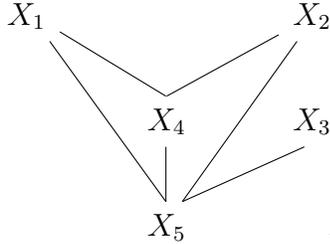
\begin{center}
 
	\begin{tikzcd}[end anchor=north]
		X_1\drar[no head, end anchor={[xshift=-2.5ex]north east}] \ar[ddr,no head, end anchor={[xshift=-2.5ex]north east}]& &X_2\ar[no head,ddl, end anchor={[xshift=-1.3ex]north east}]\ar[dl,no head] \\
		& X_{4}\dar[no head, end anchor={[xshift=-2.5ex]north east},start anchor={[xshift=-2.5ex]south east }]& X_3\ar[no head,dl, end anchor={[xshift=-1.3ex]north east},start anchor={[xshift=-3.ex]south east }] \\
		&X_5&\quad\;.
	\end{tikzcd}

\end{center}


 The composition law (\ref{group law in G}) of $N_{5,2,2}$ is given by:
 
\begin{itemize}
    \item $z_1=x_1+y_1$;
    \item $z_2=x_2+y_2$;
    \item $z_3=x_3+y_3$;
    \item $z_4=x_4+y_4+\frac{1}{2}(x_1y_2-x_2y_1)$;
    \item $z_5=x_5+y_5+\frac{1}{2}(x_1y_4-x_4y_1+x_2y_3-x_3y_2)+\frac{1}{12}(x_1-y_1)(x_1y_2-x_2y_1)$.
\end{itemize}

Since
\begin{eqnarray*}
   \mathrm{d} (L_{\bf x})_{\bf 0}= \left[\begin{matrix} 
   1 & 0 & 0& 0 & 0\\
   0 & 1 & 0& 0& 0\\
   0 & 0 &1& 0& 0\\
 -\frac{x_2}{2} & \frac{x_1}{2} &0& 1 & 0\\
 -\frac{x_4}{2}-\frac{x_1x_2}{12}& -\frac{x_3}{2}+\frac{x_1^2}{12}&\frac{x_2}{2} & \frac{x_1}{2}& 1
   \end{matrix}\right]\,,
\end{eqnarray*}
the induced left-invariant vector fields \eqref{leftinvariant vf} are: 
\begin{itemize}
\item $X_1={\partial}_{x_1}-\frac{x_2}{2}{\partial}_{ x_4}-\big(\frac{x_4}{2}+\frac{x_1x_2}{12}\big){\partial}_{ x_5}\,;$
\item $X_2={\partial}_{ x_2}+\frac{x_1}{2}{\partial}_{ x_4}+\big(\frac{x_1^2}{12}-\frac{x_3}{2}\big){\partial}_{ x_5}\,;$
\item $X_3={\partial}_{ x_3}+\frac{x_2}{2}{\partial}_{x_5}\,;$
\item $X_4={\partial}_{ x_4}+\frac{x_1}{2}{\partial}_{ x_5}\,;$
\item $X_5={\partial}_{ x_5}\,,$
\end{itemize}
and the respective left-invariant 1-forms \eqref{leftinvariant form} are:
\begin{itemize}
\item $\theta_1=dx_1$;
\item $\theta_2=dx_2$;
\item $\theta_3=dx_3$;
\item $\theta_4=dx_4-\frac{x_1}{2}dx_2+\frac{x_2}{2}dx_1$;
\item $\theta _5=dx_5-\frac{x_1}{2}dx_4-\frac{x_2}{2}dx_3+\big(\frac{x_1^2}{6}+\frac{x_3}{2}\big)dx_2+\big(\frac{x_4}{2}-\frac{x_1x_2}{6}\big)dx_1\,.$
\end{itemize}

Finally, we have\begin{eqnarray*}
   \mathrm{d} (R_{\bf x})_{\bf 0}= \left[\begin{matrix} 
   1 & 0 & 0& 0 & 0\\
   0 & 1 & 0& 0& 0\\
   0 & 0 &1& 0& 0\\
 \frac{x_2}{2} & -\frac{x_1}{2} &0& 1 & 0\\
 \frac{x_4}{2}-\frac{x_1x_2}{12}& \frac{x_3}{2}+\frac{x_1^2}{12}&-\frac{x_2}{2} & -\frac{x_1}{2}& 1
   \end{matrix}\right]\,.
\end{eqnarray*}

\subsubsection{Grading, polarizations of maximal dimension, and asymptotic cone}
The Lie algebra $N_{5,2,2}$
   is not stratifiable, but it is gradable as
   \begin{eqnarray*}
   V_1 &= & \span \{  X_1,X_2 \} \,,\\
   V_2 &= & \span \{  X_3,X_4 \} \,,\\
   V_3 &= & \span \{  X_5\} \,.
   \end{eqnarray*}

      We claim that in this Lie algebra every two complementary subspaces to the  derived subalgebra, as in \eqref{complementary_Delta}, differ by an automorphism. Indeed, when $u_i^j$ varies, 
      the matrix
$$    \left[  \begin{array}{cccccc}
     1 &  0 &  0 &  0 &  0   \\
 -u_3^4 &  1 &  0 &  0 &  0 \\
 0 &  0 &  1 &  0 &  0 \\
 u_1^4 &  u_2^4 &  u_3^4 &  1 &  0 \\
 u_1^5 &  u_2^5 & u_5^4 & u_2^4&  1 
      \end{array}\right]$$
         is a Lie algebra automorphism and sends the complementary subspace
      $\span\{X_1, X_2, X_3\}$ to $\span\lbrace X_1-u_3^4X_2+u_1^4X_4+u_1^5X_5,X_2+u_2^4X_4+u_2^5X_5,X_3+u_3^4X_4+u_3^5X_5\rbrace$, which is an arbitrary subspace as in \eqref{complementary_Delta}. In particular, every $\Delta$ as in \eqref{complementary_Delta} gives maximal dimension.
The tangent cone of each of such polarizations has Lie algebra isomorphic to $N_{5,3,2}$, see page \pageref{N532}.

 The asymptotic cone of the Lie group with Lie algebra $N_{5,2,2}$ has Lie algebra isomorphic to $N_{4,2}\times\mathbb{R}$, where $N_{4,2}$ is the filiform algebra of step 3.
 



      
   
  \subsection*{$\bf N_{5,2,3}$ } \label{free 5-dim}  The following Lie algebra is denoted as $N_{5,2,3}$ by Gong in \cite{Gong_Thesis}, as $L_{5,9}$ by de Graaf in \cite{deGraaf_classification}, as (3) by Del Barco in \cite{delBarco}, and as $\mathcal{G}_{5,4}$ by Magnin in \cite{magnin}.
  
  The non-trivial brackets are the following:
\begin{eqnarray*}\label{N523}
    [X_1,X_2]=X_3\,,\,[X_1,X_3]=X_4\,,\,[X_2,X_3]=X_5\,.
\end{eqnarray*}
This is a nilpotent Lie algebra of rank 2 and step 3 that is stratifiable, also known as the free-nilpotent Lie algebra of step 3 and 2 generators. 
Such algebra will also be studied in the section of free-nilpotent Lie algebras, see page \pageref{F23}.
The Lie brackets can be pictured with the diagram:
  
  
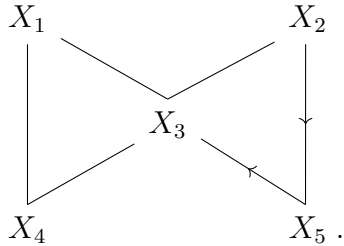
\begin{center}
 
	\begin{tikzcd}[end anchor=north]
		X_1 \ar[dr, no head]\ar[no head, dd] &  & X_2\;\;\ar[dl, no head]\ar[no head,dd, ->-=.5, end anchor={[xshift=-3.9ex]north east},start anchor={[xshift=-3.9ex]south east }] \\
		&  X_3\ar[no head, dl]\ar[no head, dr, -<-=.5, end anchor={[xshift=-3.9ex]north east}]  & &\\
		 X_4 & &  X_5\;.
	\end{tikzcd}
	\end{center}
	


 The composition law (\ref{group law in G}) of $N_{5,2,3}$ is given by:
\begin{itemize}
    \item $z_1=x_1+y_1$;
    \item $z_2=x_2+y_2$;
    \item $z_3=x_3+y_3+\frac{1}{2}(x_1y_2-x_2y_1)$;
    \item $z_4=x_4+y_4+\frac{1}{2}(x_1y_3-x_3y_1)+\frac{1}{12}(x_1-y_1)(x_1y_2-x_2y_1)$;
    \item $z_5=x_5+y_5+\frac{1}{2}(x_2y_3-x_3y_2)+\frac{1}{12}(x_2-y_2)(x_1y_2-x_2y_1)$.
\end{itemize}


Since
\begin{eqnarray*}
   \mathrm{d} (L_{\bf x})_{\bf 0}= \left[\begin{matrix} 
   1 & 0 & 0& 0 & 0\\
   0 & 1 & 0& 0& 0\\
   -\frac{x_2}{2} & \frac{x_1}{2} &1& 0& 0\\
  -\frac{x_1x_2}{12} -\frac{x_3}{2} & \frac{x_1^2}{12} &\frac{x_1}{2}& 1 & 0\\
  -\frac{x_2^2}{12}& -\frac{x_3}{2}+\frac{x_1x_2}{12}&\frac{x_2}{2} & 0& 1
   \end{matrix}\right]\,,
\end{eqnarray*}
the induced left-invariant vector fields \eqref{leftinvariant vf} are:
\begin{itemize}
\item $X_1={\partial}_{x_1}-\frac{x_2}{2}{\partial}_{ x_3}-\big(\frac{x_3}{2}+\frac{x_1x_2}{12}\big){\partial}_{x_4}-\frac{x_2^2}{12}{\partial}_{x_5}\,;$
\item $X_2={\partial}_{ x_2}+\frac{x_1}{2}{\partial}_{x_3}+\frac{x_1^2}{12}{\partial}_{x_4}-\big(\frac{x_3}{2}-\frac{x_1x_2}{12}\big){\partial}_{x_5}\,;$
\item $X_3={\partial}_{x_3}+\frac{x_1}{2}{\partial}_{x_4}+\frac{x_2}{2}{\partial}_{ x_5}\,;$
\item $X_4={\partial}_{x_4}\,;$
\item $X_5={\partial}_{x_5}\,,$
\end{itemize}
and the respective left-invariant 1-forms \eqref{leftinvariant form} are:
\begin{itemize}
\item $\theta_1=dx_1$;
\item $\theta_2=dx_2$;
\item $\theta_3=dx_3-\frac{x_1}{2} dx_2+\frac{x_2}{2}dx_1$;
\item $\theta_4=dx_4-\frac{x_1}{2}dx_3+\frac{x_1^2}{6}dx_2+\big(\frac{x_3}{2}-\frac{x_1x_2}{6}\big)dx_1$;
\item $\theta_5=dx_5-\frac{x_2}{2}dx_3+\big(\frac{x_3}{2}+\frac{x_1x_2}{6}\big)dx_2-\frac{x_2^2}{6}dx_1\,.$
\end{itemize}

Finally, we have
\begin{eqnarray*}
   \mathrm{d} (R_{\bf x})_{\bf 0}= \left[\begin{matrix} 
   1 & 0 & 0& 0 & 0\\
   0 & 1 & 0& 0& 0\\
   \frac{x_2}{2} & -\frac{x_1}{2} &1& 0& 0\\
  \frac{x_3}{2}-\frac{x_1x_2}{12}  & \frac{x_1^2}{12} &-\frac{x_1}{2}& 1 & 0\\
  -\frac{x_2^2}{12}& \frac{x_3}{2}+\frac{x_1x_2}{12}&-\frac{x_2}{2} & 0& 1
   \end{matrix}\right]\,.
\end{eqnarray*}
  
   \subsection*{$\bf N_{5,3,1}$ }    The following Lie algebra is denoted as $N_{5,3,1}$ by Gong in \cite{Gong_Thesis}, as $L_{5,4}$ by de Graaf in \cite{deGraaf_classification}, as (8) by Del Barco in \cite{delBarco}, and as $\mathcal{G}_{5,1}$ by Magnin in \cite{magnin}.
 
 The non-trivial brackets are the following:
\begin{eqnarray*}
    [X_1,X_2]=[X_3,X_4]=X_5\,.
\end{eqnarray*}
This is a nilpotent Lie algebra of rank 4 and step 2 that is stratifiable, also known as the second Heisenberg algebra. The Lie brackets can be pictured with the diagram:
  
  
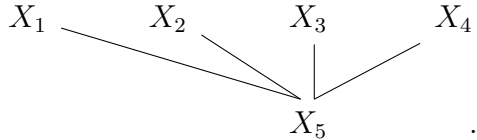
\begin{center}
 
	\begin{tikzcd}[end anchor=north]
		X_1 \ar[no head, drr,end anchor={[xshift=-3.ex]north east}]  & X_2 \ar[no head, dr,end anchor={[xshift=-3.ex]north east}] &   X_3\ar[no head, d, end anchor={[xshift=-2.ex]north east},start anchor={[xshift=-2.ex]south east}] & X_4\ar[no head, dl, end anchor={[xshift=-2.ex]north east}]\\
          & & X_5  &\quad\;.
	\end{tikzcd}
	\end{center}
	


 The composition law (\ref{group law in G}) of $N_{5,3,1}$ is given by:
\begin{itemize}
    \item $z_1=x_1+y_1$;
    \item $z_2=x_2+y_2$;
    \item $z_3=x_3+y_3$;
    \item $z_4=x_4+y_4$;
    \item $z_5=x_5+y_5+\frac{1}{2}(x_1y_2-x_2y_1+x_3y_4-x_4y_3)$.
\end{itemize}

Since
\begin{eqnarray*}
   \mathrm{d} (L_{\bf x})_{\bf 0}= \left[\begin{matrix} 
   1 & 0 & 0& 0 & 0\\
   0 & 1 & 0& 0& 0\\
  0 & 0&1& 0& 0\\
  0 & 0 & 0& 1 & 0\\
  -\frac{x_2}{2}& \frac{x_1}{2} & -\frac{x_4}{2}&\frac{x_3}{2}& 1
   \end{matrix}\right]\,,
\end{eqnarray*}
the induced left-invariant vector fields \eqref{leftinvariant vf} are:
\begin{itemize}
\item $X_1={\partial}_{x_1}-\frac{x_2}{2}{\partial}_{ x_5}\,;$
\item $X_2={\partial}_{x_2}+\frac{x_1}{2}{\partial}_{x_5}\,;$
\item $X_3={\partial}_{x_3}-\frac{x_4}{2}{\partial}_{x_5}\,;$
\item $X_4={\partial}_{x_4}+\frac{x_3}{2}{\partial}_{x_5}\,;$
\item $X_5={\partial}_{x_5}\,,$
\end{itemize}
and the respective left-invariant 1-forms \eqref{leftinvariant form} are:
\begin{itemize}
\item $\theta_1=dx_1$;
\item $\theta_2=dx_2$;
\item $\theta_3=dx_3$;
\item $\theta_4=dx_4$;
\item $\theta_5=dx_5-\frac{x_3}{2}dx_4+\frac{x_4}{2}dx_3-\frac{x_1}{2}dx_2+\frac{x_2}{2}dx_1\,.$
\end{itemize}
 
 Finally, we have
 \begin{eqnarray*}
   \mathrm{d} (R_{\bf x})_{\bf 0}= \left[\begin{matrix} 
   1 & 0 & 0& 0 & 0\\
   0 & 1 & 0& 0& 0\\
  0 & 0&1& 0& 0\\
  0 & 0 & 0& 1 & 0\\
  \frac{x_2}{2}& -\frac{x_1}{2} & \frac{x_4}{2}&-\frac{x_3}{2}& 1
   \end{matrix}\right]\,.
\end{eqnarray*}
   \subsection*{$\bf N_{5,3,2}$  } The following Lie algebra is denoted as $N_{5,3,2}$ by Gong in \cite{Gong_Thesis}, as $L_{5,8}$ by de Graaf in \cite{deGraaf_classification}, as (6) by Del Barco in \cite{delBarco}, and as $\mathcal{G}_{5,2}$ by Magnin in \cite{magnin}.  
  
  The non-trivial brackets are the following:
\begin{eqnarray*}\label{N532}
    [X_1,X_2]=X_4\,,\,[X_1,X_3]=X_5\,.
\end{eqnarray*}
This is a nilpotent Lie algebra of rank 3 and step 2 that is stratifiable. The Lie brackets can be pictured with the diagram:
  
  
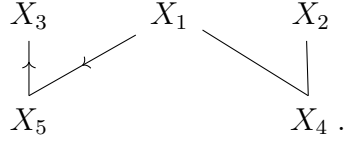
\begin{center}
 
	\begin{tikzcd}[end anchor=north]
		X_3  \ar[d, no head, -<-=.5] & X_1 \ar[dl, no head, ->-=.5] \ar[dr, no head,end anchor={[xshift=-3.8ex]north east}] & X_2\ar[d, no head, end anchor={[xshift=-3.8ex]north east}, start anchor={[xshift=-3.8ex]south east}]\;\,\\
		X_5 & & X_4\;.
	\end{tikzcd}
	\end{center}
	


 The composition law (\ref{group law in G}) of $N_{5,3,2}$ is given by:

\begin{itemize}
    \item $z_1=x_1+y_1$;
    \item $z_2=x_2+y_2$;
    \item $z_3=x_3+y_3$;
    \item $z_4=x_4+y_4+\frac{1}{2}(x_1y_2-x_2y_1)$;
    \item $z_5=x_5+y_5+\frac{1}{2}(x_1y_3-x_3y_1)$.
\end{itemize}

Since
\begin{eqnarray*}
   \mathrm{d} (L_{\bf x})_{\bf 0}= \left[\begin{matrix} 
   1 & 0 & 0& 0 & 0\\
   0 & 1 & 0& 0& 0\\
  0 & 0&1& 0& 0\\
  -\frac{x_2}{2} & \frac{x_1}{2} & 0& 1 & 0\\
  -\frac{x_3}{2}& 0 & \frac{x_1}{2}&0& 1
   \end{matrix}\right]\,,
\end{eqnarray*}
the induced left-invariant vector fields \eqref{leftinvariant vf} are:
\begin{itemize}
\item $X_1={\partial}_{x_1}-\frac{x_2}{2}{\partial}_{x_4}-\frac{x_3}{2}{\partial}_{ x_5}\,;$
\item $X_2={\partial}_{ x_2}+\frac{x_1}{2}{\partial}_{x_4}\,;$
\item $X_3={\partial}_{x_3}+\frac{x_1}{2}{\partial}_{x_5}\,;$
\item $X_4={\partial}_{x_4}\,;$
\item $X_5={\partial}_{ x_5}\,,$
\end{itemize}
and the respective left-invariant 1-forms \eqref{leftinvariant form} are:
\begin{itemize}
\item $\theta_1=dx_1$;
\item $\theta_2=dx_2$;
\item $\theta_3=dx_3$;
\item $\theta_4=dx_4-\frac{x_1}{2}dx_2+\frac{x_2}{2}dx_1$;
\item $\theta_5=dx_5-\frac{x_1}{2}dx_3+\frac{x_3}{2}dx_1\,.$
\end{itemize}

Finally, we have
\begin{eqnarray*}
   \mathrm{d} (R_{\bf x})_{\bf 0}= \left[\begin{matrix} 
   1 & 0 & 0& 0 & 0\\
   0 & 1 & 0& 0& 0\\
  0 & 0&1& 0& 0\\
  \frac{x_2}{2} & -\frac{x_1}{2} & 0& 1 & 0\\
  \frac{x_3}{2}& 0 & -\frac{x_1}{2}&0& 1
   \end{matrix}\right]\,.
\end{eqnarray*}

      \newpage
  \section{6D indecomposable    nilpotent Lie algebras}
  Among all the indecomposable nilpotent Lie algebras in dimension 6, there are 13 Carnot algebras and 11 more nilpotent Lie algebras, which are gradable.

  \subsection*{$\bf N_{6,1,1}$ non-stratifiable} The following Lie algebra is denoted as $N_{6,1,1}$ by Gong in \cite{Gong_Thesis}, as $L_{6,15}$ by de Graaf in \cite{deGraaf_classification}, as (4) by Del Barco in \cite{delBarco}, and as $\mathcal{G}_{6,19}$ by Magnin in \cite{magnin}.

The non-trivial brackets are the following:
\begin{equation*}
     [X_1,X_i]=X_{i+1}\,,\, 2 \le i\le 5\,,\,[X_2,X_3]=X_5\,,\,[X_2,X_4]=X_6\,.
\end{equation*}
This is a nilpotent Lie algebra of rank 2 and step 5 that is positively gradable, yet not stratifiable. The Lie brackets can be pictured with the  diagram:
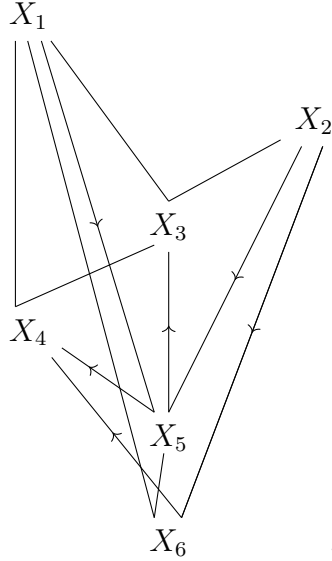
\begin{center}
 
	\begin{tikzcd}[end anchor=north]
	X_1\ar[no head, ddr,start anchor={[xshift=-.8ex]south east}]\ar[no head, ddd,start anchor={[xshift=-3.4ex]south east},end anchor={[xshift=-3.4ex]north east}]\ar[->-=.5,no head, ddddr, end anchor={[xshift=-3.5ex]north east},start anchor={[xshift=-1.5ex]south east}]\ar[no head, dddddr,end anchor={[xshift=-3.5ex]north east}, start anchor={[xshift=-2.5ex]south east}]&  & \\
	&  & X_2\ar[no head, dl]\ar[->-=.5,no head, ddddl, end anchor={[xshift=-1.5ex]north east},start anchor={[xshift=-1.8ex]south east}]\ar[->-=.5, no head, dddl]\\
	& X_3\ar[no head, dl, start anchor={[xshift=-3.5ex]south east},start anchor={[yshift=0.5ex]south east},end anchor={[xshift=-3.4ex]north east}]\ar[-<-=.5, no head, dd] & \\
	X_4 \ar[-<-=.35,no head, dr, end anchor={[xshift=-3.5ex]north east}]\ar[-<-=.5,no head, ddr,end anchor={[xshift=-1.5ex]north east}] &  &\\
	& X_5\ar[no head, d,  end anchor={[xshift=-3.5ex]north east}, start anchor={[xshift=-2.8ex]south east},start anchor={[yshift=0.7ex]south east}]&\\
	& X_6 \ar[no head, uuuur, start anchor={[xshift=-1.5ex]north east}, end anchor={[xshift=-1.8ex]south east}]&\quad\; .
	\end{tikzcd}

\end{center}

The composition law (\ref{group law in G}) of $N_{6,1,1}$ is given by:
\begin{itemize}
    \item $z_1=x_1+y_1$;
    \item $z_2=x_2+y_2$;
    \item $z_3=x_3+y_3+\frac{1}{2}(x_1y_2-x_2y_1)$;
    \item $z_4=x_4+y_4+\frac{1}{2}(x_1y_3-x_3y_1)+\frac{1}{12}(x_1-y_1)(x_1y_2-x_2y_1)$;
    \item $z_5=x_5+y_5+\frac{1}{2}(x_1y_4-x_4y_1+x_2y_3-x_3y_2)+\frac{1}{12}(x_1-y_1)(x_1y_3-x_3y_1)\salto +\frac{1}{12}(x_2-y_2)(x_1y_2-x_2y_1)-\frac{1}{24}x_1y_1(x_1y_2-x_2y_1)$;
    \item $z_6=x_6+y_6+\frac{1}{2}(x_1y_5-x_5y_1+x_2y_4-x_4y_2)+\frac{1}{12}(x_2-y_2)(x_1y_3-x_3y_1)\salto+\frac{1}{12}(x_1-y_1)(x_1y_4-x_4y_1+x_2y_3-x_3y_2)  -\frac{1}{24}x_1y_1(x_1y_3-x_3y_1)\salto-\frac{1}{24}(x_2y_1+x_1y_2)(x_1y_2-x_2y_1) +\frac{1}{720}(y_1^3-x_1^3)(x_1y_2-x_2y_1)\salto+\frac{1}{180}(x_1y_1^2-x_1^2y_1)(x_1y_2-x_2y_1)\,.$
\end{itemize}

Since
\begin{eqnarray*}
   \mathrm{d} (L_{\bf x})_{\bf 0}= \left[\begin{matrix} 
   1 & 0 & 0& 0 & 0 & 0\\
   0 & 1 & 0& 0& 0& 0\\
   -\frac{x_2}{2} & \frac{x_1}{2} &1& 0& 0& 0\\
  -\frac{x_1x_2}{12} -\frac{x_3}{2} & \frac{x_1^2}{12} &\frac{x_1}{2}& 1 & 0& 0\\
  -\frac{x_1x_3+x_2^2}{12}-\frac{x_4}{2}& \frac{x_1x_2}{12}-\frac{x_3}{2} &\frac{x_1^2}{12}+\frac{x_2}{2} & \frac{x_1}{2}& 1& 0\\ \frac{x_1^3x_2}{720}-\frac{x_1x_4+x_2x_3}{12}-\frac{x_5}{2}& -\frac{x_1x_3}{12}-\frac{x_4}{2}-\frac{x_1^4}{720} &\frac{x_1x_2}{6} &\frac{x_2}{2}+\frac{x_1^2}{12} &\frac{x_1}{2}& 1
   \end{matrix}\right]\,,
\end{eqnarray*}
the induced left-invariant vector fields \eqref{leftinvariant vf} are: 
\begin{itemize}
\item $X_1={\partial}_{x_1}-\frac{x_2}{2}{\partial}_{x_3}-\big(\frac{x_3}{2}+\frac{x_1x_2}{12}\big){\partial}_{x_4}-\big(\frac{x_4}{2}+\frac{x_1x_3+x_2^2}{12}\big){\partial}_{x_5}+\big(\frac{x_1^3x_2}{720}-\frac{x_1x_4+x_2x_3}{12}-\frac{x_5}{2}\big){\partial}_{x_6}\,;$
\item $X_2={\partial}_{x_2}+\frac{x_1}{2}{\partial}_{x_3}+\frac{x_1^2}{12}{\partial}_{x_4}+\big(\frac{x_1x_2}{12}-\frac{x_3}{2}\big){\partial}_{x_5}-\big(\frac{x_1x_3}{12}+\frac{x_4}{2}+\frac{x_1^4}{720}\big){\partial}_{x_6}\,;$
\item $X_3={\partial}_{x_3}+\frac{x_1}{2}{\partial}_{x_4}+\big(\frac{x_1^2}{12}+\frac{x_2}{2}\big){\partial}_{x_5}+\frac{x_1x_2}{6}{\partial}_{x_6}\,;$
\item $X_4={\partial}_{x_4}+\frac{x_1}{2}{\partial}_{x_5}+\big(\frac{x_2}{2}+\frac{x_1^2}{12}\big){\partial}_{x_6}\,;$
\item $X_5={\partial}_{x_5}+\frac{x_1}{2}{\partial}_{x_6}$;
\item $X_6={\partial}_{x_6}\,,$
\end{itemize}
and the respective left-invariant 1-forms \eqref{leftinvariant form} are: 
\begin{itemize}
\item $\theta_1=dx_1$;
\item $\theta_2=dx_2$;
\item $\theta_3=dx_3-\frac{x_1}{2} dx_2+\frac{x_2}{2}dx_1$;
\item $\theta_4=dx_4-\frac{x_1}{2}dx_3+\frac{x_1^2}{6}dx_2+\big(\frac{x_3}{2}-\frac{x_1x_2}{6}\big)dx_1$;
\item $\theta_5=dx_5-\frac{x_1}{2}dx_4+\big(\frac{x_1^2}{6}-\frac{x_2}{2}\big)dx_3+\big(\frac{x_1x_2}{6}+\frac{x_3}{2}-\frac{x_1^3}{24}\big)dx_2+\big(\frac{x_4}{2}-\frac{x_1x_3+x_2^2}{6}+\frac{x_1^2x_2}{24}\big)dx_1$;
\item $\theta_6=dx_6-\frac{x_1}{2}dx_5+\big(\frac{x_1^2}{6}-\frac{x_2}{2}\big)dx_4+\big(\frac{x_1x_2}{3}-\frac{x_1^3}{24}\big)dx_3+\big(\frac{x_4}{2}-\frac{x_1x_3}{6}-\frac{x_1^2x_2}{12}+\frac{x_1^4}{120}\big)dx_2$\\$  \phantom{\theta_5=}+\big(\frac{x_5}{2}-\frac{x_2x_3+x_1x_4}{6}+\frac{x_1x_2^2}{12}+\frac{x_1^2x_3}{24}-\frac{x_1^3x_2}{120}\big)dx_1\,.$
\end{itemize}

Finally, we have 

\begin{eqnarray*}
   \mathrm{d} (R_{\bf x})_{\bf 0}= \left[\begin{matrix} 
   1 & 0 & 0& 0 & 0 & 0\\
   0 & 1 & 0& 0& 0& 0\\
   \frac{x_2}{2} & -\frac{x_1}{2} &1& 0& 0& 0\\
  \frac{x_3}{2}-\frac{x_1x_2}{12}  & \frac{x_1^2}{12} &-\frac{x_1}{2}& 1 & 0& 0\\
  \frac{x_4}{2}-\frac{x_1x_3+x_2^2}{12}& \frac{x_3}{2}+\frac{x_1x_2}{12} &\frac{x_1^2}{12}-\frac{x_2}{2} & -\frac{x_1}{2}& 1& 0\\ \frac{x_5}{2}-\frac{x_1x_4+x_2x_3}{12}+\frac{x_1^3x_2}{720}& \frac{x_4}{2}-\frac{x_1x_3}{12}-\frac{x_1^4}{720} &\frac{x_1x_2}{6} &\frac{x_1^2}{12}-\frac{x_2}{2} &-\frac{x_1}{2}& 1
   \end{matrix}\right]\,.
\end{eqnarray*}

\subsubsection{Grading, polarizations of maximal dimension, and asymptotic cone}
 The Lie algebra $N_{6,1,1}$
   is not stratifiable, but it is gradable as
   \begin{eqnarray*}
   V_i &= & \span \{  X_i \} \,,\,i=1,\ldots,6\,.
   \end{eqnarray*}

      We claim that in this Lie algebra every two complementary subspaces to the  derived subalgebra as in \eqref{complementary_Delta} differ by an automorphism. Indeed, when $u_i^j$ varies, 
      the matrix
$$    \left[  \begin{array}{cccccc}
     1 &  0 &  0 &  0 &  0 &  0  \\
 0 &  1 &  0 &  0 &  0 &  0\\
 u_1^3 & u_2^3 & 1 &  0 & 0 &  0\\
 u_1^4 & u_2^4 & u_2^3 & 1 & 0 & 0\\
 u_1^5 &  u_2^5 & u_2^4-u_1^3 &  u_2^3 &  1 &  0\\
 u_1^6 &  u_2^6 &  u_2^5-u_1^4 &  u_2^4-u_1^3 &  u_2^3 &  1\\
      \end{array}\right]$$
         is a Lie algebra automorphism and   sends the complementary subspace
      $\span\{X_1, X_2\}$ to $\span\lbrace X_1+u_1^3X_3+u_1^4X_4+u_1^5X_5+u_1^6X_6,X_2+u_2^3X_3+u_2^4X_4+u_2^5X_5+u_2^6X_6\rbrace$, which is an arbitrary subspace as in \eqref{complementary_Delta}. In particular, every $\Delta$ as in \eqref{complementary_Delta} gives maximal dimension.
The tangent cone of each of such polarizations has Lie algebra isomorphic to $N_{6,2,5}$, see page \pageref{N625}.
 
 The asymptotic cone 
 of the Lie group with Lie algebra $N_{6,1,1}$ has Lie algebra isomorphic to
 the filiform algebra $N_{6,2,1}$ of step 5, see page \pageref{N621}.
      \subsection*{$\bf N_{6,1,2}$ non-stratifiable }  The following Lie algebra is denoted as $N_{6,1,2}$ by Gong in \cite{Gong_Thesis}, as $L_{6,14}$ by de Graaf in \cite{deGraaf_classification}, as (1) by Del Barco in \cite{delBarco}, and as $\mathcal{G}_{6,20}$ by Magnin in \cite{magnin}.
  
  The non-trivial brackets are the following:
\begin{equation*}
    [X_1,X_i]=X_{i+1}\,,\, 2 \le i\le 4\,,\,[X_2,X_3]=X_5\,,\,[X_2,X_5]=X_6\,,\,[X_3,X_4]=-X_6\,.
\end{equation*}
This is a nilpotent Lie algebra of rank 2 and step 5 that is positively gradable, yet not stratifiable. The Lie brackets can be pictured with the diagram:
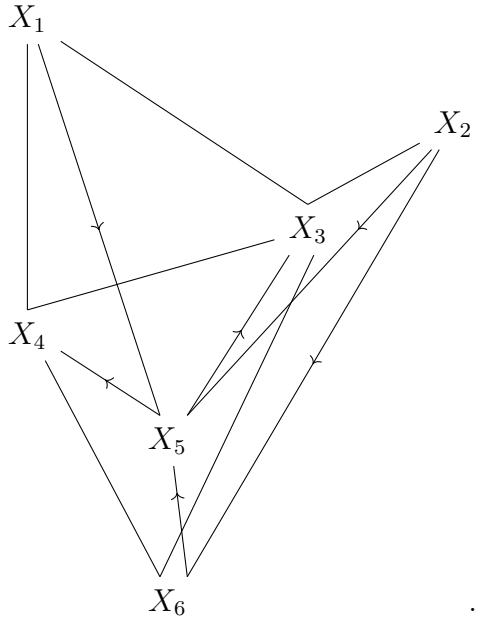
\begin{center}
 
	\begin{tikzcd}[end anchor=north]
	X_1\ar[ddrr, no head]\ar[ddddr, no head, ->-=.5, end anchor={[xshift=-3.ex]north east} ]\ar[ddd, no head] &  & &\\
	& & &X_2\ar[dddddll, no head, ->-=.5,end anchor={[xshift=-1.ex]north east} ]\ar[dl, no head]\ar[dddll, no head,->-=.3, end anchor={[xshift=-1.ex]north east} ]\\
	& & X_3\ar[ddddl, no head, end anchor={[xshift=-3.ex]north east},start anchor={[xshift=-2ex]south east} ]\ar[dll, no head]\ar[ddl, no head, -<-=.5, end anchor={[xshift=-1.ex]north east},start anchor={[xshift=-3.8ex]south east} ]&\\
	X_4\ar[dr, no head, -<-=.5, end anchor={[xshift=-3.ex]north east} ]\ar[dddr, no head, end anchor={[xshift=-3.ex]north east}] & & &\\
	& X_5 \ar[dd, no head, -<-=.3, end anchor={[xshift=-1.ex]north east} ]& &\\
	& &  & \\
	&  X_6 & & \quad\;.
	\end{tikzcd}

\end{center}

 The composition law \eqref{group law in G} of $N_{6,1,2}$ is given by:

\begin{itemize}
    \item $z_1=x_1+y_1$;
    \item $z_2=x_2+y_2$;
    \item $z_3=x_3+y_3+\frac{1}{2}(x_1y_2-x_2y_1)$;
    \item $z_4=x_4+y_4+\frac{1}{2}(x_1y_3-x_3y_1)+\frac{1}{12}(x_1-y_1)(x_1y_2-x_2y_1)$;
    \item $z_5=x_5+y_5+\frac{1}{2}(x_1y_4-x_4y_1+x_2y_3-x_3y_2)+\frac{1}{12}(x_1-y_1)(x_1y_3-x_3y_1)$\\$  \phantom{z_5=}+\frac{1}{12}(x_2-y_2)(x_1y_2-x_2y_1)-\frac{1}{24}x_1y_1(x_1y_2-x_2y_1)$;
    \item $z_6=x_6+y_6+\frac{1}{2}(x_2y_5-x_5y_2-x_3y_4+x_4y_3)+\frac{1}{12}(x_4-y_4)(x_1y_2-x_2y_1)\salto+\frac{1}{12}(x_2-y_2)(x_1y_4-x_4y_1+x_2y_3-x_3y_2)+\frac{1}{12}(y_3-x_3)(x_1y_3-x_3y_1)\salto-\frac{1}{24}[x_1y_2(x_1y_3-x_3y_1)+(x_2y_2-x_1y_3)(x_1y_2-x_2y_1)]-\frac{1}{360}y_1(x_1y_2-x_2y_1)^2\salto+\frac{1}{720}(y_1^2y_2-x_1^2x_2)(x_1y_2-x_2y_1)+\frac{1}{180}(x_1y_1y_2-x_1^2y_2)(x_1y_2-x_2y_1)]$\\$  \phantom{z_5=}-\frac{1}{120}x_1(x_1y_2-x_2y_1)^2\,.$
\end{itemize}

Since
\begin{eqnarray*}
   \mathrm{d} (L_{\bf x})_{\bf 0}= \left[\begin{matrix} 
   1 & 0 & 0& 0 & 0 & 0\\
   0 & 1 & 0& 0& 0& 0\\
   -\frac{x_2}{2} & \frac{x_1}{2} &1& 0& 0& 0\\
  -\frac{x_1x_2}{12} -\frac{x_3}{2} & \frac{x_1^2}{12} &\frac{x_1}{2}& 1 & 0& 0\\
  -\frac{x_1x_3+x_2^2}{12}-\frac{x_4}{2}& \frac{x_1x_2}{12}-\frac{x_3}{2}&\frac{x_1^2}{12}+\frac{x_2}{2} & \frac{x_1}{2}& 1& 0\\ \frac{x_3^2-2x_2x_4}{12}+\frac{x_1^2x_2^2}{720}& \frac{x_1x_4-x_2x_3}{12}-\frac{x_5}{2}-\frac{x_1^3x_2}{720} &\frac{x_4}{2}+\frac{x_2^2-x_1x_3}{12} &\frac{x_1x_2}{12}-\frac{x_3}{2} &\frac{x_2}{2}& 1
   \end{matrix}\right]\,,
\end{eqnarray*}
the induced left-invariant vector fields \eqref{leftinvariant vf} are: 
\begin{itemize}
\item $X_1={\partial}_{x_1}-\frac{x_2}{2}{\partial}_{x_3}-\big(\frac{x_3}{2}+\frac{x_1x_2}{12}\big){\partial}_{x_4}-\big(\frac{x_4}{2}+\frac{x_1x_3+x_2^2}{12}\big){\partial}_{x_5}+\big(\frac{x_3^2-2x_2x_4}{12}+\frac{x_1^2x_2^2}{720}\big){\partial}_{x_6}\,;$
\item $X_2={\partial}_{x_2}+\frac{x_1}{2}{\partial}_{x_3}+\frac{x_1^2}{12}{\partial}_{x_4}+\big(\frac{x_1x_2}{12}-\frac{x_3}{2}\big){\partial}_{x_5}+\big(\frac{x_1x_4-x_2x_3}{12}-\frac{x_5}{2}-\frac{x_1^3x_2}{720}\big){\partial}_{x_6}\,;$
\item $X_3={\partial}_{x_3}+\frac{x_1}{2}{\partial}_{x_4}+\big(\frac{x_1^2}{12}+\frac{x_2}{2}\big){\partial}_{x_5}+\big(\frac{x_4}{2}+\frac{x_2^2-x_1x_3}{12} \big){\partial}_{x_6}\,;$
\item $X_4={\partial}_{x_4}+\frac{x_1}{2}{\partial}_{x_5}+\big(\frac{x_1x_2}{12}-\frac{x_3}{2}){\partial}_{x_6}\,;$
\item $X_5={\partial}_{x_5}+\frac{x_2}{2}{\partial}_{x_6}$;
\item $X_6={\partial}_{x_6}\,,$
\end{itemize}
and the respective left-invariant 1-forms \eqref{leftinvariant form} are: 
\begin{itemize}
\item $\theta_1=dx_1$;
\item $\theta_2=dx_2$;
\item $\theta_3=dx_3-\frac{x_1}{2} dx_2+\frac{x_2}{2}dx_1$;
\item $\theta_4=dx_4-\frac{x_1}{2}dx_3+\frac{x_1^2}{6}dx_2+\big(\frac{x_3}{2}-\frac{x_1x_2}{6}\big)dx_1$;
\item $\theta_5=dx_5-\frac{x_1}{2}dx_4+\big(\frac{x_1^2}{6}-\frac{x_2}{2}\big)dx_3+\big(\frac{x_1x_2}{6}+\frac{x_3}{2}-\frac{x_1^3}{24}\big)dx_2+\big(\frac{x_4}{2}-\frac{x_1x_3+x_2^2}{6}+\frac{x_1^2x_2}{24}\big)dx_1$;
\item $\theta_6=dx_6-\frac{x_2}{2}dx_5+\big(\frac{x_3}{2}+\frac{x_1x_2}{6}\big)dx_4+\big(\frac{x_2^2-x_1x_3}{6}-\frac{x_4}{2}-\frac{x_1^2x_2}{24}\big)dx_3+\big(\frac{x_5}{2}+\frac{x_1x_4-x_2x_3}{6}$\\$  \phantom{\theta_5=}+\frac{x_1^2x_3-x_1x_2^2}{24}+\frac{x_1^3x_2}{120}\big)dx_2+\big(\frac{x_3^2-2x_2x_4}{6}+\frac{x^3_2}{24}-\frac{x_1^2x_2^2}{120}\big)dx_1\,.$
\end{itemize}

Finally, we have
\begin{eqnarray*}
   \mathrm{d} (R_{\bf x})_{\bf 0}= \left[\begin{matrix} 
   1 & 0 & 0& 0 & 0 & 0\\
   0 & 1 & 0& 0& 0& 0\\
   \frac{x_2}{2} & -\frac{x_1}{2} &1& 0& 0& 0\\
  \frac{x_3}{2}-\frac{x_1x_2}{12}  & \frac{x_1^2}{12} &-\frac{x_1}{2}& 1 & 0& 0\\
  \frac{x_4}{2}-\frac{x_1x_3+x_2^2}{12}& \frac{x_3}{2}+\frac{x_1x_2}{12}&\frac{x_1^2}{12}-\frac{x_2}{2} & -\frac{x_1}{2}& 1& 0\\ \frac{x_3^2-2x_2x_4}{12}+\frac{x_1^2x_2^2}{720}& \frac{x_5}{2}+\frac{x_1x_4-x_2x_3}{12}-\frac{x_1^3x_2}{720} &\frac{x_2^2-x_1x_3}{12}-\frac{x_4}{2} &\frac{x_3}{2}+\frac{x_1x_2}{12} &-\frac{x_2}{2}& 1
   \end{matrix}\right]\,.
\end{eqnarray*}

\subsubsection{Grading, polarizations of maximal dimension, and asymptotic cone}
 The Lie algebra $N_{6,1,2}$
   is not stratifiable, but it is gradable as
   \begin{eqnarray*}
   V_i &= & \span \{  X_i \} \,,\,i=1,\ldots,5\,,\\
   V_6&=&0\,,\\
   V_7&=&\span\{X_6\}\,.
   \end{eqnarray*}
      
        
     When $u_i^j$ varies, 
      the matrix
$$    \left[  \begin{array}{cccccc}
     1 &  0 &  0 &  0 &  0 &  0  \\
 0 &  1 &  0 &  0 &  0 &  0\\
 u_1^3 & u_2^3 & 1 &  0 & 0 &  0\\
 u_1^4 & u_2^4 & u_2^3 & 1 & 0 & 0\\
 u_1^5 &  u_2^5 & u_2^4-u_1^3 &  u_2^3 &  1 &  0\\
 u_1^6 &  u_2^6 &  u_1^4u_2^3-u_1^3u_2^4-u_1^5 &  u_1^4-u_1^3u_2^3 &  2 u_2^4-u_1^3-u_2^3u_2^3 &  1\\
      \end{array}\right]$$
         is a Lie algebra automorphism when $2u_2^4=u_2^3u_2^3$, and sends the complementary subspace $\span\lbrace X_1,X_2\rbrace$ to
      $\span\{X_1 + u_1^3 X_3+ u_1^4 X_4+u_1^5X_5+ u_1^6 X_6 , X_2 +   u_2^3 X_3+   u_2^4 X_4+ u_2^5X_5+  u_2^6 X_6\}$, which is an arbitrary subspace as in \eqref{complementary_Delta}.
The tangent cone of each of such polarizations has Lie algebra isomorphic to $N_{6,2,7}$ if $2u_2^4=u_2^3u_2^3$, to $N_{6,2,5}$ if $2u_2^4>u_2^3u_2^3$, and to $N_{6,2,5\,a}$ if $2u_2^4<u_2^3u_2^3$, see pages \pageref{N627}, \pageref{N625}, and \pageref{N625a}, respectively. In particular, each of such polarizations gives maximal Hausdorff dimension. 
 
 The asymptotic cone 
 of the Lie group with Lie algebra $N_{6,1,2}$ has Lie algebra isomorphic to
 $N_{6,2,2}$, see page \pageref{N622}.

  
   
    \subsection*{$\bf N_{6,1,3} $ non-stratifiable} The following Lie algebra is denoted as $N_{6,1,3}$ by Gong in \cite{Gong_Thesis}, as $L_{6,17}$ by de Graaf in \cite{deGraaf_classification}, as $(5)$ by Del Barco in \cite{delBarco}, and as $\mathcal{G}_{6,17}$ by Magnin in \cite{magnin}.
  
  The non-trivial brackets are the following: 
\begin{equation*}
     [X_1, X_i] = X_{i+1}\,,\, 2\leq i\leq 5\,,\, [X_2, X_3] = X_6\,. 
\end{equation*}
This is a nilpotent Lie algebra of rank 2 and step 5 that is positively gradable, yet not stratifiable. The Lie brackets can be pictured with the diagram:
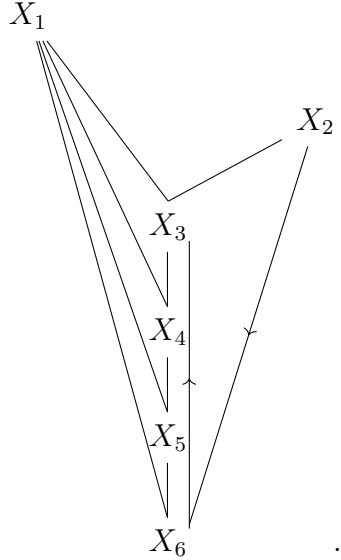
\begin{center}
 
	\begin{tikzcd}[end anchor=north]
		X_1\ar[ddr,no head] \ar[dddr,no head, end anchor={[xshift=-2.5ex]north east}]\ar[dddddr,no head, end anchor={[xshift=-2.5ex]north east}]\ar[ddddr,no head, end anchor={[xshift=-2.5ex]north east}]& & \\
		& & X_2\dlar[no head]\ar[->-=.5,no head,ddddl, end anchor={[xshift=-.9ex]north east},start anchor={[xshift=-3.ex]south east }, end anchor={[yshift=-.5ex]north east}]\\
		& X_{3}\dar[no head, end anchor={[xshift=-2.5ex]north east},start anchor={[xshift=-2.5ex]south east }]\ar[-<-=.5,no head,ddd, end anchor={[xshift=-0.9ex]north east},start anchor={[xshift=-0.9ex]south east }, end anchor={[yshift=-0.8ex]north east},start anchor={[yshift=0.8ex]south east }]&\\
		&X_4\dar[no head, end anchor={[xshift=-2.5ex]north east},start anchor={[xshift=-2.5ex]south east }]&\\
		&X_5\dar[no head,end anchor={[xshift=-2.5ex]north east},start anchor={[xshift=-2.5ex]south east }]&  \\
		& X_6 &  \phantom{X_5}\,.
	\end{tikzcd}
   
\end{center}

 
 The composition law \eqref{group law in G}  of $N_{6,1,3}$ is given by:

\begin{itemize}
    \item $z_1=x_1+y_1$;
    \item $z_2=x_2+y_2$;
    \item $z_3=x_3+y_3+\frac{1}{2}(x_1y_2-x_2y_1)$;
    \item $z_4=x_4+y_4+\frac{1}{2}(x_1y_3-x_3y_1)+\frac{1}{12}(x_1-y_1)(x_1y_2-x_2y_1)$;
    \item $z_5=x_5+y_5+\frac{1}{2}(x_1y_4-x_4y_1)+\frac{1}{12}(x_1-y_1)(x_1y_3-x_3y_1)-\frac{1}{24}x_1y_1(x_1y_2-x_2y_1)$;
    \item $z_6=x_6+y_6+\frac{1}{2}(x_1y_5-x_5y_1+x_2y_3-x_3y_2)+\frac{1}{12}(x_1-y_1)(x_1y_4-x_4y_1)\salto+\frac{1}{12}(x_2-y_2)(x_1y_2-x_2y_1)-\frac{1}{24}x_1y_1(x_1y_3-x_3y_1)+\frac{1}{720}(y_1^3-x_1^3)(x_1y_2-x_2y_1)\salto+\frac{1}{180}(x_1y_1^2-x_1^2y_1)(x_1y_2-x_2y_1)\,.$
\end{itemize}

Since 
\begin{eqnarray*}
   \mathrm{d} (L_{\bf x})_{\bf 0}= \left[\begin{matrix} 
   1 & 0 & 0& 0 & 0 & 0\\
   0 & 1 & 0& 0& 0& 0\\
   -\frac{x_2}{2} & \frac{x_1}{2} &1& 0& 0& 0\\
  -\frac{x_1x_2}{12} -\frac{x_3}{2} & \frac{x_1^2}{12} &\frac{x_1}{2}& 1 & 0& 0\\
-\frac{x_1x_3}{12}-\frac{x_4}{2}& 0 &\frac{x_1^2}{12} & \frac{x_1}{2}& 1& 0\\ 
\frac{x_1^3x_2}{720}-\frac{x_1x_4+x_2^2}{12}-\frac{x_5}{2}& \frac{x_1x_2}{12}-\frac{x_3}{2}-\frac{x_1^4}{720}&\frac{x_2}{2} &\frac{x_1^2}{12} &\frac{x_1}{2}& 1
   \end{matrix}\right]\,,
\end{eqnarray*}
the induced left-invariant vector fields \eqref{leftinvariant vf} are: 
\begin{itemize}
\item $X_1={\partial}_{x_1}-\frac{x_2}{2}{\partial}_{x_3}-\big(\frac{x_3}{2}+\frac{x_1x_2}{12}\big){\partial}_{x_4}-\big(\frac{x_4}{2}+\frac{x_1x_3}{12}\big){\partial}_{x_5}+\big(\frac{x_1^3x_2}{720}-\frac{x_1x_4+x_2^2}{12}-\frac{x_5}{2}\big){\partial}_{x_6}\,;$
\item $X_2={\partial}_{x_2}+\frac{x_1}{2}{\partial}_{x_3}+\frac{x_1^2}{12}{\partial}_{x_4}+\big(\frac{x_1x_2}{12}-\frac{x_3}{2}-\frac{x_1^4}{720}\big){\partial}_{x_6}\,;$
\item $X_3={\partial}_{x_3}+\frac{x_1}{2}{\partial}_{x_4}+\frac{x_1^2}{12}{\partial}_{x_5}+\frac{x_2}{2}{\partial}_{x_6}\,;$
\item $X_4={\partial}_{x_4}+\frac{x_1}{2}{\partial}_{x_5}+\frac{x_1^2}{12}{\partial}_{x_6}\,;$
\item $X_5={\partial}_{x_5}+\frac{x_1}{2}{\partial}_{x_6}$;
\item $X_6={\partial}_{x_6},$
\end{itemize}
and the respective left-invariant 1-forms \eqref{leftinvariant form} are:
\begin{itemize}
\item $\theta_1=dx_1$;
\item $\theta_2=dx_2$;
\item $\theta_3=dx_3-\frac{x_1}{2} dx_2+\frac{x_2}{2}dx_1$;
\item $\theta_4=dx_4-\frac{x_1}{2}dx_3+\frac{x_1^2}{6}dx_2+\big(\frac{x_3}{2}-\frac{x_1x_2}{6}\big)dx_1$;
\item $\theta_5=dx_5-\frac{x_1}{2}dx_4+\frac{x_1^2}{6}dx_3-\frac{x_1^3}{24}dx_2+\big(\frac{x_4}{2}-\frac{x_1x_3}{6}+\frac{x_1^2x_2}{24})dx_1$;
\item $\theta_6=dx_6-\frac{x_1}{2}dx_5+\frac{x_1^2}{6}dx_4-\big(\frac{x_2}{2}+\frac{x_1^3}{24}\big)dx_3+\big(\frac{x_3}{2}+\frac{x_1x_2}{6}+\frac{x_1^4}{120}\big)dx_2+\big(\frac{x_5}{2}-\frac{x_2^2+x_1x_4}{6}$\\$  \phantom{\theta_5=}+\frac{x_1^2x_3}{24}-\frac{x_1^3x_2}{120}\big)dx_1\,.$
\end{itemize}

Finally, we have
\begin{eqnarray*}
   \mathrm{d} (R_{\bf x})_{\bf 0}= \left[\begin{matrix} 
   1 & 0 & 0& 0 & 0 & 0\\
   0 & 1 & 0& 0& 0& 0\\
   \frac{x_2}{2} & -\frac{x_1}{2} &1& 0& 0& 0\\
  \frac{x_3}{2}-\frac{x_1x_2}{12}  & \frac{x_1^2}{12} &-\frac{x_1}{2}& 1 & 0& 0\\
\frac{x_4}{2}-\frac{x_1x_3}{12}& 0 &\frac{x_1^2}{12} & -\frac{x_1}{2}& 1& 0\\ 
\frac{x_5}{2}-\frac{x_1x_4+x_2^2}{12}+\frac{x_1^3x_2}{720}& \frac{x_3}{2}+\frac{x_1x_2}{12}-\frac{x_1^4}{720}&-\frac{x_2}{2} &\frac{x_1^2}{12} &-\frac{x_1}{2}& 1
   \end{matrix}\right]\,.
\end{eqnarray*}

\subsubsection{Grading, polarizations of maximal dimension, and asymptotic cone} 

  

The Lie algebra $N_{6,1,3}$
   is not stratifiable, but it is gradable as
   \begin{eqnarray*}
   V_1 &= & \span \{  X_i \} \,,\\
   V_2 &= & 0\,,\\
   V_i &= & \span \{  X_{i-1} \} \,,\,i=3,\ldots,7\,.
   \end{eqnarray*}
   
 We claim that in this Lie algebra every two complementary subspaces to the  derived subalgebra as in \eqref{complementary_Delta} differ by an automorphism. Indeed, when $u_i^j$ varies, 
      the matrix
$$    \left[  \begin{array}{cccccc}
     1 &  0 &  0 &  0 & 0& 0 \\
 0 &  1 &  0 &  0 &  0& 0 \\
 u_1^3 &  u_2^3 &  1 &  0 &  0 & 0\\
 u_1^4 &  u_2^4 &  u_2^3 &  1 &  0& 0 \\
 u_1^5 &  u_2^5 & u_2^4 & u_2^3& 1& 0\\
 u_1^6 & u_2^6 & u_2^5-u_1^3& u_2^4 & u_2^3 & 1
      \end{array}\right]$$
is a Lie algebra automorphism and sends the complementary subspace $\span\{X_1,X_2\}$ to
      $\span\{X_1+u_1^3X_3+u_1^4X_4+u_1^5X_5+u_1^6X_6, X_2+u_2^3X_3+u_2^4X_4+u_2^5X_5+u_2^6X_6\}$, which is an arbitrary subspace as in \eqref{complementary_Delta}. In particular, every $\Delta$ as in \eqref{complementary_Delta} gives maximal dimension. The tangent cone of each of such polarizations has Lie algebra isomorphic to $N_{6,2,7}$, see page \pageref{N627}.

 The asymptotic cone 
 of the Lie group with Lie algebra
 $N_{6,1,3}$ has Lie algebra isomorphic to $N_{6,2,1}$, the first-type filiform algebra of step 5, see page \pageref{N621}.
    
  \subsection*{$\bf N_{6,1,4 }$ non-stratifiable }
  The following Lie algebra is denoted as $N_{6,1,4}$ by Gong in \cite{Gong_Thesis}, as $L_{6,11}$ by de Graaf in \cite{deGraaf_classification}, as $(11)$ by Del Barco in \cite{delBarco}, and as $\mathcal{G}_{6,12}$ by Magnin in \cite{magnin}.
  
  The non-trivial brackets are the following:
\begin{equation*}
\begin{aligned}
    &[X_1,X_2]=X_{3}\,,\, [X_1,X_3]=X_4\,,\\
    [X_1,X_4]&=X_6\,,\,[X_2,X_3]=X_6\,,\,[X_2,X_5]=X_6.
\end{aligned}
\end{equation*}
This is a nilpotent Lie algebra of rank 3 and step 4 that is positively gradable, yet not stratifiable. The Lie brackets can be pictured with the diagram:
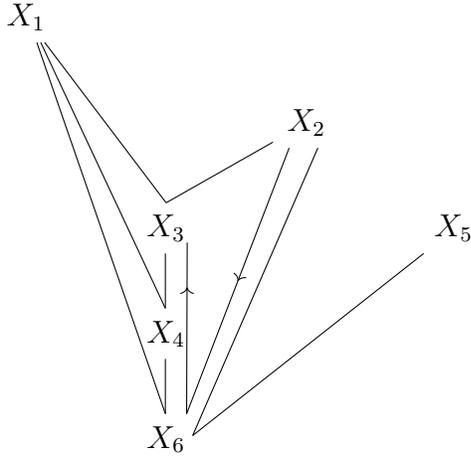
\begin{center}
 
\begin{tikzcd}[end anchor=north]
		X_1\ar[ddr,no head] 
		\ar[dddr,no head, end anchor={[xshift=-2.5ex]north east}]
		\ar[ddddr,no head, end anchor={[xshift=-2.5ex]north east}]
		& & &\\
		& & X_2\dlar[no head]\ar[->-=.5,no head,dddl,start anchor={[xshift=-3.7ex]south east }, end anchor={[xshift=1.5ex]north}]\ar[no head,dddl,end anchor={[xshift=-0.5ex]east},end anchor={[yshift=.3ex]east},start anchor={[xshift=-1.6ex]south east}] & \\
		& X_{3}
		\dar[no head, end anchor={[xshift=-2.5ex]north east},start anchor={[xshift=-2.5ex]south east }]
		\ar[-<-=.3,no head,dd, end anchor={[xshift=1.5ex]north},start anchor={[xshift=-0.9ex]south east },start anchor={[yshift=0.8ex]south east }]
		& & X_5\ar[no head,ddll, end anchor={[xshift=-0.5ex]east},end anchor={[yshift=.3ex]east}]\\
		&X_4\dar[no head, end anchor={[xshift=-2.5ex]north east},start anchor={[xshift=-2.5ex]south east }]& &\\
		&X_6&   &\phantom{X_5}\, .
	\end{tikzcd}

\end{center}

 
The composition law \eqref{group law in G} of $N_{6,1,4}$ is given by:

\begin{itemize}
    \item $z_1=x_1+y_1$;
    \item $z_2=x_2+y_2$;
    \item $z_3=x_3+y_3+\frac{1}{2}(x_1y_2-x_2y_1)$;
    \item $z_4=x_4+y_4+\frac{1}{2}(x_1y_3-x_3y_1)+\frac{1}{12}(x_1-y_1)(x_1y_2-x_2y_1)$;
    \item $z_5=x_5+y_5$;
    \item $z_6=x_6+y_6+\frac{1}{2}(x_1y_4-x_4y_1+x_2y_3-x_3y_2+x_2y_5-x_5y_2)+\frac{1}{12}(x_1-y_1)(x_1y_3-x_3y_1)\salto+\frac{1}{12}(x_2-y_2)(x_1y_2-x_2y_1)-\frac{1}{24}x_1y_1(x_1y_2-x_2y_1)$.
\end{itemize}


Since
\begin{eqnarray*}
  \mathrm{d}(L_\mathbf{x})_\mathbf{0}=
  \left[\begin{matrix} 
   1 & 0  & 0 & 0 & 0 & 0\\
   0  & 1 & 0 & 0 & 0 &0\\
   -\frac{x_2}{2} & \frac{x_1}{2} & 1 & 0 & 0 & 0\\
  -\frac{x_1x_2}{12}-\frac{x_3}{2} & \frac{x_1^2}{12} & \frac{x_1}{2} & 1& 0 & 0\\
  0 & 0 & 0 & 0& 1 & 0\\
  -\frac{x_1x_3+x_2^2}{12}-\frac{x_4}{2} & \frac{x_1x_2}{12}-\frac{x_3+x_5}{2} & \frac{x_2}{2}+\frac{x_1^2}{12} & \frac{x_1}{2}& \frac{x_2}{2} & 1
   \end{matrix}\right]\,,
\end{eqnarray*}
the induced left-invariant vector fields \eqref{leftinvariant vf} are: 
\begin{itemize}
\item $X_1={\partial}_{x_1}-\frac{x_2}{2}{\partial}_{x_3}-\big(\frac{x_3}{2}+\frac{x_1x_2}{12}\big){\partial}_{x_4}-\big(\frac{x_1x_3+x_2^2}{12}+\frac{x_4}{2}\big){\partial}_{x_6}\,;$
\item $X_2={\partial}_{x_2}+\frac{x_1}{2}{\partial}_{x_3}+\frac{x_1^2}{12}{\partial}_{x_4}+\big(\frac{x_1x_2}{12}-\frac{x_3+x_5}{2}\big){\partial}_{x_6}\,;$
\item $X_3={\partial}_{x_3}+\frac{x_1}{2}{\partial}_{x_4}+\big(\frac{x_2}{2}+\frac{x_1^2}{12}\big){\partial}_{x_6}\,;$
\item $X_4={\partial}_{x_4}+\frac{x_1}{2}{\partial}_{x_6}\,;$
\item $X_5={\partial}_{x_5}+\frac{x_2}{2}{\partial}_{x_6}$;
\item $X_6={\partial}_{x_6}\,,$
\end{itemize}
and the respective left-invariant 1-forms \eqref{leftinvariant form} are: 
\begin{itemize}
\item $\theta_1=dx_1$;
\item $\theta_2=dx_2$;
\item $\theta_3=dx_3-\frac{x_1}{2} dx_2+\frac{x_2}{2}dx_1$;
\item $\theta_4=dx_4-\frac{x_1}{2}dx_3+\frac{x_1^2}{6}dx_2+\big(\frac{x_3}{2}-\frac{x_1x_2}{6}\big)dx_1$;
\item $\theta_5=dx_5$;
\item $\theta_6=dx_6-\frac{x_2}{2}dx_5-\frac{x_1}{2}dx_4+\big(\frac{x_1^2}{6}-\frac{x_2}{2}\big)dx_3+\big(\frac{x_3+x_5}{2}+\frac{x_1x_2}{6}-\frac{x_1^3}{24}\big)dx_2+\big(\frac{x_4}{2}\saltot-\frac{x_2^2+x_1x_3}{6}+\frac{x_1^2x_2}{24}\big)dx_1\,.$
\end{itemize}

Finally, we have
\begin{eqnarray*}
  \mathrm{d}(R_\mathbf{x})_\mathbf{0}=
  \left[\begin{matrix} 
   1 & 0  & 0 & 0 & 0 & 0\\
   0  & 1 & 0 & 0 & 0 &0\\
   \frac{x_2}{2} & -\frac{x_1}{2} & 1 & 0 & 0 & 0\\
  \frac{x_3}{2}-\frac{x_1x_2}{12} & \frac{x_1^2}{12} & -\frac{x_1}{2} & 1& 0 & 0\\
  0 & 0 & 0 & 0& 1 & 0\\
  \frac{x_4}{2}-\frac{x_1x_3+x_2^2}{12} & \frac{x_3+x_5}{2}+\frac{x_1x_2}{12} & \frac{x_1^2}{12}-\frac{x_2}{2} & -\frac{x_1}{2}& -\frac{x_2}{2} & 1
   \end{matrix}\right]\,.
\end{eqnarray*}

\subsubsection{Grading, polarizations of maximal dimension, and asymptotic cone}
The Lie algebra $N_{6,1,4}$
   is not stratifiable, but it is gradable as
   \begin{eqnarray*}
   V_i &= & \span \{  X_i \} \,,\,  i=1,2,4\,,\\
   V_3 &= & \span\{X_3,X_5\} \,,\\
   V_5&=&\span\lbrace X_6\rbrace\,.
   \end{eqnarray*}




     Every   complementary subspace $\Delta$ to the  derived subalgebra is spanned by
      $X_1 +   u_1^3 X_3+ u_1^4 X_4+ u_1^6 X_6 $,
      $X_2 +   u_2^3 X_3+   u_2^4 X_4+   u_2^6 X_6$, and
       $X_5 +   u_5^3 X_3 +   u_5^4 X_4+   u_5^6 X_6$.
       Such a polarization gives maximal Hausdorff dimension if and only if 
       $u_5^3 $ is either $ -1$ or $ 0$.
       
       We claim that every two polarizations with $u_5^3 = 0$    differ by an automorphism. Likewise,  every two polarizations with $u_5^3 = -1$    differ by an automorphism.
        Indeed, when $u_i^j$ varies, 
      the matrix
$$    \left[  \begin{array}{cccccc}
     1 &  0 &  0 &  0 &  0 &  0  \\
 -u_5^4 &  1 &  0 &  0 &  0 &  0\\
 u_1^3 & u_2^3 & 1 &  0 & 0 &  0\\
 u_1^4 & u_2^4 & u_2^3 & 1 & u_5^4 & 0\\
 0 & 0 & 0 & 0 &  1 &  0\\
 u_1^6 & u_2^6 & u_2^4-u_1^3-u_5^4 u_2^3 &  u_2^3-u_5^4 & u_5^6 &  1\\
      \end{array}\right],$$
         is a Lie algebra automorphism and   sends the complementary subspace
        \begin{equation}
          \label{614_1}\span\{X_1, X_2,X_5\}
          \end{equation}
          to an arbitrary one with $u_5^3 = 0$.
      Instead, when $u_i^j$ varies, 
      the matrix
$$    \left[  \begin{array}{cccccc}
     1 &  0 &  0 &  0 &  0 &  0  \\
 -u_5^4 &  1 &  0 &  0 &  0 &  0\\
 u_1^3 & u_2^3 & 1 &  0 & -1 &  0\\
 u_1^4 & u_2^4 & u_2^3 & 1 & u_5^4 & -1\\
 0 & 0 & 0 & 0 &  1 &  0\\
 u_1^6 & u_2^6 & u_2^4-u_1^3-u_5^4 u_2^3 &  u_2^3-u_5^4 & u_5^6 &  u_5^4\\
      \end{array}\right],$$
         is a Lie algebra automorphism and sends the complementary subspace
      \begin{equation}
          \label{614_2}
  \span\{X_1, X_2,X_5-X_3\}    \end{equation}
  to an arbitrary one with $u_5^3 = -1$ .
          These two maximal-dimension polarizations \eqref{614_1} and \eqref{614_2} are not biLipschitz equivalent since they have different tangents:  $N_{6,3,3 }$ and $N_{6,2,6 }$ respectively, see pages \pageref{N633} and \pageref{N626}.

The asymptotic cone 
 of the Lie group with Lie algebra $N_{6,1,4}$ has Lie algebra isomorphic to
 $N_{5,2,1}\times\mathbb{R}$, where $N_{5,2,1}$ is the filiform algebra of step 4.
       


  \subsection*{$\bf N_{6,2,1}$ } 
  
   The following Lie algebra is denoted as $N_{6,2,1}$ by Gong in \cite{Gong_Thesis}, as $L_{6,18}$ by de Graaf in \cite{deGraaf_classification}, as (3) by Del Barco in \cite{delBarco}, and as $\mathcal{G}_{6,16}$ by Magnin in \cite{magnin}. 
   
   The non-trivial brackets are the following:
\begin{equation*}\label{N621}
   [X_1, X_i] = X_{i+1}\,,\, 2\leq i\leq 5\,.
\end{equation*}
This is a nilpotent Lie algebra of rank 2 and step 5 that is stratifiable, also known as the filiform Lie algebra of first type of dimension 6, the second type is $N_{6,2,2}$, see page \pageref{N622}. The Lie brackets can be pictured with the diagram:
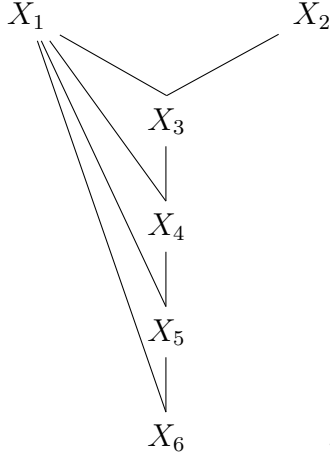
\begin{center}
 
	\begin{tikzcd}[end anchor=north]
		X_1\drar[no head] \ar[ddr,no head, end anchor={[xshift=-2.5ex]north east}]\ar[ddddr,no head, end anchor={[xshift=-2.5ex]north east}]\ar[dddr,no head, end anchor={[xshift=-2.5ex]north east}]& & X_2\dlar[no head]\\
		& X_{3}\dar[no head, end anchor={[xshift=-2.5ex]north east},start anchor={[xshift=-2.5ex]south east }]&\\
		&X_4\dar[no head, end anchor={[xshift=-2.5ex]north east},start anchor={[xshift=-2.5ex]south east }]&\\
		&X_5\dar[no head,end anchor={[xshift=-2.5ex]north east},start anchor={[xshift=-2.5ex]south east }]&  \\
		& X_6 &  \quad\;.
	\end{tikzcd}

\end{center}
 

 
The composition law \eqref{group law in G} of $N_{6,2,1}$ is given by:

\begin{itemize}
    \item $z_1=x_1+y_1$;
    \item $z_2=x_2+y_2$;
    \item $z_3=x_3+y_3+\frac{1}{2}(x_1y_2-x_2y_1)$;
    \item $z_4=x_4+y_4+\frac{1}{2}(x_1y_3-x_3y_1)+\frac{1}{12}(x_1-y_1)(x_1y_2-x_2y_1)$;
    \item $z_5=x_5+y_5+\frac{1}{2}(x_1y_4-x_4y_1)+\frac{1}{12}(x_1-y_1)(x_1y_3-x_3y_1)-\frac{1}{24}x_1y_1(x_1y_2-x_2y_1)$;
    \item $z_6=x_6+y_6+\frac{1}{2}(x_1y_5-x_5y_1)+\frac{1}{12}(x_1-y_1)(x_1y_4-x_4y_1)-\frac{1}{24}x_1y_1(x_1y_3-x_3y_1)\salto+\frac{1}{720}(y_1^3-x_1^3)(x_1y_2-x_2y_1)+\frac{1}{180}(x_1y_1^2-x_1^2y_1)(x_1y_2-x_2y_1)\,.$
\end{itemize}

Since
\begin{eqnarray*}
   \mathrm{d} (L_{\bf x})_{\bf 0}= \left[\begin{matrix} 
   1 & 0 & 0& 0 & 0 & 0\\
   0 & 1 & 0& 0& 0& 0\\
   -\frac{x_2}{2} & \frac{x_1}{2} &1& 0& 0& 0\\
  -\frac{x_1x_2}{12} -\frac{x_3}{2} & \frac{x_1^2}{12} &\frac{x_1}{2}& 1 & 0& 0\\
-\frac{x_1x_3}{12}-\frac{x_4}{2}& 0 &\frac{x_1^2}{12} & \frac{x_1}{2}& 1& 0\\ 
-\frac{x_1x_4}{12}-\frac{x_5}{2}+\frac{x_1^3x_2}{720}& -\frac{x_1^4}{720}&0 &\frac{x_1^2}{12} &\frac{x_1}{2}& 1
   \end{matrix}\right]\,,
\end{eqnarray*}
the induced left-invariant vector fields \eqref{leftinvariant vf} are: 
\begin{itemize}
\item $X_1={\partial}_{x_1}-\frac{x_2}{2}{\partial}_{x_3}-\big(\frac{x_3}{2}+\frac{x_1x_2}{12}\big){\partial}_{x_4}-\big(\frac{x_4}{2}+\frac{x_1x_3}{12}\big){\partial}_{x_5}+\big(\frac{x_1^3x_2}{720}-\frac{x_1x_4}{12}-\frac{x_5}{2}\big){\partial}_{x_6}\,;$
\item $X_2={\partial}_{x_2}+\frac{x_1}{2}{\partial}_{x_3}+\frac{x_1^2}{12}\partial_{ x_4}-\frac{x_1^4}{720}\partial_{ x_6}\,;$
\item $X_3=\partial_{ x_3}+\frac{x_1}{2}\partial_{ x_4}+\frac{x_1^2}{12}\partial_{ x_5}\,;$
\item $X_4=\partial_{ x_4}+\frac{x_1}{2}\partial_{ x_5}+\frac{x_1^2}{12}\partial_{ x_6}\,;$
\item $X_5=\partial_{ x_5}+\frac{x_1}{2}\partial_{ x_6}$;
\item $X_6=\partial_{ x_6}$,
\end{itemize}
and the respective left-invariant 1-forms \eqref{leftinvariant form} are: 
\begin{itemize}
\item $\theta_1=dx_1$;
\item $\theta_2=dx_2$;
\item $\theta_3=dx_3-\frac{x_1}{2} dx_2+\frac{x_2}{2}dx_1$;
\item $\theta_4=dx_4-\frac{x_1}{2}dx_3+\frac{x_1^2}{6}dx_2+\big(\frac{x_3}{2}-\frac{x_1x_2}{6}\big)dx_1$;
\item $\theta_5=dx_5-\frac{x_1}{2}dx_4+\frac{x_1^2}{6}dx_3-\frac{x_1^3}{24}dx_2+\big(\frac{x_4}{2}-\frac{x_1x_3}{6}+\frac{x_1^2x_2}{24})dx_1$;
\item $\theta_6=dx_6-\frac{x_1}{2}dx_5+\frac{x_1^2}{6}dx_4-\frac{x_1^3}{24}dx_3+\frac{x_1^4}{120}dx_2+\big(\frac{x_5}{2}-\frac{x_1x_4}{6}+\frac{x_1^2x_3}{24}-\frac{x_1^3x_2}{120}\big)dx_1$.
\end{itemize}

Finally, we have
\begin{eqnarray*}
   \mathrm{d} (R_{\bf x})_{\bf 0}= \left[\begin{matrix} 
   1 & 0 & 0& 0 & 0 & 0\\
   0 & 1 & 0& 0& 0& 0\\
   \frac{x_2}{2} & -\frac{x_1}{2} &1& 0& 0& 0\\
  \frac{x_3}{2}-\frac{x_1x_2}{12}  & \frac{x_1^2}{12} &-\frac{x_1}{2}& 1 & 0& 0\\
\frac{x_4}{2}-\frac{x_1x_3}{12}& 0 &\frac{x_1^2}{12} & -\frac{x_1}{2}& 1& 0\\ 
\frac{x_5}{2}-\frac{x_1x_4}{12}+\frac{x_1^3x_2}{720}& -\frac{x_1^4}{720}&0 &\frac{x_1^2}{12} &-\frac{x_1}{2}& 1
   \end{matrix}\right]\,.
\end{eqnarray*}
    \subsection*{$\bf N_{6,2,2}$  }
    The following Lie algebra is denoted as $N_{6,2,2}$ by Gong in \cite{Gong_Thesis}, as $L_{6,16}$ by de Graaf in \cite{deGraaf_classification}, as (2) by Del Barco  in \cite{delBarco}, and as $\mathcal{G}_{6,18}$ by Magnin in \cite{magnin}.
    
    The non-trivial brackets are the following:
\begin{equation*}\label{N622}
   [X_1, X_i] = X_{i+1}\,,\, 2\leq i\leq 4\,,\,[X_2,X_5]=X_6\,,\,[X_3,X_4]=-X_6\,.
\end{equation*}
This is a nilpotent Lie algebra of rank 2 and step 5 that is stratifiable, also known as the filiform Lie algebra of second type of dimension 6, the first type is $N_{6,2,1}$, see page \pageref{N621}. The Lie brackets can be pictured with the diagram:
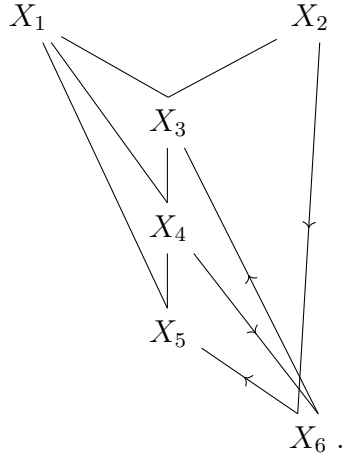
\begin{center}
 
	\begin{tikzcd}[end anchor=north]
		X_1\drar[no head] \ar[ddr,no head, end anchor={[xshift=-2.5ex]north east}]\ar[dddr,no head, end anchor={[xshift=-2.5ex]north east}]& & X_2\;\,\dlar[no head]\ar[dddd, ->-=.5, no head, end anchor={[xshift=-4.5ex]north east}] \\
		& X_{3}\dar[no head, end anchor={[xshift=-2.5ex]north east},start anchor={[xshift=-2.5ex]south east }]\ar[dddr, -<-=.5, no head, end anchor={[xshift=-3ex]north east}]& \\
		&X_4\dar[no head, end anchor={[xshift=-2.5ex]north east},start anchor={[xshift=-2.5ex]south east }]\ar[ddr, no head, ->-=.5,end anchor={[xshift=-3ex]north east}]& \\
		&X_5\drar[-<-=.5,no head,end anchor={[xshift=-4.5ex]north east}]&  \\
		&  &  X_6\; .
	\end{tikzcd}

\end{center}


 The composition law \eqref{group law in G} of $N_{6,2,2}$ is given by:

\begin{itemize}
    \item $z_1=x_1+y_1$;
    \item $z_2=x_2+y_2$;
    \item $z_3=x_3+y_3+\frac{1}{2}(x_1y_2-x_2y_1)$;
    \item $z_4=x_4+y_4+\frac{1}{2}(x_1y_3-x_3y_1)+\frac{1}{12}(x_1-y_1)(x_1y_2-x_2y_1)$;
    \item $z_5=x_5+y_5+\frac{1}{2}(x_1y_4-x_4y_1)+\frac{1}{12}(x_1-y_1)(x_1y_3-x_3y_1)-\frac{1}{24}x_1y_1(x_1y_2-x_2y_1)$;
    \item $z_6=x_6+y_6+\frac{1}{2}(x_2y_5-x_5y_2+x_4y_3-x_3y_4)+\frac{1}{12}\big(x_2-y_2)(x_1y_4-x_4y_1)\salto+\frac{1}{12}(y_3-x_3)(x_1y_3-x_3y_1)+\frac{1}{12}(x_4-y_4)(x_1y_2-x_2y_1)-\frac{1}{24}x_1y_2(x_1y_3-x_3y_1)\salto+\frac{1}{24}x_1y_3(x_1y_2-x_2y_1)+\frac{1}{720}(y_1^2y_2-x_1^2x_2)(x_1y_2-x_2y_1)-\frac{1}{360}y_1(x_1y_2-x_2y_1)^2\salto+\frac{1}{180}(x_1y_1y_
    2-x_1^2y_2)(x_1y_2-x_2y_1)-\frac{1}{120}x_1(x_1y_2-x_2y_1)^2$.
\end{itemize}

Since
\begin{eqnarray*}
   \mathrm{d} (L_{\bf x})_{\bf 0}= \left[\begin{matrix} 
   1 & 0 & 0& 0 & 0 & 0\\
   0 & 1 & 0& 0& 0& 0\\
   -\frac{x_2}{2} & \frac{x_1}{2} &1& 0& 0& 0\\
  -\frac{x_1x_2}{12} -\frac{x_3}{2} & \frac{x_1^2}{12} &\frac{x_1}{2}& 1 & 0& 0\\
-\frac{x_1x_3}{12}-\frac{x_4}{2}& 0 &\frac{x_1^2}{12} & \frac{x_1}{2}& 1& 0\\ 
 \frac{x_3^2-2x_2x_4}{12}+\frac{x_1^2x_2^2}{720}& \frac{x_1x_4}{12}-\frac{x_5}{2}-\frac{x_1^3x_2}{720}&\frac{x_4}{2}-\frac{x_1x_3}{12}& \frac{x_1x_2}{12}-\frac{x_3}{2} &\frac{x_2}{2}& 1
   \end{matrix}\right]\,,
\end{eqnarray*}
the induced left-invariant vector fields \eqref{leftinvariant vf} are:
\begin{itemize}
\item $X_1=\partial_{ x_1}-\frac{x_2}{2}\partial_{ x_3}-\big(\frac{x_3}{2}+\frac{x_1x_2}{12}\big)\partial_{ x_4}-\big(\frac{x_4}{2}+\frac{x_1x_3}{12}\big)\partial_{ x_5}+\big(\frac{x_3^2-2x_2x_4}{12}+\frac{x_1^2x_2^2}{720})\partial_{ x_6}\,;$
\item $X_2=\partial_{ x_2}+\frac{x_1}{2}\partial_{ x_3}+\frac{x_1^2}{12}\partial_{ x_4}+\big(\frac{x_1x_4}{12}-\frac{x_5}{2}-\frac{x_1^3x_2}{720}\big)\partial_{ x_6}\,;$
\item $X_3=\partial_{ x_3}+\frac{x_1}{2}\partial_{ x_4}+\frac{x_1^2}{12}\partial_{ x_5}+\big(\frac{x_4}{2}-\frac{x_1x_3}{12}\big)\partial_{ x_6}\,;$
\item $X_4=\partial_{ x_4}+\frac{x_1}{2}\partial_{ x_5}+\big(\frac{x_1x_2}{12}-\frac{x_3}{2}\big)\partial_{ x_6}\,;$
\item $X_5=\partial_{ x_5}+\frac{x_2}{2}\partial_{ x_6}$;
\item $X_6=\partial_{ x_6}$,
\end{itemize}
and the respective left-invariant 1-forms \eqref{leftinvariant vf} are: 
\begin{itemize}
\item $\theta_1=dx_1$;
\item $\theta_2=dx_2$;
\item $\theta_3=dx_3-\frac{x_1}{2} dx_2+\frac{x_2}{2}dx_1$;
\item $\theta_4=dx_4-\frac{x_1}{2}dx_3+\frac{x_1^2}{6}dx_2+\big(\frac{x_3}{2}-\frac{x_1x_2}{6}\big)dx_1$;
\item $\theta_5=dx_5-\frac{x_1}{2}dx_4+\frac{x_1^2}{6}dx_3-\frac{x_1^3}{24}dx_2+\big(\frac{x_4}{2}-\frac{x_1x_3}{6}+\frac{x_1^2x_2}{24})dx_1$;
\item $\theta_6=dx_6-\frac{x_2}{2}dx_5+\big(\frac{x_3}{2}+\frac{x_1x_2}{6}\big)dx_4-\big(\frac{x_4}{2}+\frac{x_1x_3}{6}+\frac{x_1^2x_2}{24}\big)dx_3+\big(\frac{x_5}{2}+\frac{x_1x_4}{6}+\frac{x_1^2x_3}{24}$\\$  \phantom{\theta_5=}+\frac{x_1^3x_2}{120}\big)dx_2+\big(\frac{x_3^2-2x_2x_4}{6}-\frac{x_1^2x_2^2}{120}\big)dx_1$.
\end{itemize}

Finally, we have
\begin{eqnarray*}
   \mathrm{d} (R_{\bf x})_{\bf 0}= \left[\begin{matrix} 
   1 & 0 & 0& 0 & 0 & 0\\
   0 & 1 & 0& 0& 0& 0\\
   \frac{x_2}{2} & -\frac{x_1}{2} &1& 0& 0& 0\\
  \frac{x_3}{2}-\frac{x_1x_2}{12}  & \frac{x_1^2}{12} &-\frac{x_1}{2}& 1 & 0& 0\\
\frac{x_4}{2}-\frac{x_1x_3}{12}& 0 &\frac{x_1^2}{12} & -\frac{x_1}{2}& 1& 0\\ 
 \frac{x_3^2-2x_2x_4}{12}+\frac{x_1^2x_2^2}{720}& \frac{x_5}{2}+\frac{x_1x_4}{12}-\frac{x_1^3x_2}{720}&-\frac{x_4}{2}-\frac{x_1x_3}{12}& \frac{x_3}{2}+\frac{x_1x_2}{12} &-\frac{x_2}{2}& 1
   \end{matrix}\right]\,.
\end{eqnarray*}
     \subsection*{$\bf N_{6,2,3}$ non-stratifiable }  The following Lie algebra is denoted as $N_{6,2,3}$ by Gong in \cite{Gong_Thesis}, as $L_{6,13}$ by de Graaf in \cite{deGraaf_classification}, as $(9)$ by Del Barco in \cite{delBarco}, and as $\mathcal{G}_{6,13}$ by Magnin in \cite{magnin}. 
  
  The non-trivial brackets are the following:
\begin{equation*}
   [X_1, X_2] = X_4\,, \,[X_1, X_i] = X_{i+1}\,,\, i = 4, 5\,,\, [X_2, X_3] = X_5\,,\, [X_3, X_4] = -X_6\,.
\end{equation*}
This is a nilpotent Lie algebra of rank 3 and step 4 that is positively gradable, yet not stratifiable. The Lie brackets can be pictured with the diagram:
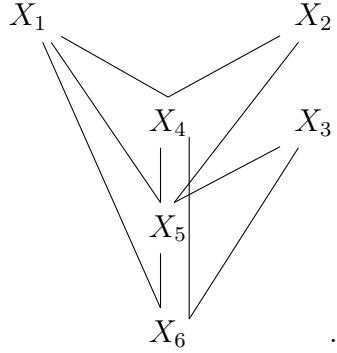
\begin{center}
 
	\begin{tikzcd}[end anchor=north]
	X_1 \ar[dr, no head]\ar[ddr, no head,end anchor={[xshift=-3.ex]north east}]\ar[dddr, no head,end anchor={[xshift=-3.ex]north east}]& & X_2\ar[no head, dl]\ar[ddl, no head,end anchor={[xshift=-2ex]north east}]   \\
	& X_4 \ar[d, no head,end anchor={[xshift=-3.ex]north east},start anchor={[xshift=-3.ex]south east}]\ar[no head,dd, end anchor={[xshift=-0.9ex]north east},start anchor={[xshift=-0.9ex]south east }, end anchor={[yshift=-0.8ex]north east},start anchor={[yshift=0.8ex]south east }]&  X_3 \ar[dl, no head,  end anchor={[xshift=-2ex]north east}]\ar[no head,ddl, end anchor={[xshift=-0.9ex]north east}, end anchor={[yshift=-0.8ex]north east}]\\
	& X_5\ar[d, no head,end anchor={[xshift=-3.ex]north east},start anchor={[xshift=-3ex]south east}] &  \\
	& X_6 & \quad\; .
	\end{tikzcd}

\end{center}

 
The composition law \eqref{group law in G} of $N_{6,2,3}$ is given by:

\begin{itemize}
    \item $z_1=x_1+y_1$;
    \item $z_2=x_2+y_2$;
    \item $z_3=x_3+y_3$;
    \item $z_4=x_4+y_4+\frac{1}{2}(x_1y_2-x_2y_1)$;
    \item $z_5=x_5+y_5+\frac{1}{2}(x_1y_4-x_4y_1+x_2y_3-x_3y_2)+\frac{1}{12}(x_1-y_1)(x_1y_2-x_2y_1)$;
    \item $z_6=x_6+y_6+\frac{1}{2}(x_1y_5-x_5y_1+x_4y_3-x_3y_4)+\frac{1}{12}(x_1-y_1)(x_1y_4-x_4y_1)\salto +\frac{1}{12}(x_1-y_1)(x_2y_3-x_3y_2)+\frac{1}{12}(y_3-x_3)(x_1y_2-x_2y_1)-\frac{1}{24}x_1y_1(x_1y_2-x_2y_1)$.
\end{itemize}




Since
\begin{eqnarray*}
 \mathrm{d}(L_\mathbf{x})_\mathbf{0}= \left[\begin{matrix} 
   1 & 0  & 0 & 0 & 0 & 0\\
   0  & 1 & 0 & 0 & 0 &0\\
  0 & 0 & 1 & 0 & 0 & 0\\
  -\frac{x_2}{2} & \frac{x_1}{2} & 0 & 1& 0 & 0\\
  -\frac{x_4}{2}-\frac{x_1x_2}{12} & \frac{x_1^2}{12}-\frac{x_3}{2} & \frac{x_2}{2} & \frac{x_1}{2}& 1 & 0\\
  \frac{x_2x_3-x_1x_4}{12}-\frac{x_5}{2} & -\frac{x_1x_3}{6} & \frac{x_4}{2}+\frac{x_1x_2}{12} & \frac{x_1^2}{12}-\frac{x_3}{2}& \frac{x_1}{2} & 1
   \end{matrix}\right]\,,
\end{eqnarray*}
the induced left-invariant vector fields \eqref{leftinvariant vf} are: 
\begin{itemize}
\item $X_1=\partial_{x_1}-\frac{x_2}{2}\partial_{x_4}-\big(\frac{x_4}{2}+\frac{x_1x_2}{12}\big)\partial_{x_5}+\big(\frac{x_2x_3-x_1x_4}{12}-\frac{x_5}{2}\big)\partial_{x_6}\,;$
\item $X_2=\partial_{x_2}+\frac{x_1}{2}\partial_{x_4}+\big(\frac{x_1^2}{12}-\frac{x_3}{2}\big)\partial_{x_5}-\frac{x_1x_3}{6}\partial_{x_6}\,;$
\item $X_3=\partial_{x_3}+\frac{x_2}{2}\partial_{x_5}+\big(\frac{x_4}{2}+\frac{x_1x_2}{12}\big)\partial_{x_6}\,;$
\item $X_4=\partial_{x_4}+\frac{x_1}{2}\partial_{x_5}+\big(\frac{x_1^2}{12}-\frac{x_3}{2}\big)\partial_{x_6}\,;$
\item $X_5=\partial_{x_5}+\frac{x_1}{2}\partial_{x_6}$;
\item $X_6=\partial_{x_6}$,
\end{itemize}
and the respective left-invariant 1-forms \eqref{leftinvariant form} are: 
\begin{itemize}
\item $\theta_1=dx_1$;
\item $\theta_2=dx_2$;
\item $\theta_3=dx_3$;
\item $\theta_4=dx_4-\frac{x_1}{2}dx_2+\frac{x_2}{2}dx_1$;
\item $\theta_5=dx_5-\frac{x_1}{2}dx_4-\frac{x_2}{2}dx_3+\big(\frac{x_3}{2}+\frac{x_1^2}{6}\big)dx_2+\big(\frac{x_4}{2}-\frac{x_1x_2}{6}\big)dx_1$;
\item $\theta_6=dx_6-\frac{x_1}{2}dx_5+\big(\frac{x_3}{2}+\frac{x_1^2}{6}\big)dx_4+\big(\frac{x_1x_2}{6}-\frac{x_4}{2}\big)dx_3-\big(\frac{x_1^3}{24}+\frac{x_1x_3}{3}\big)dx_2+\big(\frac{x_5}{2}$\\$  \phantom{z_5=}+\frac{x_2x_3-x_1x_4}{6}+\frac{x_1^2x_2}{24}\big)dx_1$.
\end{itemize}

Finally, we have
\begin{eqnarray*}
 \mathrm{d}(R_\mathbf{x})_\mathbf{0}= \left[\begin{matrix} 
   1 & 0  & 0 & 0 & 0 & 0\\
   0  & 1 & 0 & 0 & 0 &0\\
  0 & 0 & 1 & 0 & 0 & 0\\
  \frac{x_2}{2} & -\frac{x_1}{2} & 0 & 1& 0 & 0\\
  \frac{x_4}{2}-\frac{x_1x_2}{12} & \frac{x_1^2}{12}+\frac{x_3}{2} & -\frac{x_2}{2} & -\frac{x_1}{2}& 1 & 0\\
  \frac{x_2x_3-x_1x_4}{12}+\frac{x_5}{2} & -\frac{x_1x_3}{6} & \frac{x_1x_2}{12}-\frac{x_4}{2} & \frac{x_1^2}{12}+\frac{x_3}{2}& -\frac{x_1}{2} & 1
   \end{matrix}\right]\,.
\end{eqnarray*}

\subsubsection{Grading, polarizations of maximal dimension, and asymptotic cone}
The Lie algebra $N_{6,2,3}$
   is not stratifiable, but it is gradable as
   \begin{eqnarray*}
   V_1 &= & \span \{  X_1,X_2 \}\, ,\\
   V_2 &= & \span\{X_3,X_4\} \,,\\
   V_3 &= & \span \{  X_5 \}\,,\\
   V_4&=&\span\{X_6\}\,.
   \end{eqnarray*}

      Every   complementary subspace $\Delta$ to the  derived subalgebra is spanned by
      $X_1 +u_1^4X_4+u_1^5X_5+u_1^6X_6 $,
      $X_2 + u_2^4X_4+u_2^5X_5+u_2^6X_6 $, and
       $X_3 + u_3^4X_4+u_3^5X_5+u_3^6X_6 $.
       Such a polarization gives maximal Hausdorff dimension if and only if $u_3^5+u_1^4=0$.
      We claim that every two polarizations giving maximal Hausdorff dimension    differ by an automorphism. Indeed, when $u_i^j$ varies, 
      the matrix
$$    \left[  \begin{array}{cccccc}
     1 &  0 &  0 &  0 & 0& 0 \\
 -u_3^4 &  1 &  0 &  0 &  0& 0 \\
 0 & 0 &  1 &  0 &  0 & 0\\
 u_1^4 &  u_2^4 &  u_3^4 &  1 &  0& 0 \\
 u_1^5 &  u_2^5 & -u_1^4 & u_2^4& 1& 0\\
 u_1^6 & u_2^6 & u_3^6&  u_2^5 & u_2^4 & 1
      \end{array}\right]$$
      is a Lie algebra automorphism and sends the complementary subspace $\span\lbrace X_1,X_2,X_3\rbrace$ to
      $\span\{X_1-u_3^4X_2+u_1^4X_4+u_1^5X_5+u_1^6X_6, X_2+u_2^4X_4+u_2^5X_5+u_2^6X_6,X_3+u_3^4X_4-u_1^4X_5+u_3^6X_6\}$, which is an arbitrary one of maximal dimension.  
The tangent cone of each of such polarizations has Lie algebra isomorphic to $N_{6,3,1}$, see page \pageref{N631}.

 The asymptotic cone 
 of the Lie group with Lie algebra
 $N_{6,2,3}$ has Lie algebra isomorphic to $N_{5,2,1}\times\mathbb{R}$, where $N_{5,2,1}$ is the filiform algebra of step 4.


  \subsection*{$\bf N_{6,2,4}$ non-stratifiable} The following Lie algebra is denoted as $N_{6,2,4}$ by Gong in \cite{Gong_Thesis}, as $L_{6,12}$ by de Graaf in \cite{deGraaf_classification}, as $(10)$ by Del Barco in \cite{delBarco}, and as $\mathcal{G}_{6,11}$ by Magnin in \cite{magnin}.
  
  The non-trivial brackets are the following:
\begin{equation*}
   [X_1, X_2] = X_3\,, \,[X_1, X_3] = X_4\,,\, [X_1, X_4] = X_6\,,\, [X_2, X_5] = X_6\,.
\end{equation*}
This is a nilpotent Lie algebra of rank 3 and step 4 that is positively gradable, yet not stratifiable. The Lie brackets can be pictured with the diagram:
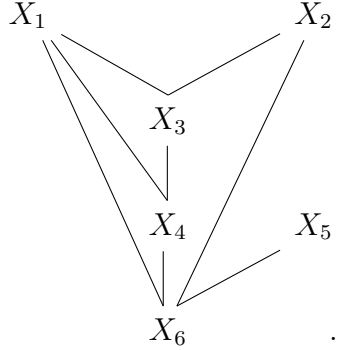
\begin{center}
 
\begin{tikzcd}[end anchor=north]
		X_1\drar[no head] \ar[ddr,no head, end anchor={[xshift=-2.5ex]north east}]\ar[dddr,no head, end anchor={[xshift=-2.8ex]north east}]& & X_2\dlar[no head]\ar[dddl, no head, end anchor={[xshift=-1.8ex]north east}] \\
		& X_{3}\dar[no head, end anchor={[xshift=-2.5ex]north east},start anchor={[xshift=-2.5ex]south east }]& \\
		&X_4\dar[no head, end anchor={[xshift=-2.8ex]north east},start anchor={[xshift=-2.8ex]south east }]&X_5\ar[dl, no head, end anchor={[xshift=-1.8ex]north east}]  \\
		&X_6&  \quad\;.
	\end{tikzcd}

\end{center}

 
The composition law \eqref{group law in G} of $N_{6,2,4}$ is given by:

\begin{itemize}
    \item $z_1=x_1+y_1$;
    \item $z_2=x_2+y_2$;
    \item $z_3=x_3+y_3+\frac{1}{2}(x_1y_2-x_2y_1)$;
    \item $z_4=x_4+y_4+\frac{1}{2}(x_1y_3-x_3y_1)+\frac{1}{12}(x_1-y_1)(x_1y_2-x_2y_1)$;
    \item $z_5=x_5+y_5$;
    \item $z_6=x_6+y_6+\frac{1}{2}(x_1y_4-x_4y_1+x_2y_5-x_5y_2)+\frac{1}{12}(x_1-y_1)(x_1y_3-x_3y_1)\salto-\frac{1}{24}x_1y_1(x_1y_2-x_2y_1)$.
\end{itemize}


Since
\begin{eqnarray*}
\mathrm{d}(L_\mathbf{x})_\mathbf{0}= \left[\begin{matrix} 
   1 & 0  & 0 & 0 & 0 & 0\\
   0  & 1 & 0 & 0 & 0 &0\\
   -\frac{x_2}{2} & \frac{x_1}{2} & 1 & 0 & 0 & 0\\
  -\frac{x_1x_2}{12}-\frac{x_3}{2} & \frac{x_1^2}{12} & \frac{x_1}{2} & 1& 0 & 0\\
  0 & 0 & 0 & 0& 1 & 0\\
  -\frac{x_1x_3}{12}-\frac{x_4}{2} & -\frac{x_5}{2} & \frac{x_1^2}{12} & \frac{x_1}{2}& \frac{x_2}{2} & 1
   \end{matrix}\right]\,,
\end{eqnarray*}
the induced left-invariant vector fields \eqref{leftinvariant vf} are: 
\begin{itemize}
\item $X_1=\partial_{ x_1}-\frac{x_2}{2}\partial_{ x_3}-\big(\frac{x_3}{2}+\frac{x_1x_2}{12}\big)\partial_{ x_4}-\big(\frac{x_1x_3}{12}+\frac{x_4}{2}\big)\partial_{ x_6}\,;$
\item $X_2=\partial_{ x_2}+\frac{x_1}{2}\partial_{ x_3}+\frac{x_1^2}{12}\partial_{ x_4}-\frac{x_5}{2}\partial_{ x_6}\,;$
\item $X_3=\partial_{ x_3}+\frac{x_1}{2}\partial_{ x_4}+\frac{x_1^2}{12}\partial_{ x_6}\,;$
\item $X_4=\partial_{ x_4}+\frac{x_1}{2}\partial_{ x_6}\,;$
\item $X_5=\partial_{ x_5}+\frac{x_2}{2}\partial_{ x_6}$;
\item $X_6=\partial_{ x_6}$,
\end{itemize}
and the respective left-invariant 1-forms \eqref{leftinvariant form} are: 
\begin{itemize}
\item $\theta_1=dx_1$;
\item $\theta_2=dx_2$;
\item $\theta_3=dx_3-\frac{x_1}{2} dx_2+\frac{x_2}{2}dx_1$;
\item $\theta_4=dx_4-\frac{x_1}{2}dx_3+\frac{x_1^2}{6}dx_2+\big(\frac{x_3}{2}-\frac{x_1x_2}{6}\big)dx_1$;
\item $\theta_5=dx_5$;
\item $\theta_6=dx_6-\frac{x_2}{2}dx_5-\frac{x_1}{2}dx_4+\frac{x_1^2}{6}dx_3+\big(\frac{x_5}{2}-\frac{x_1^3}{24}\big)dx_2+\big(\frac{x_4}{2}-\frac{x_1x_3}{6}+\frac{x_1^2x_2}{24}\big)dx_1$.
\end{itemize}

Finally, we have
 \begin{eqnarray*}
\mathrm{d}(R_\mathbf{x})_\mathbf{0}= \left[\begin{matrix} 
   1 & 0  & 0 & 0 & 0 & 0\\
   0  & 1 & 0 & 0 & 0 &0\\
   \frac{x_2}{2} & -\frac{x_1}{2} & 1 & 0 & 0 & 0\\
 \frac{x_3}{2} -\frac{x_1x_2}{12} & \frac{x_1^2}{12} & -\frac{x_1}{2} & 1& 0 & 0\\
  0 & 0 & 0 & 0& 1 & 0\\
 \frac{x_4}{2} -\frac{x_1x_3}{12} & \frac{x_5}{2} & \frac{x_1^2}{12} & -\frac{x_1}{2}& -\frac{x_2}{2} & 1
   \end{matrix}\right]\,.
\end{eqnarray*}

\subsubsection{Grading, polarizations of maximal dimension, and asymptotic cone}
The Lie algebra $N_{6,2,4}$
   is not stratifiable, but it is gradable as
   \begin{eqnarray*}
   V_1 &= & \span \{  X_1,X_2 \} \,,\\
   V_2 &= & \span\{X_3\} \,,\\
   V_3 &= & \span \{  X_4,X_5 \}\,,\\
   V_4&=&\span\{X_6\}\,.
   \end{eqnarray*}
     
      Every   complementary subspace $\Delta$ to the  derived subalgebra is spanned by
      $X_1 +u_1^3X_3+u_1^4X_4+u_1^6X_6 $,
      $X_2 +u_2^3X_3+ u_2^4X_4+u_2^6X_6 $, and
       $X_5 +u_5^3X_3+ u_5^4X_4+u_5^6X_6 $.
       Such a polarization gives maximal Hausdorff dimension if and only if $u^3_5=0$.
      We claim that every two polarizations giving maximal Hausdorff dimension    differ by an automorphism. Indeed, when $u_i^j$ varies, 
      the matrix
$$    \left[  \begin{array}{cccccc}
     1 &  0 &  0 &  0 & 0& 0 \\
 0 &  1 &  0 &  0 &  0& 0 \\
 u_1^3 & u_2^3 &  1 &  0 &  0 & 0\\
 u_1^4 &  u_2^4 &  u_2^3 &  1 &  u_5^4& 0 \\
 0 &  0 & 0 & 0& 1& 0\\
 u_1^6 & u_2^6 & u_2^4&  u_2^3 & u_5^6 & 1
      \end{array}\right]$$
      is a Lie algebra automorphism and sends the complementary subspace $\span\{X_1,X_2,X_5\}$ to
      $\span\{X_1+u_1^3X_3+u_1^4X_4+u_1^6X_6,X_2+u_2^3X_3+u_2^4X_4+u_2^6X_6, X_5+u_5^4X_4+u_5^6X_6\}$, which is an arbitrary one of maximal dimension.  
The tangent cone of each of such polarizations has Lie algebra isomorphic to $N_{6,3,3}$, see page \pageref{N633}.

 The asymptotic cone 
 of the Lie group with Lie algebra
 $N_{6,2,4}$ has Lie algebra isomorphic to $N_{5,2,1}\times\mathbb{R}$, where $N_{5,2,1}$ is the filiform algebra of step 4.

       
    \subsection*{$\bf N_{6,2,5}$  }
  
   The following Lie algebra is denoted as $N_{6,2,5}$ by Gong in \cite{Gong_Thesis}, as $L_{6,21(-1)}$ by de Graaf in \cite{deGraaf_classification}, as $(7)$ by Del Barco in \cite{delBarco}, and as $\mathcal{G}_{6,15}$ by Magnin \cite{magnin}.
   
   The non-trivial brackets are the following:
\begin{equation*}\label{N625}
   [X_1, X_i] = X_{i+1}\,,\, i = 2, 3, 5\,,\, [X_2, X_j] = X_{j+2}\,,\, j = 3, 4\,.
\end{equation*}
This is a nilpotent Lie algebra of rank 2 and step 4 that is stratifiable. The Lie brackets can be pictured with the diagram:
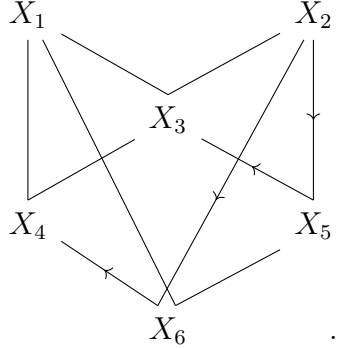
\begin{center}
 
	\begin{tikzcd}[end anchor=north]
	X_1\ar[dr, no head]\ar[dd, no head]\ar[dddr, no head,end anchor={[xshift=-1.9ex]north east}]& & X_2\ar[dl, no head]\ar[dddl, no head, ->-=.6,end anchor={[xshift=-3.2ex]north east}]\ar[dd,->-=.5, no head]\\
	& X_3\ar[dl, no head]\ar[dr, -<-=.5, no head]&\\
	X_4\ar[dr, no head, -<-=.5, end anchor={[xshift=-3.2ex]north east}]& & X_5\ar[no head, dl,  end anchor={[xshift=-1.9ex]north east}]\\
	& X_6 & \quad\;.
	\end{tikzcd}

\end{center}

 
The composition law \eqref{group law in G} of $N_{6,2,5}$ is given by:

\begin{itemize}
    \item $z_1=x_1+y_1$;
    \item $z_2=x_2+y_2$;
    \item $z_3=x_3+y_3+\frac{1}{2}(x_1y_2-x_2y_1)$;
    \item $z_4=x_4+y_4+\frac{1}{2}(x_1y_3-x_3y_1)+\frac{1}{12}(x_1-y_1)(x_1y_2-x_2y_1)$;
    \item $z_5=x_5+y_5+\frac{1}{2}(x_2y_3-x_3y_2)+\frac{1}{12}(x_2-y_2)(x_1y_2-x_2y_1)$;
    \item $z_6=x_6+y_6+\frac{1}{2}(x_1y_5-x_5y_1+x_2y_4-x_4y_2)+\frac{1}{12}(x_1-y_1)(x_2y_3-x_3y_2)\salto +\frac{1}{12}(x_2-y_2)(x_1y_3-x_3y_1)-\frac{1}{24}(x_1y_2+x_2y_1)(x_1y_2-x_2y_1)$.
\end{itemize}


Since
\begin{eqnarray*}
 \mathrm{d}(L_\mathbf{x})_\mathbf{0}=\left[\begin{matrix} 
   1 & 0  & 0 & 0 & 0 & 0\\
   0  & 1 & 0 & 0 & 0 &0\\
   -\frac{x_2}{2} & \frac{x_1}{2} & 1 & 0 & 0 & 0\\
  -\frac{x_1x_2}{12}-\frac{x_3}{2} & \frac{x_1^2}{12} & \frac{x_1}{2} & 1& 0 & 0\\
  -\frac{x_2^2}{12}& \frac{x_1x_2}{12}-\frac{x_3}{2} & \frac{x_2}{2} & 0& 1 & 0\\
  -\frac{x_2x_3}{12}-\frac{x_5}{2} & -\frac{x_4}{2}-\frac{x_1x_3}{12} & \frac{x_1x_2}{6} & \frac{x_2}{2}& \frac{x_1}{2} & 1
   \end{matrix}\right]\,,
\end{eqnarray*}
the induced left-invariant vector fields \eqref{leftinvariant vf} are: 
\begin{itemize}
\item $X_1=\partial_{x_1}- \frac{x_2}{2}\partial_{x_3}-\big(\frac{x_3}{2}+\frac{x_1x_2}{12}\big)\partial_{x_4}-\frac{x_2^2}{12}\partial_{x_5}-\big(\frac{x_2x_3}{12}+\frac{x_5}{2}\big)\partial_{x_6}\,;$
\item $X_2=\partial_{x_2}+\frac{x_1}{2}\partial_{x_3}+\frac{x_1^2}{12}\partial_{x_4}+\big(\frac{x_1x_2}{12}-\frac{x_3}{2}\big)\partial_{x_5}-\big(\frac{x_4}{2}+\frac{x_1x_3}{12}\big)\partial_{x_6}\,;$
\item $X_3=\partial_{x_3}+\frac{x_1}{2}\partial_{x_4}+\frac{x_2}{2}\partial_{x_5}+\frac{x_1x_2}{6}\partial_{x_6}\,;$
\item $X_4=\partial_{x_4}+\frac{x_2}{2}\partial_{x_6}\,;$
\item $X_5=\partial_{x_5}+\frac{x_1}{2}\partial_{x_6}$;
\item $X_6=\partial_{x_6}$,
\end{itemize}
and the respective left-invariant 1-forms \eqref{leftinvariant form} are: 
\begin{itemize}
\item $\theta_1=dx_1$;
\item $\theta_2=dx_2$;
\item $\theta_3=dx_3-\frac{x_1}{2} dx_2+\frac{x_2}{2}dx_1$;
\item $\theta_4=dx_4-\frac{x_1}{2}dx_3+\frac{x_1^2}{6}dx_2+\big(\frac{x_3}{2}-\frac{x_1x_2}{6}\big)dx_1$;
\item $\theta_5=dx_5-\frac{x_2}{2}dx_3+\big(\frac{x_3}{2}+\frac{x_1x_2}{6}\big)dx_2-\frac{x_2^2}{6}dx_1$;
\item $\theta_6=dx_6-\frac{x_1}{2}dx_5-\frac{x_2}{2}dx_4+\frac{x_1x_2}{3}dx_3+\big(\frac{x_4}{2}-\frac{x_1x_3}{6}-\frac{x_1^2x_2}{12}\big)dx_2+\big(\frac{x_5}{2}-\frac{x_2x_3}{6}+\frac{x_1x_2^2}{12}\big)dx_1$.
\end{itemize}

Finally, we have
\begin{eqnarray*}
 \mathrm{d}(R_\mathbf{x})_\mathbf{0}=\left[\begin{matrix} 
   1 & 0  & 0 & 0 & 0 & 0\\
   0  & 1 & 0 & 0 & 0 &0\\
   \frac{x_2}{2} & -\frac{x_1}{2} & 1 & 0 & 0 & 0\\
  \frac{x_3}{2}-\frac{x_1x_2}{12} & \frac{x_1^2}{12} & -\frac{x_1}{2} & 1& 0 & 0\\
  -\frac{x_2^2}{12}& \frac{x_3}{2}+\frac{x_1x_2}{12} & -\frac{x_2}{2} & 0& 1 & 0\\
  \frac{x_5}{2}-\frac{x_2x_3}{12} & \frac{x_4}{2}-\frac{x_1x_3}{12} & \frac{x_1x_2}{6} & -\frac{x_2}{2}& -\frac{x_1}{2} & 1
   \end{matrix}\right]\,.
\end{eqnarray*}
   \subsection*{$\bf N_{6,2,5\,a}$  }
  
   The following Lie algebra is denoted as $N_{6,2,5\,a}$ by Gong in \cite{Gong_Thesis}, as $L_{6,21(1)}$ by de Graaf in \cite{deGraaf_classification}, as $(8)$ by Del Barco in \cite{delBarco}, and as $\mathcal{G}_{6,15}$ by Magnin in \cite{magnin}.
   
   The non-trivial brackets are the following:
\begin{equation*}\label{N625a}
   [X_1, X_i] = X_{i+1}\,,\, i = 2, 3\,,\, [X_2, X_3] = X_{5}\,,\,[X_1,X_4]=-X_6\,,\,[X_2,X_5]=-X_6\,.
\end{equation*}
This is a nilpotent Lie algebra of rank 2 and step 4 that is stratifiable. The Lie brackets can be pictured with the diagram:
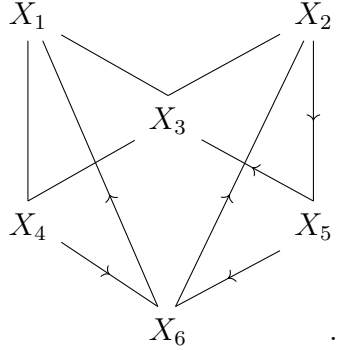
\begin{center}
 
	\begin{tikzcd}[end anchor=north]
	X_1\ar[dr, no head]\ar[dd, no head]\ar[dddr,-<-=.6, no head,end anchor={[xshift=-3.2ex]north east}]& & X_2\ar[dl, no head]\ar[dddl, no head, -<-=.6,end anchor={[xshift=-1.9ex]north east}]\ar[dd,->-=.5, no head]\\
	& X_3\ar[dl, no head]\ar[dr, -<-=.5, no head]&\\
	X_4\ar[dr, no head, ->-=.5, end anchor={[xshift=-3.2ex]north east}]& & X_5\ar[no head,->-=.5, dl,  end anchor={[xshift=-1.9ex]north east}]\\
	& X_6 & \quad\;.
	\end{tikzcd}

\end{center}

 
The composition law \eqref{group law in G} of $N_{6,2,5\,a}$ is given by:

\begin{itemize}
    \item $z_1=x_1+y_1$;
    \item $z_2=x_2+y_2$;
    \item $z_3=x_3+y_3+\frac{1}{2}(x_1y_2-x_2y_1)$;
    \item $z_4=x_4+y_4+\frac{1}{2}(x_1y_3-x_3y_1)+\frac{1}{12}(x_1-y_1)(x_1y_2-x_2y_1)$;
    \item $z_5=x_5+y_5+\frac{1}{2}(x_2y_3-x_3y_2)+\frac{1}{12}(x_2-y_2)(x_1y_2-x_2y_1)$;
    \item $z_6=x_6+y_6+\frac{1}{2}(x_4y_1-x_1y_4+x_5y_2-x_2y_5)+\frac{1}{12}(y_1-x_1)(x_1y_3-x_3y_1)\salto+\frac{1}{12}(y_2-x_2)(x_2y_3-x_3y_2)+\frac{1}{24}(x_1y_1+x_2y_2)(x_1y_2-x_2y_1)$.
\end{itemize}

Since
\begin{eqnarray*}
\mathrm{d}(L_\mathbf{x})_\mathbf{0}= \left[\begin{matrix} 
   1 & 0  & 0 & 0 & 0 & 0\\
   0  & 1 & 0 & 0 & 0 &0\\
   -\frac{x_2}{2} & \frac{x_1}{2} & 1 & 0 & 0 & 0\\
  -\frac{x_1x_2}{12}-\frac{x_3}{2} & \frac{x_1^2}{12} & \frac{x_1}{2} & 1& 0 & 0\\
  -\frac{x_2^2}{12}& \frac{x_1x_2}{12}-\frac{x_3}{2} & \frac{x_2}{2} & 0& 1 & 0\\
  \frac{x_1x_3}{12}+\frac{x_4}{2} & \frac{x_5}{2}+\frac{x_2x_3}{12} & -\frac{x_1^2+x_2^2}{12} & -\frac{x_1}{2}&- \frac{x_2}{2} & 1
   \end{matrix}\right]\,,
\end{eqnarray*}
the induced left-invariant vector fields \eqref{leftinvariant vf} are: 
\begin{itemize}
\item $X_1=\partial_{x_1} -\frac{x_2}{2}\partial_{x_3}-\big(\frac{x_3}{2}+\frac{x_1x_2}{12}\big)\partial_{x_4}-\frac{x_2^2}{12}\partial_{x_5}+\big(\frac{x_1x_3}{12}+\frac{x_4}{2}\big)\partial_{x_6}\,;$
\item $X_2=\partial_{x_2}+\frac{x_1}{2}\partial_{x_3}+\frac{x_1^2}{12}\partial_{x_4}+\big(\frac{x_1x_2}{12}-\frac{x_3}{2}\big)\partial_{x_5}+\big(\frac{x_5}{2}+\frac{x_2x_3}{12}\big)\partial_{x_6}\,;$
\item $X_3=\partial_{x_3}+\frac{x_1}{2}\partial_{x_4}+\frac{x_2}{2}\partial_{x_5}-\frac{x_1^2+x_2^2}{12}\partial_{x_6}\,;$
\item $X_4=\partial_{x_4}-\frac{x_1}{2}\partial_{x_6}\,;$
\item $X_5=\partial_{x_5}-\frac{x_2}{2}\partial_{x_6}$;
\item $X_6=\partial_{x_6}$,
\end{itemize}
and the respective left-invariant 1-forms \eqref{leftinvariant form} are: 
\begin{itemize}
\item $\theta_1=dx_1$;
\item $\theta_2=dx_2$;
\item $\theta_3=dx_3-\frac{x_1}{2} dx_2+\frac{x_2}{2}dx_1$;
\item $\theta_4=dx_4-\frac{x_1}{2}dx_3+\frac{x_1^2}{6}dx_2+\big(\frac{x_3}{2}-\frac{x_1x_2}{6}\big)dx_1$;
\item $\theta_5=dx_5-\frac{x_2}{2}dx_3+\big(\frac{x_3}{2}+\frac{x_1x_2}{6}\big)dx_2-\frac{x_2^2}{6}dx_1$;
\item $\theta_6=dx_6+\frac{x_2}{2}dx_5+\frac{x_1}{2}dx_4-\frac{x_1^2+x_2^2}{6}dx_3+\big(\frac{x_2x_3}{6}-\frac{x_5}{2}+\frac{x_1^3+x_1x_2^2}{24}\big)dx_2-\big(\frac{x_4}{2}-\frac{x_1x_3}{6}$\\$  \phantom{\theta_5=}+\frac{x_1^2x_2+x_2^3}{24}\big)dx_1$.
\end{itemize}

Finally, we have
\begin{eqnarray*}
\mathrm{d}(R_\mathbf{x})_\mathbf{0}= \left[\begin{matrix} 
   1 & 0  & 0 & 0 & 0 & 0\\
   0  & 1 & 0 & 0 & 0 &0\\
   \frac{x_2}{2} & -\frac{x_1}{2} & 1 & 0 & 0 & 0\\
  \frac{x_3}{2}-\frac{x_1x_2}{12} & \frac{x_1^2}{12} & -\frac{x_1}{2} & 1& 0 & 0\\
  -\frac{x_2^2}{12}& \frac{x_1x_2}{12}+\frac{x_3}{2} & -\frac{x_2}{2} & 0& 1 & 0\\
  \frac{x_1x_3}{12}-\frac{x_4}{2} & \frac{x_2x_3}{12}-\frac{x_5}{2} & -\frac{x_1^2+x_2^2}{12} & \frac{x_1}{2}& \frac{x_2}{2} & 1
   \end{matrix}\right]\,.
\end{eqnarray*}
   \subsection*{$\bf N_{6,2,6}$  }

   The following Lie algebra is denoted as $N_{6,2,6}$ by Gong in \cite{Gong_Thesis}, as $L_{6,20}$ by de Graaf in \cite{deGraaf_classification}, as $(14)$ by Del Barco in \cite{delBarco}, and as $\mathcal{G}_{6,10}$ by Magnin in \cite{magnin}.
   
   The non-trivial brackets are the following:
\begin{equation*}\label{N626}
  [X_1, X_2] = X_4\,,\, [X_1, X_3] = X_5\,,\, [X_1, X_4] = X_6\,,\, [X_3, X_5] = X_6\,.
\end{equation*}
This is a nilpotent Lie algebra of rank 3 and step 3 that is stratifiable. The Lie brackets can be pictured with the diagram:
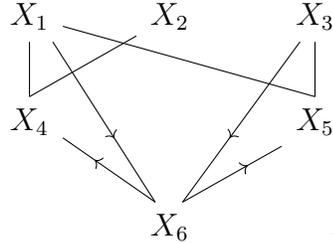
\begin{center}
 
	\begin{tikzcd}[end anchor=north]
	 X_1\ar[d, no head] \ar[drr, no head] \ar[ddr, no head, ->-=.6,end anchor={[xshift=-3.5ex]north east}]& X_2\ar[dl,no head] & X_3 \ar[d,no head]\ar[ddl, no head, ->-=.6, end anchor={[xshift=-1.5ex]north east}] \\
	X_4\ar[dr, no head, -<-=.4, end anchor={[xshift=-3.5ex]north east}]& & X_5\ar[dl, no head, -<-=.4, end anchor={[xshift=-1.5ex]north east}] \\
	 &X_6& \quad\;.
	\end{tikzcd}

\end{center}

 
The composition law \eqref{group law in G} of $N_{6,2,6}$ is given by:

\begin{itemize}
    \item $z_1=x_1+y_1$;
    \item $z_2=x_2+y_2$;
    \item $z_3=x_3+y_3$;
    \item $z_4=x_4+y_4+\frac{1}{2}(x_1y_2-x_2y_1)$;
    \item $z_5=x_5+y_5+\frac{1}{2}(x_1y_3-x_3y_1)$;
    \item $z_6=x_6+y_6+\frac{1}{2}(x_1y_4-x_4y_1+x_3y_5-x_5y_3)+\frac{1}{12}(x_1-y_1)(x_1y_2-x_2y_1)\salto+\frac{1}{12}(x_3-y_3)(x_1y_3-x_3y_1)$.
\end{itemize}


Since
\begin{eqnarray*}
 \mathrm{d}(L_\mathbf{x})_\mathbf{0}= \left[\begin{matrix} 
   1 & 0  & 0 & 0 & 0 & 0\\
   0  & 1 & 0 & 0 & 0 &0\\
   0 & 0 & 1 & 0 & 0 & 0\\
  -\frac{x_2}{2} & \frac{x_1}{2} & 0 & 1& 0 & 0\\
  -\frac{x_3}{2} & 0 & \frac{x_1}{2} & 0& 1 & 0\\
 -\frac{x_1x_2+x_3^2}{12}-\frac{x_4}{2} & \frac{x_1^2}{12} &\frac{x_1x_3}{12} -\frac{x_5}{2} & \frac{x_1}{2}& \frac{x_3}{2} & 1
   \end{matrix}\right]\,,
\end{eqnarray*}
the induced left-invariant vector fields \eqref{leftinvariant vf} are: 
\begin{itemize}
\item $X_1=\partial_{x_1} -\frac{x_2}{2}\partial_{x_4}-\frac{x_3}{2}\partial_{x_5}-\big(\frac{x_1x_2+x_3^2}{12}+\frac{x_4}{2}\big)\partial_{x_6}\,;$
\item $X_2=\partial_{x_2}+\frac{x_1}{2}\partial_{x_4}+\frac{x_1^2}{12}\partial_{x_6}\,;$
\item $X_3=\partial_{x_3}+\frac{x_1}{2}\partial_{x_5}+\big(\frac{x_1x_3}{12}-\frac{x_5}{2}\big)\partial_{x_6}\,;$
\item $X_4=\partial_{x_4}+\frac{x_1}{2}\partial_{x_6}\,;$
\item $X_5=\partial_{x_5}+\frac{x_3}{2}\partial_{x_6}$;
\item $X_6=\partial_{x_6}$,
\end{itemize}
and the respective left-invariant 1-forms \eqref{leftinvariant form} are: 
\begin{itemize}
\item $\theta_1=dx_1$;
\item $\theta_2=dx_2$;
\item $\theta_3=dx_3$;
\item $\theta_4=dx_4-\frac{x_1}{2}dx_2+\frac{x_2}{2}dx_1$;
\item $\theta_5=dx_5-\frac{x_1}{2}dx_3+\frac{x_3}{2}dx_1$;
\item $\theta_6=dx_6-\frac{x_3}{2}dx_5-\frac{x_1}{2}dx_4+\big(\frac{x_5}{2}+\frac{x_1x_3}{6}\big)dx_3+\frac{x_1^2}{6}dx_2+\big(\frac{x_4}{2}-\frac{x_3^3+x_1x_2}{6}\big)dx_1$.
\end{itemize}

Finally, we have
\begin{eqnarray*}
 \mathrm{d}(R_\mathbf{x})_\mathbf{0}= \left[\begin{matrix} 
   1 & 0  & 0 & 0 & 0 & 0\\
   0  & 1 & 0 & 0 & 0 &0\\
   0 & 0 & 1 & 0 & 0 & 0\\
  \frac{x_2}{2} & -\frac{x_1}{2} & 0 & 1& 0 & 0\\
  \frac{x_3}{2} & 0 & -\frac{x_1}{2} & 0& 1 & 0\\
 \frac{x_4}{2}-\frac{x_1x_2+x_3^2}{12} & \frac{x_1^2}{12} &\frac{x_1x_3}{12} +\frac{x_5}{2} & -\frac{x_1}{2}& -\frac{x_3}{2} & 1
   \end{matrix}\right]\,.
\end{eqnarray*}

    \subsection*{$\bf N_{6,2,7}$ }

   The following Lie algebra is denoted as $N_{6,2,7}$ by Gong in \cite{Gong_Thesis}, as $L_{6,21(0)}$ by de Graaf in \cite{deGraaf_classification}, as $(6)$ by Del Barco in \cite{delBarco}, and as $\mathcal{G}_{6,14}$ by Magnin in \cite{magnin}.
   
   The non-trivial brackets are the following:
\begin{equation*}\label{N627}
  [X_1, X_i] = X_{i+1}\,,\, i = 2, 3, 4\,,\, [X_2, X_3] = X_6\,.
\end{equation*}
This is a nilpotent Lie algebra of rank 2 and step 4 that is stratifiable. The Lie brackets can be pictured with the diagram:
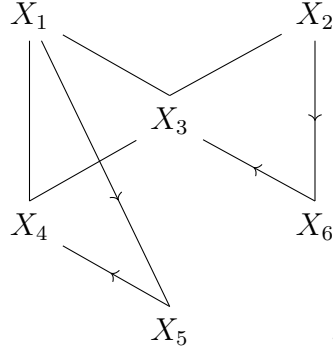
\begin{center}
 
	\begin{tikzcd}[end anchor=north]
	X_1\ar[dr, no head]\ar[dd, no head]\ar[dddr, no head, ->-=.6]& & X_2\ar[dl, no head]\ar[dd,->-=.5, no head]\\
	& X_3\ar[dl, no head]\ar[dr, -<-=.5, no head]&\\
	X_4\ar[dr, no head, -<-=.5]& & X_6\\
	& X_5 & \quad\;.
	\end{tikzcd}
   
\end{center}
 
  
The composition law \eqref{group law in G} of $N_{6,2,7}$ is given by:

\begin{itemize}
    \item $z_1=x_1+y_1$;
    \item $z_2=x_2+y_2$;
    \item $z_3=x_3+y_3+\frac{1}{2}(x_1y_2-x_2y_1)$;
    \item $z_4=x_4+y_4+\frac{1}{2}(x_1y_3-x_3y_1)+\frac{1}{12}(x_1-y_1)(x_1y_2-x_2y_1)$;
    \item $z_5=x_5+y_5+\frac{1}{2}(x_1y_4-x_4y_1)+\frac{1}{12}(x_1-y_1)(x_1y_3-x_3y_1)-\frac{1}{24}x_1y_1(x_1y_2-x_2y_1)$;
    \item $z_6=x_6+y_6+\frac{1}{2}(x_2y_3-x_3y_2)+\frac{1}{12}(x_2-y_2)(x_1y_2-x_2y_1)$.
\end{itemize}



Since
\begin{eqnarray*}
  \mathrm{d}(L_\mathbf{x})_\mathbf{0}=\left[\begin{matrix} 
   1 & 0  & 0 & 0 & 0 & 0\\
   0  & 1 & 0 & 0 & 0 &0\\
   -\frac{x_2}{2} & \frac{x_1}{2} & 1 & 0 & 0 & 0\\
 -\frac{x_3}{2}- \frac{x_1x_2}{12} & \frac{x_1^2}{12} & \frac{x_1}{2} & 1& 0 & 0\\
- \frac{x_4}{2}-\frac{x_1x_3}{12}  & 0 & \frac{x_1^2}{12} & \frac{x_1}{2}& 1 & 0\\
 -\frac{x_2^2}{12} & \frac{x_1x_2}{12}-\frac{x_3}{2} & \frac{x_2}{2} & 0 & 0 & 1
   \end{matrix}\right]\,,
\end{eqnarray*}
the induced left-invariant vector fields \eqref{leftinvariant vf} are:
\begin{itemize}
\item $X_1=\partial_{x_1} -\frac{x_2}{2}\partial_{x_3}-\big(\frac{x_3}{2}+\frac{x_1x_2}{12}\big)\partial_{x_4}-\big(\frac{x_1x_3}{12}+\frac{x_4}{2}\big)\partial_{x_5}-\frac{x_2^2}{12}\partial_{x_6}\,;$
\item $X_2=\partial_{x_2}+\frac{x_1}{2}\partial_{x_3}+\frac{x_1^2}{12}\partial_{x_4}+\big(\frac{x_1x_2}{12}-\frac{x_3}{2}\big)\partial_{x_6}\,;$
\item $X_3=\partial_{x_3}+\frac{x_1}{2}\partial_{x_4}+\frac{x_1^2}{12}\partial_{x_5}+\frac{x_2}{2}\partial_{x_6}\,;$
\item $X_4=\partial_{x_4}+\frac{x_1}{2}\partial_{x_5}\,;$
\item $X_5=\partial_{x_5}$;
\item $X_6=\partial_{x_6}$,
\end{itemize}
and the respective left-invariant 1-forms \eqref{leftinvariant form} are: 
\begin{itemize}
\item $\theta_1=dx_1$;
\item $\theta_2=dx_2$;
\item $\theta_3=dx_3-\frac{x_1}{2} dx_2+\frac{x_2}{2}dx_1$;
\item $\theta_4=dx_4-\frac{x_1}{2}dx_3+\frac{x_1^2}{6}dx_2+\big(\frac{x_3}{2}-\frac{x_1x_2}{6}\big)dx_1$;
\item $\theta_5=dx_5-\frac{x_1}{2}dx_4+\frac{x_1^2}{6}dx_3-\frac{x_1^3}{24}dx_2+\big(\frac{x_4}{2}-\frac{x_1x_3}{6}+\frac{x_1^2x_2}{24}\big)dx_1$;
\item $\theta_6=dx_6-\frac{x_2}{2}dx_3+\big(\frac{x_3}{2}+\frac{x_1x_2}{6}\big)dx_2-\frac{x_2^2}{6}dx_1$.
\end{itemize}

Finally, we have
\begin{eqnarray*}
  \mathrm{d}(R_\mathbf{x})_\mathbf{0}=\left[\begin{matrix} 
   1 & 0  & 0 & 0 & 0 & 0\\
   0  & 1 & 0 & 0 & 0 &0\\
   \frac{x_2}{2} & -\frac{x_1}{2} & 1 & 0 & 0 & 0\\
 \frac{x_3}{2}- \frac{x_1x_2}{12} & \frac{x_1^2}{12} & -\frac{x_1}{2} & 1& 0 & 0\\
 \frac{x_4}{2}-\frac{x_1x_3}{12}  & 0 & \frac{x_1^2}{12} & -\frac{x_1}{2}& 1 & 0\\
 -\frac{x_2^2}{12} & \frac{x_1x_2}{12}+\frac{x_3}{2} & -\frac{x_2}{2} & 0 & 0 & 1
   \end{matrix}\right]\,.
\end{eqnarray*}
    \subsection*{$\bf N_{6,2,8}$ non-stratifiable}


   The following Lie algebra is denoted as $N_{6,2,8}$ by Gong in \cite{Gong_Thesis}, as $L_{6,23}$ by de Graaf in \cite{deGraaf_classification}, as $(21)$ by Del Barco in \cite{delBarco}, and as $\mathcal{G}_{6,7}$ by Magnin in \cite{magnin}.
   
   The non-trivial brackets are the following:
\begin{equation*}
 [X_1, X_2] = X_3\,,\, [X_1, X_3] = X_5\,,\, [X_1, X_4] = X_6\,,\, [X_2, X_4] = X_5\,.
\end{equation*}
This is a nilpotent Lie algebra of rank 3 and step 3 that is positively gradable, yet not stratifiable. The Lie brackets can be pictured with the diagram:
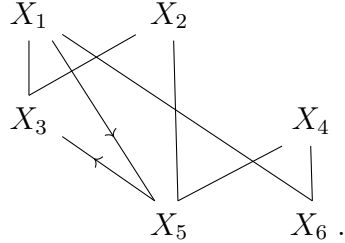
\begin{center}
 
	\begin{tikzcd}[end anchor=north]
	 X_1\ar[d, no head] \ar[ddrr, no head,end anchor={[xshift=-3.5ex]north east}] \ar[ddr, no head, ->-=.6,end anchor={[xshift=-3.5ex]north east}]& X_2\ar[dl,no head]\ar[dd, no head,  end anchor={[xshift=-1.8ex]north east}]  &  \\
	X_3\ar[dr, no head, -<-=.4, end anchor={[xshift=-3.5ex]north east}]& &X_4\;\, \ar[d,no head,end anchor={[xshift=-3.5ex]north east},start anchor={[xshift=-3.5ex]south east}]\ar[dl, no head,  end anchor={[xshift=-1.8ex]north east}] \\
	 &X_5& X_6\;.
	\end{tikzcd}

\end{center}

 
The composition law \eqref{group law in G} of $N_{6,2,8}$ is given by:

\begin{itemize}
    \item $z_1=x_1+y_1$;
    \item $z_2=x_2+y_2$;
    \item $z_3=x_3+y_3+\frac{1}{2}(x_1y_2-x_2y_1)$;
    \item $z_4=x_4+y_4$;
    \item $z_5=x_5+y_5+\frac{1}{2}(x_1y_3-x_3y_1+x_2y_4-x_4y_2)+\frac{1}{12}(x_1-y_1)(x_1y_2-x_2y_1)$;
    \item $z_6=x_6+y_6+\frac{1}{2}(x_1y_4-x_4y_1)$.
\end{itemize}


Since
\begin{eqnarray*}
 \mathrm{d}(L_\mathbf{x})_\mathbf{0}= \left[\begin{matrix} 
   1 & 0  & 0 & 0 & 0 & 0\\
   0  & 1 & 0 & 0 & 0 &0\\
   -\frac{x_2}{2} & \frac{x_1}{2} & 1 & 0 & 0 & 0\\
  0 & 0 & 0 & 1& 0 & 0\\
  -\frac{x_3}{2}-\frac{x_1x_2}{12} & \frac{x_1^2}{12}-\frac{x_4}{2} & \frac{x_1}{2} & \frac{x_2}{2}& 1 & 0\\
-\frac{x_4}{2} & 0 & 0 & \frac{x_1}{2}& 0 & 1
   \end{matrix}\right]\,,
\end{eqnarray*}
the induced left-invariant vector fields \eqref{leftinvariant vf} are: 
\begin{itemize}
\item $X_1=\partial_{x_1}- \frac{x_2}{2}\partial_{x_3}-\big(\frac{x_3}{2}+\frac{x_1x_2}{12}\big)\partial_{x_5}-\frac{x_4}{2}\partial_{x_6}\,;$
\item $X_2=\partial_{x_2}+\frac{x_1}{2}\partial_{x_3}+\big(\frac{x_1^2}{12}-\frac{x_4}{2}\big)\partial_{x_5}\,;$
\item $X_3=\partial_{x_3}+\frac{x_1}{2}\partial_{x_5}\,;$
\item $X_4=\partial_{x_4}+\frac{x_2}{2}\partial_{x_5}+\frac{x_1}{2}\partial_{x_6}\,;$
\item $X_5=\partial_{x_5}$;
\item $X_6=\partial_{x_6}$,
\end{itemize}
and the respective left-invariant 1-forms \eqref{leftinvariant form} are: 
\begin{itemize}
\item $\theta_1=dx_1$;
\item $\theta_2=dx_2$;
\item $\theta_3=dx_3-\frac{x_1}{2}dx_2+\frac{x_2}{2}dx_1$;
\item $\theta_4=dx_4$;
\item $\theta_5=dx_5-\frac{x_2}{2}dx_4-\frac{x_1}{2}dx_3+\big(\frac{x_4}{2}+\frac{x_1^2}{6}\big)dx_2+\big(\frac{x_3}{2}-\frac{x_1x_2}{6}\big)dx_1$;
\item $\theta_6=dx_6-\frac{x_1}{2}dx_4+\frac{x_4}{2}dx_1$.
\end{itemize}

Finally, we have
\begin{eqnarray*}
 \mathrm{d}(R_\mathbf{x})_\mathbf{0}= \left[\begin{matrix} 
   1 & 0  & 0 & 0 & 0 & 0\\
   0  & 1 & 0 & 0 & 0 &0\\
   \frac{x_2}{2} & -\frac{x_1}{2} & 1 & 0 & 0 & 0\\
  0 & 0 & 0 & 1& 0 & 0\\
  \frac{x_3}{2}-\frac{x_1x_2}{12} & \frac{x_1^2}{12}+\frac{x_4}{2} & -\frac{x_1}{2} & -\frac{x_2}{2}& 1 & 0\\
\frac{x_4}{2} & 0 & 0 & -\frac{x_1}{2}& 0 & 1
   \end{matrix}\right]\,.
\end{eqnarray*}
 
 \subsubsection{Grading, polarizations of maximal dimension, and asymptotic cone}
 The Lie algebra $N_{6,2,8}$
   is not stratifiable, but it is gradable as
   \begin{eqnarray*}
   V_1 &= & \span \{  X_1,X_2 \}\,,\\
   V_2 &= & \span \{  X_3,X_4 \}\,,\\
   V_3 &= & \span \{  X_5,X_6 \}\,.
   \end{eqnarray*}
      
     Every   complementary subspace $\Delta$ to the  derived subalgebra is spanned by
      $X_1 +u_1^3X_3+u_1^5X_5+u_1^6X_6 $,
      $X_2 + u_2^3X_3+u_2^5X_5+u_2^6X_6 $, and
       $X_4 + u^3_4X_3+u_4^5X_5+u_4^6X_6 $.
       
      We claim that in this Lie algebra every two complementary subspaces to the  derived subalgebra as in \eqref{complementary_Delta} differ by an automorphism. Indeed, when $u_i^j$ varies, 
      the matrix
$$    \left[  \begin{array}{cccccc}
     1 &  0 &  0 &  0 &  0 &  0  \\
 0 &  1 &  0 &  0 &  0 &  0\\
 u_1^3 & u_2^3 & 1 &  u_4^3 & 0 &  0\\
 0 & 0 & 0 & 1 & 0 & 0\\
 u_1^5 &  u_2^5 & u_2^3 &  u_4^5 &  1 &  u_4^3\\
 u_1^6 &  u_2^6 &  0 &  u_4^6 &  0 &  1\\
      \end{array}\right]$$
      sends the complementary subspace $\span\lbrace X_1,X_2,X_4\rbrace$ to
      $\span\{X_1 +u_1^3X_3+u_1^5X_5+u_1^6X_6, X_2 + u_2^3X_3+u_2^5X_5+u_2^6X_6,X_4 + u^3_4X_3+u_4^5X_5+u_4^6X_6\}$, which is an arbitrary one. In particular, every $\Delta$ as in \eqref{complementary_Delta} gives maximal dimension. The tangent cone of each  of such polarizations has Lie algebra isomorphic to $N_{6,3,6}$, see page \pageref{N636}.

 The asymptotic cone 
 of the Lie group with Lie algebra
 $N_{6,2,8}$ has Lie algebra isomorphic to $N_{6,3,4}$, see page \pageref{N634}.

    \subsection*{$\bf N_{6,2,9}$ non-stratifiable }

   The following Lie algebra is denoted as $N_{6,2,9}$ by Gong in \cite{Gong_Thesis}, as $L_{6,24(1)}$ by de Graaf in \cite{deGraaf_classification}, as $(20)$ by Del Barco in \cite{delBarco}, and as $\mathcal{G}_{6,5}$ by Magnin in \cite{magnin}.
   
   The non-trivial brackets are the following:
\begin{equation*}
   [X_1, X_2] = X_3\,,\, [X_1, X_3] = X_5
   \,,\, [X_2, X_3] = X_6\,,\, [X_2, X_4] = X_6\,.
\end{equation*}
This is a nilpotent Lie algebra of rank 3 and step 3 that is positively gradable, yet not stratifiable. The Lie brackets can be pictured with the diagram:
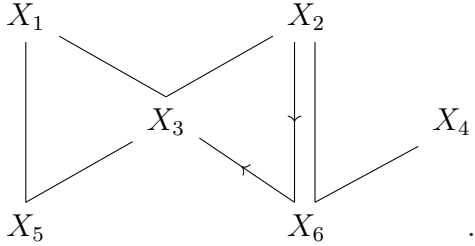
\begin{center}
 
	\begin{tikzcd}[end anchor=north]
	X_1\ar[dr, no head]\ar[dd, no head]& & X_2\ar[dl, no head]\ar[dd,->-=.5, no head,end anchor={[xshift=-3.3ex]north east},start anchor={[xshift=-3.3ex]south east}] \ar[dd,no head,end anchor={[xshift=-1.8ex]north east},start anchor={[xshift=-1.8ex]south east}]& \\
	& X_3\ar[dl, no head]\ar[dr, -<-=.5, no head,end anchor={[xshift=-3.3ex]north east}]& &X_4\ar[dl,no head, end anchor={[xshift=-1.8ex]north east}] \\
	X_5& & X_6 &\quad\; .
	\end{tikzcd}

\end{center}

 
 The composition law \eqref{group law in G} of $N_{6,2,9}$ is given by:

\begin{itemize}
    \item $z_1=x_1+y_1$;
    \item $z_2=x_2+y_2$;
    \item $z_3=x_3+y_3+\frac{1}{2}(x_1y_2-x_2y_1)$;
    \item $z_4=x_4+y_4$;
    \item $z_5=x_5+y_5+\frac{1}{2}(x_1y_3-x_3y_1)+\frac{1}{12}(x_1-y_1)(x_1y_2-x_2y_1)$;
    \item $z_6=x_6+y_6+\frac{1}{2}(x_2y_3-x_3y_2+x_2y_4-x_4y_2)+\frac{1}{12}(x_2-y_2)(x_1y_2-x_2y_1)$.
\end{itemize}



Since
\begin{eqnarray*}
 \mathrm{d}(L_\mathbf{x})_\mathbf{0}= \left[\begin{matrix} 
   1 & 0  & 0 & 0 & 0 & 0\\
   0  & 1 & 0 & 0 & 0 &0\\
   -\frac{x_2}{2} & \frac{x_1}{2} & 1 & 0 & 0 & 0\\
  0 & 0 & 0 & 1& 0 & 0\\
 - \frac{x_1x_2}{12}-\frac{x_3}{2} & \frac{x_1^2}{12} & \frac{x_1}{2}& 0& 1 & 0\\
  -\frac{x_2^2}{12} & \frac{x_1x_2}{12}-\frac{x_3+x_4}{2} & \frac{x_2}{2}  & \frac{x_2}{2}& 0 & 1
   \end{matrix}\right]\,,
\end{eqnarray*}
the induced left-invariant vector fields \eqref{leftinvariant vf} are: 
\begin{itemize}
\item $X_1=\partial_{x_1}-\frac{x_2}{2}\partial_{x_3}-\big(\frac{x_3}{2}+\frac{x_1x_2}{12}\big)\partial_{x_5}-\frac{x_2^2}{12}\partial_{x_6}\,;$
\item $X_2=\partial_{x_2}+\frac{x_1}{2}\partial_{x_3}+\frac{x_1^2}{12}\partial_{x_5}+\big(\frac{x_1x_2}{12}-\frac{x_3+x_4}{2}\big)\partial_{x_6}\,;$
\item $X_3=\partial_{x_3}+\frac{x_1}{2}\partial_{x_5}+\frac{x_2}{2}\partial_{x_6}\,;$
\item $X_4=\partial_{x_4}+\frac{x_2}{2}\partial_{x_6}\,;$
\item $X_5=\partial_{x_5}$;
\item $X_6=\partial_{x_6}$,
\end{itemize}
and the respective left-invariant 1-forms \eqref{leftinvariant form} are: 
\begin{itemize}
\item $\theta_1=dx_1$;
\item $\theta_2=dx_2$;
\item $\theta_3=dx_3-\frac{x_1}{2} dx_2+\frac{x_2}{2}dx_1$;
\item $\theta_4=dx_4$;
\item $\theta_5=dx_5-\frac{x_1}{2}dx_3+\frac{x_1^2}{6}dx_2+\big(\frac{x_3}{2}-\frac{x_1x_2}{6}\big)dx_1$;
\item $\theta_6=dx_6-\frac{x_2}{2}dx_4-\frac{x_2}{2}dx_3+\big(\frac{x_3+x_4}{2}+\frac{x_1x_2}{6}\big)dx_2-\frac{x_2^2}{6}dx_1$.
\end{itemize}

Finally, we have
\begin{eqnarray*}
 \mathrm{d}(R_\mathbf{x})_\mathbf{0}= \left[\begin{matrix} 
   1 & 0  & 0 & 0 & 0 & 0\\
   0  & 1 & 0 & 0 & 0 &0\\
   \frac{x_2}{2} & -\frac{x_1}{2} & 1 & 0 & 0 & 0\\
  0 & 0 & 0 & 1& 0 & 0\\
 \frac{x_3}{2}- \frac{x_1x_2}{12} & \frac{x_1^2}{12} & -\frac{x_1}{2}& 0& 1 & 0\\
  -\frac{x_2^2}{12} & \frac{x_1x_2}{12}+\frac{x_3+x_4}{2} & -\frac{x_2}{2}  & -\frac{x_2}{2}& 0 & 1
   \end{matrix}\right]\,.
\end{eqnarray*}

\subsubsection{Grading, polarizations of maximal dimension, and asymptotic cone}
 The Lie algebra $N_{6,2,9}$
   is not stratifiable, but it is gradable as
   \begin{eqnarray*}
   V_1 &= & \span \{  X_1,X_2 \}\,,\\ V_2 &= & \span \{  X_3,X_4 \}\,,\\
   V_3 &= & \span \{ X_5,X_6 \}\,.
   \end{eqnarray*}
      
           Every   complementary subspace $\Delta$ to the  derived subalgebra is spanned by
      $X_1 +   u_1^3 X_3+ u_1^5 X_5+ u_1^6 X_6 $,
      $X_2 +   u_2^3 X_3+   u_2^5 X_5+   u_2^6 X_6$, and
       $X_4 +   u_4^3 X_3 +   u_4^5 X_5+   u_4^6 X_6$.
       Such a polarization gives maximal Hausdorff dimension if and only if 
          $u_4^3 $ is either $ -1 $ or $0$.
      We claim that every two polarizations giving maximal Hausdorff dimension    differ by an automorphism. 
        Indeed, when $u_i^j$ varies, if we take $u_4^3=0$ 
      the matrix
$$    \left[  \begin{array}{cccccc}
     1 &  0 &  0 &  0 &  0 &  0  \\
 0 &  1 &  0 &  0 &  0 &  0\\
 u_1^3 & u_2^3 & 1 &  0 & 0 &  0\\
 0 & 0 & 0 & 1 & 0 & 0\\
 u_1^5 &  u_2^5 & u_2^3 &  u_4^5 &  1 &  0\\
 u_1^6 &  u_2^6 & -u_1^3 &  u_4^6 & 0 &  1\\
      \end{array}\right],$$
      possibly composed with the block diagonal matrix
      $$
      \left[  \begin{array}{cc} 
      0&1\\
      1&0
\end{array}\right]
\oplus
      \left[  \begin{array}{cc} 
      -1&-1\\
      0&1
\end{array}\right]
\oplus
      \left[  \begin{array}{cc} 
      0&-1\\
      -1&0
\end{array}\right]
, $$ 
         is a Lie algebra automorphism and   sends the complementary subspace $\span\lbrace X_1,X_2,X_4\rbrace$ to
      $\span\{X_1 +   u_1^3 X_3+ u_1^5 X_5+ u_1^6 X_6, X_2 +   u_2^3 X_3+   u_2^5 X_5+   u_2^6 X_6,X_4 +     u_4^5 X_5+   u_4^6 X_6\}$, which is an arbitrary one
        of maximal dimension. 
        
         If instead we take $u_4^3=-1$, as $u_i^j$ varies,
      the matrix
$$    \left[  \begin{array}{cccccc}
     1 &  0 &  0 &  0 &  0 &  0  \\
 0 &  1 &  0 &  0 &  0 &  0\\
 u_1^3 & u_2^3 & 1 &  -1 & 0 &  0\\
 0 & 0 & 0 & 1 & 0 & 0\\
 u_1^5 &  u_2^5 & u_2^3 &  u_4^5 &  -1 &  0\\
 u_1^6 &  u_2^6 & -u_1^3 &  u_4^6 & 0 &  1\\
      \end{array}\right],$$
      possibly composed with the block diagonal matrix
      $$
      \left[  \begin{array}{cc} 
      0&1\\
      1&0
\end{array}\right]
\oplus
      \left[  \begin{array}{cc} 
      -1&-1\\
      0&1
\end{array}\right]
\oplus
      \left[  \begin{array}{cc} 
      0&-1\\
      -1&0
\end{array}\right]
, $$          is a Lie algebra automorphism and  sends the complementary subspace  $\span\lbrace X_1,X_2,X_4\rbrace$ to
      $\span\{X_1 +   u_1^3 X_3+ u_1^5 X_5+ u_1^6 X_6, X_2 +   u_2^3 X_3+   u_2^5 X_5+   u_2^6 X_6,X_4 -X_3+ u_4^5 X_5+   u_4^6 X_6\}$, which is an arbitrary one
        of maximal dimension. 
        The tangent cone of each of such polarizations has Lie algebra isomorphic to $N_{6,3,3}$, see page \pageref{N633}.
        
         One should also be aware that Aut$(N_{6,2,9})$ has two connected components (see \cite{Gong_Thesis} on page 35).

 
 The asymptotic cone 
 of the Lie group with Lie algebra
 $N_{6,2,9}$ has Lie algebra isomorphic to $N_{5,2,3}\times\mathbb{R}$, where $N_{5,2,3}$ is the free nilpotent Lie algebra of step 3 and 2 generators.
    \subsection*{$\bf N_{6,2,9\,a}$ non-stratifiable  }

   The following Lie algebra is denoted as $N_{6,2,9\,a}$ by Gong in \cite{Gong_Thesis}, as $L_{6,24(-1)}$ by de Graaf in \cite{deGraaf_classification}, as $(18)$ by Del Barco in \cite{delBarco}, and as $\mathcal{G}_{6,5}$ by Magnin in \cite{magnin}.
   
   The non-trivial brackets are the following:
\begin{equation*}
\begin{aligned}
   &[X_1, X_2] = X_3\,,\, [X_1, X_3] = X_5\,,\\
   [X_2, X_3] &= -X_6\,,\, [X_2, X_4] = X_5\,,\,[X_1,X_4]=X_6\,.
\end{aligned}
\end{equation*}
This is a nilpotent Lie algebra of rank 3 and step 3 that is positively gradable, yet not stratifiable. The Lie brackets can be pictured with the diagram:

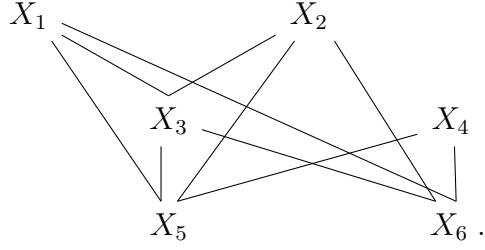
\begin{center}
 
	\begin{tikzcd}[end anchor=north]
	X_1 \ar[dr, no head]\ar[ddrrr, no head, end anchor={[xshift=-3.2ex]north east},start anchor={[yshift=1.3ex]south east}] \ar[ddr, no head,end anchor={[xshift=-3ex]north east} ]& & X_2\ar[ddl, no head, end anchor={[xshift=-1.8ex]north east}]\ar[dl, no head]\ar[ddr, no head,end anchor={[xshift=-4.7ex]north east}] & \\
	& X_3 \ar[drr, no head,end anchor={[xshift=-4.7ex]north east}]\ar[d, no head,end anchor={[xshift=-3.ex]north east},start anchor={[xshift=-3ex]south east} ]& & X_4\;\,\ar[d, no head, end anchor={[xshift=-3.2ex]north east},start anchor={[xshift=-3.2ex]south east}]\ar[dll, no head, end anchor={[xshift=-1.8ex]north east}] \\
	& X_5 & & X_6 \;.
	\end{tikzcd}

\end{center}

 


 
The composition law \eqref{group law in G} of $N_{6,2,9\,a}$ is given by:

\begin{itemize}
    \item $z_1=x_1+y_1$;
    \item $z_2=x_2+y_2$;
    \item $z_3=x_3+y_3+\frac{1}{2}(x_1y_2-x_2y_1)$;
    \item $z_4=x_4+y_4$;
    \item $z_5=x_5+y_5+\frac{1}{2}(x_1y_3-x_3y_1+x_2y_4-x_4y_2)+\frac{1}{12}(x_1-y_1)(x_1y_2-x_2y_1)$;
    \item $z_6=x_6+y_6+\frac{1}{2}(x_1y_4-x_4y_1-x_2y_3+x_3y_2)+\frac{1}{12}(y_2-x_2)(x_1y_2-x_2y_1)$.
\end{itemize}


Since
\begin{eqnarray*}
 \mathrm{d}(L_\mathbf{x})_\mathbf{0}= \left[\begin{matrix} 
   1 & 0  & 0 & 0 & 0 & 0\\
   0  & 1 & 0 & 0 & 0 &0\\
   -\frac{x_2}{2} & \frac{x_1}{2} & 1 & 0 & 0 & 0\\
  0 & 0 & 0 & 1& 0 & 0\\
  -\frac{x_1x_2}{12}-\frac{x_3}{2} & \frac{x_1^2}{12}-\frac{x_4}{2} & \frac{x_1}{2}& \frac{x_2}{2}& 1 & 0\\
  \frac{x_2^2}{12}-\frac{x_4}{2} & \frac{x_3}{2}-\frac{x_1x_2}{12} & -\frac{x_2}{2}  & \frac{x_1}{2}& 0 & 1
   \end{matrix}\right]\,,
\end{eqnarray*}
the induced left-invariant vector fields \eqref{leftinvariant vf} are: 
\begin{itemize}
\item $X_1=\partial_{x_1}-\frac{x_2}{2}\partial_{x_3}-\big(\frac{x_3}{2}+\frac{x_1x_2}{12}\big)\partial_{x_5}+\big(\frac{x_2^2}{12}-\frac{x_4}{2}\big)\partial_{x_6}\,;$
\item $X_2=\partial_{x_2}+\frac{x_1}{2}\partial_{x_3}+\big(\frac{x_1^2}{12}-\frac{x_4}{2}\big)\partial_{x_5}+\big(\frac{x_3}{2}-\frac{x_1x_2}{12}\big)\partial_{x_6}\,;$
\item $X_3=\partial_{x_3}+\frac{x_1}{2}\partial_{x_5}-\frac{x_2}{2}\partial_{x_6}\,;$
\item $X_4=\partial_{x_4}+\frac{x_2}{2}\partial_{x_5}+\frac{x_1}{2}\partial_{x_6}\,;$
\item $X_5=\partial_{x_5}$;
\item $X_6=\partial_{x_6}$,
\end{itemize}
and the respective left-invariant 1-forms \eqref{leftinvariant form} are: 
\begin{itemize}
\item $\theta_1=dx_1$;
\item $\theta_2=dx_2$;
\item $\theta_3=dx_3-\frac{x_1}{2} dx_2+\frac{x_2}{2}dx_1$;
\item $\theta_4=dx_4$;
\item $\theta_5=dx_5-\frac{x_2}{2}dx_4-\frac{x_1}{2}dx_3+\big(\frac{x_1^2}{6}+\frac{x_4}{2}\big)dx_2+\big(\frac{x_3}{2}-\frac{x_1x_2}{6}\big)dx_1$;
\item $\theta_6=dx_6-\frac{x_1}{2}dx_4+\frac{x_2}{2}dx_3-\big(\frac{x_3}{2}+\frac{x_1x_2}{6}\big)dx_2+\big(\frac{x_4}{2}+\frac{x_2^2}{6}\big)dx_1$.
\end{itemize}

Finally, we have

\begin{eqnarray*}
 \mathrm{d}(R_\mathbf{x})_\mathbf{0}= \left[\begin{matrix} 
   1 & 0  & 0 & 0 & 0 & 0\\
   0  & 1 & 0 & 0 & 0 &0\\
   \frac{x_2}{2} & -\frac{x_1}{2} & 1 & 0 & 0 & 0\\
  0 & 0 & 0 & 1& 0 & 0\\
  \frac{x_3}{2}-\frac{x_1x_2}{12} & \frac{x_1^2}{12}+\frac{x_4}{2} & -\frac{x_1}{2}& -\frac{x_2}{2}& 1 & 0\\
  \frac{x_2^2}{12}+\frac{x_4}{2} & -\frac{x_3}{2}-\frac{x_1x_2}{12} & \frac{x_2}{2}  & -\frac{x_1}{2}& 0 & 1
   \end{matrix}\right]\,.
\end{eqnarray*}

\subsubsection{Grading, polarizations of maximal dimension, and asymptotic cone}
 The Lie algebra $N_{6,2,9\, a}$
   is not stratifiable, but it is gradable as
   \begin{eqnarray*}
   V_1 &= & \span \{  X_1,X_2 \}\,,\\ V_2 &= & \span \{  X_3,X_4 \}\,,\\
   V_3 &= & \span \{ X_5,X_6 \}\,.
   \end{eqnarray*}

 Every   complementary subspace $\Delta$ to the  derived subalgebra is spanned by
      $X_1 +   u_1^3 X_3+ u_1^5 X_5+ u_1^6 X_6 $,
      $X_2 +   u_2^3 X_3+   u_2^5 X_5+   u_2^6 X_6$, and
       $X_4 +   u_4^3 X_3 +   u_4^5 X_5+   u_4^6 X_6$.
   We claim that in this Lie algebra every two complementary subspaces to the  derived subalgebra as in \eqref{complementary_Delta} differ by an automoprhism. Indeed, when $u_i^j$ varies, 
      the matrix
$$    \left[  \begin{array}{cccccc}
     1 &  0 &  0 &  0 &  0 &  0  \\
 0 &  1 &  0 &  0 &  0 &  0\\
 u_1^3 & u_2^3 & 1 &  u_4^3 & 0 &  0\\
 0 & 0 & 0 & 1 & 0 & 0\\
 u_1^5 &  u_2^5 & u_2^3 &  u_4^5 &  1 &  u_4^3\\
 u_1^6 &  u_2^6 & u_1^3 &  u_4^6 & -u_4^3 &  1\\
      \end{array}\right],$$
is a Lie algebra automorphism and sends the complementary subspace $\span\lbrace X_1,X_2,X_4\rbrace$ to
      $\span\{X_1 +   u_1^3 X_3+ u_1^5 X_5+ u_1^6 X_6, X_2 +   u_2^3 X_3+   u_2^5 X_5+   u_2^6 X_6,X_4 +   u_4^3 X_3 +   u_4^5 X_5+   u_4^6 X_6\}$, which is an arbitrary one of maximal dimension. In particular, every $\Delta$ as in \eqref{complementary_Delta} gives maximal dimension. The tangent cone of each of such polarizations has Lie algebra isomorphic to $N_{6,3,6}$, see page \pageref{N636}.

 The asymptotic cone 
 of the Lie group with Lie algebra
 $N_{6,2,9\,a}$ has Lie algebra isomorphic to $N_{5,2,3}\times\mathbb{R}$, where $N_{5,2,3}$ is the free nilpotent Lie algebra of step 3 and 2 generators.

    \subsection*{$\bf N_{6,2,10}$ non-stratifiable }


   The following Lie algebra is denoted as $N_{6,2,10}$ by Gong in \cite{Gong_Thesis}, as $L_{6,24(0)}$ by de Graaf in \cite{deGraaf_classification}, as $(19)$ by Del Barco in \cite{delBarco}, and as $\mathcal{G}_{6,8}$ by Magnin in \cite{magnin}.
   
   The non-trivial brackets are the following:
\begin{equation*}
   [X_1, X_2] = X_3\,,\, [X_1, X_3] = X_5\,,\, [X_2, X_3] = X_6\,,\, [X_2, X_4] = X_5\,.
\end{equation*}
This is a nilpotent Lie algebra of rank 3 and step 3 that is positively gradable, yet not stratifiable. The Lie brackets can be pictured with the diagram:
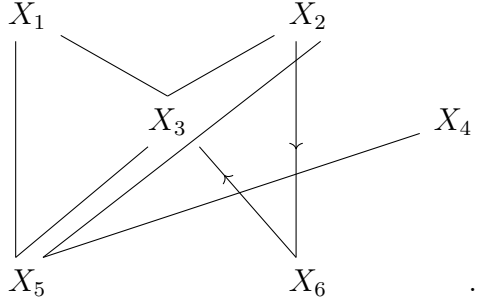
\begin{center}
 
	\begin{tikzcd}[end anchor=north]
	X_1\ar[dr, no head]\ar[ddd, no head,end anchor={[xshift=-3.3ex]north east},start anchor={[xshift=-3.3ex]south east}]& & X_2\ar[dl, no head]\ar[ddd,->-=.5, no head,end anchor={[xshift=-3.3ex]north east},start anchor={[xshift=-3.3ex]south east}] \ar[dddll, no head,start anchor={[xshift=-1.5ex]south east}, end anchor={[xshift=-1.3ex]north east}]& \\
	& X_3\ar[ddl, no head, end anchor={[xshift=-3.3ex]north east}]\ar[ddr, -<-=.3, no head,end anchor={[xshift=-3.3ex]north east}]&  &X_4\ar[ddlll, no head,end anchor={[xshift=-1.3ex]north east} ]\\
	& & &\\
	X_5& & X_6 &\quad\; .
	\end{tikzcd}

\end{center}

 
The composition law \eqref{group law in G} of $N_{6,2,10}$ is given by:

\begin{itemize}
    \item $z_1=x_1+y_1$;
    \item $z_2=x_2+y_2$;
    \item $z_3=x_3+y_3+\frac{1}{2}(x_1y_2-x_2y_1)$;
    \item $z_4=x_4+y_4$;
    \item $z_5=x_5+y_5+\frac{1}{2}(x_1y_3-x_3y_1+x_2y_4-x_4y_2)+\frac{1}{12}(x_1-y_1)(x_1y_2-x_2y_1)$;
    \item $z_6=x_6+y_6+\frac{1}{2}(x_2y_3-x_3y_2)+\frac{1}{12}(x_2-y_2)(x_1y_2-x_2y_1)$.
\end{itemize}



Since
\begin{eqnarray*}
 \mathrm{d}(L_\mathbf{x})_\mathbf{0}= \left[\begin{matrix} 
   1 & 0  & 0 & 0 & 0 & 0\\
   0  & 1 & 0 & 0 & 0 &0\\
   -\frac{x_2}{2} & \frac{x_1}{2} & 1 & 0 & 0 & 0\\
  0 & 0 & 0 & 1& 0 & 0\\
  -\frac{x_1x_2}{12}-\frac{x_3}{2} & \frac{x_1^2}{12}-\frac{x_4}{2} & \frac{x_1}{2}& \frac{x_2}{2}& 1 & 0\\
  -\frac{x_2^2}{12} & \frac{x_1x_2}{12}-\frac{x_3}{2} & \frac{x_2}{2}  & 0& 0 & 1
   \end{matrix}\right]\,,
\end{eqnarray*}
the induced left-invariant vector fields \eqref{leftinvariant vf} are: 
\begin{itemize}
\item $X_1=\partial_{x_1} -\frac{x_2}{2}\partial_{x_3}-\big(\frac{x_3}{2}+\frac{x_1x_2}{12}\big)\partial_{x_5}-\frac{x_2^2}{12}\partial_{x_6}\,;$
\item $X_2=\partial_{x_2}+\frac{x_1}{2}\partial_{x_3}+\big(\frac{x_1^2}{12}-\frac{x_4}{2}\big)\partial_{x_5}+\big(\frac{x_1x_2}{12}-\frac{x_3}{2}\big)\partial_{x_6}\,;$
\item $X_3=\partial_{x_3}+\frac{x_1}{2}\partial_{x_5}+\frac{x_2}{2}\partial_{x_6}\,;$
\item $X_4=\partial_{x_4}+\frac{x_2}{2}\partial_{x_5}\,;$
\item $X_5=\partial_{x_5}\,;$
\item $X_6=\partial_{x_6}$,
\end{itemize}
and the respective left-invariant 1-forms \eqref{leftinvariant form} are:
\begin{itemize}
\item $\theta_1=dx_1$;
\item $\theta_2=dx_2$;
\item $\theta_3=dx_3-\frac{x_1}{2} dx_2+\frac{x_2}{2}dx_1$;
\item $\theta_4=dx_4$;
\item $\theta_5=dx_5-\frac{x_2}{2}dx_4-\frac{x_1}{2}dx_3+\big(\frac{x_1^2}{6}+\frac{x_4}{2}\big)dx_2+\big(\frac{x_3}{2}-\frac{x_1x_2}{6}\big)dx_1$;
\item $\theta_6=dx_6-\frac{x_2}{2}dx_3+\big(\frac{x_3}{2}+\frac{x_1x_2}{6}\big)dx_2-\frac{x_2^2}{6}dx_1$.
\end{itemize}

Finally, we have
\begin{eqnarray*}
 \mathrm{d}(R_\mathbf{x})_\mathbf{0}= \left[\begin{matrix} 
   1 & 0  & 0 & 0 & 0 & 0\\
   0  & 1 & 0 & 0 & 0 &0\\
   \frac{x_2}{2} & -\frac{x_1}{2} & 1 & 0 & 0 & 0\\
  0 & 0 & 0 & 1& 0 & 0\\
  \frac{x_3}{2}-\frac{x_1x_2}{12} & \frac{x_1^2}{12}+\frac{x_4}{2} & -\frac{x_1}{2}& -\frac{x_2}{2}& 1 & 0\\
  -\frac{x_2^2}{12} & \frac{x_1x_2}{12}+\frac{x_3}{2} & -\frac{x_2}{2}  & 0& 0 & 1
   \end{matrix}\right]\,.
\end{eqnarray*}

\subsubsection{Grading, polarizations of maximal dimension, and asymptotic cone}
 The Lie algebra $N_{6,2,10}$
   is not stratifiable, but it is gradable as
   \begin{eqnarray*}
   V_1 &= & \span \{  X_1,X_2 \}\,,\\ V_2 &= & \span \{  X_3,X_4 \}\,,\\
   V_3 &= & \span \{ X_5,X_6 \}\,.
   \end{eqnarray*}
   
   Every   complementary subspace $\Delta$ to the  derived subalgebra is spanned by
      $X_1 +u_1^3X_3+u_1^5X_5+u_1^6X_6 $,
      $X_2 +u_2^3X_3+ u_2^5X_5+u_2^6X_6 $, and
       $X_4 +u_4^3X_3+ u^5_4X_5+u_4^6X_6 $.
       Such a polarization gives maximal Hausdorff dimension if and only if $u^3_4=0$.
      We claim that every two polarizations giving maximal Hausdorff dimension    differ by an automorphism. Indeed, when $u_i^j$ varies, 
      the matrix
$$    \left[  \begin{array}{cccccc}
     1 &  0 &  0 &  0 & 0& 0 \\
 0 &  1 &  0 &  0 &  0& 0 \\
 u_1^3 & u_2^3 &  1 &  0 &  0 & 0\\
 0 &  0 & 0 & 1& 0& 0\\
 u_1^5 &  u_2^5 &  u_2^3 &  u_4^5 &  1& 0 \\
 u_1^6 & u_2^6 & -u_1^3&  u_4^6 & 0 & 1
      \end{array}\right]$$
      is a Lie algebra automorphism and sends the complementary subspace $\span\lbrace X_1,X_2,X_4\rbrace$ to
      $\span\{X_1+u_1^3X_3+u_1^5X_5+u_1^6X_6,X_2+u_2^3X_3+u_2^5X_5+u_2^6X_6, X_4+u^5_4X_4+u_4^6X_6\}$, which is an arbitrary one of maximal dimension.  
The tangent cone of each  of such polarizations has Lie algebra isomorphic to $N_{6,3,4}$, see page \pageref{N634}.

 The asymptotic cone 
 of the Lie group with Lie algebra
 $N_{6,2,10}$ has Lie algebra isomorphic to $N_{5,2,3}\times\mathbb{R}$, where $N_{5,2,3}$ is the free nilpotent Lie algebra of step 3 and 2 generators.

   \subsection*{$\bf N_{6,3,1}$  } 

	The following Lie algebra is denoted as $N_{6,3,1}$ by Gong in \cite{Gong_Thesis}, as $L_{6,19(-1)}$ by de Graaf in \cite{deGraaf_classification}, as $(15)$ by Del Barco in \cite{delBarco}, and as $\mathcal{G}_{6,9}$ by Magnin in \cite{magnin}.
	
	The non-trivial brackets are the following:
\begin{equation*}\label{N631}
   [X_1, X_i] = X_{i+2}\,,\, i = 2, 3\,,\, [X_2, X_5] = [X_3, X_4] = X_6\,.
\end{equation*}
This is a nilpotent Lie algebra of rank 3 and step 3 that is stratifiable. The Lie brackets can be pictured with the diagram:
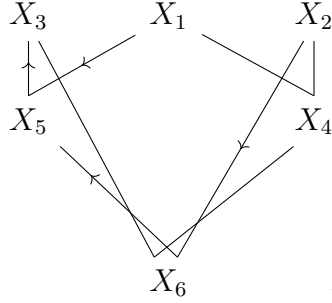
\begin{center}
 
	\begin{tikzcd}[end anchor=north]
	X_3\ar[d, no head, -<-=.5]\ar[dddr, no head,  end anchor={[xshift=-3.5ex]north east},start anchor={[xshift=-1.7ex]south east}] &X_1\ar[dr, no head,]\ar[no head, dl, ->-=.5]  &X_2\ar[no head, d]\ar[dddl, no head, ->-=.5, end anchor={[xshift=-1.8ex]north east}, start anchor={[xshift=-3.2ex]south east}]\\
	 X_5\ar[ddr, no head, -<-=.3, end anchor={[xshift=-1.8ex]north east}]& & X_4\ar[ddl, no head,  end anchor={[xshift=-3.5ex]north east}]\\
	& & \\
	 &X_6 &\quad\; .
	\end{tikzcd}
   
\end{center}

 
 The composition law \eqref{group law in G} of $N_{6,3,1}$ is given by:

\begin{itemize}
    \item $z_1=x_1+y_1$;
    \item $z_2=x_2+y_2$;
    \item $z_3=x_3+y_3$;
    \item $z_4=x_4+y_4+\frac{1}{2}(x_1y_2-x_2y_1)$;
    \item $z_5=x_5+y_5+\frac{1}{2}(x_1y_3-x_3y_1)$;
    \item $z_6=x_6+y_6+\frac{1}{2}(x_2y_5-x_5y_2+x_3y_4-x_4y_3)+\frac{1}{12}(x_3-y_3)(x_1y_2-x_2y_1)\salto+\frac{1}{12}(x_2-y_2)(x_1y_3-x_3y_1)$.
\end{itemize}



Since
\begin{eqnarray*}
  \mathrm{d}(L_\mathbf{x})_\mathbf{0}=\left[\begin{matrix} 
   1 & 0  & 0 & 0 & 0 & 0\\
   0  & 1 & 0 & 0 & 0 &0\\
   0 & 0 & 1 & 0 & 0 & 0\\
  -\frac{x_2}{2} & \frac{x_1}{2} & 0 & 1& 0 & 0\\
  -\frac{x_3}{2} & 0 & \frac{x_1}{2} & 0& 1 & 0\\
 - \frac{x_2x_3}{6} & \frac{x_1x_3}{12}-\frac{x_5}{2} & \frac{x_1x_2}{12}-\frac{x_4}{2} & \frac{x_3}{2}& \frac{x_2}{2} & 1
   \end{matrix}\right]\,,
\end{eqnarray*}
the induced left-invariant vector fields \eqref{leftinvariant vf} are: 
\begin{itemize}
\item $X_1=\partial_{x_1 }-\frac{x_2}{2}\partial_{x_4}-\frac{x_3}{2}\partial_{x_5}-\frac{x_2x_3}{6}\partial_{x_6}\,;$
\item $X_2=\partial_{x_2}+\frac{x_1}{2}\partial_{x_4}+\big(\frac{x_1x_3}{12}-\frac{x_5}{2}\big)\partial_{x_6}\,;$
\item $X_3=\partial_{x_3}+\frac{x_1}{2}\partial_{x_5}+\big(\frac{x_1x_2}{12}-\frac{x_4}{2}\big)\partial_{x_6}\,;$
\item $X_4=\partial_{x_4}+\frac{x_3}{2}\partial_{x_6}\,;$
\item $X_5=\partial_{x_5}+\frac{x_2}{2}\partial_{x_6}$;
\item $X_6=\partial_{x_6}$,
\end{itemize}
and the respective left-invariant 1-forms \eqref{leftinvariant form} are: 
\begin{itemize}
\item $\theta_1=dx_1$;
\item $\theta_2=dx_2$;
\item $\theta_3=dx_3$;
\item $\theta_4=dx_4-\frac{x_1}{2} dx_2+\frac{x_2}{2}dx_1$;
\item $\theta_5=dx_5-\frac{x_1}{2}dx_3+\frac{x_3}{2}dx_1$;
\item $\theta_6=dx_6-\frac{x_2}{2}dx_5-\frac{x_3}{2}dx_4+\big(\frac{x_4}{2}+\frac{x_1x_2}{6}\big)dx_3+\big(\frac{x_5}{2}+\frac{x_1x_3}{6}\big)dx_2-\frac{x_2x_3}{3}dx_1$.
\end{itemize}

Finally, we have
\begin{eqnarray*}
  \mathrm{d}(R_\mathbf{x})_\mathbf{0}=\left[\begin{matrix} 
   1 & 0  & 0 & 0 & 0 & 0\\
   0  & 1 & 0 & 0 & 0 &0\\
   0 & 0 & 1 & 0 & 0 & 0\\
  \frac{x_2}{2} & -\frac{x_1}{2} & 0 & 1& 0 & 0\\
  \frac{x_3}{2} & 0 & -\frac{x_1}{2} & 0& 1 & 0\\
 - \frac{x_2x_3}{6} & \frac{x_1x_3}{12}+\frac{x_5}{2} & \frac{x_1x_2}{12}+\frac{x_4}{2} & -\frac{x_3}{2}& -\frac{x_2}{2} & 1
   \end{matrix}\right]\,.
\end{eqnarray*}
   \subsection*{$\bf N_{6,3,1\,a}$  } 
	The following Lie algebra is denoted as $N_{6,3,1\,a}$ by Gong in \cite{Gong_Thesis}, as $L_{6,19(1)}$ by de Graaf in \cite{deGraaf_classification}, as $(16)$ by Del Barco in \cite{delBarco}, and as $\mathcal{G}_{6,9}$ by Magnin in \cite{magnin}.
	
	The non-trivial brackets are the following:
\begin{equation*}
   [X_1, X_i] = X_{i+2}\,,\, i = 2, 3\,,\, [X_2, X_4] = [X_3, X_5] = X_6\,.
\end{equation*}
This is a nilpotent Lie algebra of rank 3 and step 3 that is stratifiable. The Lie brackets can be pictured with the diagram:
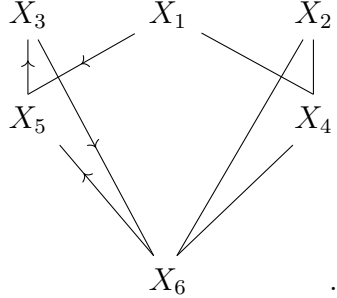
\begin{center}
 
	\begin{tikzcd}[end anchor=north]
	X_3\ar[d, no head, -<-=.5]\ar[dddr, no head,->-=.5,  end anchor={[xshift=-3.5ex]north east},start anchor={[xshift=-1.7ex]south east}] &X_1\ar[dr, no head,]\ar[no head, dl, ->-=.5]  &X_2\ar[no head, d]\ar[dddl, no head,  end anchor={[xshift=-1.8ex]north east}, start anchor={[xshift=-3.2ex]south east}]\\
	 X_5\ar[ddr, no head, -<-=.3, end anchor={[xshift=-3.5ex]north east}]& & X_4\ar[ddl, no head, end anchor={[xshift=-1.8ex]north east}]\\
	& & \\
	 &X_6 & \quad\;.
	\end{tikzcd}
   
\end{center}

 
The composition law \eqref{group law in G} of $N_{6,3,1\,a}$ is given by:

\begin{itemize}
    \item $z_1=x_1+y_1$;
    \item $z_2=x_2+y_2$;
    \item $z_3=x_3+y_3$;
    \item $z_4=x_4+y_4+\frac{1}{2}(x_1y_2-x_2y_1)$;
    \item $z_5=x_5+y_5+\frac{1}{2}(x_1y_3-x_3y_1)$;
    \item $z_6=x_6+y_6+\frac{1}{2}(x_2y_4-x_4y_2+x_3y_5-x_5y_3)+\frac{1}{12}(x_2-y_2)(x_1y_2-x_2y_1)\salto+\frac{1}{12}(x_3-y_3)(x_1y_3-x_3y_1)$.
\end{itemize}



Since
\begin{eqnarray*}
  \mathrm{d}(L_\mathbf{x})_\mathbf{0}=\left[\begin{matrix} 
   1 & 0  & 0 & 0 & 0 & 0\\
   0  & 1 & 0 & 0 & 0 &0\\
   0 & 0 & 1 & 0 & 0 & 0\\
  -\frac{x_2}{2} & \frac{x_1}{2} & 0 & 1& 0 & 0\\
  -\frac{x_3}{2} & 0 & \frac{x_1}{2} & 0& 1 & 0\\
 - \frac{x_2^2+x_3^2}{12} & \frac{x_1x_2}{12}-\frac{x_4}{2} & \frac{x_1x_3}{12}-\frac{x_5}{2} & \frac{x_2}{2}& \frac{x_3}{2} & 1
   \end{matrix}\right]\,,
\end{eqnarray*}
the induced left-invariant vector fields \eqref{leftinvariant vf} are: 
\begin{itemize}
\item $X_1=\partial_{x_1} -\frac{x_2}{2}\partial_{x_4}-\frac{x_3}{2}\partial_{x_5}-\frac{x_2^2+x_3^2}{12}\partial_{x_6}\,;$
\item $X_2=\partial_{x_2}+\frac{x_1}{2}\partial_{x_4}+\big(\frac{x_1x_2}{12}-\frac{x_4}{2}\big)\partial_{x_6}\,;$
\item $X_3=\partial_{x_3}+\frac{x_1}{2}\partial_{x_5}+\big(\frac{x_1x_3}{12}-\frac{x_5}{2}\big)\partial_{x_6}\,;$
\item $X_4=\partial_{x_4}+\frac{x_2}{2}\partial_{x_6}\,;$
\item $X_5=\partial_{x_5}+\frac{x_3}{2}\partial_{x_6}$;
\item $X_6=\partial_{x_6}$,
\end{itemize}
and the respective left-invariant 1-forms \eqref{leftinvariant form} are: 
\begin{itemize}
\item $\theta_1=dx_1$;
\item $\theta_2=dx_2$;
\item $\theta_3=dx_3$;
\item $\theta_4=dx_4-\frac{x_1}{2} dx_2+\frac{x_2}{2}dx_1$;
\item $\theta_5=dx_5-\frac{x_1}{2}dx_3+\frac{x_3}{2}dx_1$;
\item $\theta_6=dx_6-\frac{x_3}{2}dx_5-\frac{x_2}{2}dx_4+\big(\frac{x_5}{2}+\frac{x_1x_3}{6}\big)dx_3+\big(\frac{x_4}{2}+\frac{x_1x_2}{6}\big)dx_2-\frac{x_2^2+x_3^2}{6}dx_1$.
\end{itemize}

Finally, we have
\begin{eqnarray*}
  \mathrm{d}(R_\mathbf{x})_\mathbf{0}=\left[\begin{matrix} 
   1 & 0  & 0 & 0 & 0 & 0\\
   0  & 1 & 0 & 0 & 0 &0\\
   0 & 0 & 1 & 0 & 0 & 0\\
  \frac{x_2}{2} & -\frac{x_1}{2} & 0 & 1& 0 & 0\\
  \frac{x_3}{2} & 0 & -\frac{x_1}{2} & 0& 1 & 0\\
 - \frac{x_2^2+x_3^2}{12} & \frac{x_1x_2}{12}+\frac{x_4}{2} & \frac{x_1x_3}{12}+\frac{x_5}{2} & -\frac{x_2}{2}& -\frac{x_3}{2} & 1
   \end{matrix}\right]\,.
\end{eqnarray*} 
   \subsection*{$\bf N_{6,3,2}$ non-stratifiable} 
  
    The following Lie algebra is denoted as $N_{6,3,2}$ by Gong in \cite{Gong_Thesis}, as $L_{6,10}$ by de Graaf in \cite{deGraaf_classification}, as $(25)$ by Del Barco in \cite{delBarco}, and as $\mathcal{G}_{6,2}$ by Magnin in \cite{magnin}.
    
    The non-trivial brackets are the following:
\begin{equation*}
   [X_1, X_2] = X_3\,,\, [X_1, X_3] = X_6\,,\, [X_4, X_5] = X_6\,.
\end{equation*}
This is a nilpotent Lie algebra 
of rank 4 and step 3
that is positively gradable, yet not stratifiable. The Lie brackets can be pictured with the diagram:
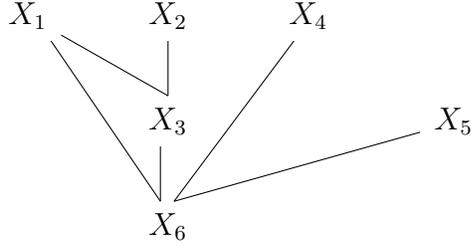
\begin{center}
 
		\begin{tikzcd}[end anchor=north]
	X_1 \ar[dr, no head]\ar[ddr, no head, end anchor={[xshift=-3ex]north east}]&X_2\ar[d, no head] & X_4\ar[ddl, no head, end anchor={[xshift=-2ex]north east}] & \\
	& X_3 \ar[d,  no head, start anchor={[xshift=-3ex]south east},end anchor={[xshift=-3ex]north east}]& &X_5\ar[dll, no head,  end anchor={[xshift=-2ex]north east}]\\
	&  X_6 & & \quad\;.
	\end{tikzcd}

\end{center}


 
 The composition law \eqref{group law in G} of $N_{6,3,2}$ is given by:

\begin{itemize}
    \item $z_1=x_1+y_1$;
    \item $z_2=x_2+y_2$;
    \item $z_3=x_3+y_3+\frac{1}{2}(x_1y_2-x_2y_1)$;
    \item $z_4=x_4+y_4$;
    \item $z_5=x_5+y_5$;
    \item $z_6=x_6+y_6+\frac{1}{2}(x_1y_3-x_3y_1+x_4y_5-x_5y_4)+\frac{1}{12}(x_1-y_1)(x_1y_2-x_2y_1)$.
\end{itemize}



Since
\begin{eqnarray*}
 \mathrm{d}(L_\mathbf{x})_\mathbf{0}= \left[\begin{matrix} 
   1 & 0  & 0 & 0 & 0 & 0\\
   0  & 1 & 0 & 0 & 0 &0\\
   -\frac{x_2}{2} & \frac{x_1}{2} & 1 & 0 & 0 & 0\\
  0 & 0 & 0 & 1& 0 & 0\\
  0 & 0 & 0& 0& 1 & 0\\
  -\frac{x_1x_2}{12}-\frac{x_3}{2} & \frac{x_1^2}{12} & \frac{x_1}{2}  & -\frac{x_5}{2}& \frac{x_4}{2} & 1
   \end{matrix}\right]\,,
\end{eqnarray*}
the induced left-invariant vector fields \eqref{leftinvariant vf} are: 
\begin{itemize}
\item $X_1=\partial_{x_1} -\frac{x_2}{2}\partial_{x_3}-\big(\frac{x_1x_2}{12}+\frac{x_3}{2}\big)\partial_{x_6}\,;$
\item $X_2=\partial_{x_2}+\frac{x_1}{2}\partial_{x_3}+\frac{x_1^2}{12}\partial_{x_6}\,;$
\item $X_3=\partial_{x_3}+\frac{x_1}{2}\partial_{x_6}\,;$
\item $X_4=\partial_{x_4}-\frac{x_5}{2}\partial_{x_6}\,;$
\item $X_5=\partial_{x_5}+\frac{x_4}{2}\partial_{x_6}$;
\item $X_6=\partial_{x_6}$,
\end{itemize}
and the respective left-invariant 1-forms \eqref{leftinvariant form} are: 
\begin{itemize}
\item $\theta_1=dx_1$;
\item $\theta_2=dx_2$;
\item $\theta_3=dx_3-\frac{x_1}{2} dx_2+\frac{x_2}{2}dx_1$;
\item $\theta_4=dx_4$;
\item $\theta_5=dx_5$;
\item $\theta_6=dx_6-\frac{x_4}{2}dx_5+\frac{x_5}{2}dx_4-\frac{x_1}{2}dx_3+\frac{x_1^2}{6}dx_2+\big(\frac{x_3}{2}-\frac{x_1x_2}{6}\big)dx_1$.
\end{itemize}

Finally, we have
\begin{eqnarray*}
 \mathrm{d}(R_\mathbf{x})_\mathbf{0}= \left[\begin{matrix} 
   1 & 0  & 0 & 0 & 0 & 0\\
   0  & 1 & 0 & 0 & 0 &0\\
   \frac{x_2}{2} & -\frac{x_1}{2} & 1 & 0 & 0 & 0\\
  0 & 0 & 0 & 1& 0 & 0\\
  0 & 0 & 0& 0& 1 & 0\\
  \frac{x_3}{2}-\frac{x_1x_2}{12} & \frac{x_1^2}{12} & -\frac{x_1}{2}  & \frac{x_5}{2}& -\frac{x_4}{2} & 1
   \end{matrix}\right]\,.
\end{eqnarray*}

\subsubsection{Grading, polarizations of maximal dimension, and asymptotic cone}
 The Lie algebra $N_{6,3,2}$
   is a not stratifiable, but it is gradable as
   \begin{eqnarray*}
   V_1 &= & \span \{  X_1, X_2, X_4 \}\,, \\
   V_2 &= & \span \{  X_3,  X_5 \}\,, \\
      V_3 &= & \span \{   X_6 \}\,. \\
   \end{eqnarray*}
   
    Every complementary subspace $\Delta$ to the  derived subalgebra is spanned by
      $X_1 +u_1^3X_3+u_1^6X_6 $,
      $X_2 +u_2^3X_3+u_2^6X_6 $,
       $X_4 +u_4^3X_3+ u_4^6X_6 $, and $X_5+u_5^3X_3+u_5^6X_6$.
      We claim that every two polarizations differ by an automorphism. Indeed, when $u_i^j$ varies, 
      the matrix
$$    \left[  \begin{array}{cccccc}
     1 &  0 &  0 &  0 &  0 &  0  \\
 0 &  1 &  0 &  0 &  0 &  0\\
 u_1^3 &  u_2^3 &  1 &  u_4^3 &  u_5^3 &  0\\
 -u_5^3 &  0 &  0 &  1 &  0 &  0\\
 u_4^3 &  0 &  0 &  0 &  1 &  0\\
 u_1^6 &  u_2^6 &  u_2^3 &  u_4^6 &  u_5^6 &  1\\
      \end{array}\right]$$
      is a Lie algebra automorphism and sends the complementary subspace $X_1,X_2,X_4,X_5\rbrace$ to
      $\span\{X_1+u_1^3X_3-u_5^3X_4+u_4^3X_5+u_1^6X_6,X_2+u_2^3X_3+u_2^6X_6, X_4+u^3_4X_3+u_4^6X_6,X_5+u_5^3X_3+u_5^6X_6\}$, which is an arbitrary one of maximal dimension.  
The tangent cone of each of such polarizations has Lie algebra isomorphic to $N_{3,2}\times N_{3,2}$, where $N_{3,2}$ is the first Heisenberg algebra.
 
 The asymptotic cone 
 of the Lie group with Lie algebra
 $N_{6,3,2}$ has Lie algebra isomorphic to $N_{4,2}\times\mathbb{R}^2$, where $N_{4,2}$ is the filiform algebra of step 3.
      
 
    \subsection*{$\bf N_{6,3,3}$  }
 	The following Lie algebra is denoted as $N_{6,3,3}$ by Gong in \cite{Gong_Thesis}, as $L_{6,19(0)}$ by de Graaf in \cite{deGraaf_classification}, as $(22)$ by Del Barco in \cite{delBarco}, and as $\mathcal{G}_{6,4}$ by Magnin in \cite{magnin}.
 	
 	The non-trivial brackets are the following:
\begin{equation*}\label{N633}
   [X_1, X_2] = X_3\,,\, [X_1, X_4] = X_6\,,\, [X_2, X_3] = X_5\,.
\end{equation*}
This is a nilpotent Lie algebra of rank 3 and step 3. The Lie brackets can be pictured with the diagram:
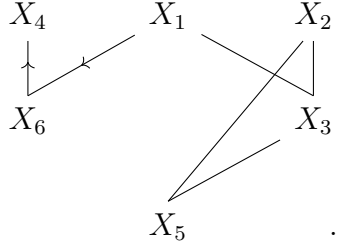
\begin{center}
 
	\begin{tikzcd}[end anchor=north]
	X_4\ar[d, no head, -<-=.5] &X_1\ar[dr, no head,]\ar[no head, dl, ->-=.5]  &X_2\ar[no head, d]\ar[ddl, no head,  start anchor={[xshift=-3.2ex]south east}]\\
	 X_6& & X_3\ar[dl, no head]\\
	 &X_5 & \quad\;.
	\end{tikzcd}
   
\end{center}

 
The composition law \eqref{group law in G} of $N_{6,3,3}$ is given by:

\begin{itemize}
    \item $z_1=x_1+y_1$;
    \item $z_2=x_2+y_2$;
    \item $z_3=x_3+y_3+\frac{1}{2}(x_1y_2-x_2y_1)$;
    \item $z_4=x_4+y_4$;
    \item $z_5=x_5+y_5+\frac{1}{2}(x_2y_3-x_3y_2)+\frac{1}{12}(x_2-y_2)(x_1y_2-x_2y_1)$;
    \item $z_6=x_6+y_6+\frac{1}{2}(x_1y_4-x_4y_1)$.
\end{itemize}



Since
\begin{eqnarray*}
 \mathrm{d}(L_\mathbf{x})_\mathbf{0}= \left[\begin{matrix} 
   1 & 0  & 0 & 0 & 0 & 0\\
   0  & 1 & 0 & 0 & 0 &0\\
   -\frac{x_2}{2}& \frac{x_1}{2} & 1 & 0 & 0 & 0\\
  0 & 0 & 0 & 1& 0 & 0\\
  -\frac{x_2^2}{12} & \frac{x_1x_2}{12}-\frac{x_3}{2} & \frac{x_2}{2} & 0& 1 & 0\\
  -\frac{x_4}{2} & 0 & 0 & \frac{x_1}{2}& 0 & 1
   \end{matrix}\right]\,,
\end{eqnarray*}
the induced left-invariant vector fields \eqref{leftinvariant vf} are: 
\begin{itemize}
\item $X_1=\partial_{x_1} -\frac{x_2}{2}\partial_{x_3}-\frac{x_2^2}{12}\partial_{x_5}-\frac{x_4}{2}\partial_{x_6}\,;$
\item $X_2=\partial_{x_2}+\frac{x_1}{2}\partial_{x_3}+\big(\frac{x_1x_2}{12}-\frac{x_3}{2}\big)\partial_{x_5}\,;$
\item $X_3=\partial_{x_3}+\frac{x_2}{2}\partial_{x_5}\,;$
\item $X_4=\partial_{x_4}+\frac{x_1}{2}\partial_{x_6}\,;$
\item $X_5=\partial_{x_5}$;
\item $X_6=\partial_{x_6}$,
\end{itemize}
and the respective left-invariant 1-forms \eqref{leftinvariant form} are:
\begin{itemize}
\item $\theta_1=dx_1$;
\item $\theta_2=dx_2$;
\item $\theta_3=dx_3-\frac{x_1}{2} dx_2+\frac{x_2}{2}dx_1$;
\item $\theta_4=dx_4$;
\item $\theta_5=dx_5-\frac{x_2}{2}dx_3+\big(\frac{x_3}{2}+\frac{x_1x_2}{6}\big)dx_2-\frac{x_2^2}{6}dx_1$;
\item $\theta_6=dx_6-\frac{x_1}{2}dx_4+\frac{x_4}{2}dx_1$.
\end{itemize}

Finally, we have
\begin{eqnarray*}
 \mathrm{d}(R_\mathbf{x})_\mathbf{0}= \left[\begin{matrix} 
   1 & 0  & 0 & 0 & 0 & 0\\
   0  & 1 & 0 & 0 & 0 &0\\
   \frac{x_2}{2}& -\frac{x_1}{2} & 1 & 0 & 0 & 0\\
  0 & 0 & 0 & 1& 0 & 0\\
  -\frac{x_2^2}{12} & \frac{x_1x_2}{12}+\frac{x_3}{2} & -\frac{x_2}{2} & 0& 1 & 0\\
  \frac{x_4}{2} & 0 & 0 & -\frac{x_1}{2}& 0 & 1
   \end{matrix}\right]\,.
\end{eqnarray*}
    \subsection*{$\bf N_{6,3,4}$  }
	The following Lie algebra is denoted as $N_{6,3,4}$ by Gong in \cite{Gong_Thesis}, as $L_{6,25}$ by de Graaf in \cite{deGraaf_classification}, as $(23)$ by Del Barco in \cite{delBarco}, and as $\mathcal{G}_{6,6}$ by Magnin in \cite{magnin}.
	
	The non-trivial brackets are the following:
\begin{equation*}\label{N634}
   [X_1, X_2] = X_3\,,\, [X_2, X_3] = X_5\,,\, [X_2, X_4] = X_6\,.
\end{equation*}
This is a nilpotent Lie algebra of rank 3 and step 3 that is stratifiable. The Lie brackets can be pictured with the diagram:
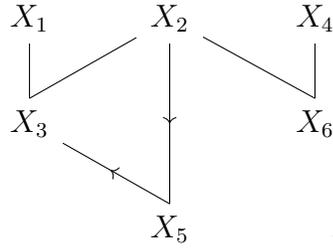
\begin{center}
 
	\begin{tikzcd}[end anchor=north]
X_1\ar[d, no head]  & X_2\ar[dl, no head] \ar[dr, no head]\ar[dd, no head, ->-=.5]& X_4\ar[d, no head] \\
 X_3 \ar[dr, no head, -<-=.5]& & X_6\\
 & X_5 & \quad\;.
	\end{tikzcd}
   
\end{center}

 
The composition law \eqref{group law in G} of $N_{6,3,4}$ is given by:

\begin{itemize}
    \item $z_1=x_1+y_1$;
    \item $z_2=x_2+y_2$;
    \item $z_3=x_3+y_3+\frac{1}{2}(x_1y_2-x_2y_1)$;
    \item $z_4=x_4+y_4$;
    \item $z_5=x_5+y_5+\frac{1}{2}(x_2y_3-x_3y_2)+\frac{1}{12}(x_2-y_2)(x_1y_2-x_2y_1)$;
    \item $z_6=x_6+y_6+\frac{1}{2}(x_2y_4-x_4y_2)$.
\end{itemize}

%


Since
\begin{eqnarray*}
 \mathrm{d}(L_\mathbf{x})_\mathbf{0}= \left[\begin{matrix} 
   1 & 0  & 0 & 0 & 0 & 0\\
   0  & 1 & 0 & 0 & 0 &0\\
   -\frac{x_2}{2}& \frac{x_1}{2} & 1 & 0 & 0 & 0\\
  0 & 0 & 0 & 1& 0 & 0\\
  -\frac{x_2^2}{12} & \frac{x_1x_2}{12}-\frac{x_3}{2} & \frac{x_2}{2} & 0& 1 & 0\\
  0 & -\frac{x_4}{2} & 0 & \frac{x_2}{2}& 0 & 1
   \end{matrix}\right]\,,
   \end{eqnarray*}
the induced left-invariant vector fields \eqref{leftinvariant vf} are: 
\begin{itemize}
\item $X_1=\partial_{x_1} -\frac{x_2}{2}\partial_{x_3} -\frac{x_2^2}{12}\partial_{x_5}\,$;
\item $X_2=\partial_{x_2}+\frac{x_1}{2}\partial_{x_3}+\big(\frac{x_1x_2}{12}-\frac{x_3}{2}\big)\partial_{x_5}-\frac{x_4}{2}\partial_{x_6}\,;$
\item $X_3=\partial_{x_3}+\frac{x_2}{2}\partial_{x_5}\,;$
\item $X_4=\partial_{x_4}+\frac{x_2}{2}\partial_{x_6}\,;$
\item $X_5=\partial_{x_5}$;
\item $X_6=\partial_{x_6}$,
\end{itemize}
and the respective left-invariant 1-forms \eqref{leftinvariant form} are: 
\begin{itemize}
\item $\theta_1=dx_1$;
\item $\theta_2=dx_2$;
\item $\theta_3=dx_3-\frac{x_1}{2} dx_2+\frac{x_2}{2}dx_1$;
\item $\theta_4=dx_4$;
\item $\theta_5=dx_5-\frac{x_2}{2}dx_3+\big(\frac{x_3}{2}+\frac{x_1x_2}{6}\big)dx_2-\frac{x_2^2}{6}dx_1$;
\item $\theta_6=dx_6-\frac{x_2}{2}dx_4+\frac{x_4}{2}dx_2$.
\end{itemize}

Finally, we have
\begin{eqnarray*}
 \mathrm{d}(R_\mathbf{x})_\mathbf{0}= \left[\begin{matrix} 
   1 & 0  & 0 & 0 & 0 & 0\\
   0  & 1 & 0 & 0 & 0 &0\\
   \frac{x_2}{2}& -\frac{x_1}{2} & 1 & 0 & 0 & 0\\
  0 & 0 & 0 & 1& 0 & 0\\
  -\frac{x_2^2}{12} & \frac{x_1x_2}{12}+\frac{x_3}{2} & -\frac{x_2}{2} & 0& 1 & 0\\
  0 & \frac{x_4}{2} & 0 & -\frac{x_2}{2}& 0 & 1
   \end{matrix}\right]\,.
   \end{eqnarray*}
   \subsection*{$\bf N_{6,3,5}$ } 

	The following Lie algebra is denoted as $N_{6,3,5}$ by Gong in \cite{Gong_Thesis}, as $L_{6,22(0)}$ by de Graaf in \cite{deGraaf_classification}, as $(29)$ by Del Barco in \cite{delBarco}, and as $\mathcal{G}_{6,1}$ by Magnin in \cite{magnin}.
	
	The non-trivial brackets are the following:
\begin{equation*}
   [X_1, X_2] = X_5\,,\, [X_1, X_4] = X_6\,,\, [X_2, X_3] = X_6\,.
\end{equation*}
This is a nilpotent Lie algebra of rank 4 and step 2 that is stratifiable. The Lie brackets can be pictured with the diagram:
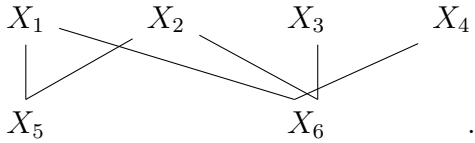
\begin{center}
 
	\begin{tikzcd}[end anchor=north]
	X_1\ar[d,no head]\ar[drr, no head, end anchor={[xshift=-3.3ex]north east }]& X_2\ar[dl, no head]\ar[dr, no head, end anchor={[xshift=-1.6ex]north east }] &X_3 \ar[d, no head, end anchor={[xshift=-1.6ex]north east },start anchor={[xshift=-1.6ex]south east }]& X_4\ar[dl, no head, end anchor={[xshift=-3.3ex]north east }]\\
	X_5 & & X_6 &\quad\;.
	\end{tikzcd}
   
\end{center}

 
The composition law \eqref{group law in G} of $N_{6,3,5}$ is given by:

\begin{itemize}
    \item $z_1=x_1+y_1$;
    \item $z_2=x_2+y_2$;
    \item $z_3=x_3+y_3$;
    \item $z_4=x_4+y_4$;
    \item $z_5=x_5+y_5+\frac{1}{2}(x_1y_2-x_2y_1)$;
    \item $z_6=x_6+y_6+\frac{1}{2}(x_1y_4-x_4y_1+x_2y_3-x_3y_2)$.
\end{itemize}

Since
\begin{eqnarray*}
 \mathrm{d}(L_\mathbf{x})_\mathbf{0}= \left[\begin{matrix} 
   1 & 0  & 0 & 0 & 0 & 0\\
   0  & 1 & 0 & 0 & 0 &0\\
   0& 0 & 1 & 0 & 0 & 0\\
  0 & 0 & 0 & 1& 0 & 0\\
  -\frac{x_2}{2} & \frac{x_1}{2} & 0 & 0& 1 & 0\\
  -\frac{x_4}{2} & -\frac{x_3}{2} & \frac{x_2}{2} & \frac{x_1}{2}& 0 & 1
   \end{matrix}\right]\,,
\end{eqnarray*}
the induced left-invariant vector fields \eqref{leftinvariant vf} are: 
\begin{itemize}
\item $X_1=\partial_{x_1} -\frac{x_2}{2}\partial_{x_5}-\frac{x_4}{2}\partial_{x_6}\,;$
\item $X_2=\partial_{x_2}+\frac{x_1}{2}\partial_{x_5}-\frac{x_3}{2}\partial_{x_6}\,;$
\item $X_3=\partial_{x_3}+\frac{x_2}{2}\partial_{x_6}\,;$
\item $X_4=\partial_{x_4}+\frac{x_1}{2}\partial_{x_6}\,;$
\item $X_5=\partial_{x_5}$;
\item $X_6=\partial_{x_6}$,
\end{itemize}
and the respective left-invariant 1-forms \eqref{leftinvariant form} are: 
\begin{itemize}
\item $\theta_1=dx_1$;
\item $\theta_2=dx_2$;
\item $\theta_3=dx_3$;
\item $\theta_4=dx_4$;
\item $\theta_5=dx_5-\frac{x_1}{2} dx_2+\frac{x_2}{2}dx_1$;
\item $\theta_6=dx_6-\frac{x_1}{2}dx_4-\frac{x_2}{2}dx_3+\frac{x_3}{2}dx_2+\frac{x_4}{2}dx_1$.
\end{itemize}

Finally, we have
\begin{eqnarray*}
 \mathrm{d}(R_\mathbf{x})_\mathbf{0}= \left[\begin{matrix} 
   1 & 0  & 0 & 0 & 0 & 0\\
   0  & 1 & 0 & 0 & 0 &0\\
   0& 0 & 1 & 0 & 0 & 0\\
  0 & 0 & 0 & 1& 0 & 0\\
  \frac{x_2}{2} & -\frac{x_1}{2} & 0 & 0& 1 & 0\\
  \frac{x_4}{2} & \frac{x_3}{2} & -\frac{x_2}{2} & -\frac{x_1}{2}& 0 & 1
   \end{matrix}\right]\,.
\end{eqnarray*}

   \subsection*{$\bf N_{6,3,6}$} 
	The following Lie algebra is denoted as $N_{6,3,6}$ by Gong in \cite{Gong_Thesis}, as $L_{6,26}$ by de Graaf in \cite{deGraaf_classification}, as $(24)$ by Del Barco in \cite{delBarco}, and as $\mathcal{G}_{6,3}$ by Magnin in \cite{magnin}.
	
	The non-trivial brackets are the following:
\begin{equation*}\label{N636}
   [X_1, X_2] = X_4\,,\, [X_1, X_3] = X_5\,,\, [X_2, X_3] = X_6\,.
\end{equation*}
This is a nilpotent Lie algebra of rank 3 and step 2 that is stratifiable, also known as the free Lie algebra of step 2 and 3 generators. The Lie brackets can be pictured with the diagram:
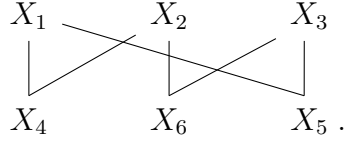
\begin{center}
 
	\begin{tikzcd}[end anchor=north]
	X_1\ar[d,no head]\ar[drr, no head, end anchor={[xshift=-4.1ex]north east }]& X_2\ar[d, no head]\ar[dl, no head] &X_3\;\;\ar[dl,no head] \ar[d, no head, end anchor={[xshift=-4.1ex]north east },start anchor={[xshift=-4.1ex]south east }]\\
	X_4 &X_6 & X_5 \;.
	\end{tikzcd}
\end{center}

 
The composition law \eqref{group law in G} of $N_{6,3,6}$ is given by:

\begin{itemize}
    \item $z_1=x_1+y_1$;
    \item $z_2=x_2+y_2$;
    \item $z_3=x_3+y_3$;
    \item $z_4=x_4+y_4+\frac{1}{2}(x_1y_2-x_2y_1)$;
    \item $z_5=x_5+y_5+\frac{1}{2}(x_1y_3-x_3y_1)$;
    \item $z_6=x_6+y_6+\frac{1}{2}(x_2y_3-x_3y_2)$.
\end{itemize}

Since
\begin{eqnarray*}
  \mathrm{d}(L_\mathbf{x})_\mathbf{0}=\left[\begin{matrix} 
   1 & 0  & 0 & 0 & 0 & 0\\
   0  & 1 & 0 & 0 & 0 &0\\
   0& 0 & 1 & 0 & 0 & 0\\
  -\frac{x_2}{2} & \frac{x_1}{2} & 0 & 1& 0 & 0\\
  -\frac{x_3}{2} & 0 & \frac{x_1}{2} & 0& 1 & 0\\
  0 & -\frac{x_3}{2} & \frac{x_2}{2} & 0& 0 & 1
   \end{matrix}\right]\,,
\end{eqnarray*}
the induced left-invariant vector fields \eqref{leftinvariant vf} are: 
\begin{itemize}
\item $X_1=\partial_{x_1 }-\frac{x_2}{2}\partial_{x_4}-\frac{x_3}{2}\partial_{x_5}\,;$
\item $X_2=\partial_{x_2}+\frac{x_1}{2}\partial_{x_4}-\frac{x_3}{2}\partial_{x_6}\,;$
\item $X_3=\partial_{x_3}+\frac{x_1}{2}\partial_{x_5}+\frac{x_2}{2}\partial_{x_6}\,;$
\item $X_4=\partial_{x_4}\,;$
\item $X_5=\partial_{x_5}$;
\item $X_6=\partial_{x_6}$,
\end{itemize}
and the respective left-invariant 1-forms \eqref{leftinvariant form} are: 
\begin{itemize}
\item $\theta_1=dx_1$;
\item $\theta_2=dx_2$;
\item $\theta_3=dx_3$;
\item $\theta_4=dx_4-\frac{x_1}{2} dx_2+\frac{x_2}{2}dx_1$;
\item $\theta_5=dx_5-\frac{x_1}{2} dx_3+\frac{x_3}{2}dx_1$;
\item $\theta_6=dx_6-\frac{x_2}{2}dx_3+\frac{x_3}{2}dx_2$.
\end{itemize}

Finally, we have
\begin{eqnarray*}
  \mathrm{d}(R_\mathbf{x})_\mathbf{0}=\left[\begin{matrix} 
   1 & 0  & 0 & 0 & 0 & 0\\
   0  & 1 & 0 & 0 & 0 &0\\
   0& 0 & 1 & 0 & 0 & 0\\
  \frac{x_2}{2} & -\frac{x_1}{2} & 0 & 1& 0 & 0\\
  \frac{x_3}{2} & 0 & -\frac{x_1}{2} & 0& 1 & 0\\
  0 & \frac{x_3}{2} & -\frac{x_2}{2} & 0& 0 & 1
   \end{matrix}\right]\,.
\end{eqnarray*}
   \subsection*{$\bf N_{6,4,4\,a}$ } 
	The following Lie algebra is denoted as $N_{6,4,4\,a}$ by Gong in \cite{Gong_Thesis}, as $L_{6,22(-1)}$ by de Graaf in \cite{deGraaf_classification}, and as $(28)$ by Del Barco in \cite{delBarco}. As a complex Lie algebra, 
	$N_{6,4,4\,a}$ is equivalent to  
	the decomposable Lie algebra $N_{3,2}\times N_{3,2}$, 
	 which is why this Lie algebra is not contained in the list produced by Magnin \cite{magnin}.
	
	The non-trivial brackets are the following:
\begin{equation*}
   [X_1, X_3] = X_5\,,\, [X_1, X_4] = X_6\,,\, [X_2, X_4] = X_5\,,\,[X_2,X_3]=-X_6\,.
\end{equation*}
This is a nilpotent Lie algebra rank 4 and step 2 that is stratifiable. The Lie brackets can be pictured with the diagram:
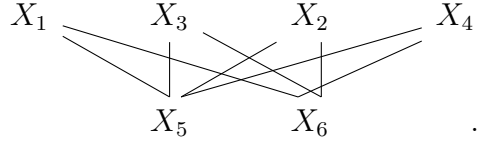
\begin{center}
 
	\begin{tikzcd}[end anchor=north]
	X_1\ar[dr,no head]\ar[drr, no head, end anchor={[xshift=-3.3ex]north east }]& X_3\ar[d, no head]\ar[dr, no head, end anchor={[xshift=-1.6ex]north east }] &X_2 \ar[dl, no head, end anchor={[xshift=-1.6ex]north east }]\ar[d, no head, end anchor={[xshift=-1.6ex]north east },start anchor={[xshift=-1.6ex]south east }]& X_4\ar[dl, no head, end anchor={[xshift=-3.3ex]north east }]\ar[dll,no head, end anchor={[xshift=-1.6ex]north east }]\\
	&X_5  & X_6 &\quad\;.
	\end{tikzcd}
   
\end{center}

 
The composition law \eqref{group law in G} of $N_{6,4,4\,a}$ is given by:

\begin{itemize}
    \item $z_1=x_1+y_1$;
    \item $z_2=x_2+y_2$;
    \item $z_3=x_3+y_3$;
    \item $z_4=x_4+y_4$;
    \item $z_5=x_5+y_5+\frac{1}{2}(x_1y_3-x_3y_1+x_2y_4-x_4y_2)$;
    \item $z_6=x_6+y_6+\frac{1}{2}(x_1y_4-x_4y_1-x_2y_3+x_3y_2)$.
\end{itemize}

Since
\begin{eqnarray*}
 \mathrm{d}(L_\mathbf{x})_\mathbf{0}= \left[\begin{matrix} 
   1 & 0  & 0 & 0 & 0 & 0\\
   0  & 1 & 0 & 0 & 0 &0\\
   0& 0 & 1 & 0 & 0 & 0\\
  0 & 0 & 0 & 1& 0 & 0\\
  -\frac{x_3}{2} & -\frac{x_4}{2} & \frac{x_1}{2} & \frac{x_2}{2}& 1 & 0\\
  -\frac{x_4}{2} & \frac{x_3}{2} & -\frac{x_2}{2} & \frac{x_1}{2}& 0 & 1
   \end{matrix}\right]\,,
\end{eqnarray*}
the induced left-invariant vector fields \eqref{leftinvariant vf} are: 
\begin{itemize}
\item $X_1=\partial_{x_1} -\frac{x_3}{2}\partial_{x_5}-\frac{x_4}{2}\partial_{x_6}\,;$
\item $X_2=\partial_{x_2}-\frac{x_4}{2}\partial_{x_5}+\frac{x_3}{2}\partial_{x_6}\,;$
\item $X_3=\partial_{x_3}+\frac{x_1}{2}\partial_{x_5}-\frac{x_2}{2}\partial_{x_6}\,;$
\item $X_4=\partial_{x_4}+\frac{x_2}{2}\partial_{x_5}+\frac{x_1}{2}\partial_{x_6}\,;$
\item $X_5=\partial_{x_5}$;
\item $X_6=\partial_{x_6}$,
\end{itemize}
and the respective left-invariant 1-forms \eqref{leftinvariant form} are: 
\begin{itemize}
\item $\theta_1=dx_1$;
\item $\theta_2=dx_2$;
\item $\theta_3=dx_3$;
\item $\theta_4=dx_4$;
\item $\theta_5=dx_5-\frac{x_2}{2}dx_4-\frac{x_1}{2}dx_3+\frac{x_4}{2} dx_2+\frac{x_3}{2}dx_1$;
\item $\theta_6=dx_6-\frac{x_1}{2}dx_4+\frac{x_2}{2}dx_3-\frac{x_3}{2}dx_2+\frac{x_4}{2}dx_1$.
\end{itemize}

Finally, we have
\begin{eqnarray*}
 \mathrm{d}(R_\mathbf{x})_\mathbf{0}= \left[\begin{matrix} 
   1 & 0  & 0 & 0 & 0 & 0\\
   0  & 1 & 0 & 0 & 0 &0\\
   0& 0 & 1 & 0 & 0 & 0\\
  0 & 0 & 0 & 1& 0 & 0\\
  \frac{x_3}{2} & \frac{x_4}{2} & -\frac{x_1}{2} & -\frac{x_2}{2}& 1 & 0\\
  \frac{x_4}{2} & -\frac{x_3}{2} & \frac{x_2}{2} & -\frac{x_1}{2}& 0 & 1
   \end{matrix}\right]\,.
\end{eqnarray*}

  \newpage
  
  \section{7D indecomposable Carnot nilpotent Lie algebras}
  The  indecomposable Carnot Lie algebras in dimension 7 are uncountable. They can be subdivided into 45 examples plus two families whose expressions depend on a real parameter $\lambda$. These two families are denoted as ($147E$) and ($147E_1$).
  In dimension 7 there are also uncountable non-stratifiable indecomposable nilpotent Lie algebras, and not all of them are gradable. For the complete list we refer to \cite{Gong_Thesis}.
  
  \subsection*{$(37A)$ }
 
 The following Lie algebra is denoted as $({37A})$ by Gong in \cite{Gong_Thesis}, and as $\mathcal{G}_{7,4,2}$ by Magnin in \cite{magnin}.
 
 The non-trivial brackets are the following:
\begin{equation*}
   [X_1, X_2] = X_5\;,\;[X_2,X_3]=X_6\;,\;[X_2,X_4]=X_7\,.
\end{equation*}
This is a nilpotent Lie algebra of rank 4 and step 2 that is stratifiable. The Lie brackets can be pictured with the diagram:
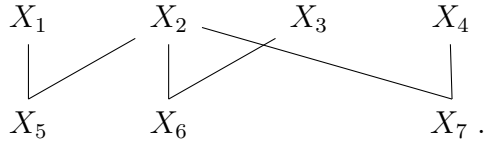
\begin{center}
 
	\begin{tikzcd}[end anchor=north]
		X_1\ar[d,no head] & X_2\ar[d, no head]\ar[dl, no head]\ar[drr, no head,end anchor={[xshift=-3.5ex]north east}] & X_3\ar[dl, no head] & X_4\;\,\ar[d, no head,,start anchor={[xshift=-3.5ex]south east},end anchor={[xshift=-3.5ex]north east}]\\
		X_5 & X_6 & & X_7\;.
	\end{tikzcd}

\end{center}

 

The composition law \eqref{group law in G} of $(37A)$ is given by:

\begin{itemize}
    \item $z_1=x_1+y_1$;
    \item $z_2=x_2+y_2$;
    \item $z_3=x_3+y_3$;
    \item $z_4=x_4+y_4$;
    \item $z_5=x_5+y_5+\frac{1}{2}(x_1y_2-x_2y_1)$;
    \item $z_6=x_6+y_6+\frac{1}{2}(x_2y_3-x_3y_2)$;
    \item $z_7=x_7+y_7+\frac{1}{2}(x_2y_4-x_4y_2)$.
\end{itemize}

Since 

\begin{eqnarray*}
 \mathrm{d}(L_\mathbf{x})_\mathbf{0}= \left[\begin{matrix} 
   1 & 0  & 0& 0& 0& 0& 0\\
   0 & 1& 0& 0& 0& 0& 0\\
  0 & 0& 1& 0& 0& 0& 0\\
   0 & 0& 0& 1& 0& 0& 0\\
   -\frac{x_2}{2} &  \frac{x_1}{2}& 0& 0& 1& 0& 0\\
 0 &  -\frac{x_3}{2}& \frac{x_2}{2}& 0& 0& 1& 0\\
 0 &  -\frac{x_4}{2}& 0& \frac{x_2}{2}& 0& 0& 1
   \end{matrix}\right]\,,
   \end{eqnarray*}
the induced left-invariant vector fields \eqref{leftinvariant vf} are: 
\begin{itemize}
\item $X_1={\partial}_{x_1}-\frac{x_2}{2}{\partial}_{x_5}\,;$
\item $X_2={\partial}_{x_2}+\frac{x_1}{2}{\partial}_{x_5}-\frac{x_3}{2}{\partial}_{x_6}-\frac{x_4}{2}{\partial}_{x_7}\,;$
\item $X_3={\partial}_{x_3}+\frac{x_2}{2}{\partial}_{x_6}\,;$
\item $X_4={\partial}_{x_4}+\frac{x_2}{2}{\partial}_{x_7}\,;$
\item $X_5={\partial}_{x_5}$;
\item $X_6={\partial}_{x_6}$;
\item $X_7={\partial}_{x_7}\,,$
\end{itemize}
and the respective left-invariant 1-forms \eqref{leftinvariant form} are: 
\begin{itemize}
\item $\theta_1=dx_1$;
\item $\theta_2=dx_2$;
\item $\theta_3=dx_3$;
\item $\theta_4=dx_4$;
\item $\theta_5=dx_5-\frac{x_1}{2}dx_2+\frac{x_2}{2}dx_1$;
\item $\theta_6=dx_6-\frac{x_2}{2}dx_3+\frac{x_3}{2}dx_2$;
\item $\theta_7=dx_7-\frac{x_2}{2}dx_4+\frac{x_4}{2}dx_2$.
\end{itemize}

Finally, we have
\begin{eqnarray*}
 \mathrm{d}(R_\mathbf{x})_\mathbf{0}= \left[\begin{matrix} 
   1 & 0  & 0& 0& 0& 0& 0\\
   0 & 1& 0& 0& 0& 0& 0\\
  0 & 0& 1& 0& 0& 0& 0\\
   0 & 0& 0& 1& 0& 0& 0\\
   \frac{x_2}{2} &  -\frac{x_1}{2}& 0& 0& 1& 0& 0\\
 0 &  \frac{x_3}{2}& -\frac{x_2}{2}& 0& 0& 1& 0\\
 0 &  \frac{x_4}{2}& 0& -\frac{x_2}{2}& 0& 0& 1
   \end{matrix}\right]\,.
   \end{eqnarray*}
  \subsection*{$(37B)$ }
 
 The following Lie algebra is denoted as $({37B})$ by Gong in \cite{Gong_Thesis}, and as $\mathcal{G}_{7,4,1}$ by Magnin in \cite{magnin}.
 
 The non-trivial brackets are the following:
\begin{equation*}
   [X_1, X_2] = X_5\,,\,[X_2,X_3]=X_6\,,\,[X_3,X_4]=X_7\,.
\end{equation*}
This is a nilpotent Lie algebra of rank 4 and step 2 that is stratifiable. The Lie brackets can be pictured with the diagram:
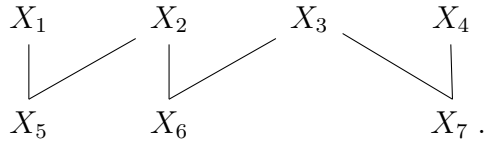
\begin{center}
 
	\begin{tikzcd}[end anchor=north]
		X_1\ar[d,no head] & X_2\ar[d, no head]\ar[dl, no head] & X_3\ar[dl, no head]\ar[dr, no head,end anchor={[xshift=-3.5ex]north east}] & X_4\;\,\ar[d, no head,end anchor={[xshift=-3.5ex]north east},start anchor={[xshift=-3.5ex]south east}]\\
		X_5 & X_6 & & X_7\;.
	\end{tikzcd}

\end{center}

 

The composition law \eqref{group law in G} of $(37B)$ is given by:

\begin{itemize}
    \item $z_1=x_1+y_1$;
    \item $z_2=x_2+y_2$;
    \item $z_3=x_3+y_3$;
    \item $z_4=x_4+y_4$;
    \item $z_5=x_5+y_5+\frac{1}{2}(x_1y_2-x_2y_1)$;
    \item $z_6=x_6+y_6+\frac{1}{2}(x_2y_3-x_3y_2)$;
    \item $z_7=x_7+y_7+\frac{1}{2}(x_3y_4-x_4y_3)$.
\end{itemize}

Since
\begin{eqnarray*}
  \mathrm{d}(L_\mathbf{x})_\mathbf{0}=\left[\begin{matrix} 
   1 & 0  & 0& 0& 0& 0& 0\\
   0 & 1& 0& 0& 0& 0& 0\\
  0 & 0& 1& 0& 0& 0& 0\\
   0 & 0& 0& 1& 0& 0& 0\\
   -\frac{x_2}{2} &  \frac{x_1}{2}& 0& 0& 1& 0& 0\\
 0 &  -\frac{x_3}{2}& \frac{x_2}{2}& 0& 0& 1& 0\\
 0 &  0& -\frac{x_4}{2}& \frac{x_3}{2}& 0& 0& 1
   \end{matrix}\right]\,,
   \end{eqnarray*}
the induced left-invariant vector fields \eqref{leftinvariant vf} are: 
\begin{itemize}
\item $X_1={\partial}_{x_1}-\frac{x_2}{2}{\partial}_{x_5}\,;$
\item $X_2={\partial}_{x_2}+\frac{x_1}{2}{\partial}_{x_5}-\frac{x_3}{2}{\partial}_{x_6}\,;$
\item $X_3={\partial}_{x_3}+\frac{x_2}{2}{\partial}_{x_6}-\frac{x_4}{2}{\partial}_{x_7}\,;$
\item $X_4={\partial}_{x_4}+\frac{x_3}{2}{\partial}_{x_7}\,;$
\item $X_5={\partial}_{x_5}$;
\item $X_6={\partial}_{x_6}$;
\item $X_7={\partial}_{x_7}\,,$
\end{itemize}
and the respective left-invariant 1-forms \eqref{leftinvariant form} are: 
\begin{itemize}
\item $\theta_1=dx_1$;
\item $\theta_2=dx_2$;
\item $\theta_3=dx_3$;
\item $\theta_4=dx_4$;
\item $\theta_5=dx_5-\frac{x_1}{2}dx_2+\frac{x_2}{2}dx_1$;
\item $\theta_6=dx_6-\frac{x_2}{2}dx_3+\frac{x_3}{2}dx_2$;
\item $\theta_7=dx_7-\frac{x_3}{2}dx_4+\frac{x_4}{2}dx_3$.
\end{itemize}

Finally, we have
\begin{eqnarray*}
  \mathrm{d}(R_\mathbf{x})_\mathbf{0}=\left[\begin{matrix} 
   1 & 0  & 0& 0& 0& 0& 0\\
   0 & 1& 0& 0& 0& 0& 0\\
  0 & 0& 1& 0& 0& 0& 0\\
   0 & 0& 0& 1& 0& 0& 0\\
   \frac{x_2}{2} &  -\frac{x_1}{2}& 0& 0& 1& 0& 0\\
 0 &  \frac{x_3}{2}& -\frac{x_2}{2}& 0& 0& 1& 0\\
 0 &  0& \frac{x_4}{2}& -\frac{x_3}{2}& 0& 0& 1
   \end{matrix}\right]\,.
   \end{eqnarray*}
  \subsection*{$(37B_1)$ }
 
 The following Lie algebra is denoted as $({37B_1})$ by Gong in \cite{Gong_Thesis}, and as $G_{7,4,1}$ by Magnin in \cite{magnin}.
 
 The non-trivial brackets are the following:
\begin{equation*}
   [X_1, X_2] = X_5\,,\,[X_1,X_3]=X_6\,,\,[X_1,X_4]=X_7\,,\,[X_2,X_4]=X_6\,,\,[X_3,X_4]=-X_5\,.
\end{equation*}
This is a nilpotent Lie algebra of rank 4 and step 2 that is stratifiable. The Lie brackets can be pictured with the diagram:
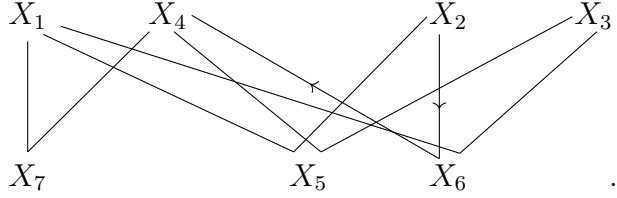
\begin{center}
 
	\begin{tikzcd}[end anchor=north]
		X_1\ar[ddrr, no head,end anchor={[xshift=-3.5ex]north east},start anchor={[xshift=-1.3ex]south east},start anchor={[yshift=.5ex]south east}]\ar[dd, no head]\ar[ddrrr, no head,end anchor={[xshift=-1.6ex]north east},end anchor={[yshift=-.1ex]north east}] & X_4\ar[ddl,no head,start anchor={[xshift=-3.8ex]south east},start anchor={[yshift=.7ex]south east}]\ar[ddr, no head, end anchor={[xshift=-1.5ex]north east},start anchor={[yshift=.7ex]south east},start anchor={[xshift=-2.ex]south east}] \ar[ddrr,-<-=.5, no head, end anchor={[yshift=-.5ex]north east},end anchor={[xshift=-3.1ex]north east}, start anchor={[xshift=-0.7ex]south east},start anchor={[yshift=1.9ex]south east}]& & X_2\ar[dd, no head,end anchor={[yshift=-.5ex]north east},end anchor={[xshift=-3.1ex]north east}, start anchor={[yshift=.5ex]south east},start anchor={[xshift=-3.1ex]south east},->-=.6]\ar[ddl, no head,end anchor={[xshift=-3.5ex]north east},start anchor={[xshift=-4.ex]south east},start anchor={[yshift=1.8ex]south east}] & X_3\ar[ddl, no head,end anchor={[xshift=-1.6ex]north east},start anchor={[xshift=-2.2ex]south east},start anchor={[yshift=.7ex]south east},end anchor={[yshift=-.1ex]north east}]\ar[ddll,no head,start anchor={[xshift=-4.ex]south east},start anchor={[yshift=1.8ex]south east},end anchor={[xshift=-1.5ex]north east}]\\
		& & & &\\
		X_7 & & X_5 & X_6 &\quad\;.
	\end{tikzcd}

\end{center}

 
  
 The composition law \eqref{group law in G} of $(37B_1)$ is given by:

\begin{itemize}
    \item $z_1=x_1+y_1$;
    \item $z_2=x_2+y_2$;
    \item $z_3=x_3+y_3$;
    \item $z_4=x_4+y_4$;
    \item $z_5=x_5+y_5+\frac{1}{2}(x_1y_2-x_2y_1+x_4y_3-x_3y_4)$;
    \item $z_6=x_6+y_6+\frac{1}{2}(x_1y_3-x_3y_1+x_2y_4-x_4y_2)$;
    \item $z_7=x_7+y_7+\frac{1}{2}(x_1y_4-x_4y_1)$.
\end{itemize}

Since

\begin{eqnarray*}
  \mathrm{d}(L_\mathbf{x})_\mathbf{0}=\left[\begin{matrix} 
   1 & 0  & 0& 0& 0& 0& 0\\
   0 & 1& 0& 0& 0& 0& 0\\
  0 & 0& 1& 0& 0& 0& 0\\
   0 & 0& 0& 1& 0& 0& 0\\
   -\frac{x_2}{2} &  \frac{x_1}{2}& \frac{x_4}{2}& -\frac{x_3}{2}& 1& 0& 0\\
 -\frac{x_3}{2} & -\frac{x_4}{2} & \frac{x_1}{2}& \frac{x_2}{2}& 0& 1& 0\\
 -\frac{x_4}{2} &  0& 0& \frac{x_1}{2}& 0& 0& 1
   \end{matrix}\right]\,,
   \end{eqnarray*}
 the induced left-invariant vector fields \eqref{leftinvariant vf} are: 
\begin{itemize}
\item $X_1={\partial}_{x_1}-\frac{x_2}{2}{\partial}_{x_5}-\frac{x_3}{2}{\partial}_{x_6}-\frac{x_4}{2}{\partial}_{x_7}\,;$
\item $X_2={\partial}_{x_2}+\frac{x_1}{2}{\partial}_{x_5}-\frac{x_4}{2}{\partial}_{x_6}\,;$
\item $X_3={\partial}_{x_3}+\frac{x_4}{2}{\partial}_{x_5}+\frac{x_1}{2}{\partial}_{x_6}\,;$
\item $X_4={\partial}_{x_4}-\frac{x_3}{2}{\partial}_{x_5}+\frac{x_2}{2}{\partial}_{x_6}+\frac{x_1}{2}{\partial}_{x_7}\,;$
\item $X_5={\partial}_{x_5}$;
\item $X_6={\partial}_{x_6}$;
\item $X_7={\partial}_{x_7}\,,$
\end{itemize}
and the respective left-invariant 1-forms \eqref{leftinvariant form} are: 
\begin{itemize}
\item $\theta_1=dx_1$;
\item $\theta_2=dx_2$;
\item $\theta_3=dx_3$;
\item $\theta_4=dx_4$;
\item $\theta_5=dx_5+\frac{x_3}{2}dx_4-\frac{x_4}{2}dx_3-\frac{x_1}{2}dx_2+\frac{x_2}{2}dx_1$;
\item $\theta_6=dx_6-\frac{x_2}{2}dx_4-\frac{x_1}{2}dx_3+\frac{x_4}{2}dx_2+\frac{x_3}{2}dx_1$;
\item $\theta_7=dx_7-\frac{x_1}{2}dx_4+\frac{x_4}{2}dx_1$.
\end{itemize}

Finally, we have
\begin{eqnarray*}
  \mathrm{d}(R_\mathbf{x})_\mathbf{0}=\left[\begin{matrix} 
   1 & 0  & 0& 0& 0& 0& 0\\
   0 & 1& 0& 0& 0& 0& 0\\
  0 & 0& 1& 0& 0& 0& 0\\
   0 & 0& 0& 1& 0& 0& 0\\
   \frac{x_2}{2} &  -\frac{x_1}{2}& -\frac{x_4}{2}& \frac{x_3}{2}& 1& 0& 0\\
 \frac{x_3}{2} & \frac{x_4}{2} & -\frac{x_1}{2}& -\frac{x_2}{2}& 0& 1& 0\\
 \frac{x_4}{2} &  0& 0& -\frac{x_1}{2}& 0& 0& 1
   \end{matrix}\right]\,.
   \end{eqnarray*}
  \subsection*{$(37C)$ }
 
 The following Lie algebra is denoted as $({37C})$ by Gong in \cite{Gong_Thesis}, and as $\mathcal{G}_{7,3,24}$ by Magnin in \cite{magnin}.
 
 The non-trivial brackets are the following:
\begin{equation*}
   [X_1, X_2] = X_5\,,\,[X_2,X_3]=X_6\,,\,[X_2,X_4]=X_7\,,\,[X_3,X_4]=X_5.
\end{equation*}
This is a nilpotent Lie algebra of rank 4 and step 2 that is stratifiable. The Lie brackets can be pictured with the diagram:
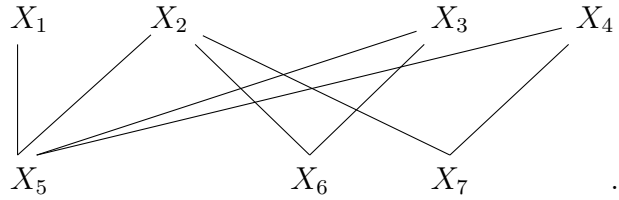
\begin{center}
 
	\begin{tikzcd}[end anchor=north]
		X_1\ar[dd, no head,end anchor={[xshift=-3.3ex]north east},start anchor={[xshift=-3.3ex]south east}]& X_2\ar[ddl,no head,start anchor={[xshift=-3.8ex]south east},start anchor={[yshift=.7ex]south east},end anchor={[xshift=-3.3ex]north east}]\ar[ddr, no head] \ar[ddrr, no head]& & X_3\ar[ddl, no head]\ar[ddlll, no head,end anchor={[xshift=-1.9ex]north east}] & X_4\ar[ddl, no head]\ar[ddllll,no head,end anchor={[xshift=-1.9ex]north east}]\\
		& & & &\\
		X_5 & & X_6 & X_7 &\quad\;.
	\end{tikzcd}

\end{center}

 
 The composition law \eqref{group law in G} of $(37C)$ is given by:

\begin{itemize}
    \item $z_1=x_1+y_1$;
    \item $z_2=x_2+y_2$;
    \item $z_3=x_3+y_3$;
    \item $z_4=x_4+y_4$;
    \item $z_5=x_5+y_5+\frac{1}{2}(x_1y_2-x_2y_1+x_3y_4-x_4y_3)$;
    \item $z_6=x_6+y_6+\frac{1}{2}(x_2y_3-x_3y_2)$;
    \item $z_7=x_7+y_7+\frac{1}{2}(x_2y_4-x_4y_2)$.
\end{itemize}

Since 

\begin{eqnarray*}
  \mathrm{d}(L_\mathbf{x})_\mathbf{0}=\left[\begin{matrix} 
   1 & 0  & 0& 0& 0& 0& 0\\
   0 & 1& 0& 0& 0& 0& 0\\
  0 & 0& 1& 0& 0& 0& 0\\
   0 & 0& 0& 1& 0& 0& 0\\
   -\frac{x_2}{2} &  \frac{x_1}{2}& -\frac{x_4}{2}& \frac{x_3}{2}& 1& 0& 0\\
 0 & -\frac{x_3}{2} & \frac{x_2}{2}& 0& 0& 1& 0\\
 0&  -\frac{x_4}{2}& 0& \frac{x_2}{2}& 0& 0& 1
   \end{matrix}\right]\,,
   \end{eqnarray*}
the induced left-invariant vector fields \eqref{leftinvariant vf} are: 
\begin{itemize}
\item $X_1={\partial}_{x_1}-\frac{x_2}{2}{\partial}_{x_5}\,;$
\item $X_2={\partial}_{x_2}+\frac{x_1}{2}{\partial}_{x_5}-\frac{x_3}{2}{\partial}_{x_6}-\frac{x_4}{2}{\partial}_{x_7}\,;$
\item $X_3={\partial}_{x_3}-\frac{x_4}{2}{\partial}_{x_5}+\frac{x_2}{2}{\partial}_{x_6}\,;$
\item $X_4={\partial}_{x_4}+\frac{x_3}{2}{\partial}_{x_5}+\frac{x_2}{2}{\partial}_{x_7}\,;$
\item $X_5={\partial}_{x_5}$;
\item $X_6={\partial}_{x_6}$;
\item $X_7={\partial}_{x_7}$,
\end{itemize}
and the respective left-invariant 1-forms \eqref{leftinvariant form} are: 
\begin{itemize}
\item $\theta_1=dx_1$;
\item $\theta_2=dx_2$;
\item $\theta_3=dx_3$;
\item $\theta_4=dx_4$;
\item $\theta_5=dx_5-\frac{x_3}{2}dx_4+\frac{x_4}{2}dx_3-\frac{x_1}{2}dx_2+\frac{x_2}{2}dx_1$;
\item $\theta_6=dx_6-\frac{x_2}{2}dx_3+\frac{x_3}{2}dx_2$;
\item $\theta_7=dx_7-\frac{x_2}{2}dx_4+\frac{x_4}{2}dx_2$.
\end{itemize}

Finally, we have
\begin{eqnarray*}
  \mathrm{d}(R_\mathbf{x})_\mathbf{0}=\left[\begin{matrix} 
   1 & 0  & 0& 0& 0& 0& 0\\
   0 & 1& 0& 0& 0& 0& 0\\
  0 & 0& 1& 0& 0& 0& 0\\
   0 & 0& 0& 1& 0& 0& 0\\
   \frac{x_2}{2} &  -\frac{x_1}{2}& \frac{x_4}{2}& -\frac{x_3}{2}& 1& 0& 0\\
 0 & \frac{x_3}{2} & -\frac{x_2}{2}& 0& 0& 1& 0\\
 0&  \frac{x_4}{2}& 0& -\frac{x_2}{2}& 0& 0& 1
   \end{matrix}\right]\,.
   \end{eqnarray*}
  \subsection*{$(37D)$ }
 
 The following Lie algebra is denoted as $({37D})$ by Gong in \cite{Gong_Thesis}, and as $\mathcal{G}_{7,3,12}$ by Magnin in \cite{magnin}.
 
 The non-trivial brackets are the following:
\begin{equation*}
   [X_1, X_2] = X_5\,,\,[X_1,X_3]=X_6\,,\,[X_2,X_4]=X_7\,,\,[X_3,X_4]=X_5\,.
\end{equation*}
This is a nilpotent Lie algebra of rank 4 and step 2 that is stratifiable. The Lie brackets can be pictured with the diagram:
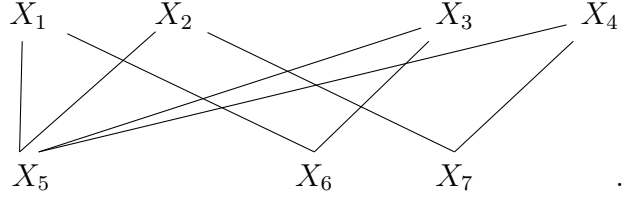
\begin{center}
 
	\begin{tikzcd}[end anchor=north]
		X_1\,\ar[dd, no head,end anchor={[xshift=-3.3ex]north east},start anchor={[xshift=-3.3ex]south east}]\ar[ddrr, no head]& X_2\ar[ddl,no head,start anchor={[xshift=-3.8ex]south east},start anchor={[yshift=.7ex]south east},end anchor={[xshift=-3.3ex]north east}] \ar[ddrr, no head]& & X_3\ar[ddl, no head]\ar[ddlll, no head,end anchor={[xshift=-1.9ex]north east}] & X_4\ar[ddl, no head]\ar[ddllll,no head,end anchor={[xshift=-1.9ex]north east}]\\
		& & & &\\
		X_5 & & X_6 & X_7 &\quad\;.
	\end{tikzcd}

\end{center}

 
 The composition law \eqref{group law in G} of $(37D)$ is given by:

\begin{itemize}
    \item $z_1=x_1+y_1$;
    \item $z_2=x_2+y_2$;
    \item $z_3=x_3+y_3$;
    \item $z_4=x_4+y_4$;
    \item $z_5=x_5+y_5+\frac{1}{2}(x_1y_2-x_2y_1+x_3y_4-x_4y_3)$;
    \item $z_6=x_6+y_6+\frac{1}{2}(x_1y_3-x_3y_1)$;
    \item $z_7=x_7+y_7+\frac{1}{2}(x_2y_4-x_4y_2)$.
\end{itemize}

Since 

\begin{eqnarray*}
  \mathrm{d}(L_\mathbf{x})_\mathbf{0}=\left[\begin{matrix} 
   1 & 0  & 0& 0& 0& 0& 0\\
   0 & 1& 0& 0& 0& 0& 0\\
  0 & 0& 1& 0& 0& 0& 0\\
   0 & 0& 0& 1& 0& 0& 0\\
   -\frac{x_2}{2} &  \frac{x_1}{2}& -\frac{x_4}{2}& \frac{x_3}{2}& 1& 0& 0\\
 -\frac{x_3}{2} & 0 & \frac{x_1}{2}& 0& 0& 1& 0\\
 0&  -\frac{x_4}{2}& 0& \frac{x_2}{2}& 0& 0& 1
   \end{matrix}\right]\,,
   \end{eqnarray*}
the induced left-invariant vector fields \eqref{leftinvariant vf} are: 
\begin{itemize}
\item $X_1={\partial}_{x_1}-\frac{x_2}{2}{\partial}_{x_5}-\frac{x_3}{2}{\partial}_{x_6}\,;$
\item $X_2={\partial}_{x_2}+\frac{x_1}{2}{\partial}_{x_5}-\frac{x_4}{2}{\partial}_{x_7}\,;$
\item $X_3={\partial}_{x_3}-\frac{x_4}{2}{\partial}_{x_5}+\frac{x_1}{2}{\partial}_{x_6}\,;$
\item $X_4={\partial}_{x_4}+\frac{x_3}{2}{\partial}_{x_5}+\frac{x_2}{2}{\partial}_{x_7}\,;$
\item $X_5={\partial}_{x_5}$;
\item $X_6={\partial}_{x_6}$;
\item $X_7={\partial}_{x_7}$,
\end{itemize}
and the respective left-invariant 1-forms \eqref{leftinvariant form} are: 
\begin{itemize}
\item $\theta_1=dx_1$;
\item $\theta_2=dx_2$;
\item $\theta_3=dx_3$;
\item $\theta_4=dx_4$;
\item $\theta_5=dx_5-\frac{x_3}{2}dx_4+\frac{x_4}{2}dx_3-\frac{x_1}{2}dx_2+\frac{x_2}{2}dx_1$;
\item $\theta_6=dx_6-\frac{x_1}{2}dx_3+\frac{x_3}{2}dx_1$;
\item $\theta_7=dx_7-\frac{x_2}{2}dx_4+\frac{x_4}{2}dx_2$.
\end{itemize}

Finally, we have
\begin{eqnarray*}
  \mathrm{d}(R_\mathbf{x})_\mathbf{0}=\left[\begin{matrix} 
   1 & 0  & 0& 0& 0& 0& 0\\
   0 & 1& 0& 0& 0& 0& 0\\
  0 & 0& 1& 0& 0& 0& 0\\
   0 & 0& 0& 1& 0& 0& 0\\
   \frac{x_2}{2} &  -\frac{x_1}{2}& \frac{x_4}{2}& -\frac{x_3}{2}& 1& 0& 0\\
 \frac{x_3}{2} & 0 & -\frac{x_1}{2}& 0& 0& 1& 0\\
 0&  \frac{x_4}{2}& 0& -\frac{x_2}{2}& 0& 0& 1
   \end{matrix}\right]\,.
   \end{eqnarray*}
  \subsection*{$(37D_1)$ }
 
 The following Lie algebra is denoted as $({37D_1})$ by Gong in \cite{Gong_Thesis}, and as $\mathcal{G}_{7,3,12}$ by Magnin in \cite{magnin}.
 
 The non-trivial brackets are the following:
\begin{equation*}
\begin{aligned}{}
    [X_1, X_2] = X_5\,,\,[X_1
    ,X_3]=&X_6\,,\,[X_1,X_4]=X_7\,,\\ [X_2,X_3]=-X_7\,,\,[X_2,X_4]=&X_6\,,\,[X_3,X_4]=-X_5.
\end{aligned}
\end{equation*}
This is a nilpotent Lie algebra of rank 4 and step 2 that is stratifiable. The Lie brackets can be pictured with the diagram:
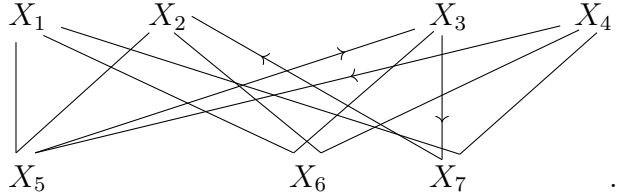
\begin{center}
 
	\begin{tikzcd}[end anchor=north]
			X_1\ar[ddrr, no head,end anchor={[xshift=-3.5ex]north east},start anchor={[xshift=-1.3ex]south east},start anchor={[yshift=.5ex]south east}]\ar[dd, no head,end anchor={[xshift=-3.3ex]north east},start anchor={[xshift=-3.3ex]south east}]\ar[ddrrr, no head,end anchor={[xshift=-1.6ex]north east},end anchor={[yshift=-.1ex]north east}] & X_2\ar[ddl,no head,start anchor={[xshift=-3.8ex]south east},start anchor={[yshift=.7ex]south east},end anchor={[xshift=-3.3ex]north east}]\ar[ddr, no head, end anchor={[xshift=-1.5ex]north east},start anchor={[yshift=.7ex]south east},start anchor={[xshift=-2.ex]south east}] \ar[ddrr,-<-=.3, no head, end anchor={[yshift=-.5ex]north east},end anchor={[xshift=-2.9ex]north east}, start anchor={[xshift=-0.7ex]south east},start anchor={[yshift=1.9ex]south east}]& & X_3\ar[dd, no head,end anchor={[yshift=-.5ex]north east},end anchor={[xshift=-2.9ex]north east}, start anchor={[yshift=.5ex]south east},start anchor={[xshift=-2.9ex]south east},->-=.7]\ar[ddl, no head,end anchor={[xshift=-3.5ex]north east},start anchor={[xshift=-3.5ex]south east},start anchor={[yshift=.8ex]south east}]\ar[ddlll,,-<-=.2 ,no head, end anchor={[xshift=-1.9ex]north east}] & X_4\ar[ddl, no head,end anchor={[xshift=-1.6ex]north east},start anchor={[xshift=-2.2ex]south east},start anchor={[yshift=.7ex]south east},end anchor={[yshift=-.1ex]north east}]\ar[ddll,no head,start anchor={[xshift=-3.5ex]south east},start anchor={[yshift=0.9ex]south east},end anchor={[xshift=-1.5ex]north east}]\ar[ddllll,->-=.4, no head,end anchor={[xshift=-1.9ex]north east}]\\
		& & & &\\
		X_5 & & X_6 & X_7 &\quad\;.
	\end{tikzcd}

\end{center}

 
 The composition law \eqref{group law in G} of $(37D_1)$ is given by:

\begin{itemize}
    \item $z_1=x_1+y_1$;
    \item $z_2=x_2+y_2$;
    \item $z_3=x_3+y_3$;
    \item $z_4=x_4+y_4$;
    \item $z_5=x_5+y_5+\frac{1}{2}(x_1y_2-x_2y_1+x_4y_3-x_3y_4)$;
    \item $z_6=x_6+y_6+\frac{1}{2}(x_1y_3-x_3y_1+x_2y_4-x_4y_2)$;
    \item $z_7=x_7+y_7+\frac{1}{2}(x_1y_4-x_4y_1+x_3y_2-x_2y_3)$.
\end{itemize}
Since 

\begin{eqnarray*}
 \mathrm{d}(L_\mathbf{x})_\mathbf{0}= \left[\begin{matrix} 
   1 & 0  & 0& 0& 0& 0& 0\\
   0 & 1& 0& 0& 0& 0& 0\\
  0 & 0& 1& 0& 0& 0& 0\\
   0 & 0& 0& 1& 0& 0& 0\\
   -\frac{x_2}{2} &  \frac{x_1}{2}& \frac{x_4}{2}& -\frac{x_3}{2}& 1& 0& 0\\
 -\frac{x_3}{2} & -\frac{x_4}{2} & \frac{x_1}{2}& \frac{x_2}{2}& 0& 1& 0\\
  -\frac{x_4}{2}&  \frac{x_3}{2}& -\frac{x_2}{2}& \frac{x_1}{2}& 0& 0& 1
   \end{matrix}\right]\,,
   \end{eqnarray*}
the induced left-invariant vector fields \eqref{leftinvariant vf} are: 
\begin{itemize}
\item $X_1={\partial}_{x_1}-\frac{x_2}{2}{\partial}_{x_5}-\frac{x_3}{2}{\partial}_{x_6}-\frac{x_4}{2}{\partial}_{x_7}\,;$
\item $X_2={\partial}_{x_2}+\frac{x_1}{2}{\partial}_{x_5}-\frac{x_4}{2}{\partial}_{x_6}+\frac{x_3}{2}{\partial}_{x_7}\,;$
\item $X_3={\partial}_{x_3}+\frac{x_4}{2}{\partial}_{x_5}+\frac{x_1}{2}{\partial}_{x_6}-\frac{x_2}{2}{\partial}_{x_7}\,;$
\item $X_4={\partial}_{x_4}-\frac{x_3}{2}{\partial}_{x_5}+\frac{x_2}{2}{\partial}_{x_6}+\frac{x_1}{2}{\partial}_{x_7}\,;$
\item $X_5={\partial}_{x_5}$;
\item $X_6={\partial}_{x_6}$;
\item $X_7={\partial}_{x_7},$
\end{itemize}
and the respective left-invariant 1-forms \eqref{leftinvariant form} are: 
\begin{itemize}
\item $\theta_1=dx_1$;
\item $\theta_2=dx_2$;
\item $\theta_3=dx_3$;
\item $\theta_4=dx_4$;
\item $\theta_5=dx_5+\frac{x_3}{2}dx_4-\frac{x_4}{2}dx_3-\frac{x_1}{2}dx_2+\frac{x_2}{2}dx_1$;
\item $\theta_6=dx_6-\frac{x_2}{2}dx_4-\frac{x_1}{2}dx_3+\frac{x_4}{2}dx_2+\frac{x_3}{2}dx_1$;
\item $\theta_7=dx_7-\frac{x_1}{2}dx_4+\frac{x_2}{2}dx_3-\frac{x_3}{2}dx_2+\frac{x_4}{2}dx_1$.
\end{itemize}

Finally, we have
\begin{eqnarray*}
 \mathrm{d}(R_\mathbf{x})_\mathbf{0}= \left[\begin{matrix} 
   1 & 0  & 0& 0& 0& 0& 0\\
   0 & 1& 0& 0& 0& 0& 0\\
  0 & 0& 1& 0& 0& 0& 0\\
   0 & 0& 0& 1& 0& 0& 0\\
   \frac{x_2}{2} &  -\frac{x_1}{2}& -\frac{x_4}{2}& \frac{x_3}{2}& 1& 0& 0\\
 \frac{x_3}{2} & \frac{x_4}{2} & -\frac{x_1}{2}& -\frac{x_2}{2}& 0& 1& 0\\
  \frac{x_4}{2}&  -\frac{x_3}{2}& \frac{x_2}{2}& -\frac{x_1}{2}& 0& 0& 1
   \end{matrix}\right]\,.
   \end{eqnarray*}
  \subsection*{$(357A)$ }
 
 The following Lie algebra is denoted as $({357A})$ by Gong in \cite{Gong_Thesis}, and as $\mathcal{G}_{7,3,6}$ by Magnin in \cite{magnin}.
 
 The non-trivial brackets are the following:
\begin{equation*}
   [X_1, X_2] = X_3\,,\,[X_1,X_3]=X_5\,,\,[X_1,X_4]=X_7\,,\,[X_2,X_4]=X_6\,.
\end{equation*}
This is a nilpotent Lie algebra of rank 3 and step 3 that is stratifiable. The Lie brackets can be pictured with the diagram:
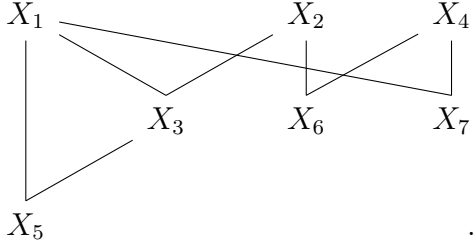
\begin{center}
 
	\begin{tikzcd}[end anchor=north]
		X_1\ar[dr,no head]\ar[dd, no head]\ar[drrr, no head] & &X_2\ar[dl, no head]\ar[d, no head] & X_4\ar[dl, no head]\ar[d, no head]\\
		&X_3  \ar[dl, no head] & X_6 &X_7\\
		X_5 & & &\quad\;.
	\end{tikzcd}

\end{center}

 
The composition law \eqref{group law in G} of $(357A)$ is given by:

\begin{itemize}
    \item $z_1=x_1+y_1$;
    \item $z_2=x_2+y_2$;
    \item $z_3=x_3+y_3+\frac{1}{2}(x_1y_2-x_2y_1)$;
    \item $z_4=x_4+y_4$;
    \item $z_5=x_5+y_5+\frac{1}{2}(x_1y_3-x_3y_1)+\frac{1}{12}(x_1-y_1)(x_1y_2-x_2y_1)$;
    \item $z_6=x_6+y_6+\frac{1}{2}(x_2y_4-x_4y_2)$;
    \item $z_7=x_7+y_7+\frac{1}{2}(x_1y_4-x_4y_1)$.
\end{itemize}

Since


\begin{eqnarray*}
  \mathrm{d}(L_\mathbf{x})_\mathbf{0}=\left[\begin{matrix} 
   1 & 0  & 0& 0& 0& 0& 0\\
   0 & 1& 0& 0& 0& 0& 0\\
   -\frac{x_2}{2} & \frac{x_1}{2}& 1& 0& 0& 0& 0\\
   0 & 0& 0& 1& 0& 0& 0\\
   -\frac{x_3}{2}-\frac{x_1x_2}{12} &  \frac{x_1^2}{12}& \frac{x_1}{2}& 0& 1& 0& 0\\
 0 &  -\frac{x_4}{2}& 0& \frac{x_2}{2}& 0& 1& 0\\
   -\frac{x_4}{2} &  0& 0& \frac{x_1}{2}& 0& 0& 1
   \end{matrix}\right]\,,
   \end{eqnarray*}
the induced left-invariant vector fields \eqref{leftinvariant vf} are: 
\begin{itemize}
\item $X_1={\partial}_{x_1}-\frac{x_2}{2}{\partial}_{x_3}-\big(\frac{x_3}{2}+\frac{x_1x_2}{12}\big){\partial}_{x_5}-\frac{x_4}{2}{\partial}_{x_7}\,;$
\item $X_2={\partial}_{x_2}+\frac{x_1}{2}{\partial}_{x_3}+\frac{x_1^2}{12}{\partial}_{x_5}-\frac{x_4}{2}{\partial}_{x_6}\,;$
\item $X_3={\partial}_{x_3}+\frac{x_1}{2}{\partial}_{x_5}\,;$
\item $X_4={\partial}_{x_4}+\frac{x_2}{2}{\partial}_{x_6}+\frac{x_1}{2}{\partial}_{x_7}\,;$
\item $X_5={\partial}_{x_5}$;
\item $X_6={\partial}_{x_6}$;
\item $X_7={\partial}_{x_7}$,
\end{itemize}
and the respective left-invariant 1-forms \eqref{leftinvariant form} are: 
\begin{itemize}
\item $\theta_1=dx_1$;
\item $\theta_2=dx_2$;
\item $\theta_3=dx_3-\frac{x_1}{2}dx_2+\frac{x_2}{2}dx_1$;
\item $\theta_4=dx_4$;
\item $\theta_5=dx_5-\frac{x_1}{2}dx_3+\frac{x_1^2}{6}dx_2+\big(\frac{x_3}{2}-\frac{x_1x_2}{6}\big)dx_1$;
\item $\theta_6=dx_6-\frac{x_2}{2}dx_4+\frac{x_4}{2}dx_2$;
\item $\theta_7=dx_7-\frac{x_1}{2}dx_4+\frac{x_4}{2}dx_1$.
\end{itemize}

Finally, we have
\begin{eqnarray*}
  \mathrm{d}(R_\mathbf{x})_\mathbf{0}=\left[\begin{matrix} 
   1 & 0  & 0& 0& 0& 0& 0\\
   0 & 1& 0& 0& 0& 0& 0\\
   \frac{x_2}{2} & -\frac{x_1}{2}& 1& 0& 0& 0& 0\\
   0 & 0& 0& 1& 0& 0& 0\\
   \frac{x_3}{2}-\frac{x_1x_2}{12} &  \frac{x_1^2}{12}& -\frac{x_1}{2}& 0& 1& 0& 0\\
 0 &  \frac{x_4}{2}& 0& -\frac{x_2}{2}& 0& 1& 0\\
   \frac{x_4}{2} &  0& 0& -\frac{x_1}{2}& 0& 0& 1
   \end{matrix}\right]\,.
   \end{eqnarray*}
  \subsection*{$(357B)$ }
 
 The following Lie algebra is denoted as $({357B})$ by Gong in \cite{Gong_Thesis}, and as $\mathcal{G}_{7,3,23}$ by Magnin in \cite{magnin}.
 
 The non-trivial brackets are the following:
\begin{equation*}
   [X_1, X_2] = X_3\,,\,[X_1,X_3]=X_5\,,\,[X_1,X_4]=X_7\,,\,[X_2,X_3]=X_6\,.
\end{equation*}
This is a nilpotent Lie algebra of rank 3 and step 3 that is stratifiable. The Lie brackets can be pictured with the diagram:
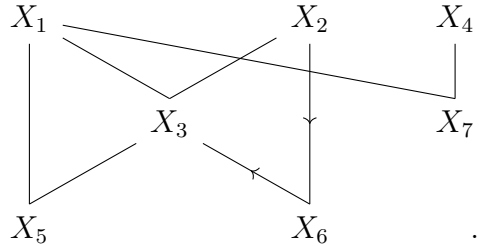
\begin{center}
 
	\begin{tikzcd}[end anchor=north]
		X_1\ar[dr,no head]\ar[dd, no head]\ar[drrr, no head] & &X_2\ar[dl, no head]\ar[dd, no head,->-=.5] & X_4\ar[d, no head]\\
		&X_3  \ar[dl, no head]\ar[dr, no head, -<-=.5] &  &X_7\\
		X_5 & &X_6 &\quad\;.
	\end{tikzcd}

\end{center}

 
 The composition law \eqref{group law in G} of $(357B)$ is given by:

\begin{itemize}
    \item $z_1=x_1+y_1$;
    \item $z_2=x_2+y_2$;
    \item $z_3=x_3+y_3+\frac{1}{2}(x_1y_2-x_2y_1)$;
    \item $z_4=x_4+y_4$;
    \item $z_5=x_5+y_5+\frac{1}{2}(x_1y_3-x_3y_1)+\frac{1}{12}(x_1-y_1)(x_1y_2-x_2y_1)$;
    \item $z_6=x_6+y_6+\frac{1}{2}(x_2y_3-x_3y_2)+\frac{1}{12}(x_2-y_2)(x_1y_2-x_2y_1)$;
    \item $z_7=x_7+y_7+\frac{1}{2}(x_1y_4-x_4y_1)$.
\end{itemize}



Since
\begin{eqnarray*}
  \mathrm{d}(L_\mathbf{x})_\mathbf{0}=\left[\begin{matrix} 
  1 & 0  & 0& 0& 0& 0& 0\\
   0 & 1& 0& 0& 0& 0& 0\\
   -\frac{x_2}{2} & \frac{x_1}{2}& 1& 0& 0& 0& 0\\
   0 & 0& 0& 1& 0& 0& 0\\
   -\frac{x_3}{2}-\frac{x_1x_2}{12} &  \frac{x_1^2}{12}& \frac{x_1}{2}& 0& 1& 0& 0\\
 -\frac{x_2^2}{12} &  -\frac{x_3}{2}+\frac{x_1x_2}{12}& \frac{x_2}{2}& 0& 0& 1& 0\\
   -\frac{x_4}{2} &  0& 0& \frac{x_1}{2}& 0& 0& 1
   \end{matrix}\right]\,,
   \end{eqnarray*}
the induced left-invariant vector fields \eqref{leftinvariant vf} are: 
\begin{itemize}
\item $X_1={\partial}_{x_1}-\frac{x_2}{2}{\partial}_{x_3}-\big(\frac{x_3}{2}+\frac{x_1x_2}{12}\big){\partial}_{x_5}-\frac{x_2^2}{12}{\partial}_{x_6}-\frac{x_4}{2}{\partial}_{x_7}\,;$
\item $X_2={\partial}_{x_2}+\frac{x_1}{2}{\partial}_{x_3}+\frac{x_1^2}{12}{\partial}_{x_5}+\big(\frac{x_1x_2}{12}-\frac{x_3}{2}\big){\partial}_{x_6}\,;$
\item $X_3={\partial}_{x_3}+\frac{x_1}{2}{\partial}_{x_5}+\frac{x_2}{2}{\partial}_{x_6}\,;$
\item $X_4={\partial}_{x_4}+\frac{x_1}{2}{\partial}_{x_7}\,;$
\item $X_5={\partial}_{x_5}$;
\item $X_6={\partial}_{x_6}$;
\item $X_7={\partial}_{x_7}$,
\end{itemize}
and the respective left-invariant 1-forms \eqref{leftinvariant form} are: 
\begin{itemize}
\item $\theta_1=dx_1$;
\item $\theta_2=dx_2$;
\item $\theta_3=dx_3-\frac{x_1}{2}dx_2+\frac{x_2}{2}dx_1$;
\item $\theta_4=dx_4$;
\item $\theta_5=dx_5-\frac{x_1}{2}dx_3+\frac{x_1^2}{6}dx_2+\big(\frac{x_3}{2}-\frac{x_1x_2}{6}\big)dx_1$;
\item $\theta_6=dx_6-\frac{x_2}{2}dx_3+\big(\frac{x_3}{2}+\frac{x_1x_2}{6}\big)dx_2-\frac{x_2^2}{6}dx_1$;
\item $\theta_7=dx_7-\frac{x_1}{2}dx_4+\frac{x_4}{2}dx_1$.
\end{itemize}

Finally, we have
\begin{eqnarray*}
  \mathrm{d}(R_\mathbf{x})_\mathbf{0}=\left[\begin{matrix} 
  1 & 0  & 0& 0& 0& 0& 0\\
   0 & 1& 0& 0& 0& 0& 0\\
 \frac{x_2}{2} & -\frac{x_1}{2}& 1& 0& 0& 0& 0\\
   0 & 0& 0& 1& 0& 0& 0\\
   \frac{x_3}{2}-\frac{x_1x_2}{12} &  \frac{x_1^2}{12}& -\frac{x_1}{2}& 0& 1& 0& 0\\
 -\frac{x_2^2}{12} &  \frac{x_3}{2}+\frac{x_1x_2}{12}&-\frac{x_2}{2}& 0& 0& 1& 0\\
   \frac{x_4}{2} &  0& 0& -\frac{x_1}{2}& 0& 0& 1
   \end{matrix}\right]\,.
   \end{eqnarray*}
  \subsection*{$(27A)$ }
 
The following Lie algebra is denoted as $({27A})$ by Gong in \cite{Gong_Thesis}, and as $\mathcal{G}_{7,4,3}$ by Magnin in \cite{magnin}.

The non-trivial brackets are the following:
\begin{equation*}
   [X_1, X_2] = X_6\,,\,[X_1,X_4]=X_7\,,\,[X_3,X_5]=X_7\,.
\end{equation*}
This is a nilpotent Lie algebra of rank 5 and step 2 that is stratifiable. The Lie brackets can be pictured with the diagram:
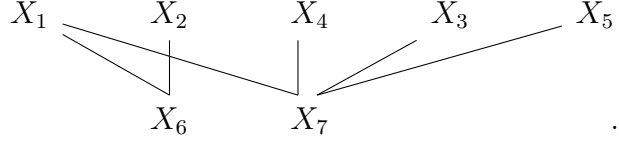
\begin{center}
 
	\begin{tikzcd}[end anchor=north]
		X_1\ar[dr,no head]\ar[drr, no head,end anchor={[xshift=-3.3ex]north east}] & X_2\ar[d, no head] & X_4\ar[d, no head,end anchor={[xshift=-3.3ex]north east},start anchor={[xshift=-3.3ex]south east}] & X_3\ar[dl, no head,end anchor={[xshift=-1.9ex]north east}] &X_5\ar[dll, no head,end anchor={[xshift=-1.9ex]north east}]\\
		 & X_6 &X_7 & &\quad\;.
	\end{tikzcd}

\end{center}

 
The composition law \eqref{group law in G} of $(27A)$ is given by;

\begin{itemize}
    \item $z_1=x_1+y_1$;
    \item $z_2=x_2+y_2$;
    \item $z_3=x_3+y_3$;
    \item $z_4=x_4+y_4$;
    \item $z_5=x_5+y_5$;
    \item $z_6=x_6+y_6+\frac{1}{2}(x_1y_2-x_2y_1)$;
    \item $z_7=x_7+y_7+\frac{1}{2}(x_1y_4-x_4y_1+x_3y_5-x_5y_3)$.
\end{itemize}

Since 

\begin{eqnarray*}
  \mathrm{d}(L_\mathbf{x})_\mathbf{0}=\left[\begin{matrix} 
   1 & 0  & 0& 0& 0& 0& 0\\
   0 & 1& 0& 0& 0& 0& 0\\
  0 & 0& 1& 0& 0& 0& 0\\
   0 & 0& 0& 1& 0& 0& 0\\
   0 &  0& 0& 0& 1& 0& 0\\
 -\frac{x_2}{2} &  \frac{x_1}{2}& 0& 0& 0& 1& 0\\
 -\frac{x_4}{2} & 0 & -\frac{x_5}{2}& \frac{x_1}{2}& \frac{x_3}{2}& 0& 1
   \end{matrix}\right]\,,
   \end{eqnarray*}
the induced left-invariant vector fields \eqref{leftinvariant vf} are: 
\begin{itemize}
\item $X_1={\partial}_{x_1}-\frac{x_2}{2}{\partial}_{x_6}-\frac{x_4}{2}{\partial}_{x_7}\,;$
\item $X_2={\partial}_{x_2}+\frac{x_1}{2}{\partial}_{x_6}\,;$
\item $X_3={\partial}_{x_3}-\frac{x_5}{2}{\partial}_{x_7}\,;$
\item $X_4={\partial}_{x_4}+\frac{x_1}{2}{\partial}_{x_7}\,;$
\item $X_5={\partial}_{x_5}+\frac{x_3}{2}{\partial}_{x_7}$;
\item $X_6={\partial}_{x_6}$;
\item $X_7={\partial}_{x_7}$,
\end{itemize}
and the respective left-invariant 1-forms \eqref{leftinvariant form} are: 
\begin{itemize}
\item $\theta_1=dx_1$;
\item $\theta_2=dx_2$;
\item $\theta_3=dx_3$;
\item $\theta_4=dx_4$;
\item $\theta_5=dx_5$;
\item $\theta_6=dx_6-\frac{x_1}{2}dx_2+\frac{x_2}{2}dx_1$;
\item $\theta_7=dx_7-\frac{x_3}{2}dx_5-\frac{x_1}{2}dx_4+\frac{x_5}{2}dx_3+\frac{x_4}{2}dx_1$.
\end{itemize}

Finally, we have
\begin{eqnarray*}
 \mathrm{d}(R_\mathbf{x})_\mathbf{0}= \left[\begin{matrix} 
   1 & 0  & 0& 0& 0& 0& 0\\
   0 & 1& 0& 0& 0& 0& 0\\
  0 & 0& 1& 0& 0& 0& 0\\
   0 & 0& 0& 1& 0& 0& 0\\
   0 &  0& 0& 0& 1& 0& 0\\
 \frac{x_2}{2} &  -\frac{x_1}{2}& 0& 0& 0& 1& 0\\
 \frac{x_4}{2} & 0 & \frac{x_5}{2}& -\frac{x_1}{2}& -\frac{x_3}{2}& 0& 1
   \end{matrix}\right]\,.
   \end{eqnarray*}
  \subsection*{$(27B)$ }
 
 The following Lie algebra is denoted as $({27B})$ by Gong in \cite{Gong_Thesis}, and as $\mathcal{G}_{7,3,19}$ by Magnin in \cite{magnin}.
 
 The non-trivial brackets are the following:
\begin{equation*}
   [X_1, X_2] = X_6\,,\,[X_1,X_5]=X_7\,,\,[X_2,X_3]=X_7\,,\,[X_3,X_4]=X_6\,.
\end{equation*}
This is a nilpotent Lie algebra of rank 5 and step 2 that is stratifiable. The Lie brackets can be pictured with the diagram:
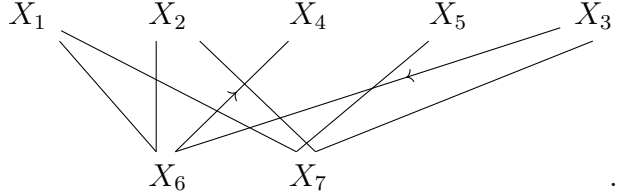
\begin{center}
 
	\begin{tikzcd}[end anchor=north]
		X_1\ar[ddr,no head,end anchor={[xshift=-3.3ex]north east}]\ar[ddrr, no head,end anchor={[xshift=-3.3ex]north east}] & X_2\ar[dd, no head,start anchor={[xshift=-3.3ex]south east},end anchor={[xshift=-3.3ex]north east}] \ar[ddr, no head,end anchor={[xshift=-1.9ex]north east}]& X_4\ar[ddl, no head,-<-=.5,end anchor={[xshift=-1.9ex]north east}] & X_5\ar[ddl, no head,end anchor={[xshift=-3.3ex]north east}] &X_3\ar[ddll, no head,end anchor={[xshift=-1.9ex]north east},start anchor={[xshift=-2.5ex]south east}]\ar[ddlll, no head,->-=.4,end anchor={[xshift=-1.9ex]north east}]\\
		& & &\\
		 & X_6 &X_7 & &\quad\;.
	\end{tikzcd}

\end{center}

 
The composition law \eqref{group law in G} of $(27B)$ is given by:

\begin{itemize}
    \item $z_1=x_1+y_1$;
    \item $z_2=x_2+y_2$;
    \item $z_3=x_3+y_3$;
    \item $z_4=x_4+y_4$;
    \item $z_5=x_5+y_5$;
    \item $z_6=x_6+y_6+\frac{1}{2}(x_1y_2-x_2y_1+x_3y_4-x_4y_3)$;
    \item $z_7=x_7+y_7+\frac{1}{2}(x_1y_5-x_5y_1+x_2y_3-x_3y_2)$.
\end{itemize}

Since 

\begin{eqnarray*}
  \mathrm{d}(L_\mathbf{x})_\mathbf{0}=\left[\begin{matrix} 
   1 & 0  & 0& 0& 0& 0& 0\\
   0 & 1& 0& 0& 0& 0& 0\\
  0 & 0& 1& 0& 0& 0& 0\\
   0 & 0& 0& 1& 0& 0& 0\\
   0 &  0& 0& 0& 1& 0& 0\\
 -\frac{x_2}{2} &  \frac{x_1}{2}& -\frac{x_4}{2}& \frac{x_3}{2}& 0& 1& 0\\
 -\frac{x_5}{2} & -\frac{x_3}{2} &\frac{x_2}{2} & 0& \frac{x_1}{2}& 0& 1
   \end{matrix}\right]\,,
   \end{eqnarray*}
the induced left-invariant vector fields \eqref{leftinvariant vf} are: 
\begin{itemize}
\item $X_1={\partial}_{x_1}-\frac{x_2}{2}{\partial}_{x_6}-\frac{x_5}{2}{\partial}_{x_7}\,;$
\item $X_2={\partial}_{x_2}+\frac{x_1}{2}{\partial}_{x_6}-\frac{x_3}{2}{\partial}_{x_7}\,;$
\item $X_3={\partial}_{x_3}-\frac{x_4}{2}{\partial}_{x_6}+\frac{x_2}{2}{\partial}_{x_7}\,;$
\item $X_4={\partial}_{x_4}+\frac{x_3}{2}{\partial}_{x_6}\,;$
\item $X_5={\partial}_{x_5}+\frac{x_1}{2}{\partial}_{x_7}$;
\item $X_6={\partial}_{x_6}$;
\item $X_7={\partial}_{x_7}$,
\end{itemize}
and the respective left-invariant 1-forms\eqref{leftinvariant form} are: 
\begin{itemize}
\item $\theta_1=dx_1$;
\item $\theta_2=dx_2$;
\item $\theta_3=dx_3$;
\item $\theta_4=dx_4$;
\item $\theta_5=dx_5$;
\item $\theta_6=dx_6-\frac{x_3}{2}dx_4+\frac{x_4}{2}dx_3-\frac{x_1}{2}dx_2+\frac{x_2}{2}dx_1$;
\item $\theta_7=dx_7-\frac{x_1}{2}dx_5-\frac{x_2}{2}dx_3+\frac{x_3}{2}dx_2+\frac{x_5}{2}dx_1$.
\end{itemize}

Finally, we have
\begin{eqnarray*}
  \mathrm{d}(R_\mathbf{x})_\mathbf{0}=\left[\begin{matrix} 
   1 & 0  & 0& 0& 0& 0& 0\\
   0 & 1& 0& 0& 0& 0& 0\\
  0 & 0& 1& 0& 0& 0& 0\\
   0 & 0& 0& 1& 0& 0& 0\\
   0 &  0& 0& 0& 1& 0& 0\\
 \frac{x_2}{2} &  -\frac{x_1}{2}& \frac{x_4}{2}& -\frac{x_3}{2}& 0& 1& 0\\
 \frac{x_5}{2} & \frac{x_3}{2} &-\frac{x_2}{2} & 0& -\frac{x_1}{2}& 0& 1
   \end{matrix}\right]\,.
   \end{eqnarray*}
  \subsection*{$(257B)$ }
 
 The following Lie algebra is denoted as $({257B})$ by Gong in \cite{Gong_Thesis}, and as $\mathcal{G}_{7,3,11}$ by Magnin in \cite{magnin}.
 
 The non-trivial brackets are the following:
\begin{equation*}
   [X_1, X_2] = X_3\,,\,[X_1, X_3] = X_6\,,\,[X_1,X_4]=X_7\,,\,[X_2,X_5]=X_7\,.
\end{equation*}
This is a nilpotent Lie algebra of rank 4 and step 3 that is stratifiable. The Lie brackets can be pictured with the diagram:
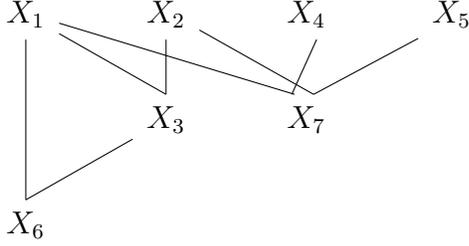
\begin{center}
 
	\begin{tikzcd}[end anchor=north]
		X_1\ar[dd,no head]\ar[drr,no head,end anchor={[xshift=-3.3ex]north east}]\ar[dr, no head] & X_2\ar[d, no head] \ar[dr, no head,end anchor={[xshift=-1.9ex]north east}]& X_4\ar[d, no head,end anchor={[xshift=-3.5ex]north east}] & X_5\ar[dl, no head,end anchor={[xshift=-1.9ex]north east}] \\
		&X_3\ar[dl, no head] &X_7 &\\
		 X_6&  & & \quad\;.
	\end{tikzcd}

\end{center}

 
The composition law \eqref{group law in G} of $(257B)$ is given by:

\begin{itemize}
    \item $z_1=x_1+y_1$;
    \item $z_2=x_2+y_2$;
    \item $z_3=x_3+y_3+\frac{1}{2}(x_1y_2-x_2y_1)$;
    \item $z_4=x_4+y_4$;
    \item $z_5=x_5+y_5$;
    \item $z_6=x_6+y_6+\frac{1}{2}(x_1y_3-x_3y_1)+\frac{1}{12}(x_1-y_1)(x_1y_2-x_2y_1)$;
    \item $z_7=x_7+y_7+\frac{1}{2}(x_1y_4-x_4y_1+x_2y_5-x_5y_2)$.
\end{itemize}

Since
\begin{eqnarray*}
  \mathrm{d}(L_\mathbf{x})_\mathbf{0}=\left[\begin{matrix} 
  1 & 0  & 0& 0& 0& 0& 0\\
   0 & 1& 0& 0& 0& 0& 0\\
  -\frac{x_2}{2} & \frac{x_1}{2}& 1& 0& 0& 0& 0\\
   0 & 0& 0& 1& 0& 0& 0\\
   0 &  0& 0& 0& 1& 0& 0\\
 -\frac{x_1x_2}{12}-\frac{x_3}{2} &  \frac{x_1^2}{12}& \frac{x_1}{2}& 0& 0& 1& 0\\
 -\frac{x_4}{2} & -\frac{x_5}{2} &0 & \frac{x_1}{2}& \frac{x_2}{2}& 0& 1
   \end{matrix}\right]\,,
   \end{eqnarray*}
the induced left-invariant vector fields \eqref{leftinvariant vf} are: 
\begin{itemize}
\item $X_1=\partial_{x_1} -\frac{x_2}{2}\partial_{x_3}-\big(\frac{x_1x_2}{12}+\frac{x_3}{2}\big)\partial_{x_6}-\frac{x_4}{2}\partial_{x_7}\,;$
\item $X_2=\partial_{x_2}+\frac{x_1}{2}\partial_{x_3}+\frac{x_1^2}{12}\partial_{x_6}-\frac{x_5}{2}\partial_{x_7}\,;$
\item $X_3=\partial_{x_3}+\frac{x_1}{2}\partial_{x_6}\,;$
\item $X_4=\partial_{x_4}+\frac{x_1}{2}\partial_{x_7}\,;$
\item $X_5=\partial_{x_5}+\frac{x_2}{2}\partial_{x_7}$;
\item $X_6=\partial_{x_6}$;
\item $X_7=\partial_{x_7}$,
\end{itemize}
and the respective left-invariant 1-forms \eqref{leftinvariant form} are: 
\begin{itemize}
\item $\theta_1=dx_1$;
\item $\theta_2=dx_2$;
\item $\theta_3=dx_3-\frac{x_1}{2}dx_2+\frac{x_2}{2}dx_1$;
\item $\theta_4=dx_4$;
\item $\theta_5=dx_5$;
\item $\theta_6=dx_6-\frac{x_1}{2}dx_3+\frac{x_1^2}{6}dx_2+\big(\frac{x_3}{2}-\frac{x_1x_2}{6}\big)dx_1$;
\item $\theta_7=dx_7-\frac{x_2}{2}dx_5-\frac{x_1}{2}dx_4+\frac{x_5}{2}dx_2+\frac{x_4}{2}dx_1$.
\end{itemize}

Finally, we have
\begin{eqnarray*}
  \mathrm{d}(R_\mathbf{x})_\mathbf{0}=\left[\begin{matrix} 
  1 & 0  & 0& 0& 0& 0& 0\\
   0 & 1& 0& 0& 0& 0& 0\\
  \frac{x_2}{2} & -\frac{x_1}{2}& 1& 0& 0& 0& 0\\
   0 & 0& 0& 1& 0& 0& 0\\
   0 &  0& 0& 0& 1& 0& 0\\
 \frac{x_3}{2}-\frac{x_1x_2}{12} &  \frac{x_1^2}{12}& -\frac{x_1}{2}& 0& 0& 1& 0\\
 \frac{x_4}{2} & \frac{x_5}{2} &0 & -\frac{x_1}{2}& -\frac{x_2}{2}& 0& 1
   \end{matrix}\right]\,.
   \end{eqnarray*}
  \subsection*{$(247A)$ }
 
 The following Lie algebra is denoted as $({247A})$ by Gong in \cite{Gong_Thesis}, and as $\mathcal{G}_{7,3,20}$ by Magnin in \cite{magnin}.
 
 The non-trivial brackets are the following:
\begin{equation*}
   [X_1, X_i] = X_{i+2}\,,\,2\le i\le 5\,.
\end{equation*}
This is a nilpotent Lie algebra of rank 3 and step 3 that is stratifiable. The Lie brackets can be pictured with the diagram:
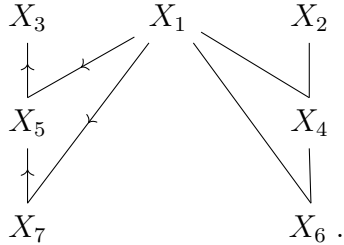
\begin{center}
 
	\begin{tikzcd}[end anchor=north]
		X_3 \ar[d,no head, -<-=.5]& X_1\ar[dr, no head,end anchor={[xshift=-3.5ex]north east}]\ar[dl, no head, ->-=.5]\ar[ddr, no head,end anchor={[xshift=-3.5ex]north east}]\ar[ddl, no head, ->-=.5] & X_2\;\,\ar[d, no head,end anchor={[xshift=-3.5ex]north east},start anchor={[xshift=-3.5ex]south east}]\\
		X_5\ar[d, no head, -<-=.5] & & X_4\;\,\ar[d, no head,end anchor={[xshift=-3.5ex]north east},start anchor={[xshift=-3.5ex]south east}]\\
		X_7 & & X_6\;.
	\end{tikzcd}

\end{center}

 
The composition law \eqref{group law in G} of $(247A)$ is given by:

\begin{itemize}
    \item $z_1=x_1+y_1$;
    \item $z_2=x_2+y_2$;
    \item $z_3=x_3+y_3$;
    \item $z_4=x_4+y_4+\frac{1}{2}(x_1y_2-x_2y_1)$;
    \item $z_5=x_5+y_5+\frac{1}{2}(x_1y_3-x_3y_1)$;
    \item $z_6=x_6+y_6+\frac{1}{2}(x_1y_4-x_4y_1)+\frac{1}{12}(x_1-y_1)(x_1y_2-x_2y_1)$;
    \item $z_7=x_7+y_7+\frac{1}{2}(x_1y_5-x_5y_1)+\frac{1}{12}(x_1-y_1)(x_1y_3-x_3y_1)$.
\end{itemize}



Since
\begin{eqnarray*}
  \mathrm{d}(L_\mathbf{x})_\mathbf{0}=\left[\begin{matrix}    1 & 0  & 0& 0& 0& 0& 0\\
   0 & 1& 0& 0& 0& 0& 0\\
   0 & 0& 1& 0& 0& 0& 0\\
   -\frac{x_2}{2} & \frac{x_1}{2}& 0& 1& 0& 0& 0\\
   -\frac{x_3}{2} &  0& \frac{x_1}{2}& 0& 1& 0& 0\\
 -\frac{x_1x_2}{12}-\frac{x_4}{2} &  \frac{x_1^2}{12}& 0& \frac{x_1}{2}& 0& 1& 0\\
 -\frac{x_1x_3}{12}-\frac{x_5}{2} & 0 &\frac{x_1^2}{12} & 0& \frac{x_1}{2}& 0& 1
   \end{matrix}\right]\,,
   \end{eqnarray*}
the induced left-invariant vector fields \eqref{leftinvariant vf} are: 
\begin{itemize}
\item $X_1=\partial_{x_1} -\frac{x_2}{2}\partial_{x_4}-\frac{x_3}{2}\partial_{x_5}-\big(\frac{x_1x_2}{12}+\frac{x_4}{2}\big)\partial_{x_6}-\big(\frac{x_1x_3}{12}+\frac{x_5}{2}\big)\partial_{x_7}\,;$
\item $X_2=\partial_{x_2}+\frac{x_1}{2}\partial_{x_4}+\frac{x_1^2}{12}\partial_{x_6}\,;$
\item $X_3=\partial_{x_3}+\frac{x_1}{2}\partial_{x_5}+\frac{x_1^2}{12}\partial_{x_7}\,;$
\item $X_4=\partial_{x_4}+\frac{x_1}{2}\partial_{x_6}\,;$
\item $X_5=\partial_{x_5}+\frac{x_1}{2}\partial_{x_7}$;
\item $X_6=\partial_{x_6}$;
\item $X_7=\partial_{x_7}$,
\end{itemize}
and the respective left-invariant 1-forms \eqref{leftinvariant form} are: 
\begin{itemize}
\item $\theta_1=dx_1$;
\item $\theta_2=dx_2$;
\item $\theta_3=dx_3$;
\item $\theta_4=dx_4-\frac{x_1}{2}dx_2+\frac{x_2}{2}dx_1$;
\item $\theta_5=dx_5-\frac{x_1}{2}dx_3+\frac{x_3}{2}dx_1$;
\item $\theta_6=dx_6-\frac{x_1}{2}dx_4+\frac{x_1^2}{6}dx_2+\big(\frac{x_4}{2}-\frac{x_1x_2}{6}\big)dx_1$;
\item $\theta_7=dx_7-\frac{x_1}{2}dx_5+\frac{x_1^2}{6}dx_3+\big(\frac{x_5}{2}-\frac{x_1x_3}{6}\big)dx_1$.
\end{itemize}

Finally, we have
\begin{eqnarray*}
  \mathrm{d}(R_\mathbf{x})_\mathbf{0}=\left[\begin{matrix}    1 & 0  & 0& 0& 0& 0& 0\\
   0 & 1& 0& 0& 0& 0& 0\\
   0 & 0& 1& 0& 0& 0& 0\\
   \frac{x_2}{2} & -\frac{x_1}{2}& 0& 1& 0& 0& 0\\
   \frac{x_3}{2} &  0& -\frac{x_1}{2}& 0& 1& 0& 0\\
 \frac{x_4}{2}-\frac{x_1x_2}{12} &  \frac{x_1^2}{12}& 0& -\frac{x_1}{2}& 0& 1& 0\\
 \frac{x_5}{2}-\frac{x_1x_3}{12} & 0 &\frac{x_1^2}{12} & 0& -\frac{x_1}{2}& 0& 1
   \end{matrix}\right]\,.
   \end{eqnarray*}
  \subsection*{$(247B)$ }
 
 The following Lie algebra is denoted as $({247B})$ by Gong in \cite{Gong_Thesis}, and as $\mathcal{G}_{7,3,21}$ by Magnin in \cite{magnin}.
 
 The non-trivial brackets are the following:
\begin{equation*}
   [X_1, X_i] = X_{i+2}\,,\,2\le i\le 4\,,\,[X_3,X_5]=X_7\,.
\end{equation*}
This is a nilpotent Lie algebra of rank 3 and step 3 that is stratifiable. The Lie brackets can be pictured with the diagram:
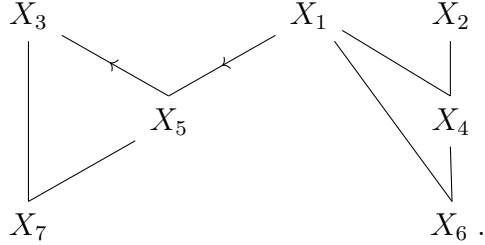
\begin{center}
 
	\begin{tikzcd}[end anchor=north]
		X_3 \ar[dr,no head, -<-=.5] \ar[dd, no head]& & X_1\ar[dr, no head,end anchor={[xshift=-3.5ex]north east}]\ar[dl, no head, ->-=.5]\ar[ddr, no head,end anchor={[xshift=-3.5ex]north east}] & X_2\;\,\ar[d, no head,end anchor={[xshift=-3.5ex]north east},start anchor={[xshift=-3.5ex]south east}]\\
		 &X_5\ar[dl, no head] & & X_4\;\,\ar[d, no head,end anchor={[xshift=-3.5ex]north east},start anchor={[xshift=-3.5ex]south east}]\\
		X_7 & & & X_6\;.
	\end{tikzcd}

\end{center}

 
The composition law \eqref{group law in G} of $(247B)$ is given by:

\begin{itemize}
    \item $z_1=x_1+y_1$;
    \item $z_2=x_2+y_2$;
    \item $z_3=x_3+y_3$;
    \item $z_4=x_4+y_4+\frac{1}{2}(x_1y_2-x_2y_1)$;
    \item $z_5=x_5+y_5+\frac{1}{2}(x_1y_3-x_3y_1)$;
    \item $z_6=x_6+y_6+\frac{1}{2}(x_1y_4-x_4y_1)+\frac{1}{12}(x_1-y_1)(x_1y_2-x_2y_1)$;
    \item $z_7=x_7+y_7+\frac{1}{2}(x_3y_5-x_5y_3)+\frac{1}{12}(x_3-y_3)(x_1y_3-x_3y_1)$.
\end{itemize}



Since
\begin{eqnarray*}
  \mathrm{d}(L_\mathbf{x})_\mathbf{0}=\left[\begin{matrix}    1 & 0  & 0& 0& 0& 0& 0\\
   0 & 1& 0& 0& 0& 0& 0\\
   0 & 0& 1& 0& 0& 0& 0\\
   -\frac{x_2}{2} & \frac{x_1}{2}& 0& 1& 0& 0& 0\\
   -\frac{x_3}{2} &  0& \frac{x_1}{2}& 0& 1& 0& 0\\
 -\frac{x_1x_2}{12}-\frac{x_4}{2} &  \frac{x_1^2}{12}& 0& \frac{x_1}{2}& 0& 1& 0\\
 -\frac{x_3^2}{12}& 0 &\frac{x_1x_3}{12}-\frac{x_5}{2} & 0& \frac{x_3}{2}& 0& 1
   \end{matrix}\right]\,,
   \end{eqnarray*}
the induces left-invariant vector fields \eqref{leftinvariant vf} are: 
\begin{itemize}
\item $X_1=\partial_{x_1} -\frac{x_2}{2}\partial_{x_4}-\frac{x_3}{2}\partial_{x_5}-\big(\frac{x_1x_2}{12}+\frac{x_4}{2}\big)\partial_{x_6}-\frac{x_3^2}{12}\partial_{x_7}\,;$
\item $X_2=\partial_{x_2}+\frac{x_1}{2}\partial_{x_4}+\frac{x_1^2}{12}\partial_{x_6}\,;$
\item $X_3=\partial_{x_3}+\frac{x_1}{2}\partial_{x_5}+\big(\frac{x_1x_3}{12}-\frac{x_5}{2}\big)\partial_{x_7}\,;$
\item $X_4=\partial_{x_4}+\frac{x_1}{2}\partial_{x_6}\,;$
\item $X_5=\partial_{x_5}+\frac{x_3}{2}\partial_{x_7}$;
\item $X_6=\partial_{x_6}$;
\item $X_7=\partial_{x_7}$,
\end{itemize}
and the respective left-invariant 1-forms \eqref{leftinvariant form} are: 
\begin{itemize}
\item $\theta_1=dx_1$;
\item $\theta_2=dx_2$;
\item $\theta_3=dx_3$;
\item $\theta_4=dx_4-\frac{x_1}{2}dx_2+\frac{x_2}{2}dx_1$;
\item $\theta_5=dx_5-\frac{x_1}{2}dx_3+\frac{x_3}{2}dx_1$;
\item $\theta_6=dx_6-\frac{x_1}{2}dx_4+\frac{x_1^2}{6}dx_2+\big(\frac{x_4}{2}-\frac{x_1x_2}{6}\big)dx_1$;
\item $\theta_7=dx_7-\frac{x_3}{2}dx_5+\big(\frac{x_5}{2}+\frac{x_1x_3}{6}\big)dx_3-\frac{x_3^2}{6}dx_1$.
\end{itemize}

Finally, we have
\begin{eqnarray*}
  \mathrm{d}(R_\mathbf{x})_\mathbf{0}=\left[\begin{matrix}    1 & 0  & 0& 0& 0& 0& 0\\
   0 & 1& 0& 0& 0& 0& 0\\
   0 & 0& 1& 0& 0& 0& 0\\
   \frac{x_2}{2} & -\frac{x_1}{2}& 0& 1& 0& 0& 0\\
   \frac{x_3}{2} &  0& -\frac{x_1}{2}& 0& 1& 0& 0\\
 \frac{x_4}{2}-\frac{x_1x_2}{12} &  \frac{x_1^2}{12}& 0& -\frac{x_1}{2}& 0& 1& 0\\
 -\frac{x_3^2}{12}& 0 &\frac{x_1x_3}{12}+\frac{x_5}{2} & 0& -\frac{x_3}{2}& 0& 1
   \end{matrix}\right]\,.
   \end{eqnarray*}
  \subsection*{$(247C)$ }
 
 The following Lie algebra is denoted as $({247C})$ by Gong in \cite{Gong_Thesis}, and as $\mathcal{G}_{7,2,43}$ by Magnin in \cite{magnin}.
 
 The non-trivial brackets are the following:
\begin{equation*}
   [X_1, X_i] = X_{i+2}\,,\,2\le i\le 4\,,\,[X_1,X_5]=X_7\,,\,[X_3,X_5]=X_6.
\end{equation*}
This is a nilpotent Lie algebra of rank 3 and step 3 that is stratifiable. The Lie brackets can be pictured with the diagram:
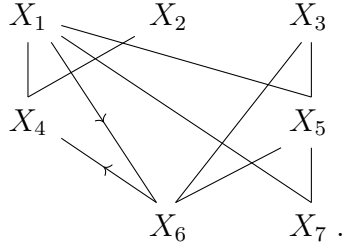
\begin{center}
 
	\begin{tikzcd}[end anchor=north]
		X_1 \ar[d, no head]\ar[ddr,no head,end anchor={[xshift=-3.3ex]north east},->-=.5]\ar[ddrr, no head,,end anchor={[xshift=-3.5ex]north east}]\ar[drr, no head,end anchor={[xshift=-3.5ex]north east}]& X_2\ar[dl, no head] &X_3\phantom{\;.}\ar[d, no head,end anchor={[xshift=-3.5ex]north east},start anchor={[xshift=-3.5ex]south east}]\ar[ddl, no head, end anchor={[xshift=-1.9ex]north east}]\\
		X_4 \ar[dr,no head, end anchor={[xshift=-3.3ex]north east},-<-=.5]& & X_5\phantom{\;.}\ar[d, no head,end anchor={[xshift=-3.5ex]north east},start anchor={[xshift=-3.5ex]south east}]\ar[dl, no head, end anchor={[xshift=-1.9ex]north east}]\\
		& X_6 & X_7\;.
	\end{tikzcd}

\end{center}

 
The composition law \eqref{group law in G} of $(247C)$ is given by:

\begin{itemize}
    \item $z_1=x_1+y_1$;
    \item $z_2=x_2+y_2$;
    \item $z_3=x_3+y_3$;
    \item $z_4=x_4+y_4+\frac{1}{2}(x_1y_2-x_2y_1)$;
    \item $z_5=x_5+y_5+\frac{1}{2}(x_1y_3-x_3y_1)$;
    \item $z_6=x_6+y_6+\frac{1}{2}(x_1y_4-x_4y_1+x_3y_5-x_5y_3)+\frac{1}{12}(x_1-y_1)(x_1y_2-x_2y_1)\salto +\frac{1}{12}(x_3-y_3)(x_1y_3-x_3y_1)$;
    \item $z_7=x_7+y_7+\frac{1}{2}(x_1y_5-x_5y_1)+\frac{1}{12}(x_1-y_1)(x_1y_3-x_3y_1)$.
\end{itemize}



Since
\begin{eqnarray*}
  \mathrm{d}(L_\mathbf{x})_\mathbf{0}=\left[\begin{matrix}    
  1 & 0  & 0& 0& 0& 0& 0\\
   0 & 1& 0& 0& 0& 0& 0\\
   0 & 0& 1& 0& 0& 0& 0\\
   -\frac{x_2}{2} & \frac{x_1}{2}& 0& 1& 0& 0& 0\\
   -\frac{x_3}{2} &  0& \frac{x_1}{2}& 0& 1& 0& 0\\
 -\frac{x_1x_2+x_3^2}{12}-\frac{x_4}{2} &  \frac{x_1^2}{12}& \frac{x_1x_3}{12}-\frac{x_5}{2}& \frac{x_1}{2}& \frac{x_3}{2}& 1& 0\\
 -\frac{x_1x_3}{12}-\frac{x_5}{2}& 0 &\frac{x_1^2}{12} & 0& \frac{x_1}{2}& 0& 1
   \end{matrix}\right]\,,
   \end{eqnarray*}
the induced left-invariant vector fields \eqref{leftinvariant vf} are: 
\begin{itemize}
\item $X_1=\partial_{x_1} -\frac{x_2}{2}\partial_{x_4}-\frac{x_3}{2}\partial_{x_5}-\big(\frac{x_1x_2+x_3^2}{12}+\frac{x_4}{2}\big)\partial_{x_6}-\big(\frac{x_1x_3}{12}+\frac{x_5}{2}\big)\partial_{x_7}\,;$
\item $X_2=\partial_{x_2}+\frac{x_1}{2}\partial_{x_4}+\frac{x_1^2}{12}\partial_{x_6}\,;$
\item $X_3=\partial_{x_3}+\frac{x_1}{2}\partial_{x_5}+\big(\frac{x_1x_3}{12}-\frac{x_5}{2}\big)\partial_{x_6}+\frac{x_1^2}{12}\partial_{x_7}\,;$
\item $X_4=\partial_{x_4}+\frac{x_1}{2}\partial_{x_6}\,;$
\item $X_5=\partial_{x_5}+\frac{x_3}{2}\partial_{x_6}+\frac{x_1}{2}\partial_{x_7}$;
\item $X_6=\partial_{x_6}$;
\item $X_7=\partial_{x_7}$,
\end{itemize}
and the respective left-invariant 1-forms \eqref{leftinvariant form} are: 
\begin{itemize}
\item $\theta_1=dx_1$;
\item $\theta_2=dx_2$;
\item $\theta_3=dx_3$;
\item $\theta_4=dx_4-\frac{x_1}{2}dx_2+\frac{x_2}{2}dx_1$;
\item $\theta_5=dx_5-\frac{x_1}{2}dx_3+\frac{x_3}{2}dx_1$;
\item $\theta_6=dx_6-\frac{x_3}{2}dx_5-\frac{x_1}{2}dx_4+\big(\frac{x_5}{2}+\frac{x_1x_3}{6}\big)dx_3+\frac{x_1^2}{6}dx_2+\big(\frac{x_4}{2}-\frac{x_1x_2+x_3^2}{6}\big)dx_1$;
\item $\theta_7=dx_7-\frac{x_1}{2}dx_5+\frac{x_1^2}{6}dx_3+\big(\frac{x_5}{2}-\frac{x_1x_3}{6}\big)dx_1$.
\end{itemize}

Finally, we have
\begin{eqnarray*}
  \mathrm{d}(R_\mathbf{x})_\mathbf{0}=\left[\begin{matrix}    
  1 & 0  & 0& 0& 0& 0& 0\\
   0 & 1& 0& 0& 0& 0& 0\\
   0 & 0& 1& 0& 0& 0& 0\\
   \frac{x_2}{2} & -\frac{x_1}{2}& 0& 1& 0& 0& 0\\
   \frac{x_3}{2} &  0& -\frac{x_1}{2}& 0& 1& 0& 0\\
 \frac{x_4}{2}-\frac{x_1x_2+x_3^2}{12} &  \frac{x_1^2}{12}& \frac{x_1x_3}{12}+\frac{x_5}{2}& -\frac{x_1}{2}& -\frac{x_3}{2}& 1& 0\\
\frac{x_5}{2} -\frac{x_1x_3}{12}& 0 &\frac{x_1^2}{12} & 0& -\frac{x_1}{2}& 0& 1
   \end{matrix}\right]\,.
   \end{eqnarray*}
  \subsection*{$(247D)$ }
 
 The following Lie algebra is denoted as $({247D})$ by Gong in \cite{Gong_Thesis}, and as $\mathcal{G}_{7,3,22}$ by Magnin in \cite{magnin}.
 
 The non-trivial brackets are the following:
\begin{equation*}
   [X_1, X_i] = X_{i+2}\,,\,2\le i\le 4\,,\,[X_2,X_5]=X_7\,,\,[X_3,X_4]=X_7\,.
\end{equation*}
This is a nilpotent Lie algebra of rank 3 and step 3 that is stratifiable. The Lie brackets can be pictured with the diagram:
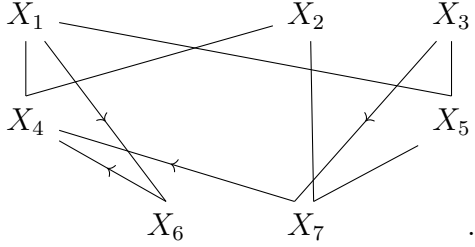
\begin{center}
 
	\begin{tikzcd}[end anchor=north]
		X_1 \ar[d, no head]\ar[ddr,no head,->-=.5]\ar[drrr, no head]& & X_2\ar[dll, no head]\ar[dd, no head, end anchor={[xshift=-1.9ex]north east}] &X_3\ar[d, no head]\ar[ddl, no head, end anchor={[xshift=-3.3ex]north east},->-=.5]\\
		X_4 \ar[dr,no head,-<-=.5]\ar[drr,no head,-<-=.5,end anchor={[xshift=-3.3ex]north east}]& & & X_5\ar[dl, no head, end anchor={[xshift=-1.9ex]north east}]\\
		& X_6&X_7 & \quad\;.
	\end{tikzcd}

\end{center}

 
The composition law \eqref{group law in G} of $(247D)$ is given by:

\begin{itemize}
    \item $z_1=x_1+y_1$;
    \item $z_2=x_2+y_2$;
    \item $z_3=x_3+y_3$;
    \item $z_4=x_4+y_4+\frac{1}{2}(x_1y_2-x_2y_1)$;
    \item $z_5=x_5+y_5+\frac{1}{2}(x_1y_3-x_3y_1)$;
    \item $z_6=x_6+y_6+\frac{1}{2}(x_1y_4-x_4y_1)+\frac{1}{12}(x_1-y_1)(x_1y_2-x_2y_1)$;
    \item $z_7=x_7+y_7+\frac{1}{2}(x_3y_4-x_4y_3+x_2y_5-x_5y_2)+\frac{1}{12}(x_3-y_3)(x_1y_2-x_2y_1)\salto+\frac{1}{12}(x_2-y_2)(x_1y_3-x_3y_1)$.
\end{itemize}



Since
\begin{eqnarray*}
  \mathrm{d}(L_\mathbf{x})_\mathbf{0}=\left[\begin{matrix}    
   1 & 0  & 0& 0& 0& 0& 0\\
   0 & 1& 0& 0& 0& 0& 0\\
   0 & 0& 1& 0& 0& 0& 0\\
   -\frac{x_2}{2} & \frac{x_1}{2}& 0& 1& 0& 0& 0\\
   -\frac{x_3}{2} &  0& \frac{x_1}{2}& 0& 1& 0& 0\\
 -\frac{x_1x_2}{12}-\frac{x_4}{2} &  \frac{x_1^2}{12}& 0& \frac{x_1}{2}& 0& 1& 0\\
 -\frac{x_2x_3}{6}& \frac{x_1x_3}{12}-\frac{x_5}{2} &\frac{x_1x_2}{12}-\frac{x_4}{2} & \frac{x_3}{2}& \frac{x_2}{2}& 0& 1
   \end{matrix}\right]\,,
   \end{eqnarray*}
the induced left-invariant vector fields \eqref{leftinvariant vf} are:
\begin{itemize}
\item $X_1=\partial_{x_1} -\frac{x_2}{2}\partial_{x_4}-\frac{x_3}{2}\partial_{x_5}-\big(\frac{x_1x_2}{12}+\frac{x_4}{2}\big)\partial_{x_6}-\frac{x_2x_3}{12}\partial_{x_7}\,;$
\item $X_2=\partial_{x_2}+\frac{x_1}{2}\partial_{x_4}+\frac{x_1^2}{12}\partial_{x_6}+\big(\frac{x_1x_3}{12}-\frac{x_5}{2}\big)\partial_{x_7}\,;$
\item $X_3=\partial_{x_3}+\frac{x_1}{2}\partial_{x_5}+\big(\frac{x_1x_2}{12}-\frac{x_4}{2}\big)\partial_{x_7}\,;$
\item $X_4=\partial_{x_4}+\frac{x_1}{2}\partial_{x_6}+\frac{x_3}{2}\partial_{x_7}\,;$
\item $X_5=\partial_{x_5}+\frac{x_2}{2}\partial_{x_7}$;
\item $X_6=\partial_{x_6}$;
\item $X_7=\partial_{x_7}$,
\end{itemize}
and the respective left-invariant 1-forms \eqref{leftinvariant form} are: 
\begin{itemize}
\item $\theta_1=dx_1$;
\item $\theta_2=dx_2$;
\item $\theta_3=dx_3$;
\item $\theta_4=dx_4-\frac{x_1}{2}dx_2+\frac{x_2}{2}dx_1$;
\item $\theta_5=dx_5-\frac{x_1}{2}dx_3+\frac{x_3}{2}dx_1$;
\item $\theta_6=dx_6-\frac{x_1}{2}dx_4+\frac{x_1^2}{6}dx_2+\big(\frac{x_4}{2}-\frac{x_1x_2}{6}\big)dx_1$;
\item $\theta_7=dx_7-\frac{x_2}{2}dx_5-\frac{x_3}{2}dx_4+\big(\frac{x_1x_2}{6}+\frac{x_4}{2}\big)dx_3+\big(\frac{x_1x_3}{6}+\frac{x_5}{2}\big)dx_2-\frac{x_2x_3}{6}dx_1$.
\end{itemize}

Finally, we have
\begin{eqnarray*}
  \mathrm{d}(R_\mathbf{x})_\mathbf{0}=\left[\begin{matrix}    
   1 & 0  & 0& 0& 0& 0& 0\\
   0 & 1& 0& 0& 0& 0& 0\\
   0 & 0& 1& 0& 0& 0& 0\\
   \frac{x_2}{2} & -\frac{x_1}{2}& 0& 1& 0& 0& 0\\
   \frac{x_3}{2} &  0& -\frac{x_1}{2}& 0& 1& 0& 0\\
 \frac{x_4}{2}-\frac{x_1x_2}{12} &  \frac{x_1^2}{12}& 0& -\frac{x_1}{2}& 0& 1& 0\\
 -\frac{x_2x_3}{6}& \frac{x_1x_3}{12}+\frac{x_5}{2} &\frac{x_1x_2}{12}+\frac{x_4}{2} & -\frac{x_3}{2}& -\frac{x_2}{2}& 0& 1
   \end{matrix}\right]\,.
   \end{eqnarray*}
  \subsection*{$(247E)$ }
 
 The following Lie algebra is denoted as $({247E})$ by Gong in \cite{Gong_Thesis}, and as $\mathcal{G}_{7,2,12}$ by Magnin in \cite{magnin}.
 
 The non-trivial brackets are the following:
\begin{equation*}
   [X_1, X_i] = X_{i+2}\,,\,2\le i\le 4\,,\,[X_1,X_5]=X_6\,,\,[X_2,X_5]=X_7\,,\,[X_3,X_4]=X_7\,.
\end{equation*}
This is a nilpotent Lie algebra of rank 3 and step 3 that is stratifiable. The Lie brackets can be pictured with the diagram:
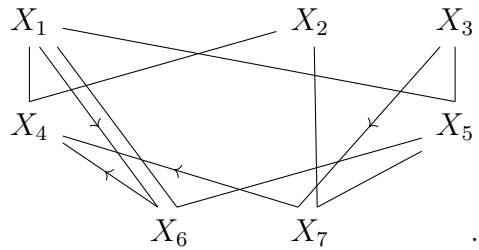
\begin{center}
 
	\begin{tikzcd}[end anchor=north]
		X_1 \ar[d, no head]\ar[ddr,no head,->-=.5,end anchor={[xshift=-3.3ex]north east},start anchor={[xshift=-1.7ex]south east}]\ar[drrr, no head]\ar[ddr,no head, end anchor={[xshift=-1.9ex]north east},start anchor={[xshift=-0.4ex]south east}]& & X_2\ar[dll, no head]\ar[dd, no head, end anchor={[xshift=-1.9ex]north east}] &X_3\ar[d, no head]\ar[ddl, no head, end anchor={[xshift=-3.3ex]north east},->-=.5]\\
		X_4 \ar[dr,no head,-<-=.5,end anchor={[xshift=-3.3ex]north east}]\ar[drr,no head,-<-=.5,end anchor={[xshift=-3.3ex]north east}]& & & X_5\ar[dl, no head, end anchor={[xshift=-1.9ex]north east}]\ar[dll,no head, end anchor={[xshift=-1.9ex]north east}]\\
		& X_6&X_7 & \quad\;.
	\end{tikzcd}

\end{center}

 
The composition law \eqref{group law in G} of $(247E)$ is given by:

\begin{itemize}
    \item $z_1=x_1+y_1$;
    \item $z_2=x_2+y_2$;
    \item $z_3=x_3+y_3$;
    \item $z_4=x_4+y_4+\frac{1}{2}(x_1y_2-x_2y_1)$;
    \item $z_5=x_5+y_5+\frac{1}{2}(x_1y_3-x_3y_1)$;
    \item $z_6=x_6+y_6+\frac{1}{2}(x_1y_4-x_4y_1+x_1y_5-x_5y_1)+\frac{1}{12}(x_1-y_1)(x_1y_2-x_2y_1)\salto+\frac{1}{12}(x_1-y_1)(x_1y_3-x_3y_1)$;
    \item $z_7=x_7+y_7+\frac{1}{2}(x_3y_4-x_4y_3+x_2y_5-x_5y_2)+\frac{1}{12}(x_3-y_3)(x_1y_2-x_2y_1)\salto+\frac{1}{12}(x_2-y_2)(x_1y_3-x_3y_1)$.
\end{itemize}



Since
\begin{eqnarray*}
  \mathrm{d}(L_\mathbf{x})_\mathbf{0}=\left[\begin{matrix}    
   1 & 0  & 0& 0& 0& 0& 0\\
   0 & 1& 0& 0& 0& 0& 0\\
   0 & 0& 1& 0& 0& 0& 0\\
   -\frac{x_2}{2} & \frac{x_1}{2}& 0& 1& 0& 0& 0\\
   -\frac{x_3}{2} &  0& \frac{x_1}{2}& 0& 1& 0& 0\\
-\frac{x_1x_2+x_1x_3}{12}-\frac{x_4+x_5}{2} &  \frac{x_1^2}{12}& \frac{x_1^2}{12}& \frac{x_1}{2}& \frac{x_1}{2}& 1& 0\\
 -\frac{x_2x_3}{6}& \frac{x_1x_3}{12}-\frac{x_5}{2} &\frac{x_1x_2}{12}-\frac{x_4}{2} & \frac{x_3}{2}& \frac{x_2}{2}& 0& 1
   \end{matrix}\right]\,,
   \end{eqnarray*}
the induced left-invariant vector fields \eqref{leftinvariant vf} are: 
\begin{itemize}
\item $X_1=\partial_{x_1} -\frac{x_2}{2}\partial_{x_4}-\frac{x_3}{2}\partial_{x_5}-\big(\frac{x_1x_2+x_1x_3}{12}+\frac{x_4+x_5}{2}\big)\partial_{x_6}-\frac{x_2x_3}{6}\partial_{x_7}\,;$
\item $X_2=\partial_{x_2}+\frac{x_1}{2}\partial_{x_4}+\frac{x_1^2}{12}\partial_{x_6}+\big(\frac{x_1x_3}{12}-\frac{x_5}{2}\big)\partial_{x_7}\,;$
\item $X_3=\partial_{x_3}+\frac{x_1}{2}\partial_{x_5}+\frac{x_1^2}{12}\partial_{x_6}+\big(\frac{x_1x_2}{12}-\frac{x_4}{2}\big)\partial_{x_7}\,;$
\item $X_4=\partial_{x_4}+\frac{x_1}{2}\partial_{x_6}+\frac{x_3}{2}\partial_{x_7}\,;$
\item $X_5=\partial_{x_5}+\frac{x_1}{2}\partial_{x_6}+\frac{x_2}{2}\partial_{x_7}$;
\item $X_6=\partial_{x_6}$;
\item $X_7=\partial_{x_7}$,
\end{itemize}
and the respective left-invariant 1-forms \eqref{leftinvariant form} are: 
\begin{itemize}
\item $\theta_1=dx_1$;
\item $\theta_2=dx_2$;
\item $\theta_3=dx_3$;
\item $\theta_4=dx_4-\frac{x_1}{2}dx_2+\frac{x_2}{2}dx_1$;
\item $\theta_5=dx_5-\frac{x_1}{2}dx_3+\frac{x_3}{2}dx_1$;
\item $\theta_6=dx_6-\frac{x_1}{2}dx_5-\frac{x_1}{2}dx_4+\frac{x_1^2}{6}dx_3+\frac{x_1^2}{6}dx_2+\big(\frac{x_4+x_5}{2}-\frac{x_1x_2+x_1x_3}{6}\big)dx_1$;
\item $\theta_7=dx_7-\frac{x_2}{2}dx_5-\frac{x_3}{2}dx_4+\big(\frac{x_1x_2}{6}+\frac{x_4}{2}\big)dx_3+\big(\frac{x_1x_3}{6}+\frac{x_5}{2}\big)dx_2-\frac{x_2x_3}{3}dx_1$.
\end{itemize}

Finally, we have
\begin{eqnarray*}
  \mathrm{d}(R_\mathbf{x})_\mathbf{0}=\left[\begin{matrix}    
   1 & 0  & 0& 0& 0& 0& 0\\
   0 & 1& 0& 0& 0& 0& 0\\
   0 & 0& 1& 0& 0& 0& 0\\
   \frac{x_2}{2} & -\frac{x_1}{2}& 0& 1& 0& 0& 0\\
   \frac{x_3}{2} &  0& -\frac{x_1}{2}& 0& 1& 0& 0\\
\frac{x_4+x_5}{2}-\frac{x_1x_2+x_1x_3}{12} &  \frac{x_1^2}{12}& \frac{x_1^2}{12}& -\frac{x_1}{2}& -\frac{x_1}{2}& 1& 0\\
 -\frac{x_2x_3}{6}& \frac{x_1x_3}{12}+\frac{x_5}{2} &\frac{x_1x_2}{12}+\frac{x_4}{2} & -\frac{x_3}{2}& -\frac{x_2}{2}& 0& 1
   \end{matrix}\right]\,.
   \end{eqnarray*}
  \subsection*{$(247E_1)$ }
 
 The following Lie algebra is denoted as $({247E_1})$ by Gong in \cite{Gong_Thesis}, and as $\mathcal{G}_{7,2,12}$ by Magnin in \cite{magnin}.
 
 The non-trivial brackets are the following:
\begin{equation*}
   [X_1, X_i] = X_{i+2}\,,\,2\le i\le 4\,,\,[X_2,X_4]=X_7\,,\,[X_3,X_4]=X_7\,.
\end{equation*}
 This is a nilpotent Lie algebra of rank 3 and step 3 that is stratifiable. The Lie brackets can be pictured with the diagram:
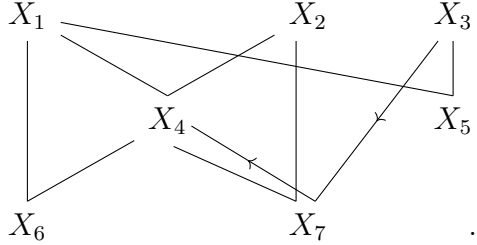
\begin{center}
 
	\begin{tikzcd}[end anchor=north]
		X_1 \ar[dr, no head]\ar[dd,no head]\ar[drrr, no head]& & X_2\ar[dl, no head]\ar[dd, no head, end anchor={[xshift=-3.3ex]north east},start anchor={[xshift=-3.3ex]south east}] &X_3\ar[d, no head]\ar[ddl, no head, end anchor={[xshift=-1.9ex]north east},->-=.5]\\
		&X_4 \ar[dl,no head]\ar[dr, no head, end anchor={[xshift=-3.3ex]north east},start anchor={[xshift=-2.ex]south east}]\ar[dr, no head, end anchor={[xshift=-1.9ex]north east},start anchor={[xshift=-0.7ex]south east},start anchor={[yshift=1.5ex]south east},-<-=.5] & & X_5\\
		X_6& &X_7 & \quad\;.
	\end{tikzcd}

\end{center}

 
The composition law \eqref{group law in G} of $(247E_1)$ is given by:

\begin{itemize}
    \item $z_1=x_1+y_1$;
    \item $z_2=x_2+y_2$;
    \item $z_3=x_3+y_3$;
    \item $z_4=x_4+y_4+\frac{1}{2}(x_1y_2-x_2y_1)$;
    \item $z_5=x_5+y_5+\frac{1}{2}(x_1y_3-x_3y_1)$;
    \item $z_6=x_6+y_6+\frac{1}{2}(x_1y_4-x_4y_1)+\frac{1}{12}(x_1-y_1)(x_1y_2-x_2y_1)$;
    \item $z_7=x_7+y_7+\frac{1}{2}(x_3y_4-x_4y_3+x_2y_4-x_4y_2)+\frac{1}{12}(x_3-y_3)(x_1y_2-x_2y_1)\salto+\frac{1}{12}(x_2-y_2)(x_1y_2-x_2y_1)$.
\end{itemize}



Since
\begin{eqnarray*}
  \mathrm{d}(L_\mathbf{x})_\mathbf{0}=\left[\begin{matrix}    
   1 & 0  & 0& 0& 0& 0& 0\\
   0 & 1& 0& 0& 0& 0& 0\\
   0 & 0& 1& 0& 0& 0& 0\\
   -\frac{x_2}{2} & \frac{x_1}{2}& 0& 1& 0& 0& 0\\
   -\frac{x_3}{2} &  0& \frac{x_1}{2}& 0& 1& 0& 0\\
 -\frac{x_1x_2}{12}-\frac{x_4}{2} &  \frac{x_1^2}{12}& 0& \frac{x_1}{2}& 0& 1& 0\\
-\frac{x_2x_3+x_2^2}{12}& \frac{x_1x_3+x_1x_2}{12}-\frac{x_4}{2} &-\frac{x_4}{2} & \frac{x_2+x_3}{2}& 0& 0& 1
   \end{matrix}\right]\,,
   \end{eqnarray*}
the left-invariant vector fields \eqref{leftinvariant vf} are: 
\begin{itemize}
\item $X_1=\partial_{x_1} -\frac{x_2}{2}\partial_{x_4}-\frac{x_3}{2}\partial_{x_5}-\big(\frac{x_1x_2}{12}+\frac{x_4}{2}\big)\partial_{x_6}-\frac{x_2x_3+x_2^2}{12}\partial_{x_7}\,;$
\item $X_2=\partial_{x_2}+\frac{x_1}{2}\partial_{x_4}+\frac{x_1^2}{12}\partial_{x_6}+\big(\frac{x_1x_3+x_1x_2}{12}-\frac{x_4}{2}\big)\partial_{x_7}\,;$
\item $X_3=\partial_{x_3}+\frac{x_1}{2}\partial_{x_5}-\frac{x_4}{2}\partial_{x_7}\,;$
\item $X_4=\partial_{x_4}+\frac{x_1}{2}\partial_{x_6}+\frac{x_2+x_3}{2}\partial_{x_7}\,;$
\item $X_5=\partial_{x_5}$;
\item $X_6=\partial_{x_6}$;
\item $X_7=\partial_{x_7}$,
\end{itemize}
and the respective left-invariant 1-forms \eqref{leftinvariant form} are: 
\begin{itemize}
\item $\theta_1=dx_1$;
\item $\theta_2=dx_2$;
\item $\theta_3=dx_3$;
\item $\theta_4=dx_4-\frac{x_1}{2}dx_2+\frac{x_2}{2}dx_1$;
\item $\theta_5=dx_5-\frac{x_1}{2}dx_3+\frac{x_3}{2}dx_1$;
\item $\theta_6=dx_6-\frac{x_1}{2}dx_4+\frac{x_1^2}{6}dx_2+\big(\frac{x_4}{2}-\frac{x_1x_2}{6}\big)dx_1$;
\item $\theta_7=dx_7-\frac{x_3+x_2}{2}dx_4+\frac{x_4}{2}dx_3+\big(\frac{x_1x_3+x_1x_2}{6}+\frac{x_4}{2}\big)dx_2-\frac{x_2x_3+x_2^2}{6}dx_1$.
\end{itemize}

Finally, we have

\begin{eqnarray*}
  \mathrm{d}(R_\mathbf{x})_\mathbf{0}=\left[\begin{matrix}    
   1 & 0  & 0& 0& 0& 0& 0\\
   0 & 1& 0& 0& 0& 0& 0\\
   0 & 0& 1& 0& 0& 0& 0\\
   \frac{x_2}{2} & -\frac{x_1}{2}& 0& 1& 0& 0& 0\\
   \frac{x_3}{2} &  0& -\frac{x_1}{2}& 0& 1& 0& 0\\
 \frac{x_4}{2}-\frac{x_1x_2}{12} &  \frac{x_1^2}{12}& 0& -\frac{x_1}{2}& 0& 1& 0\\
-\frac{x_2x_3+x_2^2}{12}& \frac{x_1x_3+x_1x_2}{12}+\frac{x_4}{2} &\frac{x_4}{2} & -\frac{x_2+x_3}{2}& 0& 0& 1
   \end{matrix}\right]\,.
   \end{eqnarray*}
  \subsection*{$(247F)$ }
 
 The following Lie algebra is denoted as $({247F})$ by Gong in \cite{Gong_Thesis}, and as $\mathcal{G}_{7,3,4}$ by Magnin in \cite{magnin}.
 
 The non-trivial brackets are the following:
\begin{equation*}
\begin{aligned}
   &[X_1, X_i] = X_{i+2}\,,\,i=2, 3\;,\,[X_2,X_4]=X_6\,,\\ [X_2 &,X_5]=X_7\,,\,[X_3,X_4]=X_7\,,\,[X_3,X_5]=X_6\,.
   \end{aligned}
\end{equation*}
 This is a nilpotent Lie algebra of rank 3 and step 3 that is stratifiable. The Lie brackets can be pictured with the diagram:
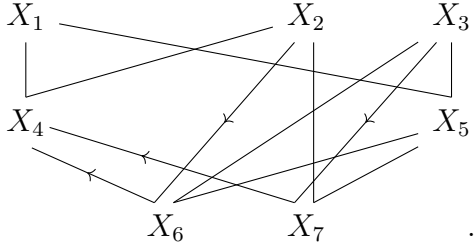
\begin{center}
 
	\begin{tikzcd}[end anchor=north]
		X_1 \ar[d, no head]\ar[drrr,no head]& & X_2\ar[dll, no head]\ar[ddl, no head, end anchor={[xshift=-3.3ex]north east},start anchor={[xshift=-3.3ex]south east},->-=.5]\ar[dd, no head, end anchor={[xshift=-1.9ex]north east},start anchor={[xshift=-1.9ex]south east}] &X_3\ar[d, no head]\ar[ddl,no head,end anchor={[xshift=-3.3ex]north east},->-=.5]\ar[ddll,no head,end anchor={[xshift=-1.9ex]north east} ]\\
		X_4 \ar[dr, no head, end anchor={[xshift=-3.3ex]north east},start anchor={[xshift=-2.ex]south east},-<-=.5]\ar[drr, no head, end anchor={[xshift=-3.3ex]north east},start anchor={[xshift=-0.7ex]south east},start anchor={[yshift=1.5ex]south east},-<-=.4]& & & X_5\ar[dl, no head, end anchor={[xshift=-1.9ex]north east}]\ar[dll,no head,end anchor={[xshift=-1.9ex]north east} ]\\
		& X_6&X_7 & \quad\;.
	\end{tikzcd}

\end{center}

 
The composition law \eqref{group law in G} of $(247F)$ is given by:

\begin{itemize}
    \item $z_1=x_1+y_1$;
    \item $z_2=x_2+y_2$;
    \item $z_3=x_3+y_3$;
    \item $z_4=x_4+y_4+\frac{1}{2}(x_1y_2-x_2y_1)$;
    \item $z_5=x_5+y_5+\frac{1}{2}(x_1y_3-x_3y_1)$;
    \item $z_6=x_6+y_6+\frac{1}{2}(x_2y_4-x_4y_2+x_3y_5-x_5y_3)+\frac{1}{12}(x_2-y_2)(x_1y_2-x_2y_1)\salto+\frac{1}{12}(x_3-y_3)(x_1y_3-x_3y_1)$;
    \item $z_7=x_7+y_7+\frac{1}{2}(x_3y_4-x_4y_3+x_2y_5-x_5y_2)+\frac{1}{12}(x_3-y_3)(x_1y_2-x_2y_1)\salto+\frac{1}{12}(x_2-y_2)(x_1y_3-x_3y_1)$.
\end{itemize}



Since
\begin{eqnarray*}
  \mathrm{d}(L_\mathbf{x})_\mathbf{0}=\left[\begin{matrix}    
    1 & 0  & 0& 0& 0& 0& 0\\
   0 & 1& 0& 0& 0& 0& 0\\
   0 & 0& 1& 0& 0& 0& 0\\
   -\frac{x_2}{2} & \frac{x_1}{2}& 0& 1& 0& 0& 0\\
   -\frac{x_3}{2} &  0& \frac{x_1}{2}& 0& 1& 0& 0\\
- \frac{x_2^2+x_3^2}{12} &  \frac{x_1x_2}{12}-\frac{x_4}{2}&\frac{x_1x_3}{12} -\frac{x_5}{2}& \frac{x_2}{2}& \frac{x_3}{2}& 1& 0\\
 -\frac{x_2x_3}{6}& \frac{x_1x_3}{12}-\frac{x_5}{2} &\frac{x_1x_2}{12}-\frac{x_4}{2} & \frac{x_3}{2}& \frac{x_2}{2}& 0& 1
   \end{matrix}\right]\,,
   \end{eqnarray*}
the induced left-invariant vector fields \eqref{leftinvariant vf} are: 
\begin{itemize}
\item $X_1=\partial_{x_1} -\frac{x_2}{2}\partial_{x_4}-\frac{x_3}{2}\partial_{x_5}-\frac{x_2^2+x_3^2}{12}\partial_{x_6}-\frac{x_2x_3}{6}\partial_{x_7}\,;$
\item $X_2=\partial_{x_2}+\frac{x_1}{2}\partial_{x_4}+\big(\frac{x_1x_2}{12}-\frac{x_4}{2}\big)\partial_{x_6}+\big(\frac{x_1x_3}{12}-\frac{x_5}{2}\big)\partial_{x_7}\,;$
\item $X_3=\partial_{x_3}+\frac{x_1}{2}\partial_{x_5}+\big(\frac{x_1x_3}{12}-\frac{x_5}{2}\big)\partial_{x_6}+\big(\frac{x_1x_2}{12}-\frac{x_4}{2}\big)\partial_{x_7}\,;$
\item $X_4=\partial_{x_4}+\frac{x_2}{2}\partial_{x_6}+\frac{x_3}{2}\partial_{x_7}\,;$
\item $X_5=\partial_{x_5}+\frac{x_3}{2}\partial_{x_6}+\frac{x_2}{2}\partial_{x_7}$;
\item $X_6=\partial_{x_6}$;
\item $X_7=\partial_{x_7}$,
\end{itemize}
and the respective left-invariant 1-forms \eqref{leftinvariant form} are: 
\begin{itemize}
\item $\theta_1=dx_1$;
\item $\theta_2=dx_2$;
\item $\theta_3=dx_3$;
\item $\theta_4=dx_4-\frac{x_1}{2}dx_2+\frac{x_2}{2}dx_1$;
\item $\theta_5=dx_5-\frac{x_1}{2}dx_3+\frac{x_3}{2}dx_1$;
\item $\theta_6=dx_6-\frac{x_3}{2}dx_5-\frac{x_2}{2}dx_4+\big(\frac{x_1x_3}{6}+\frac{x_5}{2}\big)dx_3+\big(\frac{x_1x_2}{6}+\frac{x_4}{2}\big)dx_2-\frac{x_2^2+x_3^2}{6}dx_1$;
\item $\theta_7=dx_7-\frac{x_2}{2}dx_5-\frac{x_3}{2}dx_4+\big(\frac{x_1x_2}{6}+\frac{x_4}{2}\big)dx_3+\big(\frac{x_1x_3}{6}+\frac{x_5}{2}\big)dx_2-\frac{x_2x_3}{3}dx_1$.
\end{itemize}

Finally, we have
\begin{eqnarray*}
  \mathrm{d}(R_\mathbf{x})_\mathbf{0}=\left[\begin{matrix}    
    1 & 0  & 0& 0& 0& 0& 0\\
   0 & 1& 0& 0& 0& 0& 0\\
   0 & 0& 1& 0& 0& 0& 0\\
   \frac{x_2}{2} & -\frac{x_1}{2}& 0& 1& 0& 0& 0\\
   \frac{x_3}{2} &  0& -\frac{x_1}{2}& 0& 1& 0& 0\\
- \frac{x_2^2+x_3^2}{12} &  \frac{x_1x_2}{12}+\frac{x_4}{2}&\frac{x_1x_3}{12} +\frac{x_5}{2}& -\frac{x_2}{2}& -\frac{x_3}{2}& 1& 0\\
 -\frac{x_2x_3}{6}& \frac{x_1x_3}{12}+\frac{x_5}{2} &\frac{x_1x_2}{12}+\frac{x_4}{2} & -\frac{x_3}{2}& -\frac{x_2}{2}& 0& 1
   \end{matrix}\right]\,.
   \end{eqnarray*}
  \subsection*{$(247F_1)$ }
 
 The following Lie algebra is denoted as $({247F_1})$ by Gong in \cite{Gong_Thesis}, and as $\mathcal{G}_{7,3,4}$ by Magnin in \cite{magnin}.
 
 The non-trivial brackets are the following:
\begin{equation*}
\begin{aligned}
   &[X_1, X_i] = X_{i+2}\,,\,i=2, 3\,,\,[X_2,X_4]=X_6\,,\\ [X_2&,X_5]=X_7\,,\,[X_3,X_4]=X_7\,,\,[X_3,X_5]=-X_6\,.
\end{aligned}
\end{equation*}
 This is a nilpotent Lie algebra of rank 3 and step 3 that is stratifiable. The Lie brackets can be pictured with the diagram:
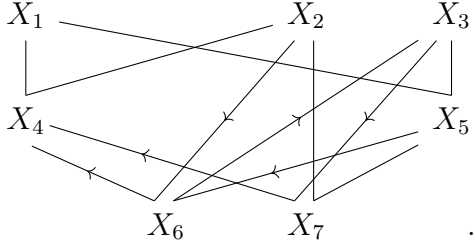
\begin{center}
 
	\begin{tikzcd}[end anchor=north]
		X_1 \ar[d, no head]\ar[drrr,no head]& & X_2\ar[dll, no head]\ar[ddl, no head, end anchor={[xshift=-3.3ex]north east},start anchor={[xshift=-3.3ex]south east},->-=.5]\ar[dd, no head, end anchor={[xshift=-1.9ex]north east},start anchor={[xshift=-1.9ex]south east}] &X_3\ar[d, no head]\ar[ddl,no head,end anchor={[xshift=-3.3ex]north east},->-=.5]\ar[ddll,no head,end anchor={[xshift=-1.9ex]north east},-<-=.5 ]\\
		X_4 \ar[dr, no head, end anchor={[xshift=-3.3ex]north east},start anchor={[xshift=-2.ex]south east},-<-=.5]\ar[drr, no head, end anchor={[xshift=-3.3ex]north east},start anchor={[xshift=-0.7ex]south east},start anchor={[yshift=1.5ex]south east},-<-=.4]& & & X_5\ar[dl, no head, end anchor={[xshift=-1.9ex]north east}]\ar[dll,no head,end anchor={[xshift=-1.9ex]north east},->-=.6 ]\\
		& X_6&X_7 & \quad\;.
	\end{tikzcd}

\end{center}

 
The composition law \eqref{group law in G} of $(247F_1)$ is given by:

\begin{itemize}
    \item $z_1=x_1+y_1$;
    \item $z_2=x_2+y_2$;
    \item $z_3=x_3+y_3$;
    \item $z_4=x_4+y_4+\frac{1}{2}(x_1y_2-x_2y_1)$;
    \item $z_5=x_5+y_5+\frac{1}{2}(x_1y_3-x_3y_1)$;
    \item $z_6=x_6+y_6+\frac{1}{2}(x_2y_4-x_4y_2-x_3y_5+x_5y_3)+\frac{1}{12}(x_2-y_2)(x_1y_2-x_2y_1)\salto+\frac{1}{12}(y_3-x_3)(x_1y_3-x_3y_1)$;
    \item $z_7=x_7+y_7+\frac{1}{2}(x_3y_4-x_4y_3+x_2y_5-x_5y_2)+\frac{1}{12}(x_3-y_3)(x_1y_2-x_2y_1)\salto +\frac{1}{12}(x_2-y_2)(x_1y_3-x_3y_1)$.
\end{itemize}



Since
\begin{eqnarray*}
  \mathrm{d}(L_\mathbf{x})_\mathbf{0}=\left[\begin{matrix}    
    1 & 0  & 0& 0& 0& 0& 0\\
   0 & 1& 0& 0& 0& 0& 0\\
   0 & 0& 1& 0& 0& 0& 0\\
   -\frac{x_2}{2} & \frac{x_1}{2}& 0& 1& 0& 0& 0\\
   -\frac{x_3}{2} &  0& \frac{x_1}{2}& 0& 1& 0& 0\\
\frac{x_3^2-x_2^2}{12} &  \frac{x_1x_2}{12}-\frac{x_4}{2}&\frac{x_5}{2}-\frac{x_1x_3}{12}& \frac{x_2}{2}& -\frac{x_3}{2}& 1& 0\\
 -\frac{x_2x_3}{6}& \frac{x_1x_3}{12}-\frac{x_5}{2} &\frac{x_1x_2}{12}-\frac{x_4}{2} & \frac{x_3}{2}& \frac{x_2}{2}& 0& 1
   \end{matrix}\right]\,,
   \end{eqnarray*}
the induced left-invariant vector fields \eqref{leftinvariant vf} are: 
\begin{itemize}
\item $X_1=\partial_{x_1} -\frac{x_2}{2}\partial_{x_4}-\frac{x_3}{2}\partial_{x_5}+\frac{x_3^2-x_2^2}{12}\partial_{x_6}-\frac{x_2x_3}{6}\partial_{x_7}\,;$
\item $X_2=\partial_{x_2}+\frac{x_1}{2}\partial_{x_4}+\big(\frac{x_1x_2}{12}-\frac{x_4}{2}\big)\partial_{x_6}+\big(\frac{x_1x_3}{12}-\frac{x_5}{2}\big)\partial_{x_7}\,;$
\item $X_3=\partial_{x_3}+\frac{x_1}{2}\partial_{x_5}+\big(\frac{x_5}{2}-\frac{x_1x_3}{12}\big)\partial_{x_6}+\big(\frac{x_1x_2}{12}-\frac{x_4}{2}\big)\partial_{x_7}\,;$
\item $X_4=\partial_{x_4}+\frac{x_2}{2}\partial_{x_6}+\frac{x_3}{2}\partial_{x_7}\,;$
\item $X_5=\partial_{x_5}-\frac{x_3}{2}\partial_{x_6}+\frac{x_2}{2}\partial_{x_7}$;
\item $X_6=\partial_{x_6}$;
\item $X_7=\partial_{x_7}$,
\end{itemize}
and the respective left-invariant 1-forms \eqref{leftinvariant form} are: 
\begin{itemize}
\item $\theta_1=dx_1$;
\item $\theta_2=dx_2$;
\item $\theta_3=dx_3$;
\item $\theta_4=dx_4-\frac{x_1}{2}dx_2+\frac{x_2}{2}dx_1$;
\item $\theta_5=dx_5-\frac{x_1}{2}dx_3+\frac{x_3}{2}dx_1$;
\item $\theta_6=dx_6+\frac{x_3}{2}dx_5-\frac{x_2}{2}dx_4-\big(\frac{x_1x_3}{6}+\frac{x_5}{2}\big)dx_3+\big(\frac{x_1x_2}{6}+\frac{x_4}{2}\big)dx_2+\frac{x_3^2-x_2^2}{6}dx_1$;
\item $\theta_7=dx_7-\frac{x_2}{2}dx_5-\frac{x_3}{2}dx_4+\big(\frac{x_1x_2}{6}+\frac{x_4}{2}\big)dx_3+\big(\frac{x_1x_3}{6}+\frac{x_5}{2}\big)dx_2-\frac{x_2x_3}{3}dx_1$.
\end{itemize}

Finally, we have
\begin{eqnarray*}
  \mathrm{d}(R_\mathbf{x})_\mathbf{0}=\left[\begin{matrix}    
    1 & 0  & 0& 0& 0& 0& 0\\
   0 & 1& 0& 0& 0& 0& 0\\
   0 & 0& 1& 0& 0& 0& 0\\
   \frac{x_2}{2} & -\frac{x_1}{2}& 0& 1& 0& 0& 0\\
   \frac{x_3}{2} &  0& -\frac{x_1}{2}& 0& 1& 0& 0\\
\frac{x_3^2-x_2^2}{12} &  \frac{x_1x_2}{12}+\frac{x_4}{2}&-\frac{x_1x_3}{12} -\frac{x_5}{2}& -\frac{x_2}{2}& \frac{x_3}{2}& 1& 0\\
 -\frac{x_2x_3}{6}& \frac{x_1x_3}{12}+\frac{x_5}{2} &\frac{x_1x_2}{12}+\frac{x_4}{2} & -\frac{x_3}{2}& -\frac{x_2}{2}& 0& 1
   \end{matrix}\right]\,.
   \end{eqnarray*}
  \subsection*{$(247G)$ }
 
 The following Lie algebra is denoted as $({247G})$ by Gong in \cite{Gong_Thesis}, and as $\mathcal{G}_{7,2,34}$ by Magnin in \cite{magnin}.
 
 The non-trivial brackets are the following:
\begin{equation*}
\begin{aligned}
   &[X_1, X_i] = X_{i+2}\,,\,2\le i\le 4\,,\,[X_1,X_5]=X_6\,,\,[X_2,X_4]=X_6\,, \\ & \phantom{===}[X_3,X_5]=X_6\,, \, [X_2,X_5]=X_7\,,\,[X_3,X_4]=X_7\,.
\end{aligned}
\end{equation*}
This is a nilpotent Lie algebra of rank 3 and step 3 that is stratifiable. The Lie brackets can be pictured with the diagram:
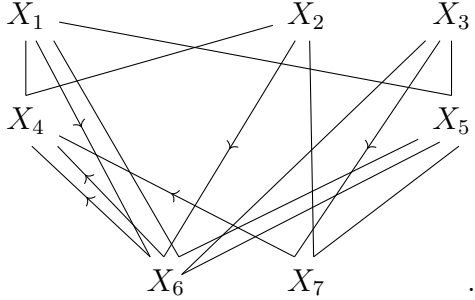
\begin{center}
 
	\begin{tikzcd}[end anchor=north]
		X_1 \ar[d, no head]\ar[dddr,no head,->-=.4,end anchor={[xshift=-3.6ex]north east},start anchor={[xshift=-1.7ex]south east}]\ar[drrr, no head]\ar[dddr,no head, end anchor={[xshift=-1.5ex]north east},start anchor={[xshift=-0.4ex]south east}]& & X_2\ar[dll, no head]\ar[ddd, no head, end anchor={[xshift=-1.9ex]north east}]\ar[dddl,no head,->-=.5,end anchor={[xshift=-2.6ex]north east}] &X_3\ar[d, no head]\ar[dddl, no head, end anchor={[xshift=-3.3ex]north east},->-=.5]\ar[dddll, no head, end anchor={[xshift=-1.3ex]north east}, end anchor={[yshift=-1.3ex]north east}]\\
		X_4 \ar[ddr,no head,-<-=.3,end anchor={[xshift=-2.6ex]north east}]\ar[ddr,no head,-<-=.5,end anchor={[xshift=-3.6ex]north east},start anchor={[xshift=-2.ex]south east}]\ar[ddrr,no head,-<-=.5,end anchor={[xshift=-3.3ex]north east}]& & & X_5\ar[ddl, no head, end anchor={[xshift=-1.9ex]north east},start anchor={[xshift=-1.9ex]south east}]\ar[ddll,no head, end anchor={[xshift=-1.5ex]north east}]\ar[ddll, no head, end anchor={[xshift=-1.3ex]north east}, end anchor={[yshift=-1.3ex]north east},start anchor={[yshift=0.3ex]south east},start anchor={[xshift=-3.3ex]south east}]\\
		& & &\\
		& X_6&X_7 & \quad\;.
	\end{tikzcd}

\end{center}

 
The composition law \eqref{group law in G} of $(247G)$ is given by:

\begin{itemize}
    \item $z_1=x_1+y_1$;
    \item $z_2=x_2+y_2$;
    \item $z_3=x_3+y_3$;
    \item $z_4=x_4+y_4+\frac{1}{2}(x_1y_2-x_2y_1)$;
    \item $z_5=x_5+y_5+\frac{1}{2}(x_1y_3-x_3y_1)$;
    \item $z_6=x_6+y_6+\frac{1}{2}(x_1y_4-x_4y_1+x_1y_5-x_5y_1+x_2y_4-x_4y_2+x_3y_5-x_5y_3)\salto+\frac{1}{12}(x_1-y_1+x_2-y_2)(x_1y_2-x_2y_1)+\frac{1}{12}(x_1-y_1+x_3-y_3)(x_1y_3-x_3y_1)$;
    \item $z_7=x_7+y_7+\frac{1}{2}(x_3y_4-x_4y_3+x_2y_5-x_5y_2)+\frac{1}{12}(x_3-y_3)(x_1y_2-x_2y_1)\salto+\frac{1}{12}(x_2-y_2)(x_1y_3-x_3y_1)$.
\end{itemize}

Since
\begin{eqnarray*}
  \mathrm{d}(L_\mathbf{x})_\mathbf{0}=\left[\begin{matrix}    
   1 & 0  & 0& 0& 0& 0& 0\\
   0 & 1& 0& 0& 0& 0& 0\\
   0 & 0& 1& 0& 0& 0& 0\\
   -\frac{x_2}{2} & \frac{x_1}{2}& 0& 1& 0& 0& 0\\
   -\frac{x_3}{2} &  0& \frac{x_1}{2}& 0& 1& 0& 0\\
-\frac{x_1x_2+x_1x_3+x_2^2+x_3^2}{12}-\frac{x_4+x_5}{2} &  \frac{x_1^2+x_1x_2}{12}-\frac{x_4}{2}& \frac{x_1^2+x_1x_3}{12}-\frac{x_5}{2}
 & \frac{x_1+x_2}{2}& \frac{x_1+x_3}{2}& 1& 0\\
 -\frac{x_2x_3}{6}& \frac{x_1x_3}{12}-\frac{x_5}{2} &\frac{x_1x_2}{12}-\frac{x_4}{2} & \frac{x_3}{2}& \frac{x_2}{2}& 0& 1
   \end{matrix}\right]\,,
   \end{eqnarray*}
the induced left-invariant vector fields \eqref{leftinvariant vf} are: 
\begin{itemize}
\item $X_1=\partial_{x_1} -\frac{x_2}{2}\partial_{x_4}-\frac{x_3}{2}\partial_{x_5}-\big(\frac{x_1x_2+x_1x_3+x_2^2+x_3^2}{12}+\frac{x_4+x_5}{2}\big)\partial_{x_6}-\frac{x_2x_3}{6}\partial_{x_7}\,;$
\item $X_2=\partial_{x_2}+\frac{x_1}{2}\partial_{x_4}+\big(\frac{x_1^2+x_1x_2}{12}-\frac{x_4}{2}\big)\partial_{x_6}+\big(\frac{x_1x_3}{12}-\frac{x_5}{2}\big)\partial_{x_7}\,;$
\item $X_3=\partial_{x_3}+\frac{x_1}{2}\partial_{x_5}+\big(\frac{x_1^2+x_1x_3}{12}-\frac{x_5}{2}\big)\partial_{x_6}+\big(\frac{x_1x_2}{12}-\frac{x_4}{2}\big)\partial_{x_7}\,;$
\item $X_4=\partial_{x_4}+\frac{x_1+x_2}{2}\partial_{x_6}+\frac{x_3}{2}\partial_{x_7}\,;$
\item $X_5=\partial_{x_5}+\frac{x_1+x_3}{2}\partial_{x_6}+\frac{x_2}{2}\partial_{x_7}$;
\item $X_6=\partial_{x_6}$;
\item $X_7=\partial_{x_7}$,
\end{itemize}
and the respective left-invariant 1-forms \eqref{leftinvariant form} are: 
\begin{itemize}
\item $\theta_1=dx_1$;
\item $\theta_2=dx_2$;
\item $\theta_3=dx_3$;
\item $\theta_4=dx_4-\frac{x_1}{2}dx_2+\frac{x_2}{2}dx_1$;
\item $\theta_5=dx_5-\frac{x_1}{2}dx_3+\frac{x_3}{2}dx_1$;
\item $\theta_6=dx_6-\frac{x_1+x_3}{2}dx_5-\frac{x_1+x_2}{2}dx_4+\big(\frac{x_1^2+x_1x_3}{6}+\frac{x_5}{2}\big)dx_3+\big(\frac{x_1^2+x_1x_2}{6}+\frac{x_4}{2}\big)dx_2\saltot+\big(\frac{x_4+x_5}{2}-\frac{x_1x_2+x_1x_3+x_2^2+x_3^2}{6}\big)dx_1$;
\item $\theta_7=dx_7-\frac{x_2}{2}dx_5-\frac{x_3}{2}dx_4+\big(\frac{x_1x_2}{6}+\frac{x_4}{2}\big)dx_3+\big(\frac{x_1x_3}{6}+\frac{x_5}{2}\big)dx_2-\frac{x_2x_3}{3}dx_1$.
\end{itemize}

Finally, we have
\begin{eqnarray*}
  \mathrm{d}(R_\mathbf{x})_\mathbf{0}=\left[\begin{matrix}    
   1 & 0  & 0& 0& 0& 0& 0\\
   0 & 1& 0& 0& 0& 0& 0\\
   0 & 0& 1& 0& 0& 0& 0\\
   \frac{x_2}{2} & -\frac{x_1}{2}& 0& 1& 0& 0& 0\\
   \frac{x_3}{2} &  0& -\frac{x_1}{2}& 0& 1& 0& 0\\
\frac{x_4+x_5}{2}-\frac{x_1x_2+x_1x_3+x_2^2+x_3^2}{12} &  \frac{x_1^2+x_1x_2}{12}+\frac{x_4}{2}& \frac{x_1^2+x_1x_3}{12}+\frac{x_5}{2}
 & -\frac{x_1+x_2}{2}& -\frac{x_1+x_3}{2}& 1& 0\\
 -\frac{x_2x_3}{6}& \frac{x_1x_3}{12}+\frac{x_5}{2} &\frac{x_1x_2}{12}+\frac{x_4}{2} & -\frac{x_3}{2}& -\frac{x_2}{2}& 0& 1
   \end{matrix}\right]\,.
   \end{eqnarray*}

  \subsection*{$(247H)$ }
 
 The following Lie algebra is denoted as $({247H})$ by Gong in \cite{Gong_Thesis}, and as $\mathcal{G}_{7,1,19}$ by Magnin in \cite{magnin}.
 
 The non-trivial brackets are the following:
\begin{equation*}
\begin{aligned}
   &[X_1, X_i] = X_{i+2}\,,\,2\le i\le 4\,,\,[X_2,X_4]=X_6\,, \\ [&X_3,X_5]=X_6\,, \, [X_2,X_5]=X_7\,,\,[X_3,X_4]=X_7\,.
\end{aligned}
\end{equation*}
This is a nilpotent Lie algebra of rank 3 and step 3 that is stratifiable. The Lie brackets can be pictured with the diagram:
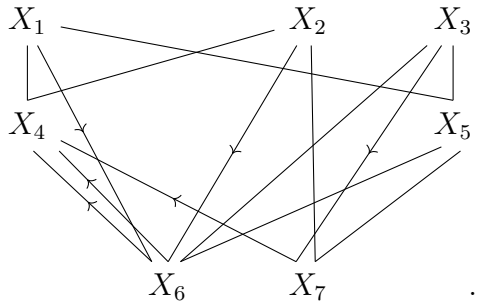
\begin{center}
 
	\begin{tikzcd}[end anchor=north]
		X_1 \ar[d, no head]\ar[dddr,no head,->-=.4,end anchor={[xshift=-3.6ex]north east},start anchor={[xshift=-1.7ex]south east}]\ar[drrr, no head]& & X_2\ar[dll, no head]\ar[ddd, no head, end anchor={[xshift=-1.9ex]north east}]\ar[dddl,no head,->-=.5,end anchor={[xshift=-2.4ex]north east}] &X_3\ar[d, no head]\ar[dddl, no head, end anchor={[xshift=-3.3ex]north east},->-=.5]\ar[dddll, no head, end anchor={[xshift=-1.5ex]north east}]\\
		X_4 \ar[ddr,no head,-<-=.3,end anchor={[xshift=-2.4ex]north east}]\ar[ddr,no head,-<-=.5,end anchor={[xshift=-3.6ex]north east},start anchor={[xshift=-2.ex]south east}]\ar[ddrr,no head,-<-=.5,end anchor={[xshift=-3.3ex]north east}]& & & X_5\ar[ddl, no head, end anchor={[xshift=-1.9ex]north east},start anchor={[xshift=-1.9ex]south east}]\ar[ddll, no head, end anchor={[xshift=-1.5ex]north east},start anchor={[yshift=0.3ex]south east},start anchor={[xshift=-3.3ex]south east}]\\
		& & &\\
		& X_6&X_7 & \quad\;.
	\end{tikzcd}

\end{center}

 
The composition law \eqref{group law in G} of $(247H)$ is given by:

\begin{itemize}
    \item $z_1=x_1+y_1$;
    \item $z_2=x_2+y_2$;
    \item $z_3=x_3+y_3$;
    \item $z_4=x_4+y_4+\frac{1}{2}(x_1y_2-x_2y_1)$;
    \item $z_5=x_5+y_5+\frac{1}{2}(x_1y_3-x_3y_1)$;
    \item $z_6=x_6+y_6+\frac{1}{2}(x_1y_4-x_4y_1+x_2y_4-x_4y_2+x_3y_5-x_5y_3)+\frac{1}{12}(x_3-y_3)(x_1y_3-x_3y_1)\salto+\frac{1}{12}(x_1-y_1+x_2-y_2)(x_1y_2-x_2y_1)$;
    \item $z_7=x_7+y_7+\frac{1}{2}(x_3y_4-x_4y_3+x_2y_5-x_5y_2)+\frac{1}{12}(x_3-y_3)(x_1y_2-x_2y_1)\salto+\frac{1}{12}(x_2-y_2)(x_1y_3-x_3y_1)$.
\end{itemize}

Since
\begin{eqnarray*}
  \mathrm{d}(L_\mathbf{x})_\mathbf{0}=\left[\begin{matrix}    
   1 & 0  & 0& 0& 0& 0& 0\\
   0 & 1& 0& 0& 0& 0& 0\\
   0 & 0& 1& 0& 0& 0& 0\\
   -\frac{x_2}{2} & \frac{x_1}{2}& 0& 1& 0& 0& 0\\
   -\frac{x_3}{2} &  0& \frac{x_1}{2}& 0& 1& 0& 0\\
-\frac{x_1x_2+x_2^2+x_3^2}{12}-\frac{x_4}{2} &  \frac{x_1^2+x_1x_2}{12}-\frac{x_4}{2}& \frac{x_1x_3}{12}-\frac{x_5}{2}
 & \frac{x_1+x_2}{2}& \frac{x_3}{2}& 1& 0\\
 -\frac{x_2x_3}{6}& \frac{x_1x_3}{12}-\frac{x_5}{2} &\frac{x_1x_2}{12}-\frac{x_4}{2} & \frac{x_3}{2}& \frac{x_2}{2}& 0& 1
   \end{matrix}\right]\,,
   \end{eqnarray*}
the induced left-invariant vector fields \eqref{leftinvariant vf} are: 
\begin{itemize}
\item $X_1=\partial_{x_1} -\frac{x_2}{2}\partial_{x_4}-\frac{x_3}{2}\partial_{x_5}-\big(\frac{x_1x_2+x_2^2+x_3^2}{12}+\frac{x_4}{2}\big)\partial_{x_6}-\frac{x_2x_3}{6}\partial_{x_7}\,;$
\item $X_2=\partial_{x_2}+\frac{x_1}{2}\partial_{x_4}+\big(\frac{x_1^2+x_1x_2}{12}-\frac{x_4}{2}\big)\partial_{x_6}+\big(\frac{x_1x_3}{12}-\frac{x_5}{2}\big)\partial_{x_7}\,;$
\item $X_3=\partial_{x_3}+\frac{x_1}{2}\partial_{x_5}+\big(\frac{x_1x_3}{12}-\frac{x_5}{2}\big)\partial_{x_6}+\big(\frac{x_1x_2}{12}-\frac{x_4}{2}\big)\partial_{x_7}\,;$
\item $X_4=\partial_{x_4}+\frac{x_1+x_2}{2}\partial_{x_6}+\frac{x_3}{2}\partial_{x_7}\,;$
\item $X_5=\partial_{x_5}+\frac{x_3}{2}\partial_{x_6}+\frac{x_2}{2}\partial_{x_7}$;
\item $X_6=\partial_{x_6}$;
\item $X_7=\partial_{x_7}$,
\end{itemize}
and the respective left-invariant 1-forms \eqref{leftinvariant form} are: 
\begin{itemize}
\item $\theta_1=dx_1$;
\item $\theta_2=dx_2$;
\item $\theta_3=dx_3$;
\item $\theta_4=dx_4-\frac{x_1}{2}dx_2+\frac{x_2}{2}dx_1$;
\item $\theta_5=dx_5-\frac{x_1}{2}dx_3+\frac{x_3}{2}dx_1$;
\item $\theta_6=dx_6-\frac{x_3}{2}dx_5-\frac{x_1+x_2}{2}dx_4+\big(\frac{x_1x_3}{6}+\frac{x_5}{2}\big)dx_3+\big(\frac{x_1^2+x_1x_2}{6}+\frac{x_4}{2}\big)dx_2+\big(\frac{x_4}{2}\saltot-\frac{x_1x_2+x_2^2+x_3^2}{6}\big)dx_1$;
\item $\theta_7=dx_7-\frac{x_2}{2}dx_5-\frac{x_3}{2}dx_4+\big(\frac{x_1x_2}{6}+\frac{x_4}{2}\big)dx_3+\big(\frac{x_1x_3}{6}+\frac{x_5}{2}\big)dx_2-\frac{x_2x_3}{3}dx_1$.
\end{itemize}

Finally, we have
\begin{eqnarray*}
  \mathrm{d}(R_\mathbf{x})_\mathbf{0}=\left[\begin{matrix}    
   1 & 0  & 0& 0& 0& 0& 0\\
   0 & 1& 0& 0& 0& 0& 0\\
   0 & 0& 1& 0& 0& 0& 0\\
   \frac{x_2}{2} & -\frac{x_1}{2}& 0& 1& 0& 0& 0\\
   \frac{x_3}{2} &  0& -\frac{x_1}{2}& 0& 1& 0& 0\\
\frac{x_4}{2}-\frac{x_1x_2+x_2^2+x_3^2}{12} &  \frac{x_1^2+x_1x_2}{12}+\frac{x_4}{2}& \frac{x_1x_3}{12}+\frac{x_5}{2}
 & -\frac{x_1+x_2}{2}& -\frac{x_3}{2}& 1& 0\\
 -\frac{x_2x_3}{6}& \frac{x_1x_3}{12}+\frac{x_5}{2} &\frac{x_1x_2}{12}+\frac{x_4}{2} & -\frac{x_3}{2}& -\frac{x_2}{2}& 0& 1
   \end{matrix}\right]\,.
   \end{eqnarray*}
  \subsection*{$(247H_1)$ }
 
 The following Lie algebra is denoted as $({247H_1})$ by Gong in \cite{Gong_Thesis}, and as $\mathcal{G}_{7,1,19}$ by Magnin in \cite{magnin}.
 
 The non-trivial brackets are the following:
\begin{equation*}
\begin{aligned}
   &[X_1, X_i] = X_{i+2}\;,\,2\le i\le 4\,,\,[X_2,X_4]=X_6\,, \\ [&X_3,X_5]=-X_6\,, \, [X_2,X_5]=X_7\,,\,[X_3,X_4]=X_7\,.
\end{aligned}
\end{equation*}
This is a nilpotent Lie algebra of rank 3 and step 3 that is stratifiable. The Lie brackets can be pictured with the diagram:
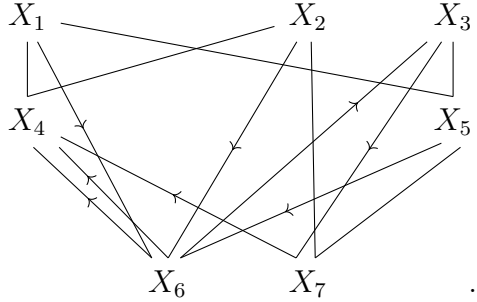
\begin{center}
 
	\begin{tikzcd}[end anchor=north]
		X_1 \ar[d, no head]\ar[dddr,no head,->-=.4,end anchor={[xshift=-3.6ex]north east},start anchor={[xshift=-1.7ex]south east}]\ar[drrr, no head]& & X_2\ar[dll, no head]\ar[ddd, no head, end anchor={[xshift=-1.9ex]north east}]\ar[dddl,no head,->-=.5,end anchor={[xshift=-2.4ex]north east}] &X_3\ar[d, no head]\ar[dddl, no head, end anchor={[xshift=-3.3ex]north east},->-=.5]\ar[dddll, no head, end anchor={[xshift=-1.5ex]north east},-<-=.3]\\
		X_4 \ar[ddr,no head,-<-=.3,end anchor={[xshift=-2.4ex]north east}]\ar[ddr,no head,-<-=.5,end anchor={[xshift=-3.6ex]north east},start anchor={[xshift=-2.ex]south east}]\ar[ddrr,no head,-<-=.5,end anchor={[xshift=-3.3ex]north east}]& & & X_5\ar[ddl, no head, end anchor={[xshift=-1.9ex]north east},start anchor={[xshift=-1.9ex]south east}]\ar[ddll, no head,->-=.6, end anchor={[xshift=-1.5ex]north east},start anchor={[yshift=0.3ex]south east},start anchor={[xshift=-3.3ex]south east}]\\
		& & &\\
		& X_6&X_7 & \quad\;.
	\end{tikzcd}

\end{center}

 
The composition law \eqref{group law in G} of $(247H_1)$ is given by:

\begin{itemize}
    \item $z_1=x_1+y_1$;
    \item $z_2=x_2+y_2$;
    \item $z_3=x_3+y_3$;
    \item $z_4=x_4+y_4+\frac{1}{2}(x_1y_2-x_2y_1)$;
    \item $z_5=x_5+y_5+\frac{1}{2}(x_1y_3-x_3y_1)$;
    \item $z_6=x_6+y_6+\frac{1}{2}(x_1y_4-x_4y_1+x_2y_4-x_4y_2-x_3y_5+x_5y_3)+\frac{1}{12}(y_3-x_3)(x_1y_3-x_3y_1)\salto+\frac{1}{12}(x_1-y_1+x_2-y_2)(x_1y_2-x_2y_1)$;
    \item $z_7=x_7+y_7+\frac{1}{2}(x_3y_4-x_4y_3+x_2y_5-x_5y_2)+\frac{1}{12}(x_3-y_3)(x_1y_2-x_2y_1)\salto+\frac{1}{12}(x_2-y_2)(x_1y_3-x_3y_1)$.
\end{itemize}

Since
\begin{eqnarray*}
 \mathrm{d}(L_\mathbf{x})_\mathbf{0}= \left[\begin{matrix}    
   1 & 0  & 0& 0& 0& 0& 0\\
   0 & 1& 0& 0& 0& 0& 0\\
   0 & 0& 1& 0& 0& 0& 0\\
   -\frac{x_2}{2} & \frac{x_1}{2}& 0& 1& 0& 0& 0\\
   -\frac{x_3}{2} &  0& \frac{x_1}{2}& 0& 1& 0& 0\\
\frac{x_3^2-x_1x_2-x_2^2}{12}-\frac{x_4}{2} &  \frac{x_1^2+x_1x_2}{12}-\frac{x_4}{2}& \frac{x_5}{2}-\frac{x_1x_3}{12}
 & \frac{x_1+x_2}{2}& -\frac{x_3}{2}& 1& 0\\
 -\frac{x_2x_3}{6}& \frac{x_1x_3}{12}-\frac{x_5}{2} &\frac{x_1x_2}{12}-\frac{x_4}{2} & \frac{x_3}{2}& \frac{x_2}{2}& 0& 1
   \end{matrix}\right]\,,
   \end{eqnarray*}
the induced left-invariant vector fields \eqref{leftinvariant vf} are: 
\begin{itemize}
\item $X_1=\partial_{x_1} -\frac{x_2}{2}\partial_{x_4}-\frac{x_3}{2}\partial_{x_5}+\big(\frac{x_3^2-x_1x_2-x_2^2}{12}-\frac{x_4}{2}\big)\partial_{x_6}-\frac{x_2x_3}{6}\partial_{x_7}\,;$
\item $X_2=\partial_{x_2}+\frac{x_1}{2}\partial_{x_4}+\big(\frac{x_1^2+x_1x_2}{12}-\frac{x_4}{2}\big)\partial_{x_6}+\big(\frac{x_1x_3}{12}-\frac{x_5}{2}\big)\partial_{x_7}\,;$
\item $X_3=\partial_{x_3}+\frac{x_1}{2}\partial_{x_5}+\big(\frac{x_5}{2}-\frac{x_1x_3}{12}\big)\partial_{x_6}+\big(\frac{x_1x_2}{12}-\frac{x_4}{2}\big)\partial_{x_7}\,;$
\item $X_4=\partial_{x_4}+\frac{x_1+x_2}{2}\partial_{x_6}+\frac{x_3}{2}\partial_{x_7}\,;$
\item $X_5=\partial_{x_5}-\frac{x_3}{2}\partial_{x_6}+\frac{x_2}{2}\partial_{x_7}$;
\item $X_6=\partial_{x_6}$;
\item $X_7=\partial_{x_7}$,
\end{itemize}
and the respective left-invariant 1-forms \eqref{leftinvariant form} are: 
\begin{itemize}
\item $\theta_1=dx_1$;
\item $\theta_2=dx_2$;
\item $\theta_3=dx_3$;
\item $\theta_4=dx_4-\frac{x_1}{2}dx_2+\frac{x_2}{2}dx_1$;
\item $\theta_5=dx_5-\frac{x_1}{2}dx_3+\frac{x_3}{2}dx_1$;
\item $\theta_6=dx_6+\frac{x_3}{2}dx_5-\frac{x_1+x_2}{2}dx_4-\big(\frac{x_1x_3}{6}+\frac{x_5}{2}\big)dx_3+\big(\frac{x_1^2+x_1x_2}{6}+\frac{x_4}{2}\big)dx_2+\big(\frac{x_4}{2}\saltot+\frac{x_3^2-x_1x_2-x_2^2}{6}\big)dx_1$;
\item $\theta_7=dx_7-\frac{x_2}{2}dx_5-\frac{x_3}{2}dx_4+\big(\frac{x_1x_2}{6}+\frac{x_4}{2}\big)dx_3+\big(\frac{x_1x_3}{6}+\frac{x_5}{2}\big)dx_2-\frac{x_2x_3}{3}dx_1$.
\end{itemize}

Finally, we have
\begin{eqnarray*}
 \mathrm{d}(R_\mathbf{x})_\mathbf{0}= \left[\begin{matrix}    
   1 & 0  & 0& 0& 0& 0& 0\\
   0 & 1& 0& 0& 0& 0& 0\\
   0 & 0& 1& 0& 0& 0& 0\\
   \frac{x_2}{2} & -\frac{x_1}{2}& 0& 1& 0& 0& 0\\
   \frac{x_3}{2} &  0& -\frac{x_1}{2}& 0& 1& 0& 0\\
\frac{x_3^2-x_1x_2-x_2^2}{12}+\frac{x_4}{2} &  \frac{x_1^2+x_1x_2}{12}+\frac{x_4}{2}& -\frac{x_5}{2}-\frac{x_1x_3}{12}
 & -\frac{x_1+x_2}{2}& \frac{x_3}{2}& 1& 0\\
 -\frac{x_2x_3}{6}& \frac{x_1x_3}{12}+\frac{x_5}{2} &\frac{x_1x_2}{12}+\frac{x_4}{2} & -\frac{x_3}{2}& -\frac{x_2}{2}& 0& 1
   \end{matrix}\right]\,.
   \end{eqnarray*}
  \subsection*{$(247I)$ }
 
 The following Lie algebra is denoted as $({247I})$ by Gong in \cite{Gong_Thesis}, and as $\mathcal{G}_{7,3,5}$ by Magnin in \cite{magnin}.
 
 The non-trivial brackets are the following:
\begin{equation*}
   [X_1, X_i] = X_{i+2}\,,\,i=2, 3\,,\,[X_2,X_5]=X_6\,,\,[X_3,X_4]=X_6\,,\,[X_3,X_5]=X_7\,.
\end{equation*}
 This is a nilpotent Lie algebra of rank 3 and step 3 that is stratifiable. The Lie brackets can be pictured with the diagram:
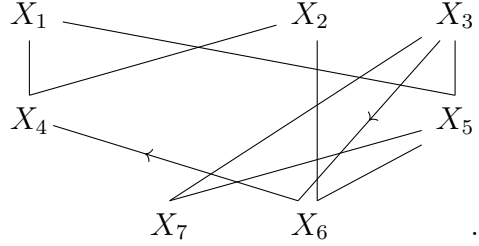
\begin{center}
 
	\begin{tikzcd}[end anchor=north]
		X_1 \ar[d, no head]\ar[drrr,no head]& & X_2\ar[dll, no head]\ar[dd, no head, end anchor={[xshift=-1.9ex]north east},start anchor={[xshift=-1.9ex]south east}] &X_3\ar[d, no head]\ar[ddl,no head,end anchor={[xshift=-3.3ex]north east},->-=.5]\ar[ddll,no head ]\\
		X_4 \ar[drr, no head, end anchor={[xshift=-3.3ex]north east},start anchor={[xshift=-0.7ex]south east},start anchor={[yshift=1.5ex]south east},-<-=.4]& & & X_5\ar[dl, no head, end anchor={[xshift=-1.9ex]north east}]\ar[dll,no head ]\\
		& X_7&X_6 & \quad\;.
	\end{tikzcd}

\end{center}

 
The composition law \eqref{group law in G} of $(247I)$ is given by:

\begin{itemize}
    \item $z_1=x_1+y_1$;
    \item $z_2=x_2+y_2$;
    \item $z_3=x_3+y_3$;
    \item $z_4=x_4+y_4+\frac{1}{2}(x_1y_2-x_2y_1)$;
    \item $z_5=x_5+y_5+\frac{1}{2}(x_1y_3-x_3y_1)$;
    \item $z_6=x_6+y_6+\frac{1}{2}(x_3y_4-x_4y_3+x_2y_5-x_5y_2)+\frac{1}{12}(x_3-y_3)(x_1y_2-x_2y_1)\salto+\frac{1}{12}(x_2-y_2)(x_1y_3-x_3y_1)$;
    \item $z_7=x_7+y_7+\frac{1}{2}(x_3y_5-x_5y_3)+\frac{1}{12}(x_3-y_3)(x_1y_3-x_3y_1)$.
\end{itemize}


 

Since
\begin{eqnarray*}
  \mathrm{d}(L_\mathbf{x})_\mathbf{0}=\left[\begin{matrix}    
     1 & 0  & 0& 0& 0& 0& 0\\
   0 & 1& 0& 0& 0& 0& 0\\
   0 & 0& 1& 0& 0& 0& 0\\
   -\frac{x_2}{2} & \frac{x_1}{2}& 0& 1& 0& 0& 0\\
   -\frac{x_3}{2} &  0& \frac{x_1}{2}& 0& 1& 0& 0\\
   -\frac{x_2x_3}{6}& \frac{x_1x_3}{12}-\frac{x_5}{2} &\frac{x_1x_2}{12}-\frac{x_4}{2} & \frac{x_3}{2}& \frac{x_2}{2}& 1& 0\\
 -\frac{x_3^2}{12} &  0&\frac{x_1x_3}{12} -\frac{x_5}{2}& 0& \frac{x_3}{2}& 0& 1
   \end{matrix}\right]\,,
   \end{eqnarray*}
the induced left-invariant vector fields \eqref{leftinvariant vf} are: 
\begin{itemize}
\item $X_1=\partial_{x_1} -\frac{x_2}{2}\partial_{x_4}-\frac{x_3}{2}\partial_{x_5}-\frac{x_2x_3}{6}\partial_{x_6}-\frac{x_3^2}{12}\partial_{x_7}\,;$
\item $X_2=\partial_{x_2}+\frac{x_1}{2}\partial_{x_4}+\big(\frac{x_1x_3}{12}-\frac{x_5}{2}\big)\partial_{x_6}\,;$
\item $X_3=\partial_{x_3}+\frac{x_1}{2}\partial_{x_5}+\big(\frac{x_1x_2}{12}-\frac{x_4}{2}\big)\partial_{x_6}+\big(\frac{x_1x_3}{12}-\frac{x_5}{2}\big)\partial_{x_7}\,;$
\item $X_4=\partial_{x_4}+\frac{x_3}{2}\partial_{x_6}\,;$
\item $X_5=\partial_{x_5}+\frac{x_2}{2}\partial_{x_6}+\frac{x_3}{2}\partial_{x_7}$;
\item $X_6=\partial_{x_6}$;
\item $X_7=\partial_{x_7}$,
\end{itemize}
and the respective left-invariant 1-forms \eqref{leftinvariant form} are: 
\begin{itemize}
\item $\theta_1=dx_1$;
\item $\theta_2=dx_2$;
\item $\theta_3=dx_3$;
\item $\theta_4=dx_4-\frac{x_1}{2}dx_2+\frac{x_2}{2}dx_1$;
\item $\theta_5=dx_5-\frac{x_1}{2}dx_3+\frac{x_3}{2}dx_1$;
\item $\theta_6=dx_6-\frac{x_2}{2}dx_5-\frac{x_3}{2}dx_4+\big(\frac{x_1x_2}{6}+\frac{x_4}{2}\big)dx_3+\big(\frac{x_1x_3}{6}+\frac{x_5}{2}\big)dx_2-\frac{x_2x_3}{3}dx_1$;
\item $\theta_7=dx_7-\frac{x_3}{2}dx_5+\big(\frac{x_1x_3}{6}+\frac{x_5}{2}\big)dx_3-\frac{x_3^2}{6}dx_1$.
\end{itemize}

Finally, we have
\begin{eqnarray*}
  \mathrm{d}(R_\mathbf{x})_\mathbf{0}=\left[\begin{matrix}    
     1 & 0  & 0& 0& 0& 0& 0\\
   0 & 1& 0& 0& 0& 0& 0\\
   0 & 0& 1& 0& 0& 0& 0\\
   \frac{x_2}{2} & -\frac{x_1}{2}& 0& 1& 0& 0& 0\\
   \frac{x_3}{2} &  0& -\frac{x_1}{2}& 0& 1& 0& 0\\
   -\frac{x_2x_3}{6}& \frac{x_1x_3}{12}+\frac{x_5}{2} &\frac{x_1x_2}{12}+\frac{x_4}{2} & -\frac{x_3}{2}& -\frac{x_2}{2}& 1& 0\\
 -\frac{x_3^2}{12} &  0&\frac{x_1x_3}{12} +\frac{x_5}{2}& 0& -\frac{x_3}{2}& 0& 1
   \end{matrix}\right]\,.
   \end{eqnarray*}
  \subsection*{$(247J)$ }
 
 The following Lie algebra is denoted as $({247J})$ by Gong in \cite{Gong_Thesis}, and as $\mathcal{G}_{7,2,26}$ by Magnin in \cite{magnin}.
 
 The non-trivial brackets are the following:
\begin{equation*}
\begin{aligned}
    &[X_1, X_i] = X_{i+2}\,,\,i=2, 3\,,\,[X_1,X_5]=X_6\,,\\
    [X_2&,X_5]=X_7\,,\,[X_3,X_4]=X_7\,,\,[X_3,X_5]=X_6\,.
\end{aligned}
\end{equation*}
 This is a nilpotent Lie algebra of rank 3 and step 3 that is stratifiable. The non-trivial Lie brackets can be pictured with the diagram:
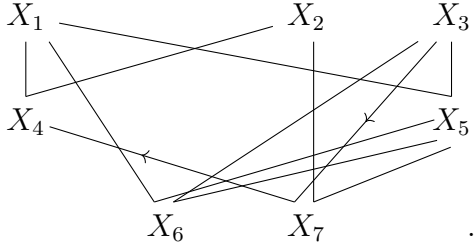
\begin{center}
 
	\begin{tikzcd}[end anchor=north]
		X_1 \ar[ddr, no head, end anchor={[xshift=-3.3ex]north east}]\ar[d, no head]\ar[drrr,no head]& & X_2\ar[dll, no head]\ar[dd, no head, end anchor={[xshift=-1.9ex]north east},start anchor={[xshift=-1.9ex]south east}] &X_3\ar[d, no head]\ar[ddl,no head,end anchor={[xshift=-3.3ex]north east},->-=.5]\ar[ddll,no head,end anchor={[xshift=-1.9ex]north east} ]\\
		X_4 \ar[drr, no head, end anchor={[xshift=-3.3ex]north east},start anchor={[xshift=-0.7ex]south east},start anchor={[yshift=1.5ex]south east},-<-=.4]& & & X_5\ar[dl, no head, end anchor={[xshift=-1.9ex]north east},start anchor={[xshift=-2.5ex]south east}]\ar[dll, no head, end anchor={[xshift=-3.3ex]north east},start anchor={[xshift=-3.7ex]south east} ,start anchor={[yshift=2ex]south east}]\ar[dll,no head,end anchor={[xshift=-1.9ex]north east},start anchor={[xshift=-3.5ex]south east} ,start anchor={[yshift=.5ex]south east}]\\
		& X_6&X_7 & \quad\;.
	\end{tikzcd}

\end{center}

 
The composition law \eqref{group law in G} of $(247J)$ is given by:

\begin{itemize}
    \item $z_1=x_1+y_1$;
    \item $z_2=x_2+y_2$;
    \item $z_3=x_3+y_3$;
    \item $z_4=x_4+y_4+\frac{1}{2}(x_1y_2-x_2y_1)$;
    \item $z_5=x_5+y_5+\frac{1}{2}(x_1y_3-x_3y_1)$;
    \item $z_6=x_6+y_6+\frac{1}{2}(x_1y_5-x_5y_1+x_3y_5-x_5y_3)+\frac{1}{12}(x_1-y_1)(x_1y_3-x_3y_1)\salto+\frac{1}{12}(x_3-y_3)(x_1y_3-x_3y_1)$;
    \item $z_7=x_7+y_7+\frac{1}{2}(x_3y_4-x_4y_3+x_2y_5-x_5y_2)+\frac{1}{12}(x_3-y_3)(x_1y_2-x_2y_1)\salto+\frac{1}{12}(x_2-y_2)(x_1y_3-x_3y_1)$.
\end{itemize}



Since
\begin{eqnarray*}
  \mathrm{d}(L_\mathbf{x})_\mathbf{0}=\left[\begin{matrix}    
    1 & 0  & 0& 0& 0& 0& 0\\
   0 & 1& 0& 0& 0& 0& 0\\
   0 & 0& 1& 0& 0& 0& 0\\
   -\frac{x_2}{2} & \frac{x_1}{2}& 0& 1& 0& 0& 0\\
   -\frac{x_3}{2} &  0& \frac{x_1}{2}& 0& 1& 0& 0\\
-\frac{x_1x_3+x_3^2}{12}-\frac{x_5}{2} &  0&\frac{x_1^2+x_1x_3}{12} -\frac{x_5}{2}& 0& \frac{x_1+x_3}{2}& 1& 0\\
 -\frac{x_2x_3}{6}& \frac{x_1x_3}{12}-\frac{x_5}{2} &\frac{x_1x_2}{12}-\frac{x_4}{2} & \frac{x_3}{2}& \frac{x_2}{2}& 0& 1
   \end{matrix}\right]\,,
   \end{eqnarray*}
the induced left-invariant vector fields \eqref{leftinvariant vf} are:
\begin{itemize}
\item $X_1=\partial_{x_1} -\frac{x_2}{2}\partial_{x_4}-\frac{x_3}{2}\partial_{x_5}-\big(\frac{x_1x_3+x_3^2}{12}+\frac{x_5}{2}\big)\partial_{x_6}-\frac{x_2x_3}{6}\partial_{x_7}\,;$
\item $X_2=\partial_{x_2}+\frac{x_1}{2}\partial_{x_4}+\big(\frac{x_1x_3}{12}-\frac{x_5}{2}\big)\partial_{x_7}\,;$
\item $X_3=\partial_{x_3}+\frac{x_1}{2}\partial_{x_5}+\big(\frac{x_1^2+x_1x_3}{12}-\frac{x_5}{2}\big)\partial_{x_6}+\big(\frac{x_1x_2}{12}-\frac{x_4}{2}\big)\partial_{x_7}\,;$
\item $X_4=\partial_{x_4}+\frac{x_3}{2}\partial_{x_7}\,;$
\item $X_5=\partial_{x_5}+\frac{x_1+x_3}{2}\partial_{x_6}+\frac{x_2}{2}\partial_{x_7}$;
\item $X_6=\partial_{x_6}$;
\item $X_7=\partial_{x_7}$,
\end{itemize}
and the respective left-invariant 1-forms \eqref{leftinvariant form} are: 
\begin{itemize}
\item $\theta_1=dx_1$;
\item $\theta_2=dx_2$;
\item $\theta_3=dx_3$;
\item $\theta_4=dx_4-\frac{x_1}{2}dx_2+\frac{x_2}{2}dx_1$;
\item $\theta_5=dx_5-\frac{x_1}{2}dx_3+\frac{x_3}{2}dx_1$;
\item $\theta_6=dx_6-\frac{x_1+x_3}{2}dx_5+\big(\frac{x_1^2+x_1x_3}{6}+\frac{x_5}{2}\big)dx_3+\big(\frac{x_5}{2}-\frac{x_1x_3+x_3^2}{6}\big)dx_1$;
\item $\theta_7=dx_7-\frac{x_2}{2}dx_5-\frac{x_3}{2}dx_4+\big(\frac{x_1x_2}{6}+\frac{x_4}{2}\big)dx_3+\big(\frac{x_1x_3}{6}+\frac{x_5}{2}\big)dx_2-\frac{x_2x_3}{3}dx_1$.
\end{itemize}

Finally, we have
\begin{eqnarray*}
  \mathrm{d}(R_\mathbf{x})_\mathbf{0}=\left[\begin{matrix}    
    1 & 0  & 0& 0& 0& 0& 0\\
   0 & 1& 0& 0& 0& 0& 0\\
   0 & 0& 1& 0& 0& 0& 0\\
   \frac{x_2}{2} & -\frac{x_1}{2}& 0& 1& 0& 0& 0\\
   \frac{x_3}{2} &  0& -\frac{x_1}{2}& 0& 1& 0& 0\\
\frac{x_5}{2}-\frac{x_1x_3+x_3^2}{12} &  0&\frac{x_1^2+x_1x_3}{12}+ \frac{x_5}{2}& 0& -\frac{x_1+x_3}{2}& 1& 0\\
 -\frac{x_2x_3}{6}& \frac{x_1x_3}{12}+\frac{x_5}{2} &\frac{x_1x_2}{12}+\frac{x_4}{2} & -\frac{x_3}{2}& -\frac{x_2}{2}& 0& 1
   \end{matrix}\right]\,.
   \end{eqnarray*}
  \subsection*{$(247K)$ }
 
 The following Lie algebra is denoted as $({247K})$ by Gong in \cite{Gong_Thesis}, and as $\mathcal{G}_{7,2,35}$ by Magnin in \cite{magnin}.
 
 The non-trivial brackets are the following:
\begin{equation*}
   [X_1, X_i] = X_{i+2}\,,\,2\le i\le 4\,,\,[X_2,X_5]=X_7\,,\,[X_3,X_4]=X_7\,,\,[X_3,X_5]=X_6\,.
\end{equation*}
 This is a nilpotent Lie algebra of rank 3 and step 3 that is stratifiable. The Lie brackets can be pictured with the diagram:
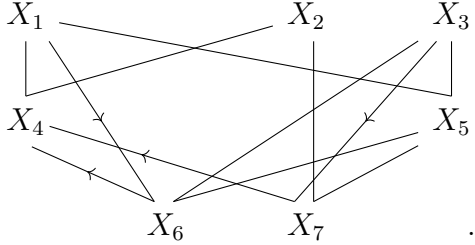
\begin{center}
 
	\begin{tikzcd}[end anchor=north]
		X_1 \ar[d, no head]\ar[drrr,no head]\ar[ddr, no head, end anchor={[xshift=-3.3ex]north east},->-=.5]& & X_2\ar[dll, no head]\ar[dd, no head, end anchor={[xshift=-1.9ex]north east},start anchor={[xshift=-1.9ex]south east}] &X_3\ar[d, no head]\ar[ddl,no head,end anchor={[xshift=-3.3ex]north east},->-=.5]\ar[ddll,no head,end anchor={[xshift=-1.9ex]north east} ]\\
		X_4 \ar[dr, no head, end anchor={[xshift=-3.3ex]north east},start anchor={[xshift=-2.ex]south east},-<-=.5]\ar[drr, no head, end anchor={[xshift=-3.3ex]north east},start anchor={[xshift=-0.7ex]south east},start anchor={[yshift=1.5ex]south east},-<-=.4]& & & X_5\ar[dl, no head, end anchor={[xshift=-1.9ex]north east}]\ar[dll,no head,end anchor={[xshift=-1.9ex]north east} ]\\
		& X_6&X_7 & \quad\;.
	\end{tikzcd}

\end{center}

 
The composition law \eqref{group law in G} of $(247K)$ is given by:

\begin{itemize}
    \item $z_1=x_1+y_1$;
    \item $z_2=x_2+y_2$;
    \item $z_3=x_3+y_3$;
    \item $z_4=x_4+y_4+\frac{1}{2}(x_1y_2-x_2y_1)$;
    \item $z_5=x_5+y_5+\frac{1}{2}(x_1y_3-x_3y_1)$;
    \item $z_6=x_6+y_6+\frac{1}{2}(x_1y_4-x_4y_1+x_3y_5-x_5y_3)+\frac{1}{12}(x_1-y_1)(x_1y_2-x_2y_1)\salto+\frac{1}{12}(x_3-y_3)(x_1y_3-x_3y_1)$;
    \item $z_7=x_7+y_7+\frac{1}{2}(x_3y_4-x_4y_3+x_2y_5-x_5y_2)+\frac{1}{12}(x_3-y_3)(x_1y_2-x_2y_1)\salto+\frac{1}{12}(x_2-y_2)(x_1y_3-x_3y_1)$.
\end{itemize}



Since
\begin{eqnarray*}
  \mathrm{d}(L_\mathbf{x})_\mathbf{0}=\left[\begin{matrix}    
    1 & 0  & 0& 0& 0& 0& 0\\
   0 & 1& 0& 0& 0& 0& 0\\
   0 & 0& 1& 0& 0& 0& 0\\
   -\frac{x_2}{2} & \frac{x_1}{2}& 0& 1& 0& 0& 0\\
   -\frac{x_3}{2} &  0& \frac{x_1}{2}& 0& 1& 0& 0\\
-\frac{x_1x_2+x_3^2}{12}-\frac{x_4}{2} &  \frac{x_1^2}{12}&\frac{x_1x_3}{12} -\frac{x_5}{2}& \frac{x_1}{2}& \frac{x_3}{2}& 1& 0\\
 -\frac{x_2x_3}{6}& \frac{x_1x_3}{12}-\frac{x_5}{2} &\frac{x_1x_2}{12}-\frac{x_4}{2} & \frac{x_3}{2}& \frac{x_2}{2}& 0& 1
   \end{matrix}\right]\,,
   \end{eqnarray*}
the induced left-invariant vector fields \eqref{leftinvariant vf} are:
\begin{itemize}
\item $X_1=\partial_{x_1} -\frac{x_2}{2}\partial_{x_4}-\frac{x_3}{2}\partial_{x_5}-\big(\frac{x_1x_2+x_3^2}{12}+\frac{x_4}{2}\big)\partial_{x_6}-\frac{x_2x_3}{6}\partial_{x_7}\,;$
\item $X_2=\partial_{x_2}+\frac{x_1}{2}\partial_{x_4}+\frac{x_1^2}{12}\partial_{x_6}+\big(\frac{x_1x_3}{12}-\frac{x_5}{2}\big)\partial_{x_7}\,;$
\item $X_3=\partial_{x_3}+\frac{x_1}{2}\partial_{x_5}+\big(\frac{x_1x_3}{12}-\frac{x_5}{2}\big)\partial_{x_6}+\big(\frac{x_1x_2}{12}-\frac{x_4}{2}\big)\partial_{x_7}\,;$
\item $X_4=\partial_{x_4}+\frac{x_1}{2}\partial_{x_6}+\frac{x_3}{2}\partial_{x_7}\,;$
\item $X_5=\partial_{x_5}+\frac{x_3}{2}\partial_{x_6}+\frac{x_2}{2}\partial_{x_7}$;
\item $X_6=\partial_{x_6}$;
\item $X_7=\partial_{x_7}$,
\end{itemize}
and the respective left-invariant 1-forms \eqref{leftinvariant form} are: 
\begin{itemize}
\item $\theta_1=dx_1$;
\item $\theta_2=dx_2$;
\item $\theta_3=dx_3$;
\item $\theta_4=dx_4-\frac{x_1}{2}dx_2+\frac{x_2}{2}dx_1$;
\item $\theta_5=dx_5-\frac{x_1}{2}dx_3+\frac{x_3}{2}dx_1$;
\item $\theta_6=dx_6-\frac{x_3}{2}dx_5-\frac{x_1}{2}dx_4+\big(\frac{x_1x_3}{6}+\frac{x_5}{2}\big)dx_3+\frac{x_1^2}{6}dx_2+\big(\frac{x_4}{2}-\frac{x_1x_2+x_3^2}{6}\big)dx_1$;
\item $\theta_7=dx_7-\frac{x_2}{2}dx_5-\frac{x_3}{2}dx_4+\big(\frac{x_1x_2}{6}+\frac{x_4}{2}\big)dx_3+\big(\frac{x_1x_3}{6}+\frac{x_5}{2}\big)dx_2-\frac{x_2x_3}{3}dx_1$.
\end{itemize}

Finally, we have
\begin{eqnarray*}
  \mathrm{d}(R_\mathbf{x})_\mathbf{0}=\left[\begin{matrix}    
    1 & 0  & 0& 0& 0& 0& 0\\
   0 & 1& 0& 0& 0& 0& 0\\
   0 & 0& 1& 0& 0& 0& 0\\
   \frac{x_2}{2} & -\frac{x_1}{2}& 0& 1& 0& 0& 0\\
   \frac{x_3}{2} &  0& -\frac{x_1}{2}& 0& 1& 0& 0\\
\frac{x_4}{2}-\frac{x_1x_2+x_3^2}{12} &  \frac{x_1^2}{12}&\frac{x_1x_3}{12} +\frac{x_5}{2}& -\frac{x_1}{2}& -\frac{x_3}{2}& 1& 0\\
 -\frac{x_2x_3}{6}& \frac{x_1x_3}{12}+\frac{x_5}{2} &\frac{x_1x_2}{12}+\frac{x_4}{2} & -\frac{x_3}{2}& -\frac{x_2}{2}& 0& 1
   \end{matrix}\right]\,.
   \end{eqnarray*}
  \subsection*{$(247N)$ }
 
 The following Lie algebra is denoted as $({247N})$ by Gong in \cite{Gong_Thesis}, and as $\mathcal{G}_{7,2,44}$ by Magnin in \cite{magnin}.
 
The non-trivial brackets are the following:
\begin{equation*}
   [X_1, X_i] = X_{i+2}\,,\,i=2,3\,,\,[X_1,X_5]=X_6\,,\,[X_2,X_3]=X_7\,,\,[X_2,X_4]=X_6\,.
\end{equation*}
 This is a nilpotent Lie algebra of rank 3 and step 3 that is stratifiable. The Lie brackets can be pictured with the diagram:
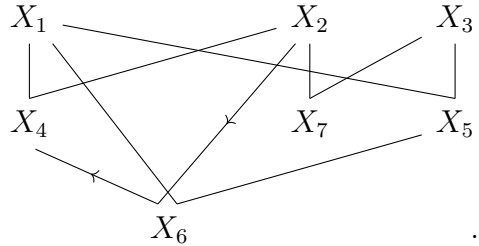
\begin{center}
 
	\begin{tikzcd}[end anchor=north]
		X_1 \ar[d, no head]\ar[drrr,no head]\ar[ddr, no head, end anchor={[xshift=-1.9ex]north east}]& & X_2\ar[ddl, no head, end anchor={[xshift=-3.3ex]north east},->-=.5]\ar[dll, no head]\ar[d, no head] &X_3\ar[d, no head]\ar[dl,no head]\\
		X_4 \ar[dr, no head, end anchor={[xshift=-3.3ex]north east},start anchor={[xshift=-2.ex]south east},-<-=.5]& &X_7 & X_5\ar[dll,no head,end anchor={[xshift=-1.9ex]north east} ]\\
		& X_6& & \quad\;.
	\end{tikzcd}

\end{center}

 
The composition law \eqref{group law in G} of $(247N)$ is given by:

\begin{itemize}
    \item $z_1=x_1+y_1$;
    \item $z_2=x_2+y_2$;
    \item $z_3=x_3+y_3$;
    \item $z_4=x_4+y_4+\frac{1}{2}(x_1y_2-x_2y_1)$;
    \item $z_5=x_5+y_5+\frac{1}{2}(x_1y_3-x_3y_1)$;
    \item $z_6=x_6+y_6+\frac{1}{2}(x_2y_4-x_4y_2+x_1y_5-x_5y_1)+\frac{1}{12}(x_2-y_2)(x_1y_2-x_2y_1)\salto+\frac{1}{12}(x_1-y_1)(x_1y_3-x_3y_1)]$;
    \item $z_7=x_7+y_7+\frac{1}{2}(x_2y_3-x_3y_2)$.
\end{itemize}



Since
\begin{eqnarray*}
  \mathrm{d}(L_\mathbf{x})_\mathbf{0}=\left[\begin{matrix}    
      
   1 & 0  & 0& 0& 0& 0& 0\\
   0 & 1& 0& 0& 0& 0& 0\\
   0 & 0& 1& 0& 0& 0& 0\\
   -\frac{x_2}{2} & \frac{x_1}{2}& 0& 1& 0& 0& 0\\
   -\frac{x_3}{2} &  0& \frac{x_1}{2}& 0& 1& 0& 0\\
 -\frac{x_2^2+x_1x_3}{12}-\frac{x_5}{2} &  \frac{x_1x_2}{12}-\frac{x_4}{2}&\frac{x_1^2}{12}& \frac{x_2}{2}& \frac{x_1}{2}& 1& 0\\
 0& -\frac{x_3}{2}&\frac{x_2}{2} & 0& 0& 0& 1
   \end{matrix}\right]\,,
   \end{eqnarray*}
the induced left-invariant vector fields \eqref{leftinvariant vf} are:
\begin{itemize}
\item $X_1=\partial_{x_1} -\frac{x_2}{2}\partial_{x_4}-\frac{x_3}{2}\partial_{x_5}-\big(\frac{x_1x_3+x_2^2}{12}+\frac{x_5}{2}\big)\partial_{x_6}\,;$
\item $X_2=\partial_{x_2}+\frac{x_1}{2}\partial_{x_4}+\big(\frac{x_1x_2}{12}-\frac{x_4}{2}\big)\partial_{x_6}-\frac{x_3}{2}\partial_{x_7}\,;$
\item $X_3=\partial_{x_3}+\frac{x_1}{2}\partial_{x_5}+\frac{x_1^2}{12}\partial_{x_6}+\frac{x_2}{2}\partial_{x_7}\,;$
\item $X_4=\partial_{x_4}+\frac{x_2}{2}\partial_{x_6}\,;$
\item $X_5=\partial_{x_5}+\frac{x_1}{2}\partial_{x_6}$;
\item $X_6=\partial_{x_6}$;
\item $X_7=\partial_{x_7}$,
\end{itemize}
and the respective left-invariant 1-forms \eqref{leftinvariant form} are: 
\begin{itemize}
\item $\theta_1=dx_1$;
\item $\theta_2=dx_2$;
\item $\theta_3=dx_3$;
\item $\theta_4=dx_4-\frac{x_1}{2}dx_2+\frac{x_2}{2}dx_1$;
\item $\theta_5=dx_5-\frac{x_1}{2}dx_3+\frac{x_3}{2}dx_1$;
\item $\theta_6=dx_6-\frac{x_1}{2}dx_5-\frac{x_2}{2}dx_4+\frac{x_1^2}{6}dx_3+\big(\frac{x_4}{2}+\frac{x_1x_2}{6}\big)dx_2+\big(\frac{x_5}{2}-\frac{x_1x_3+x_2^2}{6}\big)dx_1$;
\item $\theta_7=dx_7-\frac{x_2}{2}dx_3+\frac{x_3}{2}dx_2$.
\end{itemize}

Finally, we have
\begin{eqnarray*}
  \mathrm{d}(R_\mathbf{x})_\mathbf{0}=\left[\begin{matrix}    
      
   1 & 0  & 0& 0& 0& 0& 0\\
   0 & 1& 0& 0& 0& 0& 0\\
   0 & 0& 1& 0& 0& 0& 0\\
   \frac{x_2}{2} & -\frac{x_1}{2}& 0& 1& 0& 0& 0\\
   \frac{x_3}{2} &  0& -\frac{x_1}{2}& 0& 1& 0& 0\\
 \frac{x_5}{2}-\frac{x_2^2+x_1x_3}{12} &  \frac{x_1x_2}{12}+\frac{x_4}{2}&\frac{x_1^2}{12}& -\frac{x_2}{2}& -\frac{x_1}{2}& 1& 0\\
 0& \frac{x_3}{2}&-\frac{x_2}{2} & 0& 0& 0& 1
   \end{matrix}\right]\,.
   \end{eqnarray*}
  \subsection*{$(247P)$ }
 
 The following Lie algebra is denoted as $({247P})$ by Gong in \cite{Gong_Thesis}, and as $\mathcal{G}_{7,3,1(i_{\lambda})}$ with $\lambda=0,1$ by Magnin in \cite{magnin}.
 
 The non-trivial brackets are the following:
\begin{equation*}
   [X_1, X_i] = X_{i+2}\,,\,i=2,3\,,\,[X_2,X_5]=X_7\,,\,[X_3,X_4]=X_7\,,\,[X_2,X_3]=X_6\,.
\end{equation*}
 This is a nilpotent Lie algebra of rank 3 and step 3 that is stratifiable. The Lie brackets can be pictured with the diagram:
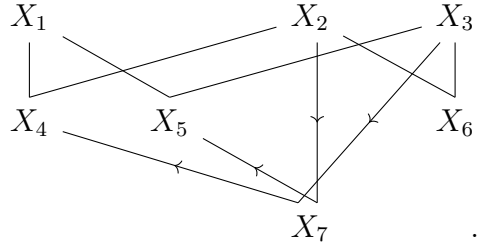
\begin{center}
 
	\begin{tikzcd}[end anchor=north]
		X_1\ar[d,no head]\ar[dr, no head]& & X_2\ar[dll,no head]\ar[dd,no head, end anchor={[xshift=-1.9ex]north east},start anchor={[xshift=-1.9ex]south east},->-=.5] \ar[dr, no head]& X_3\ar[d, no head]\ar[dll, no head]\ar[ddl,no head, end anchor={[xshift=-3.3ex]north east},->-=.5]\\
		X_4\ar[drr,no head, end anchor={[xshift=-3.3ex]north east},-<-=.5]  & X_5\ar[dr,no head, end anchor={[xshift=-1.9ex]north east},-<-=.5] & & X_6\\
		&  & X_7&\quad\;.
	\end{tikzcd}

\end{center}

 
The composition law \eqref{group law in G} of $(247P)$ is given by:

\begin{itemize}
    \item $z_1=x_1+y_1$;
    \item $z_2=x_2+y_2$;
    \item $z_3=x_3+y_3$;
    \item $z_4=x_4+y_4+\frac{1}{2}(x_1y_2-x_2y_1)$;
    \item $z_5=x_5+y_5+\frac{1}{2}(x_1y_3-x_3y_1)$;
    \item $z_6=x_6+y_6+\frac{1}{2}(x_2y_3-x_3y_2)$;
    \item $z_7=x_7+y_7+\frac{1}{2}(x_3y_4-x_4y_3+x_2y_5-x_5y_2)+\frac{1}{12}(x_3-y_3)(x_1y_2-x_2y_1)\salto+\frac{1}{12}(x_2-y_2)(x_1y_3-x_3y_1)$.
\end{itemize}



Since
\begin{eqnarray*}
  \mathrm{d}(L_\mathbf{x})_\mathbf{0}=\left[\begin{matrix}    
    1 & 0  & 0& 0& 0& 0& 0\\
   0 & 1& 0& 0& 0& 0& 0\\
   0 & 0& 1& 0& 0& 0& 0\\
  - \frac{x_2}{2} & \frac{x_1}{2}& 0& 1& 0& 0& 0\\
   -\frac{x_3}{2} &  0& \frac{x_1}{2}& 0& 1& 0& 0\\
0 &  -\frac{x_3}{2}& \frac{x_2}{2}& 0& 0& 1& 0\\
 -\frac{x_2x_3}{6}& \frac{x_1x_3}{12}-\frac{x_5}{2} &\frac{x_1x_2}{12}-\frac{x_4}{2} & \frac{x_3}{2}& \frac{x_2}{2}& 0& 1
   \end{matrix}\right]\,,
   \end{eqnarray*}
the induced left-invariant vector fields \eqref{leftinvariant vf} are:
\begin{itemize}
\item $X_1=\partial_{x_1} -\frac{x_2}{2}\partial_{x_4}-\frac{x_3}{2}\partial_{x_5}-\frac{x_2x_3}{6}\partial_{x_7}\,;$
\item $X_2=\partial_{x_2}+\frac{x_1}{2}\partial_{x_4}-\frac{x_3}{2}\partial_{x_6}+\big(\frac{x_1x_3}{12}-\frac{x_5}{2}\big)\partial_{x_7}\,;$
\item $X_3=\partial_{x_3}+\frac{x_1}{2}\partial_{x_5}+\frac{x_2}{2}\partial_{x_6}+\big(\frac{x_1x_2}{12}-\frac{x_4}{2}\big)\partial_{x_7}\,;$
\item $X_4=\partial_{x_4}+\frac{x_3}{2}\partial_{x_7}\,;$
\item $X_5=\partial_{x_5}+\frac{x_2}{2}\partial_{x_7}$;
\item $X_6=\partial_{x_6}$;
\item $X_7=\partial_{x_7}$,
\end{itemize}
and the respective left-invariant 1-forms \eqref{leftinvariant form} are: 
\begin{itemize}
\item $\theta_1=dx_1$;
\item $\theta_2=dx_2$;
\item $\theta_3=dx_3$;
\item $\theta_4=dx_4-\frac{x_1}{2}dx_2+\frac{x_2}{2}dx_1$;
\item $\theta_5=dx_5-\frac{x_1}{2}dx_3+\frac{x_3}{2}dx_1$;
\item $\theta_6=dx_6-\frac{x_2}{2}dx_3+\frac{x_3}{2}dx_2$;
\item $\theta_7=dx_7-\frac{x_2}{2}dx_5-\frac{x_3}{2}dx_4+\big(\frac{x_1x_2}{6}+\frac{x_4}{2}\big)dx_3+\big(\frac{x_1x_3}{6}+\frac{x_5}{2}\big)dx_2-\frac{x_2x_3}{3}dx_1$.
\end{itemize}

Finally, we have 
\begin{eqnarray*}
  \mathrm{d}(R_\mathbf{x})_\mathbf{0}=\left[\begin{matrix}    
    1 & 0  & 0& 0& 0& 0& 0\\
   0 & 1& 0& 0& 0& 0& 0\\
   0 & 0& 1& 0& 0& 0& 0\\
   \frac{x_2}{2} & -\frac{x_1}{2}& 0& 1& 0& 0& 0\\
   \frac{x_3}{2} &  0& -\frac{x_1}{2}& 0& 1& 0& 0\\
0 &  \frac{x_3}{2}& -\frac{x_2}{2}& 0& 0& 1& 0\\
 -\frac{x_2x_3}{6}& \frac{x_1x_3}{12}+\frac{x_5}{2} &\frac{x_1x_2}{12}+\frac{x_4}{2} & -\frac{x_3}{2}& -\frac{x_2}{2}& 0& 1
   \end{matrix}\right]\,.
   \end{eqnarray*}
  \subsection*{$(247P_1)$ }
 
 The following Lie algebra is denoted as $({247P_1})$ by Gong in \cite{Gong_Thesis}, and as $\mathcal{G}_{7,3,1(i_{\lambda})}$ with $\lambda=0,1$ by Magnin in \cite{magnin}.
 
 The non-trivial brackets are the following:
\begin{equation*}
   [X_1, X_i] = X_{i+2}\;,\,i=2,3\,,\,[X_2,X_4]=X_7\,,\,[X_3,X_5]=X_7\,,\,[X_2,X_3]=X_6\,.
\end{equation*}
 This is a nilpotent Lie algebra of rank 3 and step 3 that is stratifiable. The Lie brackets can be pictured with the diagram:
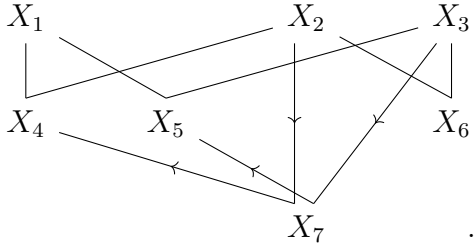
\begin{center}
 
	\begin{tikzcd}[end anchor=north]
		X_1\ar[d,no head]\ar[dr, no head]& & X_2\ar[dll,no head]\ar[dd,no head, end anchor={[xshift=-3.3ex]north east},start anchor={[xshift=-3.3ex]south east},->-=.5] \ar[dr, no head]& X_3\ar[d, no head]\ar[dll, no head]\ar[ddl,no head, end anchor={[xshift=-1.9ex]north east},->-=.5]\\
		X_4\ar[drr,no head, end anchor={[xshift=-3.3ex]north east},-<-=.5]  & X_5\ar[dr,no head, end anchor={[xshift=-1.9ex]north east},-<-=.5] & & X_6\\
		&  & X_7&\quad\;.
	\end{tikzcd}

\end{center}

 


The composition law \eqref{group law in G} of $(247P_1)$ is given by:
\begin{itemize}
    \item $z_1=x_1+y_1$;
    \item $z_2=x_2+y_2$;
    \item $z_3=x_3+y_3$;
    \item $z_4=x_4+y_4+\frac{1}{2}(x_1y_2-x_2y_1)$;
    \item $z_5=x_5+y_5+\frac{1}{2}(x_1y_3-x_3y_1)$;
    \item $z_6=x_6+y_6+\frac{1}{2}(x_2y_3-x_3y_2)$;
    \item $z_7=x_7+y_7+\frac{1}{2}(x_3y_5-x_5y_3+x_2y_4-x_4y_2)+\frac{1}{12}(x_2-y_2)(x_1y_2-x_2y_1)\salto+\frac{1}{12}(x_3-y_3)(x_1y_3-x_3y_1)$.
\end{itemize}



Since
\begin{eqnarray*}
  \mathrm{d}(L_\mathbf{x})_\mathbf{0}=\left[\begin{matrix}    
    1 & 0  & 0& 0& 0& 0& 0\\
   0 & 1& 0& 0& 0& 0& 0\\
   0 & 0& 1& 0& 0& 0& 0\\
   -\frac{x_2}{2} & \frac{x_1}{2}& 0& 1& 0& 0& 0\\
   -\frac{x_3}{2} &  0& \frac{x_1}{2}& 0& 1& 0& 0\\
0 &  -\frac{x_3}{2}& \frac{x_2}{2}& 0& 0& 1& 0\\
 -\frac{x_2^2+x_3^2}{12}& \frac{x_1x_2}{12}-\frac{x_4}{2} &\frac{x_1x_3}{12}-\frac{x_5}{2} & \frac{x_2}{2}& \frac{x_3}{2}& 0& 1
   \end{matrix}\right]\,,
   \end{eqnarray*}
the induced left-invariant vector fields \eqref{leftinvariant vf} are: 
\begin{itemize}
\item $X_1=\partial_{x_1}-\frac{x_2}{2}\partial_{x_4}-\frac{x_3}{2}\partial_{x_5}-\frac{x_2^2+x_3^2}{12}\partial_{x_7}\,;$
\item $X_2=\partial_{x_2}+\frac{x_1}{2}\partial_{x_4}-\frac{x_3}{2}\partial_{x_6}+\big(\frac{x_1x_2}{12}-\frac{x_4}{2}\big)\partial_{x_7}\,;$
\item $X_3=\partial_{x_3}+\frac{x_1}{2}\partial_{x_5}+\frac{x_2}{2}\partial_{x_6}+\big(\frac{x_1x_3}{12}-\frac{x_5}{2}\big)\partial_{x_7}\,;$
\item $X_4=\partial_{x_4}+\frac{x_2}{2}\partial_{x_7}\,;$
\item $X_5=\partial_{x_5}+\frac{x_3}{2}\partial_{x_7}$;
\item $X_6=\partial_{x_6}$;
\item $X_7=\partial_{x_7}$,
\end{itemize}
and the respective left-invariant 1-forms \eqref{leftinvariant form} are: 
\begin{itemize}
\item $\theta_1=dx_1$;
\item $\theta_2=dx_2$;
\item $\theta_3=dx_3$;
\item $\theta_4=dx_4-\frac{x_1}{2}dx_2+\frac{x_2}{2}dx_1$;
\item $\theta_5=dx_5-\frac{x_1}{2}dx_3+\frac{x_3}{2}dx_1$;
\item $\theta_6=dx_6-\frac{x_2}{2}dx_3+\frac{x_3}{2}dx_2$;
\item $\theta_7=dx_7-\frac{x_3}{2}dx_5-\frac{x_2}{2}dx_4+\big(\frac{x_1x_3}{6}+\frac{x_5}{2}\big)dx_3+\big(\frac{x_1x_2}{6}+\frac{x_4}{2}\big)dx_2-\frac{x_2^2+x_3^2}{6}dx_1$.
\end{itemize}

Finally, we have
\begin{eqnarray*}
  \mathrm{d}(R_\mathbf{x})_\mathbf{0}=\left[\begin{matrix}    
    1 & 0  & 0& 0& 0& 0& 0\\
   0 & 1& 0& 0& 0& 0& 0\\
   0 & 0& 1& 0& 0& 0& 0\\
   \frac{x_2}{2} & -\frac{x_1}{2}& 0& 1& 0& 0& 0\\
   \frac{x_3}{2} &  0& -\frac{x_1}{2}& 0& 1& 0& 0\\
0 &  \frac{x_3}{2}& -\frac{x_2}{2}& 0& 0& 1& 0\\
 -\frac{x_2^2+x_3^2}{12}& \frac{x_1x_2}{12}+\frac{x_4}{2} &\frac{x_1x_3}{12}+\frac{x_5}{2} & -\frac{x_2}{2}& -\frac{x_3}{2}& 0& 1
   \end{matrix}\right]\,.
   \end{eqnarray*}
  \subsection*{$(2457A)$ }
 
 The following Lie algebra is denoted as $({2457A})$ by Gong in \cite{Gong_Thesis}, and as $\mathcal{G}_{7,3,2}$ by Magnin in \cite{magnin}.

 The non-trivial brackets are the following:
\begin{equation*}
   [X_1, X_i] = X_{i+1}\,,\,i=2,3
   \,,\,[X_1, X_i] = X_{i+2}\,,\,i=4,5\,.
\end{equation*}
 This is a nilpotent Lie algebra of rank 3 and step 4 that is stratifiable. The Lie brackets can be pictured with the diagram:
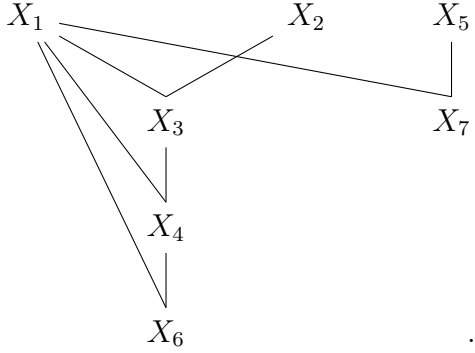
\begin{center}
 
	\begin{tikzcd}[end anchor=north]
		X_1\ar[dddr, no head]\ar[dr,no head]\ar[ddr, no head]\ar[drrr,no head]& & X_2\ar[dl,no head] & X_5\ar[d, no head]\\
		  &X_3\ar[d, no head]  & & X_7\\
		& X_ 4\ar[d, no head]& &\\
		& X_6 & &\quad\;.
	\end{tikzcd}

\end{center}

 


The composition law \eqref{group law in G} of $(2457A)$ is given by:
\begin{itemize}
    \item $z_1=x_1+y_1$;
    \item $z_2=x_2+y_2$;
    \item $z_3=x_3+y_3+\frac{1}{2}(x_1y_2-x_2y_1)$;
    \item $z_4=x_4+y_4+\frac{1}{2}(x_1y_3-x_3y_1)+\frac{1}{12}(x_1-y_1)(x_1y_2-x_2y_1)$;
    \item $z_5=x_5+y_5$;
    \item $z_6=x_6+y_6+\frac{1}{2}(x_1y_4-x_4y_1)+\frac{1}{12}(x_1-y_1)(x_1y_3-x_3y_1)-\frac{1}{24}x_1y_1(x_1y_2-x_2y_1)$;
    \item $z_7=x_7+y_7+\frac{1}{2}(x_1y_5-x_5y_1)$.
\end{itemize}



Since
\begin{eqnarray*}
  \mathrm{d}(L_\mathbf{x})_\mathbf{0}=\left[\begin{matrix}    
     1 & 0  & 0& 0& 0& 0& 0\\
   0 & 1& 0& 0& 0& 0& 0\\
   -\frac{x_2}{2} & \frac{x_1}{2}& 1& 0& 0& 0& 0\\
   -\frac{x_1x_2}{12}-\frac{x_3}{2} &\frac{x_1^2}{12} & \frac{x_1}{2}& 1& 0& 0& 0\\
   0 &  0& 0& 0& 1& 0& 0\\
 -\frac{x_1x_3}{12}-\frac{x_4}{2} &  0& \frac{x_1^2}{12}& \frac{x_1}{2}& 0& 1& 0\\
 -\frac{x_5}{2}& 0 &0 & 0& \frac{x_1}{2}& 0& 1
   \end{matrix}\right]\,,
   \end{eqnarray*}
the induced left-invariant vector fields \eqref{leftinvariant vf} are: 
\begin{itemize}
\item $X_1=\partial_{x_1}-\frac{x_2}{2}\partial_{x_3}-\big(\frac{x_3}{2}+\frac{x_1x_2}{12}\big)\partial_{x_4}-\big(\frac{x_1x_3}{12}+\frac{x_4}{2}\big)\partial_{x_6}-\frac{x_5}{2}\partial_{x_7}\,;$
\item $X_2=\partial_{x_2}+\frac{x_1}{2}\partial_{x_3}+\frac{x_1^2}{12}\partial_{x_4}\,;$
\item $X_3=\partial_{x_3}+\frac{x_1}{2}\partial_{x_4}+\frac{x_1^2}{12}\partial_{x_6}\,;$
\item $X_4=\partial_{x_4}+\frac{x_1}{2}\partial_{x_6}\,;$
\item $X_5=\partial_{x_5}+\frac{x_1}{2}\partial_{x_7}$;
\item $X_6=\partial_{x_6}$;
\item $X_7=\partial_{x_7}$,
\end{itemize}
and the respective left-invariant 1-forms \eqref{leftinvariant form} are: 
\begin{itemize}
\item $\theta_1=dx_1$;
\item $\theta_2=dx_2$;
\item $\theta_3=dx_3-\frac{x_1}{2}dx_2+\frac{x_2}{2}dx_1$;
\item $\theta_4=dx_4-\frac{x_1}{2}dx_3+\frac{x_1^2}{6}dx_2+\big(\frac{x_3}{2}-\frac{x_1x_2}{6}\big)dx_1$;
\item $\theta_5=dx_5$;
\item $\theta_6=dx_6-\frac{x_1}{2}dx_4+\frac{x_1^2}{6}dx_3-\frac{x_1^3}{24}dx_2+\big(\frac{x_1^2x_2}{24}-\frac{x_1x_3}{6}+\frac{x_4}{2}\big)dx_1$;
\item $\theta_7=dx_7-\frac{x_1}{2}dx_5+\frac{x_5}{2}dx_1$.
\end{itemize}

Finally, we have
\begin{eqnarray*}
  \mathrm{d}(R_\mathbf{x})_\mathbf{0}=\left[\begin{matrix}    
     1 & 0  & 0& 0& 0& 0& 0\\
   0 & 1& 0& 0& 0& 0& 0\\
   \frac{x_2}{2} & -\frac{x_1}{2}& 1& 0& 0& 0& 0\\
   \frac{x_3}{2}-\frac{x_1x_2}{12} &\frac{x_1^2}{12} & -\frac{x_1}{2}& 1& 0& 0& 0\\
   0 &  0& 0& 0& 1& 0& 0\\
 \frac{x_4}{2}-\frac{x_1x_3}{12} &  0& \frac{x_1^2}{12}& -\frac{x_1}{2}& 0& 1& 0\\
 \frac{x_5}{2}& 0 &0 & 0& -\frac{x_1}{2}& 0& 1
   \end{matrix}\right]\,.
   \end{eqnarray*}
  \subsection*{$(2457B)$ }
 
 The following Lie algebra is denoted as $({2457B})$ by Gong in \cite{Gong_Thesis}, and as $\mathcal{G}_{7,3,3}$ by Magnin in \cite{magnin}. 
 
 The non-trivial brackets are the following:
\begin{equation*}
   [X_1, X_i] = X_{i+1}\,,\,i=2,3
   \,,\,[X_1, X_4] = X_7\,,\,[X_2, X_5] = X_6\,.
\end{equation*}
 This is a nilpotent Lie algebra of rank 3 and step 4 that is stratifiable. The Lie brackets can be pictured with the diagram:
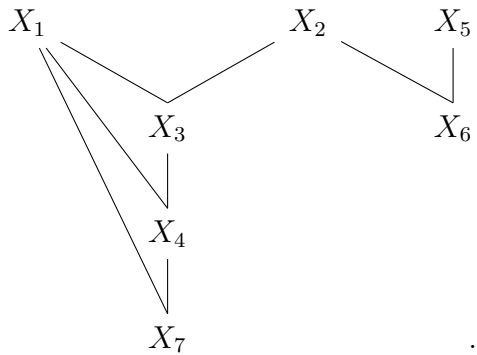
\begin{center}
 
	\begin{tikzcd}[end anchor=north]
		X_1\ar[dddr, no head]\ar[dr,no head]\ar[ddr, no head]& & X_2\ar[dl,no head]\ar[dr,no head] & X_5\ar[d, no head]\\
		  &X_3\ar[d, no head]  & & X_6\\
		& X_ 4\ar[d, no head]& &\\
		& X_7 & &\quad\;.
	\end{tikzcd}

\end{center}

 
 The composition law \eqref{group law in G} of $(2457B)$ is given by:

\begin{itemize}
    \item $z_1=x_1+y_1$;
    \item $z_2=x_2+y_2$;
    \item $z_3=x_3+y_3+\frac{1}{2}(x_1y_2-x_2y_1)$;
    \item $z_4=x_4+y_4+\frac{1}{2}(x_1y_3-x_3y_1)+\frac{1}{12}(x_1-y_1)(x_1y_2-x_2y_1)$;
    \item $z_5=x_5+y_5$;
    \item $z_6=x_6+y_6+\frac{1}{2}(x_2y_5-x_5y_2)$;
    \item $z_7=x_7+y_7+\frac{1}{2}(x_1y_4-x_4y_1)+\frac{1}{12}(x_1-y_1)(x_1y_3-x_3y_1)-\frac{1}{24}x_1y_1(x_1y_2-x_2y_1)$.
\end{itemize}



Since
\begin{eqnarray*}
  \mathrm{d}(L_\mathbf{x})_\mathbf{0}=\left[\begin{matrix}    
   1 & 0  & 0& 0& 0& 0& 0\\
   0 & 1& 0& 0& 0& 0& 0\\
   -\frac{x_2}{2} & \frac{x_1}{2}& 1& 0& 0& 0& 0\\
   -\frac{x_1x_2}{12}-\frac{x_3}{2} &\frac{x_1^2}{12} & \frac{x_1}{2}& 1& 0& 0& 0\\
   0 &  0& 0& 0& 1& 0& 0\\
 0& -\frac{x_5}{2} &0 & 0& \frac{x_2}{2}& 1& 0\\
 -\frac{x_1x_3}{12}-\frac{x_4}{2} &  0& \frac{x_1^2}{12}& \frac{x_1}{2}& 0& 0& 1\\
   \end{matrix}\right]\,,
   \end{eqnarray*}
the induced left-invariant vector fields \eqref{leftinvariant vf} are: 
\begin{itemize}
\item $X_1=\partial_{x_1}-\frac{x_2}{2}\partial_{x_3}-\big(\frac{x_3}{2}+\frac{x_1x_2}{12}\big)\partial_{x_4}-\big(\frac{x_1x_3}{12}+\frac{x_4}{2}\big)\partial_{x_7}\,;$
\item $X_2=\partial_{x_2}+\frac{x_1}{2}\partial_{x_3}+\frac{x_1^2}{12}\partial_{x_4}-\frac{x_5}{2}\partial_{x_6}\,;$
\item $X_3=\partial_{x_3}+\frac{x_1}{2}\partial_{x_4}+\frac{x_1^2}{12}\partial_{x_7}\,;$
\item $X_4=\partial_{x_4}+\frac{x_1}{2}\partial_{x_7}\,;$
\item $X_5=\partial_{x_5}+\frac{x_2}{2}\partial_{x_6}$;
\item $X_6=\partial_{x_6}$;
\item $X_7=\partial_{x_7}$,
\end{itemize}
and the respective left-invariant 1-forms \eqref{leftinvariant form} are: 
\begin{itemize}
\item $\theta_1=dx_1$;
\item $\theta_2=dx_2$;
\item $\theta_3=dx_3-\frac{x_1}{2}dx_2+\frac{x_2}{2}dx_1$;
\item $\theta_4=dx_4-\frac{x_1}{2}dx_3+\frac{x_1^2}{6}dx_2+\big(\frac{x_3}{2}-\frac{x_1x_2}{6}\big)dx_1$;
\item $\theta_5=dx_5$;
\item $\theta_6=dx_6-\frac{x_2}{2}dx_5+\frac{x_5}{2}dx_2$;
\item $\theta_7=dx_7-\frac{x_1}{2}dx_4+\frac{x_1^2}{6}dx_3-\frac{x_1^3}{24}dx_2+\big(\frac{x_4}{2}-\frac{x_1x_3}{6}+\frac{x_1^2x_2}{24}\big)dx_1$.
\end{itemize}

Finally, we have

\begin{eqnarray*}
  \mathrm{d}(R_\mathbf{x})_\mathbf{0}=\left[\begin{matrix}    
   1 & 0  & 0& 0& 0& 0& 0\\
   0 & 1& 0& 0& 0& 0& 0\\
   \frac{x_2}{2} & -\frac{x_1}{2}& 1& 0& 0& 0& 0\\
   \frac{x_3}{2}-\frac{x_1x_2}{12} &\frac{x_1^2}{12} & -\frac{x_1}{2}& 1& 0& 0& 0\\
   0 &  0& 0& 0& 1& 0& 0\\
 0& \frac{x_5}{2} &0 & 0& -\frac{x_2}{2}& 1& 0\\
 \frac{x_4}{2}-\frac{x_1x_3}{12} &  0& \frac{x_1^2}{12}& -\frac{x_1}{2}& 0& 0& 1\\
   \end{matrix}\right]\,.
   \end{eqnarray*}
  \subsection*{$(2457L)$ }
 
 The following Lie algebra is denoted as $({2457L})$ by Gong in \cite{Gong_Thesis}, and as $\mathcal{G}_{7,2,9}$ by Magnin in \cite{magnin}.
 
 The non-trivial brackets are the following:
\begin{equation*}
\begin{aligned}
  & [X_1, X_i] = X_{i+1}\,,\,i=2,3
   \,,\,[X_1, X_4] = X_6\,,\\ [X_1, X_5] = &X_7\,, \, [X_2, X_3] = X_5\,,\,[X_2, X_4] = X_7\,,\,[X_2, X_5] = X_6\,.
   \end{aligned}
\end{equation*}
 This is a nilpotent Lie algebra of rank 2 and step 4 that is stratifiable. The Lie brackets can be pictured with the diagram:
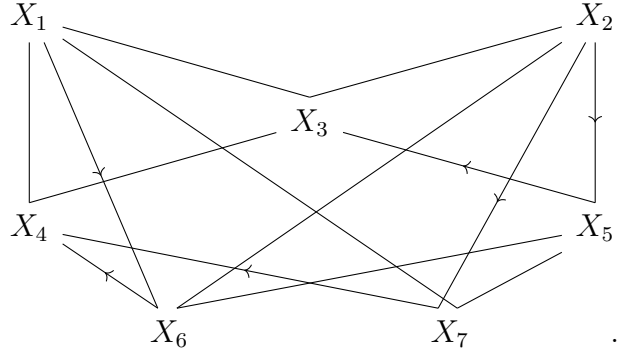
\begin{center}
 
	\begin{tikzcd}[end anchor=north]
		X_1\ar[dddrrr, no head, end anchor={[xshift=-1.9ex]north east}]\ar[dddr, no head,end anchor={[xshift=-3.3ex]north east},->-=.5]\ar[drr,no head]\ar[dd, no head]& & & & X_2\ar[dddl, no head,end anchor={[xshift=-3.3ex]north east},->-=.6]\ar[dddlll, no head ,end anchor={[xshift=-1.9ex]north east}]\ar[dll,no head]\ar[dd, no head,->-=.5] \\
		 & &X_3\ar[dll, no head] \ar[drr, no head,-<-=.5]& & \\
		X_4\ar[drrr, no head,end anchor={[xshift=-3.3ex]north east},-<-=.5]\ar[dr, no head,end anchor={[xshift=-3.3ex]north east},-<-=.5]& & & & X_5\ar[dlll, no head, ,end anchor={[xshift=-1.9ex]north east}]\ar[dl, no head, ,end anchor={[xshift=-1.9ex]north east}]\\
		& X_6& &X_7 & \quad\;.
	\end{tikzcd}

\end{center}

 
 The composition law \eqref{group law in G} of $(2457L)$ is given by:

\begin{itemize}
    \item $z_1=x_1+y_1$;
    \item $z_2=x_2+y_2$;
    \item $z_3=x_3+y_3+\frac{1}{2}(x_1y_2-x_2y_1)$;
    \item $z_4=x_4+y_4+\frac{1}{2}(x_1y_3-x_3y_1)+\frac{1}{12}(x_1-y_1)(x_1y_2-x_2y_1)$;
    \item $z_5=x_5+y_5+\frac{1}{2}(x_2y_3-x_3y_2)+\frac{1}{12}(x_2-y_2)(x_1y_2-x_2y_1)$;
    \item $z_6=x_6+y_6+\frac{1}{2}(x_1y_4-x_4y_1+x_2y_5-x_5y_2)+\frac{1}{12}(x_1-y_1)(x_1y_3-x_3y_1)\salto+\frac{1}{12}(x_2-y_2)(x_2y_3-x_3y_2)-\frac{1}{24}(x_1y_1+x_2y_2)(x_1y_2-x_2y_1)$;
    \item $z_7=x_7+y_7+\frac{1}{2}(x_1y_5-x_5y_1+x_2y_4-x_4y_2)+\frac{1}{12}(x_2-y_2)(x_1y_3-x_3y_1)\salto+\frac{1}{12}(x_1-y_1)(x_2y_3-x_3y_2)-\frac{1}{24}(x_1y_2+x_2y_1)(x_1y_2-x_2y_1)$.
\end{itemize}



Since
\begin{eqnarray*}
  \mathrm{d}(L_\mathbf{x})_\mathbf{0}=\left[\begin{matrix}    
      1 & 0  & 0& 0& 0& 0& 0\\
   0 & 1& 0& 0& 0& 0& 0\\
   -\frac{x_2}{2} & \frac{x_1}{2}& 1& 0& 0& 0& 0\\
   -\frac{x_1x_2}{12}-\frac{x_3}{2} &\frac{x_1^2}{12} & \frac{x_1}{2}& 1& 0& 0& 0\\
   -\frac{x_2^2}{12} &  \frac{x_1x_2}{12}-\frac{x_3}{2}& \frac{x_2}{2}& 0& 1& 0& 0\\
 -\frac{x_1x_3}{12}-\frac{x_4}{2} &  -\frac{x_2x_3}{12}-\frac{x_5}{2}& \frac{x_1^2+x_2^2}{12}& \frac{x_1}{2}& \frac{x_2}{2}& 1& 0\\
 -\frac{x_2x_3}{12}-\frac{x_5}{2}& -\frac{x_1x_3}{12}-\frac{x_4}{2} &\frac{x_1x_2}{6} & \frac{x_2}{2}& \frac{x_1}{2}& 0& 1
   \end{matrix}\right]\,,
   \end{eqnarray*}
the induced left-invariant vector fields \eqref{leftinvariant vf} are: 
\begin{itemize}
\item $X_1=\partial_{x_1} -\frac{x_2}{2}\partial_{x_3}-\big(\frac{x_3}{2}+\frac{x_1x_2}{12}\big)\partial_{x_4}-\frac{x_2^2}{12}\partial_{x_5}-\big(\frac{x_1x_3}{12}+\frac{x_4}{2}\big)\partial_{x_6}-\big(\frac{x_2x_3}{12}+\frac{x_5}{2}\big)\partial_{x_7}\,;$
\item $X_2=\partial_{x_2}+\frac{x_1}{2}\partial_{x_3}+\frac{x_1^2}{12}\partial_{x_4}+\big(\frac{x_1x_2}{12}-\frac{x_3}{2}\big)\partial_{x_5}-\big(\frac{x_2x_3}{12}+\frac{x_5}{2}\big)\partial_{x_6}-\big(\frac{x_1x_3}{12}+\frac{x_4}{2}\big)\partial_{x_7}\,;$
\item $X_3=\partial_{x_3}+\frac{x_1}{2}\partial_{x_4}+\frac{x_2}{2}\partial_{x_5}+\frac{x_1^2+x_2^2}{12}\partial_{x_6}+\frac{x_1x_2}{6}\partial_{x_7}\,;$
\item $X_4=\partial_{x_4}+\frac{x_1}{2}\partial_{x_6}+\frac{x_2}{2}\partial_{x_7}\,;$
\item $X_5=\partial_{x_5}+\frac{x_2}{2}\partial_{x_6}+\frac{x_1}{2}\partial_{x_7}$;
\item $X_6=\partial_{x_6}$;
\item $X_7=\partial_{x_7}$,
\end{itemize}
and the respective left-invariant 1-forms \eqref{leftinvariant form} are: 
\begin{itemize}
\item $\theta_1=dx_1$;
\item $\theta_2=dx_2$;
\item $\theta_3=dx_3-\frac{x_1}{2}dx_2+\frac{x_2}{2}dx_1$;
\item $\theta_4=dx_4-\frac{x_1}{2}dx_3+\frac{x_1^2}{6}dx_2+\big(\frac{x_3}{2}-\frac{x_1x_2}{6}\big)dx_1$;
\item $\theta_5=dx_5-\frac{x_2}{2}dx_3+\big(\frac{x_1x_2}{6}+\frac{x_3}{2}\big)dx_2-\frac{x_2^2}{6}dx_1$;
\item $\theta_6=dx_6-\frac{x_2}{2}dx_5-\frac{x_1}{2}dx_4+\frac{x_1^2+x_2^2}{6}dx_3+\big(\frac{x_5}{2}-\frac{x_2x_3}{6}-\frac{x_1^3+x_1x_2^2}{24}\big)dx_2+\big(\frac{x_4}{2}-\frac{x_1x_3}{6}\saltot+\frac{x_1^2x_2+x_2^3}{24}\big)dx_1$;
\item $\theta_7=dx_7-\frac{x_1}{2}dx_5-\frac{x_2}{2}dx_4+\frac{x_1x_2}{3}dx_3+\big(\frac{x_4}{2}-\frac{x_1x_3}{6}-\frac{x_1^2x_2}{12}\big)dx_2+\big(\frac{x_5}{2}-\frac{x_2x_3}{6}\saltot+\frac{x_1x_2^2}{12}\big) dx_1$.
\end{itemize}

Finally, we have
\begin{eqnarray*}
  \mathrm{d}(R_\mathbf{x})_\mathbf{0}=\left[\begin{matrix}    
      1 & 0  & 0& 0& 0& 0& 0\\
   0 & 1& 0& 0& 0& 0& 0\\
   \frac{x_2}{2} & -\frac{x_1}{2}& 1& 0& 0& 0& 0\\
   \frac{x_3}{2}-\frac{x_1x_2}{12} &\frac{x_1^2}{12} & -\frac{x_1}{2}& 1& 0& 0& 0\\
   -\frac{x_2^2}{12} &  \frac{x_1x_2}{12}+\frac{x_3}{2}& -\frac{x_2}{2}& 0& 1& 0& 0\\
 \frac{x_4}{2}-\frac{x_1x_3}{12} &  \frac{x_5}{2}-\frac{x_2x_3}{12}& \frac{x_1^2+x_2^2}{12}& -\frac{x_1}{2}& -\frac{x_2}{2}& 1& 0\\
 \frac{x_5}{2}-\frac{x_2x_3}{12}& \frac{x_4}{2}-\frac{x_1x_3}{12} &\frac{x_1x_2}{6} & -\frac{x_2}{2}& -\frac{x_1}{2}& 0& 1
   \end{matrix}\right]\,.
   \end{eqnarray*}

  \subsection*{$(2457L_1)$ }
 
 The following Lie algebra is denoted as $({2457L_1})$ by Gong in \cite{Gong_Thesis}, and as $\mathcal{G}_{7,2,9}$ by Magnin in \cite{magnin}.
 
 The non-trivial brackets are the following:
\begin{equation*}
\begin{aligned}
  & [X_1, X_i] = X_{i+1}\,,\,i=2,3
   \,,\,[X_1, X_4] = X_6\,,\\ [X_1, X_5] = &X_7\,,\,  [X_2, X_3] = X_5\,,\,[X_2, X_4] = X_7\,,\,[X_2, X_5] = -X_6\,.
   \end{aligned}
\end{equation*}
 This is a nilpotent Lie algebra of rank 2 and step 4 that is stratifiable. The Lie brackets can be pictured with the diagram:
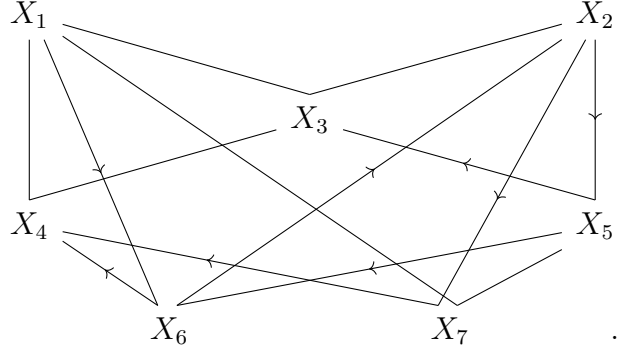
\begin{center}
 
	\begin{tikzcd}[end anchor=north]
		X_1\ar[dddrrr, no head, end anchor={[xshift=-1.9ex]north east}]\ar[dddr, no head,end anchor={[xshift=-3.3ex]north east},->-=.5]\ar[drr,no head]\ar[dd, no head]& & & & X_2\ar[dddl, no head,end anchor={[xshift=-3.3ex]north east},->-=.6]\ar[dddlll, no head ,end anchor={[xshift=-1.9ex]north east},-<-=.5]\ar[dll,no head]\ar[dd, no head,->-=.5] \\
		 & &X_3\ar[dll, no head] \ar[drr, no head,-<-=.5]& & \\
		X_4\ar[drrr, no head,end anchor={[xshift=-3.3ex]north east},-<-=.4]\ar[dr, no head,end anchor={[xshift=-3.3ex]north east},-<-=.5]& & & & X_5\ar[dlll, no head, ,end anchor={[xshift=-1.9ex]north east},->-=.5]\ar[dl, no head, ,end anchor={[xshift=-1.9ex]north east}]\\
		& X_6& &X_7 & \quad\;.
	\end{tikzcd}

\end{center}

 
 The composition law \eqref{group law in G} of $(2457L_1)$ is given by:

\begin{itemize}
    \item $z_1=x_1+y_1$;
    \item $z_2=x_2+y_2$;
    \item $z_3=x_3+y_3+\frac{1}{2}(x_1y_2-x_2y_1)$;
    \item $z_4=x_4+y_4+\frac{1}{2}(x_1y_3-x_3y_1)+\frac{1}{12}(x_1-y_1)(x_1y_2-x_2y_1)$;
    \item $z_5=x_5+y_5+\frac{1}{2}(x_2y_3-x_3y_2)+\frac{1}{12}(x_2-y_2)(x_1y_2-x_2y_1)$;
    \item $z_6=x_6+y_6+\frac{1}{2}(x_1y_4-x_4y_1-x_2y_5+x_5y_2)+\frac{1}{12}(x_1-y_1)(x_1y_3-x_3y_1)\salto\frac{1}{12}(y_2-x_2)(x_2y_3-x_3y_2)+\frac{1}{24}(x_2y_2-x_1y_1)(x_1y_2-x_2y_1)$;
    \item $z_7=x_7+y_7+\frac{1}{2}(x_1y_5-x_5y_1+x_2y_4-x_4y_2)+\frac{1}{12}(x_2-y_2)(x_1y_3-x_3y_1)\salto+\frac{1}{12}(x_1-y_1)(x_2y_3-x_3y_2)]-\frac{1}{24}(x_1y_2+x_2y_1)(x_1y_2-x_2y_1)$.
\end{itemize}



Since
\begin{eqnarray*}
  \mathrm{d}(L_\mathbf{x})_\mathbf{0}=\left[\begin{matrix}    
      1 & 0  & 0& 0& 0& 0& 0\\
   0 & 1& 0& 0& 0& 0& 0\\
   -\frac{x_2}{2} & \frac{x_1}{2}& 1& 0& 0& 0& 0\\
   -\frac{x_1x_2}{12}-\frac{x_3}{2} &\frac{x_1^2}{12} & \frac{x_1}{2}& 1& 0& 0& 0\\
   -\frac{x_2^2}{12} &  \frac{x_1x_2}{12}-\frac{x_3}{2}& \frac{x_2}{2}& 0& 1& 0& 0\\
 -\frac{x_1x_3}{12}-\frac{x_4}{2} &  \frac{x_2x_3}{12}+\frac{x_5}{2}& \frac{x_1^2-x_2^2}{12}& \frac{x_1}{2}& -\frac{x_2}{2}& 1& 0\\
 -\frac{x_2x_3}{12}-\frac{x_5}{2}& -\frac{x_1x_3}{12}-\frac{x_4}{2} &\frac{x_1x_2}{6} & \frac{x_2}{2}& \frac{x_1}{2}& 0& 1
   \end{matrix}\right]\,,
   \end{eqnarray*}
the induced left-invariant vector fields \eqref{leftinvariant vf} are: 
\begin{itemize}
\item $X_1=\partial_{x_1} -\frac{x_2}{2}\partial_{x_3}-\big(\frac{x_3}{2}+\frac{x_1x_2}{12}\big)\partial_{x_4}-\frac{x_2^2}{12}\partial_{x_5}-\big(\frac{x_1x_3}{12}+\frac{x_4}{2}\big)\partial_{x_6}-\big(\frac{x_5}{2}+\frac{x_2x_3}{12}\big)\partial_{x_7}\,;$
\item $X_2=\partial_{x_2}+\frac{x_1}{2}\partial_{x_3}+\frac{x_1^2}{12}\partial_{x_4}+\big(\frac{x_1x_2}{12}-\frac{x_3}{2}\big)\partial_{x_5}+\big(\frac{x_2x_3}{12}+\frac{x_5}{2}\big)\partial_{x_6}-\big(\frac{x_1x_3}{12}+\frac{x_4}{2}\big)\partial_{x_7}\,;$
\item $X_3=\partial_{x_3}+\frac{x_1}{2}\partial_{x_4}+\frac{x_2}{2}\partial_{x_5}+\frac{x_1^2-x_2^2}{12}\partial_{x_6}+\frac{x_1x_2}{6}\partial_{x_7}\,;$
\item $X_4=\partial_{x_4}+\frac{x_1}{2}\partial_{x_6}+\frac{x_2}{2}\partial_{x_7}\,;$
\item $X_5=\partial_{x_5}-\frac{x_2}{2}\partial_{x_6}+\frac{x_1}{2}\partial_{x_7}$;
\item $X_6=\partial_{x_6}$;
\item $X_7=\partial_{x_7}$,
\end{itemize}
and the respective left-invariant 1-forms \eqref{leftinvariant form} are: 
\begin{itemize}
\item $\theta_1=dx_1$;
\item $\theta_2=dx_2$;
\item $\theta_3=dx_3-\frac{x_1}{2}dx_2+\frac{x_2}{2}dx_1$;
\item $\theta_4=dx_4-\frac{x_1}{2}dx_3+\frac{x_1^2}{6}dx_2+\big(\frac{x_3}{2}-\frac{x_1x_2}{6}\big)dx_1$;
\item $\theta_5=dx_5-\frac{x_2}{2}dx_3+\big(\frac{x_3}{2}+\frac{x_1x_2}{6}\big)dx_2-\frac{x_2^2}{6}dx_1$;
\item $\theta_6=dx_6+\frac{x_2}{2}dx_5-\frac{x_1}{2}dx_4+\frac{x_1^2-x_2^2}{6}dx_3+\big(\frac{x_2x_3}{6}-\frac{x_5}{2}+\frac{x_1x_2^2-x_1^3}{24}\big)dx_2+\big(\frac{x_4}{2}-\frac{x_1x_3}{6}\saltot+\frac{x_1^2x_2-x_2^3}{24}\big)dx_1$;
\item $\theta_7=dx_7-\frac{x_1}{2}dx_5-\frac{x_2}{2}dx_4+\frac{x_1x_2}{3}dx_3+\big(\frac{x_4}{2}-\frac{x_1x_3}{6}-\frac{x_1^2x_2}{12}\big)dx_2+\big(\frac{x_5}{2}-\frac{x_2x_3}{6}+\frac{x_1x_2^2}{12}\big)dx_1$.
\end{itemize}

Finally, we have
\begin{eqnarray*}
  \mathrm{d}(R_\mathbf{x})_\mathbf{0}=\left[\begin{matrix}    
      1 & 0  & 0& 0& 0& 0& 0\\
   0 & 1& 0& 0& 0& 0& 0\\
   \frac{x_2}{2} & -\frac{x_1}{2}& 1& 0& 0& 0& 0\\
   \frac{x_3}{2}-\frac{x_1x_2}{12} &\frac{x_1^2}{12} & -\frac{x_1}{2}& 1& 0& 0& 0\\
   -\frac{x_2^2}{12} &  \frac{x_1x_2}{12}+\frac{x_3}{2}& -\frac{x_2}{2}& 0& 1& 0& 0\\
 \frac{x_4}{2}-\frac{x_1x_3}{12} &  \frac{x_2x_3}{12}-\frac{x_5}{2}& \frac{x_1^2-x_2^2}{12}& -\frac{x_1}{2}& \frac{x_2}{2}& 1& 0\\
 \frac{x_5}{2}-\frac{x_2x_3}{12}& \frac{x_4}{2}-\frac{x_1x_3}{12} &\frac{x_1x_2}{6} & -\frac{x_2}{2}& -\frac{x_1}{2}& 0& 1
   \end{matrix}\right]\,.
   \end{eqnarray*}
  \subsection*{$(2457M)$ }
 
 The following Lie algebra is denoted as $({2457M})$ by Gong in \cite{Gong_Thesis}, and as $\mathcal{G}_{7,2,8}$ by Magnin in \cite{magnin}.
 
 The non-trivial brackets are the following:
\begin{equation*}
\begin{aligned}
  & \;[X_1, X_i] = X_{i+1}\,,\,i=2,3
   \,,\,[X_1, X_4] = X_7\,,\\ [&X_1, X_5] = X_6\,, \, [X_2, X_3] = X_5\,,\,[X_2, X_4] = X_6\,.
   \end{aligned}
\end{equation*}
 This is a nilpotent Lie algebra of rank 2 and step 4 that is stratifiable. The Lie brackets can be pictured with the diagram:
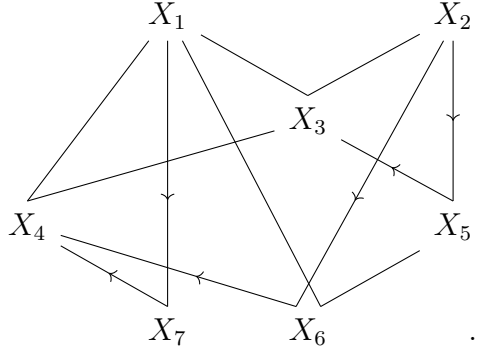
\begin{center}

\begin{tikzcd}[end anchor=north]
& X_1\ar[dr, no head]\ar[ddd, no head,->-=.6]\ar[ddl, no head] \ar[dddr, no head,  end anchor={[xshift=-1.5ex]north east}]&  & X_2\ar[dl, no head]\ar[dd,no head, ->-=.5]\ar[dddl, no head, ->-=.6, end anchor={[xshift=-3.3ex]north east}]\\
 & & X_3\ar[dll, no head]\ar[dr, no head, -<-=.5] &\\
 X_4\ar[dr, no head,-<-=.5] \ar[drr, no head, -<-=.6, end anchor={[xshift=-3.3ex]north east}]& & & X_5\ar[dl, no head,  end anchor={[xshift=-1.5ex]north east}]\\
 & X_7 & X_6 &\quad\;.
 
\end{tikzcd}

\end{center}

 
 The composition law \eqref{group law in G} of $(2457M)$ is given by:

\begin{itemize}
    \item $z_1=x_1+y_1$;
    \item $z_2=x_2+y_2$;
    \item $z_3=x_3+y_3+\frac{1}{2}(x_1y_2-x_2y_1)$;
    \item $z_4=x_4+y_4+\frac{1}{2}(x_1y_3-x_3y_1)+\frac{1}{12}(x_1-y_1)(x_1y_2-x_2y_1)$;
    \item $z_5=x_5+y_5+\frac{1}{2}(x_2y_3-x_3y_2)+\frac{1}{12}(x_2-y_2)(x_1y_2-x_2y_1)$;
    \item $z_6=x_6+y_6+\frac{1}{2}(x_2y_4-x_4y_2+x_1y_5-x_5y_1)+\frac{1}{12}(x_1-y_1)(x_2y_3-x_3y_2)\salto+\frac{1}{12}(x_2-y_2)(x_1y_3-x_3y_1)-\frac{1}{24}(x_1y_2+x_2y_1)(x_1y_2-x_2y_1)$;
    \item $z_7=x_7+y_7+\frac{1}{2}(x_1y_4-x_4y_1)+\frac{1}{12}(x_1-y_1)(x_1y_3-x_3y_1)-\frac{1}{24}x_1y_1(x_1y_2-x_2y_1)$.
\end{itemize}



Since
\begin{eqnarray*}
  \mathrm{d}(L_\mathbf{x})_\mathbf{0}=\left[\begin{matrix}    
      1 & 0  & 0& 0& 0& 0& 0\\
   0 & 1& 0& 0& 0& 0& 0\\
   -\frac{x_2}{2} & \frac{x_1}{2}& 1& 0& 0& 0& 0\\
   -\frac{x_1x_2}{12}-\frac{x_3}{2} &\frac{x_1^2}{12} & \frac{x_1}{2}& 1& 0& 0& 0\\
   -\frac{x_2^2}{12} &  \frac{x_1x_2}{12}-\frac{x_3}{2}& \frac{x_2}{2}& 0& 1& 0& 0\\
 -\frac{x_2x_3}{12}-\frac{x_5}{2} &  -\frac{x_1x_3}{12}-\frac{x_4}{2}& \frac{x_1x_2}{6}& \frac{x_2}{2}& \frac{x_1}{2}& 1& 0\\
 -\frac{x_1x_3}{12}-\frac{x_4}{2}& 0 &\frac{x_1^2}{12} & \frac{x_1}{2}& 0& 0& 1
   \end{matrix}\right]\,,
   \end{eqnarray*}
the induced left-invariant vector fields \eqref{leftinvariant vf} are: 
\begin{itemize}
\item $X_1=\partial_{x_1} -\frac{x_2}{2}\partial_{x_3}-\big(\frac{x_3}{2}+\frac{x_1x_2}{12}\big)\partial_{x_4}-\frac{x_2^2}{12}\partial_{x_5}-\big(\frac{x_2x_3}{12}+\frac{x_5}{2}\big)\partial_{x_6}-\big(\frac{x_1x_3}{12}+\frac{x_4}{2}\big)\partial_{x_7}\,;$
\item $X_2=\partial_{x_2}+\frac{x_1}{2}\partial_{x_3}+\frac{x_1^2}{12}\partial_{x_4}+\big(\frac{x_1x_2}{12}-\frac{x_3}{2}\big)\partial_{x_5}-\big(\frac{x_1x_3}{12}+\frac{x_4}{2}\big)\partial_{x_6}\,;$
\item $X_3=\partial_{x_3}+\frac{x_1}{2}\partial_{x_4}+\frac{x_2}{2}\partial_{x_5}+\frac{x_1x_2}{6}\partial_{x_6}+\frac{x_1^2}{12}\partial_{x_7}\,;$
\item $X_4=\partial_{x_4}+\frac{x_2}{2}\partial_{x_6}+\frac{x_1}{2}\partial_{x_7}\,;$
\item $X_5=\partial_{x_5}+\frac{x_1}{2}\partial_{x_6}$;
\item $X_6=\partial_{x_6}$;
\item $X_7=\partial_{x_7}$,
\end{itemize}
and the respective left-invariant 1-forms \eqref{leftinvariant form} are: 
\begin{itemize}
\item $\theta_1=dx_1$;
\item $\theta_2=dx_2$;
\item $\theta_3=dx_3-\frac{x_1}{2}dx_2+\frac{x_2}{2}dx_1$;
\item $\theta_4=dx_4-\frac{x_1}{2}dx_3+\frac{x_1^2}{6}dx_2+\big(\frac{x_3}{2}-\frac{x_1x_2}{6}\big)dx_1$;
\item $\theta_5=dx_5-\frac{x_2}{2}dx_3+\big(\frac{x_3}{2}+\frac{x_1x_2}{6}\big)dx_2-\frac{x_2^2}{6}dx_1$;
\item $\theta_6=dx_6-\frac{x_1}{2}dx_5-\frac{x_2}{2}dx_4+\frac{x_1x_2}{3}dx_3+\big(\frac{x_4}{2}-\frac{x_1x_3}{6}-\frac{x_1^2x_2}{12}\big)dx_2+\big(\frac{x_5}{2}-\frac{x_2x_3}{6}+\frac{x_1x_2^2}{12}\big)dx_1$;
\item $\theta_7=dx_7-\frac{x_1}{2}dx_4+\frac{x_1^2}{6}dx_3-\frac{x_1^3}{24}dx_2+\big(\frac{x_4}{2}-\frac{x_1x_3}{6}+\frac{x_1^2x_2}{24}\big)dx_1$.
\end{itemize}

Finally, we have
\begin{eqnarray*}
  \mathrm{d}(R_\mathbf{x})_\mathbf{0}=\left[\begin{matrix}    
      1 & 0  & 0& 0& 0& 0& 0\\
   0 & 1& 0& 0& 0& 0& 0\\
   \frac{x_2}{2} & -\frac{x_1}{2}& 1& 0& 0& 0& 0\\
   \frac{x_3}{2}-\frac{x_1x_2}{12} &\frac{x_1^2}{12} & -\frac{x_1}{2}& 1& 0& 0& 0\\
   -\frac{x_2^2}{12} &  \frac{x_1x_2}{12}+\frac{x_3}{2}& -\frac{x_2}{2}& 0& 1& 0& 0\\
 \frac{x_5}{2}-\frac{x_2x_3}{12} &  \frac{x_4}{2}-\frac{x_1x_3}{12}& \frac{x_1x_2}{6}& -\frac{x_2}{2}& -\frac{x_1}{2}& 1& 0\\
 \frac{x_4}{2}-\frac{x_1x_3}{12}& 0 &\frac{x_1^2}{12} & -\frac{x_1}{2}& 0& 0& 1
   \end{matrix}\right]\,.
   \end{eqnarray*}
 \subsection*{$(23457A)$}
 
 The following Lie algebra is denoted as $(23457A)$ by Gong in \cite{Gong_Thesis}, and as $\mathcal{G}_{7,2,7}$ by Magnin in \cite{magnin}.
 
 The non-trivial brackets are the following
 \begin{equation*}
     [X_1,X_i]=X_{i+1}\,,\,i=2,3,4\,,\,[X_1,X_5]=X_6\,,\,[X_2,X_3]=X_7\,.
 \end{equation*}
This is a nilpotent Lie algebra of rank 2 and step 5 that is stratifiable. The Lie brackets can be pictured with the diagram:
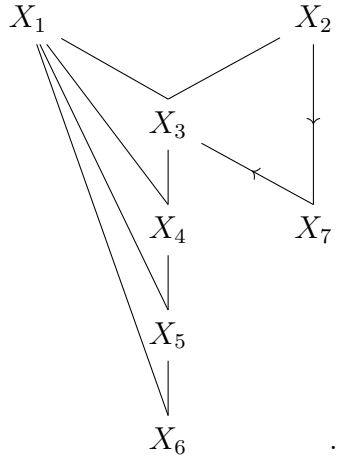
\begin{center}
 
	\begin{tikzcd}[end anchor=north]
		X_1\ar[dr,no head] \ar[ddr, no head]\ar[ dddr, no head]\ar[ddddr, no head]& & X_2\ar[dd, no head, ->-=.5]\ar[dl, no head]\\
		& X_3\ar[d, no head]\ar[dr, no head, -<-=.5] & \\
		& X_4\ar[d, no head] & X_7\\
		& X_5\ar[d, no head] & \\
		& X_6 &\quad\;. \\
	\end{tikzcd}
	\end{center}
	
The composition law \eqref{group law in G} of $(23457A)$ is given by:
\begin{itemize}
    \item $z_1=x_1+y_1$;
    \item $z_2=x_2+y_2$;
    \item $z_3=x_3+y_3+\frac{1}{2}(x_1y_2-x_2y_1)$;
    \item $z_4=x_4+y_4+\frac{1}{2}(x_1y_3-x_3y_1)+\frac{1}{12}(x_1-y_1)(x_1y_2-x_2y_1)$;
    \item $z_5=x_5+y_5+\frac{1}{2}(x_1y_4-x_4y_1)+\frac{1}{12}(x_1-y_1)(x_1y_3-x_3y_1)-\frac{1}{24}x_1y_1(x_1y_2-x_2y_1)$;
    \item $z_6=x_6+y_6+\frac{1}{2}(x_1y_5-x_5y_1)+\frac{1}{12}(x_1-y_1)(x_1y_4-x_4y_1)-\frac{1}{24}x_1y_1(x_1y_3-x_3y_1)\salto+\frac{1}{180}(x_1y_1^2-x_1^2y_1)(x_1y_2-x_2y_1)+\frac{1}{720}(y_1^3-x_1^3)(x_1y_2-x_2y_1)$;
    \item $z_7=x_7+y_7+\frac{1}{2}(x_2y_3-x_3y_2)+\frac{1}{12}(x_2-y_2)(x_1y_2-x_2y_1)$.
\end{itemize}

Since
\begin{eqnarray*}
  \mathrm{d}(L_\mathbf{x})_\mathbf{0}=\left[\begin{matrix}    
      1 & 0  & 0& 0& 0& 0& 0\\
   0 & 1& 0& 0& 0& 0& 0\\
   -\frac{x_2}{2} & \frac{x_1}{2}& 1& 0& 0& 0& 0\\
   -\frac{x_1x_2}{12}-\frac{x_3}{2} &\frac{x_1^2}{12} & \frac{x_1}{2}& 1& 0& 0& 0\\
   -\frac{x_1x_3}{12}-\frac{x_4}{2} &  0& \frac{x_1^2}{12}& \frac{x_1}{2}& 1& 0& 0\\
 \frac{x_1^3x_2}{720}-\frac{x_1x_4}{12}-\frac{x_5}{2} &  -\frac{x_1^4}{720}& 0& \frac{x_1^2}{12}& \frac{x_1}{2}& 1& 0\\
 -\frac{x_2^2}{12}& \frac{x_1x_2}{12}-\frac{x_3}{2} &\frac{x_2}{2} & 0& 0& 0& 1
   \end{matrix}\right]\,,
   \end{eqnarray*}
the induced left-invariant vector fields \eqref{leftinvariant vf} are: 
\begin{itemize}
\item $X_1=\partial_{x_1} -\frac{x_2}{2}\partial_{x_3}-\big(\frac{x_3}{2}+\frac{x_1x_2}{12}\big)\partial_{x_4}-\big(\frac{x_1x_3}{12}+\frac{x_4}{2}\big)\partial_{x_5}+\big(\frac{x_1^3x_2}{720}-\frac{x_1x_4}{12}-\frac{x_5}{2}\big)\partial_{x_6}-\frac{x_2^2}{12}\partial_{x_7}\,;$
\item $X_2=\partial_{x_2}+\frac{x_1}{2}\partial_{x_3}+\frac{x_1^2}{12}\partial_{x_4}-\frac{x_1^4}{720}\partial_{x_6}+\big(\frac{x_1x_2}{12}-\frac{x_3}{2}\big)\partial_{x_7}\,;$
\item $X_3=\partial_{x_3}+\frac{x_1}{2}\partial_{x_4}+\frac{x_1^2}{12}\partial_{x_5}+\frac{x_2}{2}\partial_{x_7}\,;$
\item $X_4=\partial_{x_4}+\frac{x_1}{2}\partial_{x_5}+\frac{x_1^2}{12}\partial_{x_6}\,;$
\item $X_5=\partial_{x_5}+\frac{x_1}{2}\partial_{x_6}$;
\item $X_6=\partial_{x_6}$;
\item $X_7=\partial_{x_7}$,
\end{itemize}
and the respective left-invariant 1-forms \eqref{leftinvariant form} are: 
\begin{itemize}
\item $\theta_1=dx_1$;
\item $\theta_2=dx_2$;
\item $\theta_3=dx_3-\frac{x_1}{2}dx_2+\frac{x_2}{2}dx_1$;
\item $\theta_4=dx_4-\frac{x_1}{2}dx_3+\frac{x_1^2}{6}dx_2+\big(\frac{x_3}{2}-\frac{x_1x_2}{6}\big)dx_1$;
\item $\theta_5=dx_5-\frac{x_1}{2}dx_4+\frac{x_1^2}{6}dx_3-\frac{x_1^3}{24}dx_2+\big(\frac{x_4}{2}-\frac{x_1x_3}{6}+\frac{x_1^2x_2}{24}\big)dx_1$;
\item $\theta_6=dx_6-\frac{x_1}{2}dx_5+\frac{x_1^2}{6}dx_4-\frac{x_1^3}{24}dx_3+\frac{x_1^4}{120}dx_2+\big(\frac{x_5}{2}-\frac{x_1x_4}{6}+\frac{x_1^2x_3}{24}-\frac{x_1^3x_2}{120}\big)dx_1$;
\item $\theta_7=dx_7-\frac{x_2}{2}dx_3+\big(\frac{x_3}{2}+\frac{x_1x_2}{6}\big)dx_2-\frac{x_2^2}{6}dx_1$.
\end{itemize}

Finally, we have
	\begin{eqnarray*}
  \mathrm{d}(R_\mathbf{x})_\mathbf{0}=\left[\begin{matrix}    
      1 & 0  & 0& 0& 0& 0& 0\\
   0 & 1& 0& 0& 0& 0& 0\\
   \frac{x_2}{2} & -\frac{x_1}{2}& 1& 0& 0& 0& 0\\
   \frac{x_3}{2}-\frac{x_1x_2}{12} &\frac{x_1^2}{12} & -\frac{x_1}{2}& 1& 0& 0& 0\\
   \frac{x_4}{2}-\frac{x_1x_3}{12} &  0& \frac{x_1^2}{12}& -\frac{x_1}{2}& 1& 0& 0\\
 \frac{x_1^3x_2}{720}-\frac{x_1x_4}{12}+\frac{x_5}{2} &  -\frac{x_1^4}{720}& 0& \frac{x_1^2}{12}& -\frac{x_1}{2}& 1& 0\\
 -\frac{x_2^2}{12}& \frac{x_1x_2}{12}+\frac{x_3}{2} &-\frac{x_2}{2} & 0& 0& 0& 1
   \end{matrix}\right]\,.
   \end{eqnarray*}

 \subsection*{$(23457B)$}
 
 The following Lie algebra is denoted as $(23457B)$ by Gong in \cite{Gong_Thesis}, and as $\mathcal{G}_{7,2,6}$ by Magnin in \cite{magnin}.
 
 The non-trivial brackets are the following
 \begin{equation*}
     [X_1,X_i]=X_{i+1}\,,\,i=2,3,4\,,\,[X_2,X_5]=X_6\,,\,[X_2,X_3]=X_7\,,\,[X_3,X_4]=-X_6\,.
 \end{equation*}
This is a nilpotent Lie algebra of rank 2 and step 5 that is stratifiable. The Lie brackets can be pictured with the diagram:
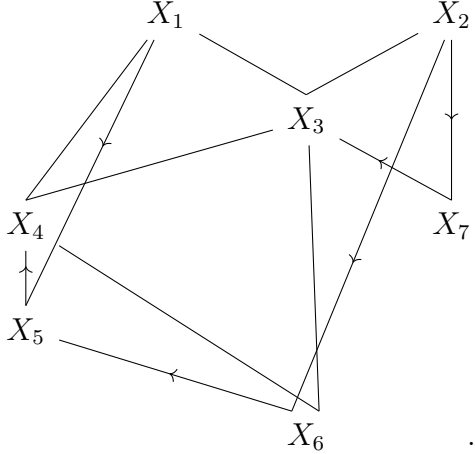
\begin{center}
 
	\begin{tikzcd}[end anchor=north]
		& X_1\ar[dr,no head]\ar[ddl, no head]\ar[dddl, no head, ->-=.4] & & X_2\ar[dd, no head, ->-=.5]\ar[dl, no head]\ar[ddddl, no head, ->-=.6,end anchor={[xshift=-3.5ex]north east}]\\
		& & X_3\ar[dll, no head]\ar[dr, no head, -<-=.4]\ar[ddd, no head, end anchor={[xshift=-1.5ex]north east}] & \\
		X_4\ar[d, no head, -<-=.4]\ar[ddrr, no head, end anchor={[xshift=-1.5ex]north east}]& &  & X_7\\
		X_5\ar[drr, no head, -<-=.5, end anchor={[xshift=-3.5ex]north east}]& &  & \\
		& & X_6 &\quad\;. \\
	\end{tikzcd}
	\end{center}
	
The composition law \eqref{group law in G} of $(23457B)$ is given by:
\begin{itemize}
    \item $z_1=x_1+y_1$;
    \item $z_2=x_2+y_2$;
    \item $z_3=x_3+y_3+\frac{1}{2}(x_1y_2-x_2y_1)$;
    \item $z_4=x_4+y_4+\frac{1}{2}(x_1y_3-x_3y_1)+\frac{1}{12}(x_1-y_1)(x_1y_2-x_2y_1)$;
    \item $z_5=x_5+y_5+\frac{1}{2}(x_1y_4-x_4y_1)+\frac{1}{12}(x_1-y_1)(x_1y_3-x_3y_1)-\frac{1}{24}x_1y_1(x_1y_2-x_2y_1)$;
    \item $z_6=x_6+y_6+\frac{1}{2}(x_2y_5-x_5y_2-x_3y_4+x_4y_3)+\frac{1}{12}(x_2-y_2)(x_1y_4-x_4y_1)\salto+\frac{1}{12}(y_3-x_3)(x_1y_3-x_3y_1)+\frac{1}{12}(x_4-y_4)(x_1y_2-x_2y_1)-\frac{1}{120}x_1(x_1y_2-x_2y_1)^2\salto-\frac{1}{24}[x_1y_2(x_1y_3-x_3y_1)-x_1y_3(x_1y_2-x_2y_1)]+\frac{1}{720}(y_1^2y_2-x_1^2x_2)(x_1y_2-x_2y_1)\salto+\frac{1}{180}x_1y_1y_2(x_1y_2-x_2y_1)-\frac{1}{180}x_1^2y_2(x_1y_2-x_2y_1)-\frac{1}{360}y_1(x_1y_2-x_2y_1)^2$;
    \item $z_7=x_7+y_7+\frac{1}{2}(x_2y_3-x_3y_2)+\frac{1}{12}(x_2-y_2)(x_1y_2-x_2y_1)$.
\end{itemize}

Since
\begin{eqnarray*}
  \mathrm{d}(L_\mathbf{x})_\mathbf{0}=\left[\begin{matrix}    
      1 & 0  & 0& 0& 0& 0& 0\\
   0 & 1& 0& 0& 0& 0& 0\\
   -\frac{x_2}{2} & \frac{x_1}{2}& 1& 0& 0& 0& 0\\
   -\frac{x_1x_2}{12}-\frac{x_3}{2} &\frac{x_1^2}{12} & \frac{x_1}{2}& 1& 0& 0& 0\\
   -\frac{x_1x_3}{12}-\frac{x_4}{2} &  0& \frac{x_1^2}{12}& \frac{x_1}{2}& 1& 0& 0\\
 \frac{x_1^2x_2^2}{720}+\frac{x_3^2-2x_2x_4}{12} &  \frac{x_1x_4}{12}-\frac{x_1^3x_2}{720}-\frac{x_5}{2}& \frac{x_4}{2}-\frac{x_1x_3}{12}& \frac{x_1x_2}{12}-\frac{x_3}{2}& \frac{x_2}{2}& 1& 0\\
 -\frac{x_2^2}{12}& \frac{x_1x_2}{12}-\frac{x_3}{2} &\frac{x_2}{2} & 0& 0& 0& 1
   \end{matrix}\right]\,,
   \end{eqnarray*}
the induced left-invariant vector fields \eqref{leftinvariant vf} are: 
\begin{itemize}
\item $X_1=\partial_{x_1} -\frac{x_2}{2}\partial_{x_3}-\big(\frac{x_3}{2}+\frac{x_1x_2}{12}\big)\partial_{x_4}-\big(\frac{x_1x_3}{12}+\frac{x_4}{2}\big)\partial_{x_5}+\big(\frac{x_1^2x_2^2}{720}+\frac{x_3^2-2x_2x_4}{12}\big)\partial_{x_6}-\frac{x_2^2}{12}\partial_{x_7}\,;$
\item $X_2=\partial_{x_2}+\frac{x_1}{2}\partial_{x_3}+\frac{x_1^2}{12}\partial_{x_4}+\big(\frac{x_1x_4}{12}-\frac{x_1^3x_2}{720}-\frac{x_5}{2}\big)\partial_{x_6}+\big(\frac{x_1x_2}{12}-\frac{x_3}{2}\big)\partial_{x_7}\,;$
\item $X_3=\partial_{x_3}+\frac{x_1}{2}\partial_{x_4}+\frac{x_1^2}{12}\partial_{x_5}+\big(\frac{x_4}{2}-\frac{x_1x_3}{12}\big)\partial_{x_6}+\frac{x_2}{2}\partial_{x_7}\,;$
\item $X_4=\partial_{x_4}+\frac{x_1}{2}\partial_{x_5}+\big(\frac{x_1x_2}{12}-\frac{x_3}{2}\big)\partial_{x_6}\,;$
\item $X_5=\partial_{x_5}+\frac{x_2}{2}\partial_{x_6}$;
\item $X_6=\partial_{x_6}$;
\item $X_7=\partial_{x_7}$,
\end{itemize}
and the respective left-invariant 1-forms \eqref{leftinvariant form} are: 
\begin{itemize}
\item $\theta_1=dx_1$;
\item $\theta_2=dx_2$;
\item $\theta_3=dx_3-\frac{x_1}{2}dx_2+\frac{x_2}{2}dx_1$;
\item $\theta_4=dx_4-\frac{x_1}{2}dx_3+\frac{x_1^2}{6}dx_2+\big(\frac{x_3}{2}-\frac{x_1x_2}{6}\big)dx_1$;
\item $\theta_5=dx_5-\frac{x_1}{2}dx_4+\frac{x_1^2}{6}dx_3-\frac{x_1^3}{24}dx_2+\big(\frac{x_4}{2}-\frac{x_1x_3}{6}+\frac{x_1^2x_2}{24}\big)dx_1$;
\item $\theta_6=dx_6-\frac{x_2}{2}dx_5+\big(\frac{x_1x_2}{6}+\frac{x_3}{2}\big)dx_4-\big(\frac{x_4}{2}+\frac{x_1x_3}{6}+\frac{x_1^2x_2}{24}\big)dx_3+\big(\frac{x_5}{2}+\frac{x_1x_4}{6}+\frac{x_1^2x_3}{24}\saltot+\frac{x_1^3x_2}{120}\big)dx_2+\big(\frac{x_3^2-2x_2x_4}{6}-\frac{x_1^2x_2^2}{120}\big)dx_1$;
\item $\theta_7=dx_7-\frac{x_2}{2}dx_3+\big(\frac{x_3}{2}+\frac{x_1x_2}{6}\big)dx_2-\frac{x_2^2}{6}dx_1$.
\end{itemize}

Finally, we have
\begin{eqnarray*}
  \mathrm{d}(R_\mathbf{x})_\mathbf{0}=\left[\begin{matrix}    
      1 & 0  & 0& 0& 0& 0& 0\\
   0 & 1& 0& 0& 0& 0& 0\\
   \frac{x_2}{2} & -\frac{x_1}{2}& 1& 0& 0& 0& 0\\
   \frac{x_3}{2}-\frac{x_1x_2}{12} &\frac{x_1^2}{12} & -\frac{x_1}{2}& 1& 0& 0& 0\\
   \frac{x_4}{2}-\frac{x_1x_3}{12} &  0& \frac{x_1^2}{12}& -\frac{x_1}{2}& 1& 0& 0\\
 \frac{x_1^2x_2^2}{720}+\frac{x_3^2-2x_2x_4}{12} &  \frac{x_1x_4}{12}-\frac{x_1^3x_2}{720}+\frac{x_5}{2}& -\frac{x_4}{2}-\frac{x_1x_3}{12}& \frac{x_1x_2}{12}+\frac{x_3}{2}& -\frac{x_2}{2}& 1& 0\\
 -\frac{x_2^2}{12}& \frac{x_1x_2}{12}+\frac{x_3}{2} &-\frac{x_2}{2} & 0& 0& 0& 1
   \end{matrix}\right]\,.
   \end{eqnarray*}

 \subsection*{$(23457C)$}
 
 The following Lie algebra is denoted as $(23457C)$ by Gong in \cite{Gong_Thesis}, and as $\mathcal{G}_{7,2,4}$ by Magnin in \cite{magnin}.
 
 The non-trivial brackets are the following
 \begin{equation*}
     [X_1,X_i]=X_{i+1}\,,\,i=2,3,4\,,\,[X_1,X_5]=X_6\,,\,[X_2,X_5]=X_7\,,\,[X_3,X_4]=-X_7\,.
 \end{equation*}
This is a nilpotent Lie algebra of rank 2 and step 5 that is stratifiable. The Lie brackets can be pictured with the diagram:
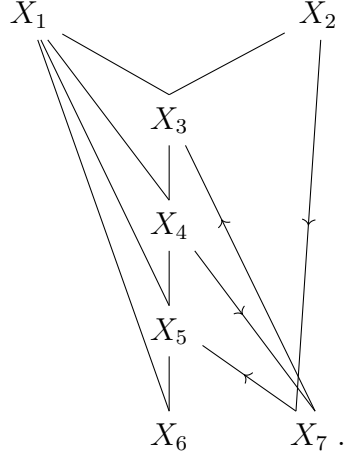
\begin{center}
 
	\begin{tikzcd}[end anchor=north]
		 X_1\ar[dr,no head]\ar[ddr, no head]\ar[dddr, no head]\ar[ddddr, no head] &  & X_2\ar[dl, no head]\ar[dddd, no head, ->-=.5, end anchor={[xshift=-4.7ex]north east}]\\
		 & X_3\ar[d, no head]\ar[dddr, no head, -<-=.3, end anchor={[xshift=-3.3ex]north east}] & \\
		 & X_4\ar[d, no head]\ar[ddr, no head, ->-=.4, end anchor={[xshift=-3.3ex]north east}]  & \\
		 & X_5 \ar[d, no head]\ar[dr, no head, -<-=.5, end anchor={[xshift=-4.7ex]north east}]& \\
		 & X_6 & X_7\;. \\
	\end{tikzcd}
	\end{center}
	
The composition law \eqref{group law in G} of $(23457C)$ is given by:
\begin{itemize}
    \item $z_1=x_1+y_1$;
    \item $z_2=x_2+y_2$;
    \item $z_3=x_3+y_3+\frac{1}{2}(x_1y_2-x_2y_1)$;
    \item $z_4=x_4+y_4+\frac{1}{2}(x_1y_3-x_3y_1)+\frac{1}{12}(x_1-y_1)(x_1y_2-x_2y_1)$;
    \item $z_5=x_5+y_5+\frac{1}{2}(x_1y_4-x_4y_1)+\frac{1}{12}(x_1-y_1)(x_1y_3-x_3y_1)-\frac{1}{24}x_1y_1(x_1y_2-x_2y_1)$;
    \item $z_6=x_6+y_6+\frac{1}{2}(x_1y_5-x_5y_1)+\frac{1}{12}(x_1-y_1)(x_1y_4-x_4y_1)-\frac{1}{24}x_1y_1(x_1y_3-x_3y_1)\salto+\frac{1}{180}(x_1y_1^2-x_1^2y_1)(x_1y_2-x_2y_1)+\frac{1}{720}(y_1^3-x_1^3)(x_1y_2-x_2y_1)$;
    \item $z_7=x_7+y_7+\frac{1}{2}(x_2y_5-x_5y_2-x_3y_4+x_4y_3)+\frac{1}{12}(x_2-y_2)(x_1y_4-x_4y_1)\salto\frac{1}{12}(y_3-x_3)(x_1y_3-x_3y_1)+\frac{1}{12}(x_4-y_4)(x_1y_2-x_2y_1)-\frac{1}{24}x_1y_2(x_1y_3-x_3y_1)\salto-\frac{1}{24}x_1y_3(x_1y_2-x_2y_1)+\frac{1}{180}(x_1y_1y_2-x_1^2y_2)(x_1y_2-x_2y_1)\salto+\frac{1}{720}(y_1^2y_1-x_1^2x_2)(x_1y_2-x_2y_1)-\frac{1}{360}y_1(x_1y_2-x_2y_1)^2-\frac{1}{120}x_1(x_1y_2-x_2y_1)^2$.
\end{itemize}

Since
\begin{eqnarray*}
  \mathrm{d}(L_\mathbf{x})_\mathbf{0}=\left[\begin{matrix}    
      1 & 0  & 0& 0& 0& 0& 0\\
   0 & 1& 0& 0& 0& 0& 0\\
   -\frac{x_2}{2} & \frac{x_1}{2}& 1& 0& 0& 0& 0\\
   -\frac{x_1x_2}{12}-\frac{x_3}{2} &\frac{x_1^2}{12} & \frac{x_1}{2}& 1& 0& 0& 0\\
   -\frac{x_1x_3}{12}-\frac{x_4}{2} &  0& \frac{x_1^2}{12}& \frac{x_1}{2}& 1& 0& 0\\
 \frac{x_1^3x_2}{720}-\frac{x_1x_4}{12}-\frac{x_5}{2} &  -\frac{x_1^4}{720}& 0& \frac{x_1^2}{12}& \frac{x_1}{2}& 1& 0\\
 \frac{x_1^2x_2^2}{720}+\frac{x_3^2-2x_2x_4}{12}& \frac{x_1x_4}{12}-\frac{x_1^3x_2}{720}-\frac{x_5}{2} &\frac{x_4}{2}-\frac{x_1x_3}{12} & \frac{x_1x_2}{12}-\frac{x_3}{2}& \frac{x_2}{2}& 0& 1
   \end{matrix}\right]\,,
   \end{eqnarray*}
the induced left-invariant vector fields \eqref{leftinvariant vf} are: 
\begin{itemize}
\item $X_1=\partial_{x_1} -\frac{x_2}{2}\partial_{x_3}-\big(\frac{x_3}{2}+\frac{x_1x_2}{12}\big)\partial_{x_4}-\big(\frac{x_1x_3}{12}+\frac{x_4}{2}\big)\partial_{x_5}+\big(\frac{x_1^3x_2}{720}-\frac{x_1x_4}{12}-\frac{x_5}{2}\big)\partial_{x_6}\saltox+\big(\frac{x_1^2x_2^2}{720}+\frac{x_3^2-2x_2x_4}{12}\big)\partial_{x_7}\,;$
\item $X_2=\partial_{x_2}+\frac{x_1}{2}\partial_{x_3}+\frac{x_1^2}{12}\partial_{x_4}-\frac{x_1^4}{720}\partial_{x_6}+\big(\frac{x_1x_4}{12}-\frac{x_5}{2}-\frac{x_1^3x_2}{720}\big)\partial_{x_7}\,;$
\item $X_3=\partial_{x_3}+\frac{x_1}{2}\partial_{x_4}+\frac{x_1^2}{12}\partial_{x_5}+\big(\frac{x_4}{2}-\frac{x_1x_3}{12}\big)\partial_{x_7}\,;$
\item $X_4=\partial_{x_4}+\frac{x_1}{2}\partial_{x_5}+\frac{x_1^2}{12}\partial_{x_6}+\big(\frac{x_1x_2}{12}-\frac{x_3}{2}\big)\partial_{x_7}\,;$
\item $X_5=\partial_{x_5}+\frac{x_1}{2}\partial_{x_6}+\frac{x_2}{2}\partial_{x_7}$;
\item $X_6=\partial_{x_6}$;
\item $X_7=\partial_{x_7}$,
\end{itemize}
and the respective left-invariant 1-forms \eqref{leftinvariant form} are: 
\begin{itemize}
\item $\theta_1=dx_1$;
\item $\theta_2=dx_2$;
\item $\theta_3=dx_3-\frac{x_1}{2}dx_2+\frac{x_2}{2}dx_1$;
\item $\theta_4=dx_4-\frac{x_1}{2}dx_3+\frac{x_1^2}{6}dx_2+\big(\frac{x_3}{2}-\frac{x_1x_2}{6}\big)dx_1$;
\item $\theta_5=dx_5-\frac{x_1}{2}dx_4+\frac{x_1^2}{6}dx_3-\frac{x_1^3}{24}dx_2+\big(\frac{x_4}{2}-\frac{x_1x_3}{6}+\frac{x_1^2x_2}{24}\big)dx_1$;
\item $\theta_6=dx_6-\frac{x_1}{2}dx_5+\frac{x_1^2}{6}dx_4-\frac{x_1^3}{24}dx_3+\frac{x_1^4}{120}dx_2+\big(\frac{x_5}{2}-\frac{x_1x_4}{6}+\frac{x_1^2x_3}{24}-\frac{x_1^3x_2}{120}\big)dx_1$;
\item $\theta_7=dx_7-\frac{x_2}{2}dx_5+\big(\frac{x_3}{2}+\frac{x_1x_2}{6}\big)dx_4-\big(\frac{x_4}{2}+\frac{x_1x_3}{6}+\frac{x_1^2x_2}{24}\big)dx_3+\big(\frac{x_5}{2}+\frac{x_1x_4}{6}+\frac{x_1^2x_3}{24}\saltot+\frac{x_1^3x_2}{120}\big)dx_2+\big(\frac{x_3^2-2x_2x_4}{6}-\frac{x_1^2x_2^2}{120}\big)dx_1$.
\end{itemize}

Finally, we have
\begin{eqnarray*}
  \mathrm{d}(R_\mathbf{x})_\mathbf{0}=\left[\begin{matrix}    
      1 & 0  & 0& 0& 0& 0& 0\\
   0 & 1& 0& 0& 0& 0& 0\\
   \frac{x_2}{2} & -\frac{x_1}{2}& 1& 0& 0& 0& 0\\
   \frac{x_3}{2}-\frac{x_1x_2}{12} &\frac{x_1^2}{12} & -\frac{x_1}{2}& 1& 0& 0& 0\\
   \frac{x_4}{2}-\frac{x_1x_3}{12} &  0& \frac{x_1^2}{12}& -\frac{x_1}{2}& 1& 0& 0\\
 \frac{x_1^3x_2}{720}-\frac{x_1x_4}{12}+\frac{x_5}{2} &  -\frac{x_1^4}{720}& 0& \frac{x_1^2}{12}& -\frac{x_1}{2}& 1& 0\\
 \frac{x_1^2x_2^2}{720}+\frac{x_3^2-2x_2x_4}{12}& \frac{x_1x_4}{12}-\frac{x_1^3x_2}{720}+\frac{x_5}{2} &-\frac{x_4}{2}-\frac{x_1x_3}{12} & \frac{x_1x_2}{12}+\frac{x_3}{2}& -\frac{x_2}{2}& 0& 1
   \end{matrix}\right]\,.
   \end{eqnarray*}

 \subsection*{$(17)$ }
 
 The following Lie algebra is denoted as $({17})$ by Gong in \cite{Gong_Thesis}, and as $\mathcal{G}_{7,4,4}$ by Magnin in \cite{magnin}. This is the 7-dimensional Heisenberg Lie algebra.
 
 The non-trivial brackets are the following:
\begin{equation*}
   [X_1, X_2] = [X_3,X_4]=[X_5,X_6]=X_7\,.
\end{equation*}
This is a nilpotent Lie algebra of rank 6 and step 2 that is stratifiable. The Lie brackets can be pictured with the diagram:
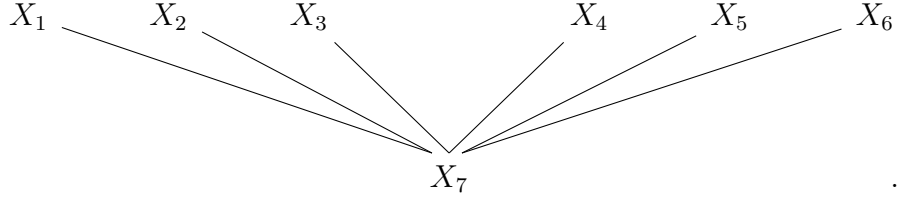
\begin{center}
 
	\begin{tikzcd}[end anchor=north]
		X_1\ar[ddrrr, no head, end anchor={[xshift=-3.7ex]north east}] & X_2\ar[ddrr, no head, end anchor={[xshift=-3.7ex]north east}] & X_3\ar[ddr, no head] & &X_4\ar[ddl, no head] & X_5\ar[ddll, no head, end anchor={[xshift=-1.5ex]north east}] & X_6\ar[ddlll, no head, end anchor={[xshift=-1.5ex]north east}] \\
		& & & & & & \\
		 &  & & X_7 & & &\quad\;.
	\end{tikzcd}

\end{center}

 

The composition law \eqref{group law in G} of $(17)$ is given by:

\begin{itemize}
    \item $z_1=x_1+y_1$;
    \item $z_2=x_2+y_2$;
    \item $z_3=x_3+y_3$;
    \item $z_4=x_4+y_4$;
    \item $z_5=x_5+y_5$;
    \item $z_6=x_6+y_6$;
    \item $z_7=x_7+y_7+\frac{1}{2}(x_1y_2-x_2y_1+x_3y_4-x_4y_3+x_5y_6-x_6y_5)$.
\end{itemize}

Since 

\begin{eqnarray*}
 \mathrm{d}(L_\mathbf{x})_\mathbf{0}= \left[\begin{matrix} 
   1 & 0  & 0& 0& 0& 0& 0\\
   0 & 1& 0& 0& 0& 0& 0\\
  0 & 0& 1& 0& 0& 0& 0\\
   0 & 0& 0& 1& 0& 0& 0\\
   0 &  0& 0& 0& 1& 0& 0\\
 0 &  0& 0& 0& 0& 1& 0\\
 -\frac{x_2}{2} &  \frac{x_1}{2}& -\frac{x_4}{2}& \frac{x_3}{2}& -\frac{x_6}{2}& \frac{x_5}{2}& 1
   \end{matrix}\right]\,,
   \end{eqnarray*}
the induced left-invariant vector fields \eqref{leftinvariant vf} are: 
\begin{itemize}
\item $X_1={\partial}_{x_1}-\frac{x_2}{2}{\partial}_{x_7}\,;$
\item $X_2={\partial}_{x_2}+\frac{x_1}{2}{\partial}_{x_7}\,;$
\item $X_3={\partial}_{x_3}-\frac{x_4}{2}{\partial}_{x_7}\,;$
\item $X_4={\partial}_{x_4}+\frac{x_3}{2}{\partial}_{x_7}\,;$
\item $X_5={\partial}_{x_5}-\frac{x_6}{2}\partial_{x_7}$;
\item $X_6={\partial}_{x_6}+\frac{x_5}{2}\partial_{x_7}$;
\item $X_7={\partial}_{x_7}\,,$
\end{itemize}
and the respective left-invariant 1-forms \eqref{leftinvariant form} are: 
\begin{itemize}
\item $\theta_1=dx_1$;
\item $\theta_2=dx_2$;
\item $\theta_3=dx_3$;
\item $\theta_4=dx_4$;
\item $\theta_5=dx_5$;
\item $\theta_6=dx_6$;
\item $\theta_7=dx_7-\frac{x_5}{2}dx_6+\frac{x_6}{2}dx_5-\frac{x_3}{2}dx_4+\frac{x_4}{2}dx_3-\frac{x_1}{2}dx_2+\frac{x_2}{2}dx_1$.
\end{itemize}

Finally, we have
\begin{eqnarray*}
 \mathrm{d}(R_\mathbf{x})_\mathbf{0}= \left[\begin{matrix} 
   1 & 0  & 0& 0& 0& 0& 0\\
   0 & 1& 0& 0& 0& 0& 0\\
  0 & 0& 1& 0& 0& 0& 0\\
   0 & 0& 0& 1& 0& 0& 0\\
   0 &  0& 0& 0& 1& 0& 0\\
 0 &  0& 0& 0& 0& 1& 0\\
 \frac{x_2}{2} &  -\frac{x_1}{2}& \frac{x_4}{2}& -\frac{x_3}{2}& \frac{x_6}{2}& -\frac{x_5}{2}& 1
   \end{matrix}\right]\,.
   \end{eqnarray*}
 \subsection*{$(147D)$ }
 
 The following Lie algebra is denoted as $({147D})$ by Gong in \cite{Gong_Thesis}, and as $\mathcal{G}_{7,2,2}$ by Magnin in \cite{magnin}.
 
The non-trivial brackets are the following:
\begin{equation*}
\begin{aligned}
   &[X_1, X_2] = X_{4}\,,\,[X_1,X_3]=-X_6\,,\,[X_2,X_3]=X_5\,,\\ [X_1&,X_5]=[X_1,X_6]=[X_2,X_6]=X_7\,,\,[X_3,X_4]=-2X_7\,.
\end{aligned}
\end{equation*}
 This is a nilpotent Lie algebra of rank 3 and step 3 that is stratifiable. The Lie brackets can be pictured with the diagram:
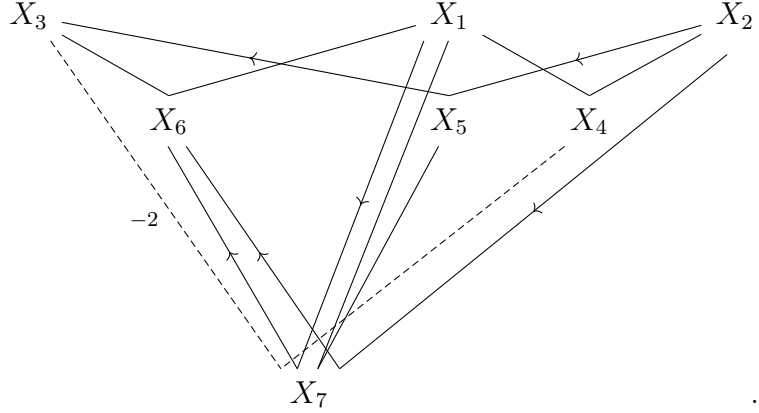
\begin{center}
 
	\begin{tikzcd}[end anchor=north]
		X_3 \ar[dr, no head]\ar[drrr, no head, -<-=.5]\ar[dddddrr, no head,dashed, "-2" ',end anchor={[xshift=-4.5ex]north east}]& & & X_1\ar[dll, no head]\ar[dddddl, no head, ->-=.5, end anchor={[xshift=-3.3ex]north east},start anchor={[xshift=-4.3ex]south east}]\ar[dddddl, no head,  end anchor={[xshift=-1.8ex]north east},start anchor={[xshift=-2.5ex]south east}]\ar[dr, no head]& & X_2\ar[dl, no head]\ar[dll, no head, ->-=.5]\ar[dddddlll, no head, ->-=.5,end anchor={[xshift=-0.2ex]north east},start anchor={[xshift=-3.ex]south east},start anchor={[yshift=-1.ex]south east}]\\
		& X_6\ar[ddddr, no head, -<-=.5, end anchor={[xshift=-3.3ex]north east},start anchor={[xshift=-2.5ex]south east}]\ar[ddddr, no head, -<-=.5, end anchor={[xshift=-0.2ex]north east}] & & X_5\ar[ddddl, no head, end anchor={[xshift=-1.8ex]north east}] & X_4\ar[ddddll, no head, dashed, end anchor={[xshift=-4.5ex]north east}] &\\
		& & & & &\\
		& & & & &\\
		& & & & &\\
		& & X_7 & & & \quad\;.
	\end{tikzcd}
\end{center}

 
The composition law \eqref{group law in G} of $(147D)$ is given by:

\begin{itemize}
    \item $z_1=x_1+y_1$;
    \item $z_2=x_2+y_2$;
    \item $z_3=x_3+y_3$;
    \item $z_4=x_4+y_4+\frac{1}{2}(x_1y_2-x_2y_1)$;
    \item $z_5=x_5+y_5+\frac{1}{2}(x_2y_3-x_3y_2)$;
    \item $z_6=x_6+y_6+\frac{1}{2}(x_3y_1-x_1y_3)$;
    \item $z_7=x_7+y_7+\frac{1}{2}(x_1y_5-x_5y_1+x_1y_6-x_6y_1+x_2y_6-x_6y_2-2x_3y_4+2x_4y_3)\salto+\frac{1}{12}[x_1(x_2y_3-x_3y_2)-x_2(x_1y_3-x_3y_1)-x_1(x_1y_3-x_3y_1)-2x_3(x_1y_2-x_2y_1)]\salto-\frac{1}{12}[y_1(x_2y_3-x_3y_2)-y_2(x_1y_3-x_3y_1)-y_1(x_1y_3-x_3y_1)-2y_3(x_1y_2-x_2y_1)]$.
\end{itemize}



Since
\begin{eqnarray*}
  \mathrm{d}(L_\mathbf{x})_\mathbf{0}=\left[\begin{matrix}    
      
   1 & 0  & 0& 0& 0& 0& 0\\
   0 & 1& 0& 0& 0& 0& 0\\
   0 & 0& 1& 0& 0& 0& 0\\
   -\frac{x_2}{2} & \frac{x_1}{2}& 0& 1& 0& 0& 0\\
  0 &  -\frac{x_3}{2}& \frac{x_2}{2}& 0& 1& 0& 0\\
 \frac{x_3}{2} &  0&-\frac{x_1}{2}& 0& 0& 1& 0\\
 \frac{x_1x_3+3x_2x_3}{12}-\frac{x_5+x_6}{2}& -\frac{x_6}{2}-\frac{x_1x_3}{4}&{x_4}-\frac{x_1^2}{12} & -x_3& \frac{x_1}{2}& \frac{x_1+x_2}{2}& 1
   \end{matrix}\right]\,,
   \end{eqnarray*}
the induced left-invariant vector fields \eqref{leftinvariant vf} are:
\begin{itemize}
\item $X_1=\partial_{x_1} -\frac{x_2}{2}\partial_{x_4}+\frac{x_3}{2}\partial_{x_6}+\big(\frac{x_1x_3+3x_2x_3}{12}-\frac{x_5+x_6}{2}\big)\partial_{x_7}\,;$
\item $X_2=\partial_{x_2}+\frac{x_1}{2}\partial_{x_4}-\frac{x_3}{2}\partial_{x_5}-\big(\frac{x_6}{2}+\frac{x_1x_3}{4}\big)\partial_{x_7}\,;$
\item $X_3=\partial_{x_3}+\frac{x_2}{2}\partial_{x_5}-\frac{x_1}{2}\partial_{x_6}+\big({x_4}-\frac{x_1^2}{12}\big)\partial_{x_7}\,;$
\item $X_4=\partial_{x_4}-x_3\partial_{x_7}\,;$
\item $X_5=\partial_{x_5}+\frac{x_1}{2}\partial_{x_7}$;
\item $X_6=\partial_{x_6}+\frac{x_1+x_2}{2}\partial_{x_7}$;
\item $X_7=\partial_{x_7}$,
\end{itemize}
and the respective left-invariant 1-forms \eqref{leftinvariant form} are: 
\begin{itemize}
\item $\theta_1=dx_1$;
\item $\theta_2=dx_2$;
\item $\theta_3=dx_3$;
\item $\theta_4=dx_4-\frac{x_1}{2}dx_2+\frac{x_2}{2}dx_1$;
\item $\theta_5=dx_5-\frac{x_2}{2}dx_3+\frac{x_3}{2}dx_2$;
\item $\theta_6=dx_6+\frac{x_1}{2}dx_3-\frac{x_3}{2}dx_1$;
\item $\theta_7=dx_7-\frac{x_1+x_2}{2}dx_6-\frac{x_1}{2}dx_5+x_3dx_4-\big(x_4+\frac{x_1^2}{6}\big)dx_3+\big(\frac{x_6}{2}-\frac{x_1x_3}{2}\big)dx_2+\big(\frac{x_5+x_6}{2}\saltot+\frac{x_1x_3+3x_2x_3}{6}\big)dx_1$.
\end{itemize}

Finally, we have
\begin{eqnarray*}
  \mathrm{d}(R_\mathbf{x})_\mathbf{0}=\left[\begin{matrix}    
   1 & 0  & 0& 0& 0& 0& 0\\
   0 & 1& 0& 0& 0& 0& 0\\
   0 & 0& 1& 0& 0& 0& 0\\
   \frac{x_2}{2} & -\frac{x_1}{2}& 0& 1& 0& 0& 0\\
  0 &  \frac{x_3}{2}& -\frac{x_2}{2}& 0& 1& 0& 0\\
 -\frac{x_3}{2} &  0& \frac{x_1}{2}& 0& 0& 1& 0\\
 \frac{x_1x_3+3x_2x_3}{12}+\frac{x_5+x_6}{2}& \frac{x_6}{2}-\frac{x_1x_3}{4}&-{x_4}-\frac{x_1^2}{12} & x_3& -\frac{x_1}{2}& -\frac{x_1+x_2}{2}& 1
   \end{matrix}\right]\,.
   \end{eqnarray*}
 \subsection*{$(147E)$ }
 
 The following one-parameter family of Lie algebras is denoted as $({147E})$ by Gong in \cite{Gong_Thesis}, and as $\mathcal{G}_{7,3,1(i_{\lambda})}$ with $\lambda\in ( 0,1)$ by Magnin in \cite{magnin}.
 
For $\lambda\in( 0,1)$, the non-trivial brackets are the following:
\begin{equation*}
\begin{aligned}
   &[X_1, X_2] = X_{4}\,,\,[X_1,X_3]=-X_6\,,[X_2,X_3]=X_5\,,\\ [X_1&,X_5]=-X_7\,,[X_2,X_6]=\lambda X_7\,,\,[X_3,X_4]=(1-\lambda)X_7\,.
\end{aligned}
\end{equation*}
Let us stress that when $\lambda=0,1$, the Lie algebra we obtain is isomorphic to $(247P)$.
 This is a nilpotent Lie algebra of rank 3 and step 3 that is stratifiable. The Lie brackets can be pictured with the diagram:
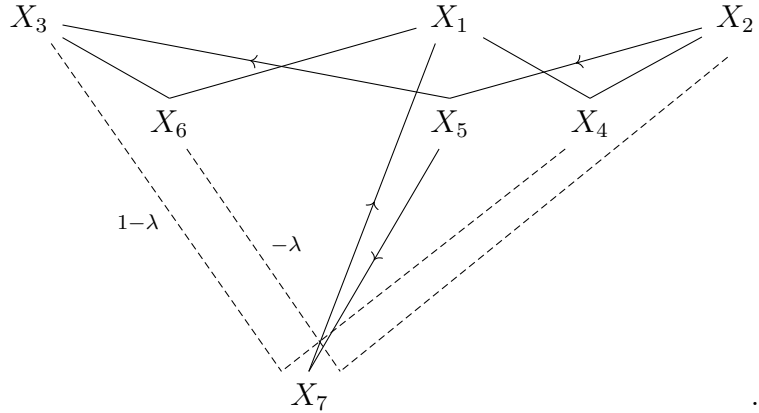
\begin{center}
 
	\begin{tikzcd}[end anchor=north]
		X_3 \ar[dr, no head]\ar[drrr, no head, -<-=.5]\ar[dddddrr, no head,dashed, "1-\lambda"', end anchor={[xshift=-4.5ex]north east}]& & & X_1\ar[dll, no head]\ar[dddddl, no head, -<-=.5, end anchor={[xshift=-2.5ex]north east},start anchor={[xshift=-3.5ex]south east}]\ar[dr, no head]& & X_2\ar[dl, no head]\ar[dll, no head, ->-=.5]\ar[dddddlll, dashed, no head,  end anchor={[xshift=-0.2ex]north east},start anchor={[xshift=-3.ex]south east},start anchor={[yshift=-1.ex]south east}]\\
		& X_6\ar[ddddr, dashed, "-\lambda", no head,   end anchor={[xshift=-0.2ex]north east}] & & X_5\ar[ddddl, no head, end anchor={[xshift=-2.5ex]north east},->-=.5] & X_4\ar[ddddll, dashed, no head, end anchor={[xshift=-4.5ex]north east}] &\\
		& & & & &\\
		& & & & &\\
		& & & & &\\
		& & X_7 & & & \quad\;.
	\end{tikzcd}
\end{center}

 
The composition law \eqref{group law in G} of $(147E)$ is given by:

\begin{itemize}
    \item $z_1=x_1+y_1$;
    \item $z_2=x_2+y_2$;
    \item $z_3=x_3+y_3$;
    \item $z_4=x_4+y_4+\frac{1}{2}(x_1y_2-x_2y_1)$;
    \item $z_5=x_5+y_5+\frac{1}{2}(x_2y_3-x_3y_2)$;
    \item $z_6=x_6+y_6+\frac{1}{2}(x_3y_1-x_1y_3)$;
    \item $z_7=x_7+y_7+\frac{1}{2}[\lambda(x_2y_6-x_6y_2)+(1-\lambda)(x_3y_4-x_4y_3)-x_1y_5+x_5y_1]\salto+\frac{1}{12}[(1-\lambda)x_3(x_1y_2-x_2y_1)-\lambda x_2(x_1y_3-x_3y_1)-x_1(x_2y_3-x_3y_2)]\salto-\frac{1}{12}[(1-\lambda)y_3(x_1y_2-x_2y_1)-\lambda y_2(x_1y_3-x_3y_1)-y_1(x_2y_3-x_3y_2)]$.
\end{itemize}



Since
\begin{eqnarray*}
  \mathrm{d}(L_\mathbf{x})_\mathbf{0}=\left[\begin{matrix}    
      
   1 & 0  & 0& 0& 0& 0& 0\\
   0 & 1& 0& 0& 0& 0& 0\\
   0 & 0& 1& 0& 0& 0& 0\\
   -\frac{x_2}{2} & \frac{x_1}{2}& 0& 1& 0& 0& 0\\
  0 &  -\frac{x_3}{2}& \frac{x_2}{2}& 0& 1& 0& 0\\
 \frac{x_3}{2} &  0&-\frac{x_1}{2}& 0& 0& 1& 0\\
 \frac{x_5}{2}+\frac{(2\lambda-1)x_2x_3}{12}& \frac{(2-\lambda)x_1x_3}{12}-\frac{\lambda x_6}{2}&-\frac{(1-\lambda)x_4}{2}-\frac{(1+\lambda)x_1x_2}{12} & \frac{(1-\lambda)x_3}{2}& -\frac{x_1}{2}& \frac{\lambda x_2}{2}& 1
   \end{matrix}\right]\,,
   \end{eqnarray*}
the induced left-invariant vector fields \eqref{leftinvariant vf} are:
\begin{itemize}
\item $X_1=\partial_{x_1} -\frac{x_2}{2}\partial_{x_4}+\frac{x_3}{2}\partial_{x_6}+\big(\frac{x_5}{2}+\frac{(2\lambda-1)x_2x_3}{12}\big)\partial_{x_7}\,;$
\item $X_2=\partial_{x_2}+\frac{x_1}{2}\partial_{x_4}-\frac{x_3}{2}\partial_{x_5}+\big(\frac{(2-\lambda)x_1x_3}{12}-\frac{\lambda x_6}{2}\big)\partial_{x_7}\,;$
\item $X_3=\partial_{x_3}+\frac{x_2}{2}\partial_{x_5}-\frac{x_1}{2}\partial_{x_6}-\big(\frac{(1-\lambda)x_4}{2}+\frac{(1+\lambda)x_1x_2}{12}\big)\partial_{x_7}\,;$
\item $X_4=\partial_{x_4}+\frac{(1-\lambda)x_3}{2}\partial_{x_7}\,;$
\item $X_5=\partial_{x_5}-\frac{x_1}{2}\partial_{x_7}$;
\item $X_6=\partial_{x_6}+\frac{\lambda x_2}{2}\partial_{x_7}$;
\item $X_7=\partial_{x_7}$,
\end{itemize}
and the respective left-invariant 1-forms \eqref{leftinvariant form} are: 
\begin{itemize}
\item $\theta_1=dx_1$;
\item $\theta_2=dx_2$;
\item $\theta_3=dx_3$;
\item $\theta_4=dx_4-\frac{x_1}{2}dx_2+\frac{x_2}{2}dx_1$;
\item $\theta_5=dx_5-\frac{x_2}{2}dx_3+\frac{x_3}{2}dx_2$;
\item $\theta_6=dx_6+\frac{x_1}{2}dx_3-\frac{x_3}{2}dx_1$;
\item $\theta_7=dx_7-\frac{\lambda x_2}{2}dx_6+\frac{x_1}{2}dx_5-\frac{(1-\lambda)x_3}{2}dx_4+\big(\frac{(1-\lambda)x_4}{2}-\frac{(1+\lambda)x_1x_2}{6}\big)dx_3+\big(\frac{\lambda x_6}{2}\saltot+\frac{(2-\lambda)x_1x_3}{6}\big)dx_2+\big(\frac{(2\lambda-1)x_2x_3}{6}-\frac{x_5}{2}\big)dx_1$.
\end{itemize}

Finally, we have
\begin{eqnarray*}
  \mathrm{d}(R_\mathbf{x})_\mathbf{0}=\left[\begin{matrix}    
      
   1 & 0  & 0& 0& 0& 0& 0\\
   0 & 1& 0& 0& 0& 0& 0\\
   0 & 0& 1& 0& 0& 0& 0\\
   \frac{x_2}{2} & -\frac{x_1}{2}& 0& 1& 0& 0& 0\\
  0 &  \frac{x_3}{2}& -\frac{x_2}{2}& 0& 1& 0& 0\\
 -\frac{x_3}{2} &  0&\frac{x_1}{2}& 0& 0& 1& 0\\
 \frac{(2\lambda-1)x_2x_3}{12}-\frac{x_5}{2}& \frac{(2-\lambda)x_1x_3}{12}+\frac{\lambda x_6}{2}&\frac{(1-\lambda)x_4}{2}-\frac{(1+\lambda)x_1x_2}{12} & -\frac{(1-\lambda)x_3}{2}& \frac{x_1}{2}& -\frac{\lambda x_2}{2}& 1
   \end{matrix}\right]\,.
   \end{eqnarray*}
 \subsection*{$(147E_1)$ }
 
 The following one-parameter family of Lie algebras is denoted as $({147E_1})$ by Gong in \cite{Gong_Thesis}, and as $\mathcal{G}_{7,3,1(i_{\lambda})}$ with $\lambda >1$ by Magnin in \cite{magnin}.
 
For $\lambda>1$, the non-trivial brackets are the following:
\begin{equation*}
\begin{aligned}
   &[X_1, X_2] = X_{4}\,,\,[X_1,X_3]=-X_6\,,\,[X_2,X_3]=X_5\,,\\  [X_2,X_6]=2&X_7\,,\,[X_1,X_6]=-\lambda X_7\,,\,[X_2,X_5]=\lambda X_7\,,\,[X_3,X_4]=-2X_7\,.
\end{aligned}
\end{equation*}
Let us stress that when $\lambda=1$, the Lie algebra we obtain is isomorphic to $(247P_1)$.
 This is a nilpotent Lie algebra of rank 3 and step 3 that is stratifiable. The Lie brackets can be pictured with the diagram:
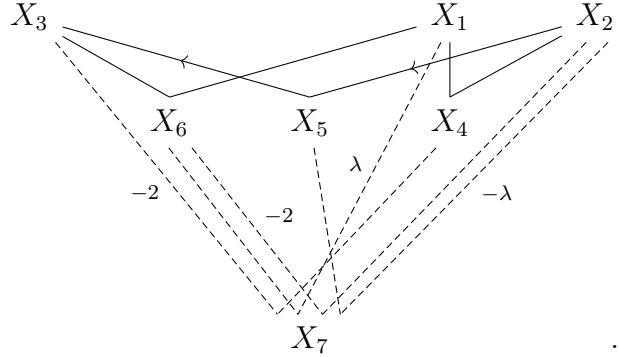
\begin{center}
 
	\begin{tikzcd}[end anchor=north]
		X_3 \ar[dr, no head]\ar[drr, no head,-<-=.5]\ar[ddddrr, no head, dashed, "-2"',end anchor={[xshift=-4.8ex]north east}]& & & X_1\ar[ddddl,dashed, "\lambda"', no head, end anchor={[xshift=-3.3ex]north east}]\ar[d, no head]\ar[dll, no head] & X_2\ar[dl, no head]\ar[dll, no head, ->-=.6]\ar[ddddll, dashed, "-\lambda", no head,   end anchor={[xshift=-.2ex]north east},start anchor={[xshift=-1.5ex]south east}]\ar[ddddll, no head, dashed , end anchor={[xshift=-1.5ex]north east},start anchor={[xshift=-3.1ex]south east}]\\
		& X_6\ar[dddr, dashed, no head,  end anchor={[xshift=-3.3ex]north east},start anchor={[xshift=-2.5ex]south east}]\ar[dddr, no head, dashed, "-2", end anchor={[xshift=-1.5ex]north east}, start anchor={[xshift=-.8ex]south east}] & X_5\ar[ddd, dashed, no head,  end anchor={[xshift=-.2ex]north east}] & X_4\ar[dddl, dashed, no head, end anchor={[xshift=-4.8ex]north east}] &\\
		& & & &\\
		& & & &\\
		& & X_7 & &\quad\;.
	\end{tikzcd}
\end{center}

 
The composition law \eqref{group law in G} of $(147E_1)$ is given by:

\begin{itemize}
    \item $z_1=x_1+y_1$;
    \item $z_2=x_2+y_2$;
    \item $z_3=x_3+y_3$;
    \item $z_4=x_4+y_4+\frac{1}{2}(x_1y_2-x_2y_1)$;
    \item $z_5=x_5+y_5+\frac{1}{2}(x_2y_3-x_3y_2)$;
    \item $z_6=x_6+y_6+\frac{1}{2}(x_3y_1-x_1y_3)$;
    \item $z_7=x_7+y_7+x_2y_6-x_6y_2-x_3y_4+x_4y_3+\frac{1}{2}[\lambda(x_2y_5-x_5y_2)-\lambda(x_1y_6-x_6y_1)]\salto+\frac{1}{12}(\lambda x_2-\lambda y_2)(x_2y_3-x_3y_2)+\frac{1}{12}(\lambda x_1-\lambda y_1)(x_1y_3-x_3y_1)\salto +\frac{1}{6}(y_2-x_2)(x_1y_3-x_3y_1)+\frac{1}{6}(y_3-x_3)(x_1y_2-x_2y_1)$.
\end{itemize}



Since
\begin{eqnarray*}
  \mathrm{d}(L_\mathbf{x})_\mathbf{0}=\left[\begin{matrix}    
      
   1 & 0  & 0& 0& 0& 0& 0\\
   0 & 1& 0& 0& 0& 0& 0\\
   0 & 0& 1& 0& 0& 0& 0\\
   -\frac{x_2}{2} & \frac{x_1}{2}& 0& 1& 0& 0& 0\\
  0 &  -\frac{x_3}{2}& \frac{x_2}{2}& 0& 1& 0& 0\\
 \frac{x_3}{2} &  0&-\frac{x_1}{2}& 0& 0& 1& 0\\
 a_1& a_2& x_4+\frac{\lambda x_1^2-2 x_1x_2+\lambda x_2^2}{12} & -x_3& \frac{\lambda x_2}{2}& \frac{2x_2-\lambda x_1}{2}& 1
   \end{matrix}\right]\,,
   \end{eqnarray*}
    where
   \begin{eqnarray*}
   a_1&=&\frac{\lambda x_6}{2}+\frac{4x_2x_3-\lambda x_1x_3}{12}\,;\\ 
   a_2&=&-\frac{\lambda x_5+2x_6}{2}-\frac{2x_1x_3+\lambda x_2x_3}{12}\,,
   \end{eqnarray*}
the induced left-invariant vector fields \eqref{leftinvariant vf} are:
\begin{itemize}
\item $X_1=\partial_{x_1} -\frac{x_2}{2}\partial_{x_4}+\frac{x_3}{2}\partial_{x_6}+\big(\frac{\lambda x_6}{2}+\frac{4x_2x_3-\lambda x_1x_3}{12}\big)\partial_{x_7}\,;$
\item $X_2=\partial_{x_2}+\frac{x_1}{2}\partial_{x_4}-\frac{x_3}{2}\partial_{x_5}-\big(\frac{\lambda x_5+2x_6}{2}+\frac{2x_1x_3+\lambda x_2x_3}{12}\big)\partial_{x_7}\,;$
\item $X_3=\partial_{x_3}+\frac{x_2}{2}\partial_{x_5}-\frac{x_1}{2}\partial_{x_6}+\big(x_4+\frac{\lambda x_1^2-2x_1x_2+\lambda x_2^2}{12}\big)\partial_{x_7}\,;$
\item $X_4=\partial_{x_4}-x_3\partial_{x_7}\,;$
\item $X_5=\partial_{x_5}+\frac{\lambda x_2}{2}\partial_{x_7}$;
\item $X_6=\partial_{x_6}+\frac{2 x_2-\lambda x_1}{2}\partial_{x_7}$;
\item $X_7=\partial_{x_7}$,
\end{itemize}
and the respective left-invariant 1-forms \eqref{leftinvariant form} are: 
\begin{itemize}
\item $\theta_1=dx_1$;
\item $\theta_2=dx_2$;
\item $\theta_3=dx_3$;
\item $\theta_4=dx_4-\frac{x_1}{2}dx_2+\frac{x_2}{2}dx_1$;
\item $\theta_5=dx_5-\frac{x_2}{2}dx_3+\frac{x_3}{2}dx_2$;
\item $\theta_6=dx_6+\frac{x_1}{2}dx_3-\frac{x_3}{2}dx_1$;
\item $\theta_7=dx_7+\frac{\lambda x_1-2x_2}{2}dx_6-\frac{\lambda x_2}{2}dx_5+x_3dx_4+\big(\frac{\lambda x_1^2-2x_1x_2+\lambda x_2^2}{6}-x_4\big)dx_3+\big(\frac{\lambda x_5+2x_6}{2}\saltot-\frac{\lambda x_2x_3+2x_1x_3}{6}\big)dx_2-\big(\frac{\lambda x_6}{2}+\frac{\lambda x_1x_3-4x_2x_3}{6}\big)dx_1$.
\end{itemize}

Finally, we have
\begin{eqnarray*}
  \mathrm{d}(R_\mathbf{x})_\mathbf{0}=\left[\begin{matrix}    
   1 & 0  & 0& 0& 0& 0& 0\\
   0 & 1& 0& 0& 0& 0& 0\\
   0 & 0& 1& 0& 0& 0& 0\\
   \frac{x_2}{2} & -\frac{x_1}{2}& 0& 1& 0& 0& 0\\
  0 &  \frac{x_3}{2}& -\frac{x_2}{2}& 0& 1& 0& 0\\
 -\frac{x_3}{2} &  0&\frac{x_1}{2}& 0& 0& 1& 0\\
 a_1& a_2& \frac{\lambda x_1^2-2 x_1x_2+\lambda x_2^2}{12}-x_4 & x_3& -\frac{\lambda x_2}{2}& \frac{\lambda x_1-2x_2}{2}& 1
   \end{matrix}\right]\,,
      \end{eqnarray*}
     where
   \begin{eqnarray*}
   a_1&=&\frac{4x_2x_3-\lambda x_1x_3}{12}-\frac{\lambda x_6}{2}\,;\\ 
   a_2&=&\frac{\lambda x_5+2x_6}{2}-\frac{2x_1x_3+\lambda x_2x_3}{12}\,,
   \end{eqnarray*}

 \subsection*{$(137A)$ }
 
 The following Lie algebra is denoted as $({137A})$ by Gong in \cite{Gong_Thesis}, and as $\mathcal{G}_{7,3,16}$ by Magnin in \cite{magnin}.
 
The non-trivial brackets are the following:
\begin{equation*}
   [X_1, X_2] = X_{5}\,,\,[X_1,X_5]=X_7\,,[X_3,X_4]=X_6\,,\,[X_3,X_6]=X_7\,.
\end{equation*}
 This is a nilpotent Lie algebra of rank 4 and step 3 that is stratifiable. The Lie brackets can be pictured with the diagram:
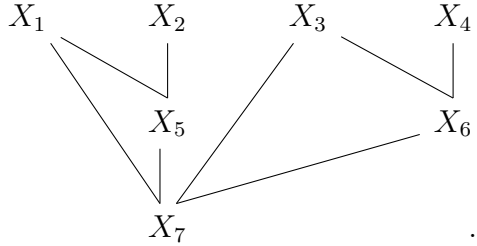
\begin{center}
 
	\begin{tikzcd}[end anchor=north]
	X_1\ar[dr, no head]\ar[ddr, no head,end anchor={[xshift=-3.ex]north east}]& X_2\ar[d, no head] & X_3 \ar[dr, no head]\ar[ddl, no head,end anchor={[xshift=-1.8ex]north east}]& X_4\ar[d, no head]\\
	 & X_5\ar[d, no head,end anchor={[xshift=-3.ex]north east},start anchor={[xshift=-3.ex]south east}] & & X_6\ar[dll, no head,end anchor={[xshift=-1.8ex]north east}]\\
	 & X_7 & &\quad\;.
	\end{tikzcd}
\end{center}

 
The composition law \eqref{group law in G} of $(137A)$ is given by:

\begin{itemize}
    \item $z_1=x_1+y_1$;
    \item $z_2=x_2+y_2$;
    \item $z_3=x_3+y_3$;
    \item $z_4=x_4+y_4$;
    \item $z_5=x_5+y_5+\frac{1}{2}(x_1y_2-x_2y_1)$;
    \item $z_6=x_6+y_6+\frac{1}{2}(x_3y_4-x_4y_3)$;
    \item $z_7=x_7+y_7+\frac{1}{2}(x_1y_5-x_5y_1+x_3y_6-x_6y_3)+\frac{1}{12}(x_1-y_1)(x_1y_2-x_2y_1)\salto+\frac{1}{12}(x_3-y_3)(x_3y_4-x_4y_3)$.
\end{itemize}



Since
\begin{eqnarray*}
  \mathrm{d}(L_\mathbf{x})_\mathbf{0}=\left[\begin{matrix}    
      
   1 & 0  & 0& 0& 0& 0& 0\\
   0 & 1& 0& 0& 0& 0& 0\\
   0 & 0& 1& 0& 0& 0& 0\\
   0 & 0& 0& 1& 0& 0& 0\\
  -\frac{x_2}{2} &  \frac{x_1}{2}& 0& 0& 1& 0& 0\\
 0 &  0&-\frac{x_4}{2}& \frac{x_3}{2}& 0& 1& 0\\
 -\frac{x_1x_2}{12}-\frac{x_5}{2}& \frac{x_1^2}{12}&-\frac{x_6}{2}-\frac{x_3x_4}{12} & \frac{x_3^2}{12}& \frac{x_1}{2}& \frac{x_3}{2}& 1
   \end{matrix}\right]\,,
   \end{eqnarray*}
the induced left-invariant vector fields \eqref{leftinvariant vf} are:
\begin{itemize}
\item $X_1=\partial_{x_1} -\frac{x_2}{2}\partial_{x_5}-\big(\frac{x_1x_2}{12}+\frac{x_5}{2}\big)\partial_{x_7}\,;$
\item $X_2=\partial_{x_2}+\frac{x_1}{2}\partial_{x_5}+\frac{x_1^2}{12}\partial_{x_7}\,;$
\item $X_3=\partial_{x_3}-\frac{x_4}{2}\partial_{x_6}-\big(\frac{x_6}{2}+\frac{x_3x_4}{12}\big)\partial_{x_7}\,;$
\item $X_4=\partial_{x_4}+\frac{x_3}{2}\partial_{x_6}+\frac{x_3^2}{12}\partial_{x_7}\,;$
\item $X_5=\partial_{x_5}+\frac{x_1}{2}\partial_{x_7}$;
\item $X_6=\partial_{x_6}+\frac{x_3}{2}\partial_{x_7}$;
\item $X_7=\partial_{x_7}$,
\end{itemize}
and the respective left-invariant 1-forms \eqref{leftinvariant form} are: 
\begin{itemize}
\item $\theta_1=dx_1$;
\item $\theta_2=dx_2$;
\item $\theta_3=dx_3$;
\item $\theta_4=dx_4$;
\item $\theta_5=dx_5-\frac{x_1}{2}dx_2+\frac{x_2}{2}dx_1$;
\item $\theta_6=dx_6-\frac{x_3}{2}dx_4+\frac{x_4}{2}dx_3$;
\item $\theta_7=dx_7-\frac{x_3}{2}dx_6-\frac{x_1}{2}dx_5+\frac{x_3^2}{6}dx_4+\big(\frac{x_6}{2}-\frac{x_3x_4}{6}\big)dx_3+\frac{x_1^2}{6}dx_2+\big(\frac{x_5}{2}-\frac{x_1x_2}{6}\big)dx_1$.
\end{itemize}

Finally, we have
\begin{eqnarray*}
  \mathrm{d}(R_\mathbf{x})_\mathbf{0}=\left[\begin{matrix}    
    1 & 0  & 0& 0& 0& 0& 0\\
   0 & 1& 0& 0& 0& 0& 0\\
   0 & 0& 1& 0& 0& 0& 0\\
   0 & 0& 0& 1& 0& 0& 0\\
  \frac{x_2}{2} &  -\frac{x_1}{2}& 0& 0& 1& 0& 0\\
 0 &  0&\frac{x_4}{2}& -\frac{x_3}{2}& 0& 1& 0\\
 \frac{x_5}{2}-\frac{x_1x_2}{12}& \frac{x_1^2}{12}&\frac{x_6}{2}-\frac{x_3x_4}{12} & \frac{x_3^2}{12}& -\frac{x_1}{2}& -\frac{x_3}{2}& 1
   \end{matrix}\right]\,.
   \end{eqnarray*}
 \subsection*{$(137A_1)$ }
 
 The following Lie algebra is denoted as $({137A_1})$ by Gong in \cite{Gong_Thesis}, and as $\mathcal{G}_{7,3,16}$ by Magnin in \cite{magnin}.
 
The non-trivial brackets are the following:
\begin{equation*}
\begin{aligned}
   &\;[X_1, X_3] = X_{5}\,,\,[X_1,X_4]=X_6\,,[X_1,X_5]=X_7\,,\\ [&X_2,X_3]=-X_6\,,\,[X_2,X_4]=X_5\,,\,[X_2,X_6]=X_7\,.
\end{aligned}
\end{equation*}
 This is a nilpotent Lie algebra of rank 4 and step 3 that is stratifiable. The Lie brackets can be pictured with the diagram:
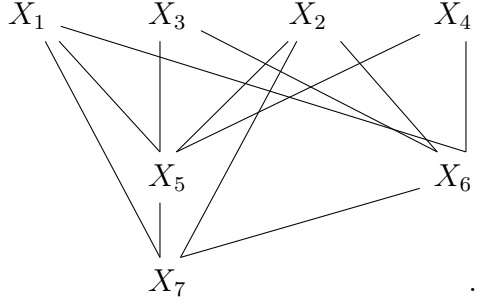
\begin{center}
 
	\begin{tikzcd}[end anchor=north]
	X_1\ar[ddr, no head,end anchor={[xshift=-3.ex]north east}]\ar[dddr, no head, end anchor={[xshift=-3.ex]north east}]\ar[ddrrr, no head, end anchor={[xshift=-1.5ex]north east}]& X_3\ar[dd, no head,end anchor={[xshift=-3.ex]north east},start anchor={[xshift=-3.ex]south east}]\ar[ddrr, no head, end anchor={[xshift=-3.6ex]north east}]& X_2\ar[dddl, no head, end anchor={[xshift=-1.5ex]north east}]\ar[ddl, no head, end anchor={[xshift=-1.8ex]north east}] \ar[ddr, no head, end anchor={[xshift=-3.6ex]north east}]& X_4\ar[dd, no head, end anchor={[xshift=-1.5ex]north east},start anchor={[xshift=-1.5ex]south east}]\ar[ddll, no head,end anchor={[xshift=-1.8ex]north east}]\\
	& & & \\
	& X_5\ar[d, no head, end anchor={[xshift=-3.ex]north east},start anchor={[xshift=-3.ex]south east}] & & X_6\ar[dll, no head, end anchor={[xshift=-1.5ex]north east}]\\
	& X_7 & &\quad\;.
	\end{tikzcd}
\end{center}

 
The composition law \eqref{group law in G} of $(137A_1)$ is given by:

\begin{itemize}
    \item $z_1=x_1+y_1$;
    \item $z_2=x_2+y_2$;
    \item $z_3=x_3+y_3$;
    \item $z_4=x_4+y_4$;
    \item $z_5=x_5+y_5+\frac{1}{2}(x_1y_3-x_3y_1+x_2y_4-x_4y_2)$;
    \item $z_6=x_6+y_6+\frac{1}{2}(x_1y_4-x_4y_1-x_2y_3+x_3y_2)$;
    \item $z_7=x_7+y_7+\frac{1}{2}(x_1y_5-x_5y_1+x_2y_6-x_6y_2)+\frac{1}{12}(x_1-y_1)(x_1y_3-x_3y_1+x_2y_4-x_4y_2)\salto+\frac{1}{12}(x_2-y_2)(x_1y_4-x_4y_1-x_2y_3+x_3y_2)]$.
\end{itemize}



Since
\begin{eqnarray*}
  \mathrm{d}(L_\mathbf{x})_\mathbf{0}=\left[\begin{matrix}    
      
   1 & 0  & 0& 0& 0& 0& 0\\
   0 & 1& 0& 0& 0& 0& 0\\
   0 & 0& 1& 0& 0& 0& 0\\
   0 & 0& 0& 1& 0& 0& 0\\
  -\frac{x_3}{2} &  -\frac{x_4}{2}& \frac{x_1}{2}& \frac{x_2}{2}& 1& 0& 0\\
 -\frac{x_4}{2} &  \frac{x_3}{2}&-\frac{x_2}{2}& \frac{x_1}{2}& 0& 1& 0\\
 -\frac{x_1x_3+x_2x_4}{12}-\frac{x_5}{2}& \frac{x_2x_3-x_1x_4}{12}-\frac{x_6}{2}&\frac{x_1^2-x_2^2}{12} & \frac{x_1x_2}{6}& \frac{x_1}{2}& \frac{x_2}{2}& 1
   \end{matrix}\right]\,,
   \end{eqnarray*}
the induced left-invariant vector fields \eqref{leftinvariant vf} are:
\begin{itemize}
\item $X_1=\partial_{x_1} -\frac{x_3}{2}\partial_{x_5}-\frac{x_4}{2}\partial_{x_6}-\big(\frac{x_5}{2}+\frac{x_1x_3+x_2x_4}{12}\big)\partial_{x_7}\,;$
\item $X_2=\partial_{x_2}-\frac{x_4}{2}\partial_{x_5}+\frac{x_3}{2}\partial_{x_6}+\big(\frac{x_2x_3-x_1x_4}{12}-\frac{x_6}{2}\big)\partial_{x_7}\,;$
\item $X_3=\partial_{x_3}+\frac{x_1}{2}\partial_{x_5}-\frac{x_2}{2}\partial_{x_6}+\frac{x_1^2-x_2^2}{12}\partial_{x_7}\,;$
\item $X_4=\partial_{x_4}+\frac{x_2}{2}\partial_{x_5}+\frac{x_1}{2}\partial_{x_6}+\frac{x_1x_2}{6}\partial_{x_7}\,;$
\item $X_5=\partial_{x_5}+\frac{x_1}{2}\partial_{x_7}$;
\item $X_6=\partial_{x_6}+\frac{x_2}{2}\partial_{x_7}$;
\item $X_7=\partial_{x_7}$,
\end{itemize}
and the respective left-invariant 1-forms \eqref{leftinvariant form} are: 
\begin{itemize}
\item $\theta_1=dx_1$;
\item $\theta_2=dx_2$;
\item $\theta_3=dx_3$;
\item $\theta_4=dx_4$;
\item $\theta_5=dx_5-\frac{x_2}{2}dx_4-\frac{x_1}{2}dx_3+\frac{x_4}{2}dx_2+\frac{x_3}{2}dx_1$;
\item $\theta_6=dx_6-\frac{x_1}{2}dx_4+\frac{x_2}{2}dx_3-\frac{x_3}{2}dx_2+\frac{x_4}{2}dx_1$;
\item $\theta_7=dx_7-\frac{x_2}{2}dx_6-\frac{x_1}{2}dx_5+\frac{x_1x_2}{3}dx_4+\frac{x_1^2-x_2^2}{6}dx_3+\big(\frac{x_6}{2}+\frac{x_2x_3-x_1x_4}{6}\big)dx_2+\big(\frac{x_5}{2}\saltot-\frac{x_1x_3+x_2x_4}{6}\big)dx_1$.
\end{itemize}

Finally, we have
\begin{eqnarray*}
  \mathrm{d}(R_\mathbf{x})_\mathbf{0}=\left[\begin{matrix}    
    1 & 0  & 0& 0& 0& 0& 0\\
   0 & 1& 0& 0& 0& 0& 0\\
   0 & 0& 1& 0& 0& 0& 0\\
   0 & 0& 0& 1& 0& 0& 0\\
  \frac{x_3}{2} &  \frac{x_4}{2}& -\frac{x_1}{2}& -\frac{x_2}{2}& 1& 0& 0\\
 \frac{x_4}{2} &  -\frac{x_3}{2}& \frac{x_2}{2}& -\frac{x_1}{2}& 0& 1& 0\\
 \frac{x_5}{2}-\frac{x_1x_3+x_2x_4}{12}& \frac{x_2x_3-x_1x_4}{12}+\frac{x_6}{2}&\frac{x_1^2-x_2^2}{12} & \frac{x_1x_2}{6}& -\frac{x_1}{2}& -\frac{x_2}{2}& 1
   \end{matrix}\right]\,.
   \end{eqnarray*}
 \subsection*{$(137C)$ }
 
 The following Lie algebra is denoted as $({137C})$ by Gong in \cite{Gong_Thesis}, and as $\mathcal{G}_{7,3,10}$ by Magnin in \cite{magnin}.
 
The non-trivial brackets are the following:
\begin{equation*}
   [X_1, X_2] = X_{5}\,,\,[X_1,X_4]=X_6\,,[X_1,X_6]=X_7\,,\,[X_2,X_3]=X_6\,,\,[X_3,X_5]=-X_7\,.
\end{equation*}
 This is a nilpotent Lie algebra of rank 4 and step 3 that is stratifiable. The Lie brackets can be pictured with the diagram:
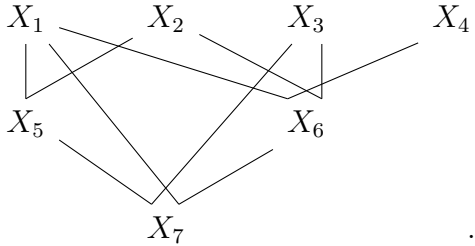
\begin{center}
 
	\begin{tikzcd}[end anchor=north]
	X_1\ar[ddr, no head, end anchor={[xshift=-1.5ex]north east}]\ar[d, no head]\ar[drr, no head, end anchor={[xshift=-3.8ex]north east}]& X_2\ar[dl, no head] \ar[dr, no head, end anchor={[xshift=-1.3ex]north east}]& X_3\ar[ddl, no head, end anchor={[xshift=-3.5ex]north east}]\ar[d, no head, end anchor={[xshift=-1.3ex]north east},start anchor={[xshift=-1.3ex]south east}] & X_4\ar[dl, no head, end anchor={[xshift=-3.8ex]north east}]\\
	X_5 \ar[dr, no head, end anchor={[xshift=-3.5ex]north east}]& & X_6\ar[dl, no head, end anchor={[xshift=-1.5ex]north east}] &\\
	&X_7 &  &\quad\;.
	\end{tikzcd}
\end{center}

 
The composition law \eqref{group law in G} of $(137C)$ is given by:

\begin{itemize}
    \item $z_1=x_1+y_1$;
    \item $z_2=x_2+y_2$;
    \item $z_3=x_3+y_3$;
    \item $z_4=x_4+y_4$;
    \item $z_5=x_5+y_5+\frac{1}{2}(x_1y_2-x_2y_1)$;
    \item $z_6=x_6+y_6+\frac{1}{2}(x_1y_4-x_4y_1+x_2y_3-x_3y_2)$;
    \item $z_7=x_7+y_7+\frac{1}{2}(x_1y_6-x_6y_1-x_3y_5+x_5y_3)+\frac{1}{12}(x_1-y_1)(x_1y_4-x_4y_1+x_2y_3-x_3y_2)\salto+\frac{1}{12}(y_3-x_3)(x_1y_2-x_2y_1)]$.
\end{itemize}



Since
\begin{eqnarray*}
  \mathrm{d}(L_\mathbf{x})_\mathbf{0}=\left[\begin{matrix}    
      
   1 & 0  & 0& 0& 0& 0& 0\\
   0 & 1& 0& 0& 0& 0& 0\\
   0 & 0& 1& 0& 0& 0& 0\\
   0 & 0& 0& 1& 0& 0& 0\\
  -\frac{x_2}{2} &  \frac{x_1}{2}& 0& 0& 1& 0& 0\\
 -\frac{x_4}{2} &  -\frac{x_3}{2}& \frac{x_2}{2}&\frac{x_1}{2} & 0& 1& 0\\
 \frac{x_2x_3-x_1x_4}{12}-\frac{x_6}{2}& -\frac{x_1x_3}{6}&\frac{x_1x_2}{12}+\frac{x_5}{2} & \frac{x_1^2}{12}& -\frac{x_3}{2}& \frac{x_1}{2}& 1
   \end{matrix}\right]\,,
   \end{eqnarray*}
the induced left-invariant vector fields \eqref{leftinvariant vf} are:
\begin{itemize}
\item $X_1=\partial_{x_1} -\frac{x_2}{2}\partial_{x_5}-\frac{x_4}{2}\partial_{x_6}+\big(\frac{x_2x_3-x_1x_4}{12}-\frac{x_6}{2}\big)\partial_{x_7}\,;$
\item $X_2=\partial_{x_2}+\frac{x_1}{2}\partial_{x_5}-\frac{x_3}{2}\partial_{x_6}-\frac{x_1x_3}{6}\partial_{x_7}\,;$
\item $X_3=\partial_{x_3}+\frac{x_2}{2}\partial_{x_6}+\big(\frac{x_1x_2}{12}+\frac{x_5}{2}\big)\partial_{x_7}\,;$
\item $X_4=\partial_{x_4}+\frac{x_1}{2}\partial_{x_6}+\frac{x_1^2}{12}\partial_{x_7}\,;$
\item $X_5=\partial_{x_5}-\frac{x_3}{2}\partial_{x_7}$;
\item $X_6=\partial_{x_6}+\frac{x_1}{2}\partial_{x_7}$;
\item $X_7=\partial_{x_7}$,
\end{itemize}
and the respective left-invariant 1-forms \eqref{leftinvariant form} are: 
\begin{itemize}
\item $\theta_1=dx_1$;
\item $\theta_2=dx_2$;
\item $\theta_3=dx_3$;
\item $\theta_4=dx_4$;
\item $\theta_5=dx_5-\frac{x_1}{2}dx_2+\frac{x_2}{2}dx_1$;
\item $\theta_6=dx_6-\frac{x_1}{2}dx_4-\frac{x_2}{2}dx_3+\frac{x_3}{2}dx_2+\frac{x_4}{2}dx_1$;
\item $\theta_7=dx_7-\frac{x_1}{2}dx_6+\frac{x_3}{2}dx_5+\frac{x_1^2}{6}dx_4+\big(\frac{x_1x_2}{6}-\frac{x_5}{2}\big)dx_3-\frac{x_1x_3}{3}dx_2+\big(\frac{x_6}{2}+\frac{x_2x_3-x_1x_4}{6}\big)dx_1$.
\end{itemize}

Finally, we have
\begin{eqnarray*}
  \mathrm{d}(R_\mathbf{x})_\mathbf{0}=\left[\begin{matrix}    
     1 & 0  & 0& 0& 0& 0& 0\\
   0 & 1& 0& 0& 0& 0& 0\\
   0 & 0& 1& 0& 0& 0& 0\\
   0 & 0& 0& 1& 0& 0& 0\\
  \frac{x_2}{2} &  -\frac{x_1}{2}& 0& 0& 1& 0& 0\\
 \frac{x_4}{2} &  \frac{x_3}{2}& -\frac{x_2}{2}&-\frac{x_1}{2} & 0& 1& 0\\
 \frac{x_2x_3-x_1x_4}{12}+\frac{x_6}{2}& -\frac{x_1x_3}{6}&\frac{x_1x_2}{12}-\frac{x_5}{2} & \frac{x_1^2}{12}& \frac{x_3}{2}& -\frac{x_1}{2}& 1
   \end{matrix}\right]\,.
   \end{eqnarray*}
 \subsection*{$(12457H)$}
 
 The following Lie algebra is denoted as $(12457H)$ by Gong in \cite{Gong_Thesis}, and as $\mathcal{G}_{7,2,5}$ by Magnin in \cite{magnin}.
 
 The non-trivial brackets are the following:
 \begin{equation*}
     [X_1,X_i]=X_{i+1}\,,\,i=2,3,5,6\,,\,[X_2,X_j]=X_{j+2}\,,\,j=3,4\,,\,[X_3,X_4]=X_7\,.
 \end{equation*}
This is a nilpotent Lie algebra of rank 2 and step 5 that is stratifiable. The Lie brackets can be pictured with the diagram:
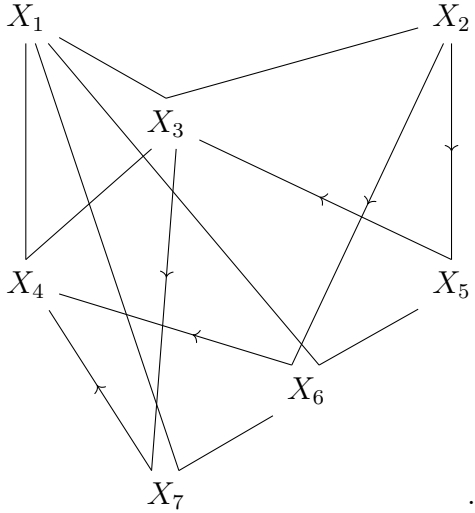
\begin{center}
 
	\begin{tikzcd}[end anchor=north]
		 X_1\ar[dr, no head]\ar[ddd, no head]\ar[ddddrr, no head,end anchor={[xshift=-1.5ex]north east}]\ar[dddddr, no head, end anchor={[xshift=-1.5ex]north east}]& & & X_2\ar[dll, no head]\ar[ddd, no head, ->-=.5]\ar[ddddl, no head, ->-=.5, end anchor={[xshift=-3.5ex]north east}]\\
		 & X_3\ar[ddl, no head,start anchor={[xshift=-3.5ex]south east}]\ar[ddrr, no head, -<-=.5]\ar[dddd, no head, ->-=.4, end anchor={[xshift=-3.5ex]north east},start anchor={[xshift=-1.7ex]south east}]& &\\
		 & & & \\
		 X_4\ar[drr, no head, -<-=.6, end anchor={[xshift=-3.5ex]north east}]\ar[ddr, no head, -<-=.5,end anchor={[xshift=-3.5ex]north east}] & &  &X_5\ar[dl, no head, end anchor={[xshift=-1.5ex]north east}]\\
		 &  &X_6 \ar[dl, no head,end anchor={[xshift=-1.5ex]north east}]& \\
		 & X_7 & &\quad\;.
	\end{tikzcd}
	\end{center}
	
The composition law \eqref{group law in G} of $(12457H)$ is given by:
\begin{itemize}
    \item $z_1=x_1+y_1$;
    \item $z_2=x_2+y_2$;
    \item $z_3=x_3+y_3+\frac{1}{2}(x_1y_2-x_2y_1)$;
    \item $z_4=x_4+y_4+\frac{1}{2}(x_1y_3-x_3y_1)+\frac{1}{12}(x_1-y_1)(x_1y_2-x_2y_1)$;
    \item $z_5=x_5+y_5+\frac{1}{2}(x_2y_3-x_3y_2)+\frac{1}{12}(x_2-y_2)(x_1y_2-x_2y_1)$;
    \item $z_6=x_6+y_6+\frac{1}{2}(x_1y_5-x_5y_1+x_2y_4-x_4y_2)+\frac{1}{12}(x_1-y_1)(x_2y_3-x_3y_2)\salto+\frac{1}{12}(x_2-y_2)(x_1y_3-x_3y_1)-\frac{1}{24}(x_1y_2+x_2y_1)(x_1y_2-x_2y_1)$;
    \item $z_7=x_7+y_7+\frac{1}{2}(x_1y_6-x_6y_1+x_3y_4-x_4y_3)+\frac{1}{12}(x_1-y_1)(x_1y_5-x_5y_1+x_2y_4-x_4y_2)\salto+\frac{1}{12}(x_3-y_3)(x_1y_3-x_3y_1)+\frac{1}{12}(y_4-x_4)(x_1y_2-x_2y_1)-\frac{1}{24}[x_1y_3(x_1y_2-x_2y_1)\salto+x_2y_1(x_1y_3-x_3y_1)+x_1y_1(x_2y_3-x_3y_2)]+\frac{1}{360}(y_1+3x_1)(x_1y_2-x_2y_1)^2\salto+\frac{1}{180}(x_2y_1^2+x_1y_1y_2-2x_1x_2y_1)(x_1y_2-x_2y_1)+\frac{2}{720}(y_1^2y_2-x_1^2x_2)(x_1y_2-x_2y_1)$.
\end{itemize}

Since
\begin{eqnarray*}
  \mathrm{d}(L_\mathbf{x})_\mathbf{0}=\left[\begin{matrix}    
      1 & 0  & 0& 0& 0& 0& 0\\
   0 & 1& 0& 0& 0& 0& 0\\
   -\frac{x_2}{2} & \frac{x_1}{2}& 1& 0& 0& 0& 0\\
   -\frac{x_1x_2}{12}-\frac{x_3}{2} &\frac{x_1^2}{12} & \frac{x_1}{2}& 1& 0& 0& 0\\
   -\frac{x_2^2}{12} &  \frac{x_1x_2}{12}-\frac{x_3}{2}& \frac{x_2}{2}& 0& 1& 0& 0\\
 -\frac{x_2x_3}{12}-\frac{x_5}{2} &  -\frac{x_1x_3}{12}-\frac{x_4}{2}& \frac{x_1x_2}{6}& \frac{x_2}{2}& \frac{x_1}{2}& 1& 0\\
 \frac{x_1^2x_2^2}{360}+\frac{x_2x_4-x_3^2-x_1x_5}{12}-\frac{x_6}{2}& -\frac{x_1x_4}{6}-\frac{x_1^3x_2}{360} &\frac{x_1x_3}{12}-\frac{x_4}{2} & \frac{x_1x_2}{12}+\frac{x_3}{2}& \frac{x_1^2}{12}& \frac{x_1}{2}& 1
   \end{matrix}\right]\,,
   \end{eqnarray*}
the induced left-invariant vector fields \eqref{leftinvariant vf} are: 
\begin{itemize}
\item $X_1=\partial_{x_1} -\frac{x_2}{2}\partial_{x_3}-\big(\frac{x_3}{2}+\frac{x_1x_2}{12}\big)\partial_{x_4}-\frac{x_2^2}{12}\partial_{x_5}-\big(\frac{x_2x_3}{12}+\frac{x_5}{2}\big)\partial_{x_6}+\big(\frac{x_1^2x_2^2}{360}+\frac{x_2x_4-x_3^2-x_1x_5}{12}\saltox-\frac{x_6}{2}\big)\partial_{x_7}\,;$
\item $X_2=\partial_{x_2}+\frac{x_1}{2}\partial_{x_3}+\frac{x_1^2}{12}\partial_{x_4}+\big(\frac{x_1x_2}{12}-\frac{x_3}{2}\big)\partial_{x_5}-\big(\frac{x_4}{2}+\frac{x_1x_3}{12}\big)\partial_{x_6}-\big(\frac{x_1x_4}{6}+\frac{x_1^3x_2}{360}\big)\partial_{x_7}\,;$
\item $X_3=\partial_{x_3}+\frac{x_1}{2}\partial_{x_4}+\frac{x_2}{2}\partial_{x_5}+\frac{x_1x_2}{6}\partial_{x_6}+\big(\frac{x_1x_3}{12}-\frac{x_4}{2}\big)\partial_{x_7}\,;$
\item $X_4=\partial_{x_4}+\frac{x_2}{2}\partial_{x_6}+\big(\frac{x_1x_2}{12}+\frac{x_3}{2}\big)\partial_{x_7}\,;$
\item $X_5=\partial_{x_5}+\frac{x_1}{2}\partial_{x_6}+\frac{x_1^2}{12}\partial_{x_7}$;
\item $X_6=\partial_{x_6}+\frac{x_1}{2}\partial_{x_7}$;
\item $X_7=\partial_{x_7}$,
\end{itemize}
and the respective left-invariant 1-forms \eqref{leftinvariant form} are: 
\begin{itemize}
\item $\theta_1=dx_1$;
\item $\theta_2=dx_2$;
\item $\theta_3=dx_3-\frac{x_1}{2}dx_2+\frac{x_2}{2}dx_1$;
\item $\theta_4=dx_4-\frac{x_1}{2}dx_3+\frac{x_1^2}{6}dx_2+\big(\frac{x_3}{2}-\frac{x_1x_2}{6}\big)dx_1$;
\item $\theta_5=dx_5-\frac{x_2}{2}dx_3+\big(\frac{x_1x_2}{6}+\frac{x_3}{2}\big)dx_2-\frac{x_2^2}{6}dx_1$;
\item $\theta_6=dx_6-\frac{x_1}{2}dx_5-\frac{x_2}{2}dx_4+\frac{x_1x_2}{3}dx_3+\big(\frac{x_4}{2}-\frac{x_1x_3}{6}-\frac{x_1^2x_2}{12}\big)dx_2+\big(\frac{x_5}{2}-\frac{x_2x_3}{6}+\frac{x_1x_2^2}{12}\big)dx_1$;
\item $\theta_7=dx_7-\frac{x_1}{2}dx_6+\frac{x_1^2}{6}dx_5+\big(\frac{x_1x_2}{6}-\frac{x_3}{2}\big)dx_4+\big(\frac{x_4}{2}+\frac{x_1x_3}{6}-\frac{x_1^2x_2}{12}\big)dx_3+\big(\frac{x_1^3x_2}{60}\saltot-\frac{x_1x_4}{3}\big)dx_2+\big(\frac{x_6}{2}+\frac{x_2x_4-x_3^2-x_1x_5}{6}+\frac{x_1x_2x_3}{12}-\frac{x_1^2x_2^2}{60}\big)dx_1$.
\end{itemize}

Finally, we have
\begin{eqnarray*}
  \mathrm{d}(R_\mathbf{x})_\mathbf{0}=\left[\begin{matrix}    
      1 & 0  & 0& 0& 0& 0& 0\\
   0 & 1& 0& 0& 0& 0& 0\\
   \frac{x_2}{2} & -\frac{x_1}{2}& 1& 0& 0& 0& 0\\
   \frac{x_3}{2}-\frac{x_1x_2}{12} &\frac{x_1^2}{12} & -\frac{x_1}{2}& 1& 0& 0& 0\\
   -\frac{x_2^2}{12} &  \frac{x_1x_2}{12}+\frac{x_3}{2}& -\frac{x_2}{2}& 0& 1& 0& 0\\
 \frac{x_5}{2}-\frac{x_2x_3}{12} &  \frac{x_4}{2}-\frac{x_1x_3}{12}& \frac{x_1x_2}{6}& -\frac{x_2}{2}& -\frac{x_1}{2}& 1& 0\\
 \frac{x_1^2x_2^2}{360}+\frac{x_2x_4-x_3^2-x_1x_5}{12}+\frac{x_6}{2}& -\frac{x_1x_4}{6}-\frac{x_1^3x_2}{360} &\frac{x_1x_3}{12}+\frac{x_4}{2} & \frac{x_1x_2}{12}-\frac{x_3}{2}& \frac{x_1^2}{12}& -\frac{x_1}{2}& 1
   \end{matrix}\right]\,.
   \end{eqnarray*}

 \subsection*{$(12457L)$}
 
 The following Lie algebra is denoted as $(12457L)$ by Gong in \cite{Gong_Thesis}, and as $\mathcal{G}_{7,1,17}$ by Magnin in \cite{magnin}.
 
 The non-trivial brackets are the following:
 \begin{equation*}
 \begin{aligned}
     &[X_1,X_i]=X_{i+1}\,,\,i=2,3,5,6\,,\,[X_2,X_j]=X_{j+2}\,,\,j=3,4\,,\,\\
     &\qquad [X_2,X_6]=X_7\,,\,[X_3,X_4]=X_7\,,\,[X_3,X_5]=-X_7\,.
     \end{aligned}
 \end{equation*}
This is a nilpotent Lie algebra of rank 2 and step 5 that is stratifiable. The Lie brackets can be pictured with the diagram:
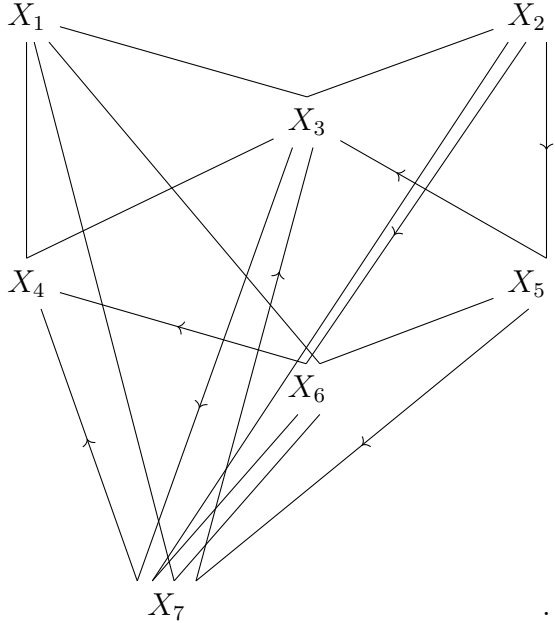
\begin{center}
 
	\begin{tikzcd}[end anchor=north]
		 X_1\ar[drr, no head]\ar[ddd, no head]\ar[ddddrr, no head, end anchor={[xshift=-1.5ex]north east}]\ar[dddddddr, no head, end anchor={[xshift=-1.9ex]north east}] & &  & & X_2\ar[dddddddlll, no head,end anchor={[xshift=-3.5ex]north east},start anchor={[xshift=-3.8ex]south east}]\ar[dll, no head]\ar[ddd, no head, ->-=.5,end anchor={[xshift=-1.ex]north east},start anchor={[xshift=-1.ex]south east}]\ar[ddddll, no head, ->-=.6,end anchor={[xshift=-2.5ex]north east},start anchor={[xshift=-2.5ex]south east}]\\
		 & & X_3\ar[ddll, no head]\ar[ddrr, no head,-<-=.3,end anchor={[xshift=-1.ex]north east}] \ar[ddddddl, no head, ->-=.6, end anchor={[xshift=-4.6ex]north east},start anchor={[xshift=-3.5ex]south east}]\ar[ddddddl, no head, start anchor={[xshift=-2.ex]south east},end anchor={[xshift=-.3ex]north east},-<-=.3
		 ]& &\\
		 & & & & \\
		 X_4\ar[ddddr, no head, -<-=.5, end anchor={[xshift=-4.6ex]north east}]\ar[drr, no head, -<-=.5, end anchor={[xshift=-2.5ex]north east}] & & & & X_5\ar[dll, no head, end anchor={[xshift=-1.5ex]north east}]\ar[ddddlll, no head, end anchor={[xshift=-.3ex]north east},->-=.5,start anchor={[xshift=-2.3ex]south east}]\\
		 & &X_6 \ar[dddl, no head, end anchor={[xshift=-3.5ex]north east},start anchor={[xshift=-3.1ex]south east}]\ar[dddl, no head, end anchor={[xshift=-1.9ex]north east},start anchor={[xshift=-1.5ex]south east}]&  &\\
		 & & & &\\
		 & & & &\\
		 & X_7 & & &\quad\;.
		 \end{tikzcd}
	\end{center}
	
The composition law \eqref{group law in G} of $(12457L)$ is given by:
\begin{itemize}
    \item $z_1=x_1+y_1$;
    \item $z_2=x_2+y_2$;
    \item $z_3=x_3+y_3+\frac{1}{2}(x_1y_2-x_2y_1)$;
    \item $z_4=x_4+y_4+\frac{1}{2}(x_1y_3-x_3y_1)+\frac{1}{12}(x_1-y_1)(x_1y_2-x_2y_1)$;
    \item $z_5=x_5+y_5+\frac{1}{2}(x_2y_3-x_3y_2)+\frac{1}{12}(x_2-y_2)(x_1y_2-x_2y_1)$;
    \item $z_6=x_6+y_6+\frac{1}{2}(x_1y_5-x_5y_1+x_2y_4-x_4y_2)+\frac{1}{12}(x_1-y_1)(x_2y_3-x_3y_2)\salto+\frac{1}{12}(x_2-y_2)(x_1y_3-x_3y_1)-\frac{1}{24}(x_1y_2+x_2y_1)(x_1y_2-x_2y_1)$;
    \item $z_7=x_7+y_7+\frac{1}{2}(x_1y_6-x_6y_1+x_2y_6-x_6y_2+x_3y_4-x_4y_3-x_3y_5+x_5y_3)\salto+\frac{1}{12}(x_1-y_1)(x_1y_5-x_5y_1+x_2y_4-x_4y_2)+\frac{1}{12}(x_2-y_2)(x_1y_5-x_5y_1+x_2y_4-x_4y_2)\salto+\frac{1}{12}(x_3-y_3)(x_1y_3-x_3y_1-x_2y_3+x_3y_2)+\frac{1}{12}(x_5-y_5-x_4+y_4)(x_1y_2-x_2y_1)\salto+\frac{1}{24}(x_2y_3-x_1y_3)(x_1y_2-x_2y_1)-\frac{1}{24}(x_2y_1+x_2y_2)(x_1y_3-x_3y_1)\salto-\frac{1}{24}(x_1y_1+x_1y_2)(x_2y_3-x_3y_1)+\frac{1}{180}(x_2y_1^2+x_1y_2^2+x_1y_1y_2)(x_1y_2-x_2y_1)\salto+\frac{1}{180}(x_2y_1y_2-2x_1x_2y_1-2x_1x_2y_2)(x_1y_2-x_2y_1)+\frac{1}{360}(y_1-y_2)(x_1y_2-x_2y_1)^2\salto+\frac{1}{360}(y_1^2y_2+y_1y_2^2-x_1^2x_2-x_1x_2^2)(x_1y_2-x_2y_1)+\frac{1}{120}(x_1-y_1)(x_1y_2-x_2y_1)^2$.
\end{itemize}

Since
\begin{eqnarray*}
  \mathrm{d}(L_\mathbf{x})_\mathbf{0}=\left[\begin{matrix}    
      1 & 0  & 0& 0& 0& 0& 0\\
   0 & 1& 0& 0& 0& 0& 0\\
   -\frac{x_2}{2} & \frac{x_1}{2}& 1& 0& 0& 0& 0\\
   -\frac{x_1x_2}{12}-\frac{x_3}{2} &\frac{x_1^2}{12} & \frac{x_1}{2}& 1& 0& 0& 0\\
   -\frac{x_2^2}{12} &  \frac{x_1x_2}{12}-\frac{x_3}{2}& \frac{x_2}{2}& 0& 1& 0& 0\\
 -\frac{x_2x_3}{12}-\frac{x_5}{2} &  -\frac{x_1x_3}{12}-\frac{x_4}{2}& \frac{x_1x_2}{6}& \frac{x_2}{2}& \frac{x_1}{2}& 1& 0\\
      a_1 & a_2  & a_3& \frac{x_1x_2+x_2^2}{12}+\frac{x_3}{2}& \frac{x_1^2+x_1x_2}{12}-\frac{x_3}{2}& \frac{x_1+x_2}{2}& 1
   \end{matrix}\right],
   \end{eqnarray*}
   where
         \begin{eqnarray*}
         a_1 &=&\frac{x_1^2x_2^2+x_1x_2^3}{360}+\frac{x_2x_4-x_3^2-x_1x_5-2x_2x_5}{12}-\frac{x_6}{2}
        ,\\
         a_2 &=&\frac{x_1x_5-x_2x_4-2x_1x_4+x_3^2}{12}-\frac{x_1^3x_2+x_1^2x_2^2}{360}-\frac{x_6}{2}   ,\\ 
         a_3&=&\frac{x_1x_3-x_2x_3}{12}+\frac{x_5-x_4}{2} ,
        \end{eqnarray*}
the induced left-invariant vector fields \eqref{leftinvariant vf} are: 
\begin{itemize}
\item $X_1=\partial_{x_1} -\frac{x_2}{2}\partial_{x_3}-\big(\frac{x_3}{2}+\frac{x_1x_2}{12}\big)\partial_{x_4}-\frac{x_2^2}{12}\partial_{x_5}-\big(\frac{x_2x_3}{12}+\frac{x_5}{2}\big)\partial_{x_6}+\big(\frac{x_1^2x_2^2+x_1x_2^3}{360}\saltox+\frac{x_2x_4-x_3^2-x_1x_5-2x_2x_5}{12}-\frac{x_6}{2}\big)\partial_{x_7}\,;$
\item $X_2=\partial_{x_2}+\frac{x_1}{2}\partial_{x_3}+\frac{x_1^2}{12}\partial_{x_4}+\big(\frac{x_1x_2}{12}-\frac{x_3}{2}\big)\partial_{x_5}-\big(\frac{x_4}{2}+\frac{x_1x_3}{12}\big)\partial_{x_6}+\big(\frac{x_1x_5-x_2x_4-2x_1x_4+x_3^2}{12}\saltox-\frac{x_1^3x_2+x_1^2x_2^2}{360}-\frac{x_6}{2}\big)\partial_{x_7}\,;$
\item $X_3=\partial_{x_3}+\frac{x_1}{2}\partial_{x_4}+\frac{x_2}{2}\partial_{x_5}+\frac{x_1x_2}{6}\partial_{x_6}+\big(\frac{x_1x_3-x_2x_3}{12}+\frac{x_5-x_4}{2}\big)\partial_{x_7}\,;$
\item $X_4=\partial_{x_4}+\frac{x_2}{2}\partial_{x_6}+\big(\frac{x_1x_2+x_2^2}{12}+\frac{x_3}{2}\big)\partial_{x_7}\,;$
\item $X_5=\partial_{x_5}+\frac{x_1}{2}\partial_{x_6}+\big(\frac{x_1^2+x_1x_2}{12}-\frac{x_3}{2}\big)\partial_{x_7}$;
\item $X_6=\partial_{x_6}+\frac{x_1+x_2}{2}\partial_{x_7}$;
\item $X_7=\partial_{x_7}$,
\end{itemize}
and the respective left-invariant 1-forms \eqref{leftinvariant form} are: 
\begin{itemize}
\item $\theta_1=dx_1$;
\item $\theta_2=dx_2$;
\item $\theta_3=dx_3-\frac{x_1}{2}dx_2+\frac{x_2}{2}dx_1$;
\item $\theta_4=dx_4-\frac{x_1}{2}dx_3+\frac{x_1^2}{6}dx_2+\big(\frac{x_3}{2}-\frac{x_1x_2}{6}\big)dx_1$;
\item $\theta_5=dx_5-\frac{x_2}{2}dx_3+\big(\frac{x_1x_2}{6}+\frac{x_3}{2}\big)dx_2-\frac{x_2^2}{6}dx_1$;
\item $\theta_6=dx_6-\frac{x_1}{2}dx_5-\frac{x_2}{2}dx_4+\frac{x_1x_2}{3}dx_3+\big(\frac{x_4}{2}-\frac{x_1x_3}{6}-\frac{x_1^2x_2}{12}\big)dx_2+\big(\frac{x_5}{2}-\frac{x_2x_3}{6}+\frac{x_1x_2^2}{12}\big)dx_1$;
\item $\theta_7=dx_7-\frac{x_1+x_2}{2}dx_6+\big(\frac{x_3}{2}+\frac{x_1^2+x_1x_2}{6}\big)dx_5+\big(\frac{x_1x_2+x_2^2}{6}-\frac{x_3}{2}\big)dx_4+\big(\frac{x_4-x_5}{2}\saltot+\frac{x_1x_3-x_2x_3}{6}-\frac{x_1^2x_2+x_1x_2^2}{12}\big)dx_3+\big(\frac{x_1^2x_2^2+x_1^3x_2}{60}+\frac{x_1x_2x_3}{12}+\frac{x_3^2-2x_1x_4-x_2x_4+x_1x_5}{6}\saltot+\frac{x_6}{2}\big)dx_2+\big(\frac{x_6}{2}+\frac{x_2x_4-x_3^2-x_1x_5-2x_2x_5}{6}+\frac{x_1x_2x_3}{12}-\frac{x_1^2x_2^2+x_1x_2^3}{60}\big)dx_1$.
\end{itemize}

Finally, we have
\begin{eqnarray*}
  \mathrm{d}(R_\mathbf{x})_\mathbf{0}=\left[\begin{matrix}    
      1 & 0  & 0& 0& 0& 0& 0\\
   0 & 1& 0& 0& 0& 0& 0\\
   \frac{x_2}{2} & -\frac{x_1}{2}& 1& 0& 0& 0& 0\\
   \frac{x_3}{2}-\frac{x_1x_2}{12} &\frac{x_1^2}{12} & -\frac{x_1}{2}& 1& 0& 0& 0\\
   -\frac{x_2^2}{12} &  \frac{x_1x_2}{12}+\frac{x_3}{2}& -\frac{x_2}{2}& 0& 1& 0& 0\\
 \frac{x_5}{2}-\frac{x_2x_3}{12} &  \frac{x_4}{2}-\frac{x_1x_3}{12}& \frac{x_1x_2}{6}& -\frac{x_2}{2}& -\frac{x_1}{2}& 1& 0\\
 a_1& a_2 &a_3 & \frac{x_1x_2+x_2^2}{12}-\frac{x_3}{2}& \frac{x_1^2+x_1x_2}{12}+\frac{x_3}{2}& -\frac{x_1+x_2}{2}& 1
   \end{matrix}\right]\,,
   \end{eqnarray*}
   where
    \begin{eqnarray*}
         a_1 &=&\frac{x_1^2x_2^2+x_1x_2^3}{360}+\frac{x_2x_4-x_3^2-x_1x_5-2x_2x_5}{12}+\frac{x_6}{2}
        ,\\
         a_2 &=&\frac{x_1x_5-x_2x_4-2x_1x_4+x_3^2}{12}-\frac{x_1^3x_2+x_1^2x_2^2}{360}+\frac{x_6}{2} ,\\ 
         a_3&=&\frac{x_1x_3-x_2x_3}{12}+\frac{x_4-x_5}{2} .
        \end{eqnarray*}

 \subsection*{$(12457L_1)$}
 
 The following Lie algebra is denoted as $(12457L_1)$ by Gong in \cite{Gong_Thesis}, and as $\mathcal{G}_{7,1,17}$ by Magnin in \cite{magnin}.
 
 The non-trivial brackets are the following:
 \begin{equation*}
 \begin{aligned}
     &[X_1,X_i]=X_{i+1}\,,\,i=2,3\,,\,[X_1,X_4]=-X_{6}\,,[X_1,X_6]=X_7\,,\,\\
     &\qquad [X_2,X_3]=X_5\,,\,[X_2,X_5]=-X_6\,,\,[X_3,X_5]=-X_7\,.
     \end{aligned}
 \end{equation*}
This is a nilpotent Lie algebra of rank 2 and step 5 that is stratifiable. The Lie brackets can be pictured with the diagram:
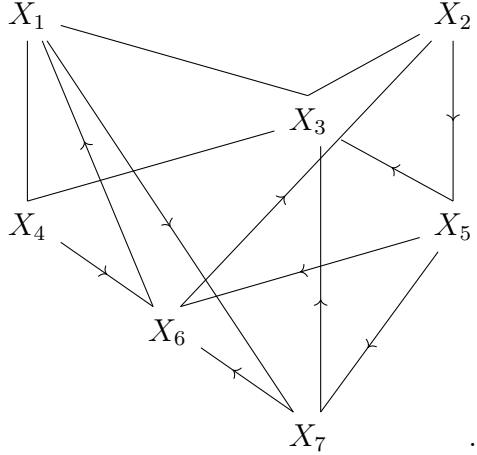
\begin{center}
 
	\begin{tikzcd}[end anchor=north]
		X_1 \ar[drr, no head]\ar[dd, no head]\ar[dddr, no head, -<-=.4,end anchor={[xshift=-3.5ex]north east}]\ar[ddddrr, no head, ->-=.5,end anchor={[xshift=-3.5ex]north east}]& & & X_2\ar[dddll, no head,end anchor={[xshift=-1.5ex]north east},-<-=.6] \ar[dl, no head]\ar[dd, no head, ->-=.5]\\
		& &X_3\ar[dll, no head]\ar[dr, no head ,-<-=.5]\ar[ddd, no head, -<-=.6,end anchor={[xshift=-1.5ex]north east},start anchor={[xshift=-1.5ex]south east}] &  \\
		X_4\ar[dr, no head, ->-=.5,end anchor={[xshift=-3.5ex]north east}]& & & X_5\ar[dll, no head,end anchor={[xshift=-1.5ex]north east}, ->-=.5]\ar[ddl, no head, ->-=.6,end anchor={[xshift=-1.5ex]north east}] \\
		& X_6\ar[dr, no head, -<-=.4,end anchor={[xshift=-3.5ex]north east}]& & \\
		& & X_7 &\quad\;.
		 \end{tikzcd}
	\end{center}
	
The composition law \eqref{group law in G} of $(12457L_1)$ is given by:
\begin{itemize}
    \item $z_1=x_1+y_1$;
    \item $z_2=x_2+y_2$;
    \item $z_3=x_3+y_3+\frac{1}{2}(x_1y_2-x_2y_1)$;
    \item $z_4=x_4+y_4+\frac{1}{2}(x_1y_3-x_3y_1)+\frac{1}{12}(x_1-y_1)(x_1y_2-x_2y_1)$;
    \item $z_5=x_5+y_5+\frac{1}{2}(x_2y_3-x_3y_2)+\frac{1}{12}(x_2-y_2)(x_1y_2-x_2y_1)$;
    \item $z_6=x_6+y_6+\frac{1}{2}(x_5y_2-x_2y_5+x_4y_1-x_1y_4)+\frac{1}{12}(y_1-x_1)(x_1y_3-x_3y_1)\salto+\frac{1}{12}(y_2-x_2)(x_2y_3-x_3y_2)]+\frac{1}{24}(x_1y_1+x_2y_2)(x_1y_2-x_2y_1)$;
    \item $z_7=x_7+y_7+\frac{1}{2}(x_1y_6-x_6y_1-x_3y_5+x_5y_3)+\frac{1}{12}(x_5-y_5)(x_1y_2-x_2y_1)\salto+\frac{1}{12}(y_1-x_1)(x_2y_5-x_5y_2+x_1y_4-x_4y_1)+\frac{1}{12}(y_3-x_3)(x_2y_3-x_3y_2)]\salto+\frac{1}{24}[x_2y_3(x_1y_2-x_2y_1)+x_1y_1(x_1y_3-x_3y_1)+x_2y_1(x_2y_3-x_3y_2)]\salto+\frac{1}{180}(x_1^2y_1+x_2^2y_1-x_1y_1^2-x_2y_1y_2)(x_1y_2-x_2y_1)-\frac{1}{360}y_2(x_1y_2-x_2y_1)^2\salto+\frac{1}{720}(x_1^3+x_1x_2^2-y_1^3-y_1y_2^2)(x_1y_2-x_2y_1)-\frac{1}{120}x_2(x_1y_2-x_2y_1)^2$.
\end{itemize}

Since
\begin{eqnarray*}
  \mathrm{d}(L_\mathbf{x})_\mathbf{0}=\left[\begin{matrix}    
      1 & 0  & 0& 0& 0& 0& 0\\
   0 & 1& 0& 0& 0& 0& 0\\
   -\frac{x_2}{2} & \frac{x_1}{2}& 1& 0& 0& 0& 0\\
   -\frac{x_1x_2}{12}-\frac{x_3}{2} &\frac{x_1^2}{12} & \frac{x_1}{2}& 1& 0& 0& 0\\
   -\frac{x_2^2}{12} &  \frac{x_1x_2}{12}-\frac{x_3}{2}& \frac{x_2}{2}& 0& 1& 0& 0\\
 \frac{x_1x_3}{12}+\frac{x_4}{2} &  \frac{x_2x_3}{12}+\frac{x_5}{2}& -\frac{x_1^2+x_2^2}{12}& -\frac{x_1}{2}& -\frac{x_2}{2}& 1& 0\\
 a_1& a_2 &\frac{x_5}{2}-\frac{x_2x_3}{12} & -\frac{x_1^2}{12}& -\frac{x_1x_2}{12}-\frac{x_3}{2}& \frac{x_1}{2}& 1
   \end{matrix}\right]\,,
   \end{eqnarray*}
   where
    \begin{eqnarray*}
         a_1 &=&\frac{x_1x_4-x_2x_5}{12}-\frac{x_6}{2}-\frac{x_1^3x_2+x_1x_2^3}{720}
        ,\\
         a_2 &=&\frac{x_3^2+2x_1x_5}{12}+\frac{x_1^4+x_1^2x_2^2}{720}\,,
        \end{eqnarray*}
the induced left-invariant vector fields \eqref{leftinvariant vf} are: 
\begin{itemize}
\item $X_1=\partial_{x_1} -\frac{x_2}{2}\partial_{x_3}-\big(\frac{x_3}{2}+\frac{x_1x_2}{12}\big)\partial_{x_4}-\frac{x_2^2}{12}\partial_{x_5}+\big(\frac{x_1x_3}{12}+\frac{x_4}{2}\big)\partial_{x_6}+\big(\frac{x_1x_4-x_2x_5}{12}-\frac{x_6}{2}\saltox-\frac{x_1^3x_2+x_1x_2^3}{720}\big)\partial_{x_7}\,;$
\item $X_2=\partial_{x_2}+\frac{x_1}{2}\partial_{x_3}+\frac{x_1^2}{12}\partial_{x_4}+\big(\frac{x_1x_2}{12}-\frac{x_3}{2}\big)\partial_{x_5}+\big(\frac{x_5}{2}+\frac{x_2x_3}{12}\big)\partial_{x_6}+\big(\frac{x_3^2+2x_1x_5}{12}\saltox+\frac{x_1^4+x_1^2x_2^2}{720}\big)\partial_{x_7}\,;$
\item $X_3=\partial_{x_3}+\frac{x_1}{2}\partial_{x_4}+\frac{x_2}{2}\partial_{x_5}-\frac{x_1^2+x_2^2}{12}\partial_{x_6}+\big(\frac{x_5}{2}-\frac{x_2x_3}{12}\big)\partial_{x_7}\,;$
\item $X_4=\partial_{x_4}-\frac{x_1}{2}\partial_{x_6}-\frac{x_1^2}{12}\partial_{x_7}\,;$
\item $X_5=\partial_{x_5}-\frac{x_2}{2}\partial_{x_6}-\big(\frac{x_1x_2}{12}+\frac{x_3}{2}\big)\partial_{x_7}$;
\item $X_6=\partial_{x_6}+\frac{x_1}{2}\partial_{x_7}$;
\item $X_7=\partial_{x_7}$,
\end{itemize}
and the respective left-invariant 1-forms \eqref{leftinvariant form} are: 
\begin{itemize}
\item $\theta_1=dx_1$;
\item $\theta_2=dx_2$;
\item $\theta_3=dx_3-\frac{x_1}{2}dx_2+\frac{x_2}{2}dx_1$;
\item $\theta_4=dx_4-\frac{x_1}{2}dx_3+\frac{x_1^2}{6}dx_2+\big(\frac{x_3}{2}-\frac{x_1x_2}{6}\big)dx_1$;
\item $\theta_5=dx_5-\frac{x_2}{2}dx_3+\big(\frac{x_1x_2}{6}+\frac{x_3}{2}\big)dx_2-\frac{x_2^2}{6}dx_1$;
\item $\theta_6=dx_6+\frac{x_2}{2}dx_5+\frac{x_1}{2}dx_4-\frac{x_1^2+x_2^2}{6}dx_3+\big(\frac{x_1^3+x_1x_2^2}{24}+\frac{x_2x_3}{6}-\frac{x_5}{2}\big)dx_2+\big(\frac{x_1x_3}{6}\saltot-\frac{x_1^2x_2+x_2^3}{24}-\frac{x_4}{2}\big)dx_1$;
\item $\theta_7=dx_7-\frac{x_1}{2}dx_6+\big(\frac{x_3}{2}-\frac{x_1x_2}{6}\big)dx_5-\frac{x_1^2}{6}dx_4+\big(\frac{x_1^3+x_1x_2^2}{24}-\frac{x_5}{2}-\frac{x_2x_3}{6}\big)dx_3+\big(\frac{x_3^2+2x_1x_5}{6}\saltot-\frac{x_1^4+x_1^2x_2^2}{120}\big)dx_2+\big(\frac{x_6}{2}+\frac{x_1x_4-x_2x_5}{6}-\frac{x_1^2x_3+x_2^2x_3}{24}+\frac{x_1^3x_2+x_1x_2^3}{120}\big)dx_1$.
\end{itemize}

Finally, we have
\begin{eqnarray*}
  \mathrm{d}(R_\mathbf{x})_\mathbf{0}=\left[\begin{matrix}    
            1 & 0  & 0& 0& 0& 0& 0\\
   0 & 1& 0& 0& 0& 0& 0\\
   \frac{x_2}{2} & -\frac{x_1}{2}& 1& 0& 0& 0& 0\\
   \frac{x_3}{2}-\frac{x_1x_2}{12} &\frac{x_1^2}{12} & -\frac{x_1}{2}& 1& 0& 0& 0\\
   -\frac{x_2^2}{12} &  \frac{x_1x_2}{12}+\frac{x_3}{2}& -\frac{x_2}{2}& 0& 1& 0& 0\\
 \frac{x_1x_3}{12}-\frac{x_4}{2} &  \frac{x_2x_3}{12}-\frac{x_5}{2}& -\frac{x_1^2+x_2^2}{12}& \frac{x_1}{2}& \frac{x_2}{2}& 1& 0\\
 a_1& a_2 &-\frac{x_5}{2}-\frac{x_2x_3}{12} & -\frac{x_1^2}{12}& \frac{x_3}{2}-\frac{x_1x_2}{12}& -\frac{x_1}{2}& 1
   \end{matrix}\right]\,,
      \end{eqnarray*}
    where
    \begin{eqnarray*}
         a_1 &=&\frac{x_1x_4-x_2x_5}{12}+\frac{x_6}{2}-\frac{x_1^3x_2+x_1x_2^3}{720}
        ,\\
         a_2 &=&\frac{x_3^2+2x_1x_5}{12}+\frac{x_1^4+x_1^2x_2^2}{720}\,.
        \end{eqnarray*}

 \subsection*{$(123457A)$}
 
 The following Lie algebra is denoted as $(123457A)$ by Gong in \cite{Gong_Thesis}, and as $\mathcal{G}_{7,2,3}$ by Magnin in \cite{magnin}.
 
 The non-trivial brackets are the following:
 \begin{equation*}
     [X_1,X_i]=X_{i+1}\,,\,2\le i\le 6\,.
 \end{equation*}
This is a nilpotent Lie algebra of rank 2 and step 6 that is stratifiable, also known as the filiform Lie algebra of dimension 7. The Lie brackets can be pictured with the diagram:
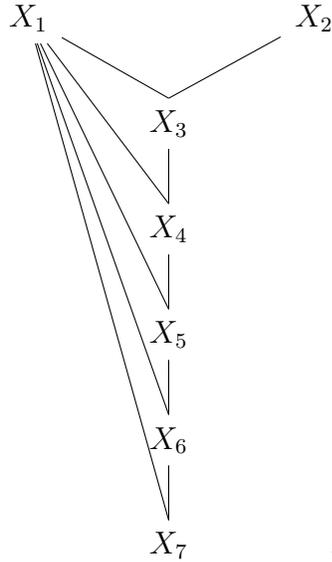
\begin{center}
 
	\begin{tikzcd}[end anchor=north]
		 X_1\ar[dr, no head]\ar[ddr, no head]\ar[dddr, no head]\ar[ddddr, no head]\ar[dddddr, no head] & & X_2\ar[dl, no head]\\
		 & X_3\ar[d, no head] & \\
		 & X_4\ar[d, no head] &\\
		 & X_5\ar[d, no head]&\\
		 & X_6\ar[d, no head] & \\
		 & X_7 &\quad\;.
		 \end{tikzcd}
	\end{center}
	
The composition law \eqref{group law in G} of $(123457A)$ is given by:
\begin{itemize}
    \item $z_1=x_1+y_1$;
    \item $z_2=x_2+y_2$;
    \item $z_3=x_3+y_3+\frac{1}{2}(x_1y_2-x_2y_1)$;
    \item $z_4=x_4+y_4+\frac{1}{2}(x_1y_3-x_3y_1)+\frac{1}{12}(x_1-y_1)(x_1y_2-x_2y_1)$;
    \item $z_5=x_5+y_5+\frac{1}{2}(x_1y_4-x_4y_1)+\frac{1}{12}(x_1-y_1)(x_1y_3-x_3y_1)-\frac{1}{24}x_1y_1(x_1y_2-x_2y_1)$;
    \item $z_6=x_6+y_6+\frac{1}{2}(x_1y_5-x_5y_1)+\frac{1}{12}(x_1-y_1)(x_1y_4-x_4y_1)-\frac{1}{24}x_1y_1(x_1y_3-x_3y_1)\salto+\frac{1}{720}(y_1^3-x_1^3)(x_1y_2-x_2y_1)+\frac{1}{180}(x_1y_1^2-x_1^2y_1)(x_1y_2-x_2y_1)$;
    \item $z_7=x_7+y_7+\frac{1}{2}(x_1y_6-x_6y_1)+\frac{1}{12}(x_1-y_1)(x_1y_5-x_5y_1)-\frac{1}{24}x_1y_1(x_1y_4-x_4y_1)\salto+\frac{1}{180}(x_1y_1^2-x_1^2y_1)(x_1y_3-x_3y_1)+\frac{1}{720}(y_1^3-x_1^3)(x_1y_3-x_3y_1)\salto+\frac{1}{360}x_1^2y_1^2(x_1y_2-x_2y_1)+\frac{1}{1440}(x_1^3y_1+x_1y_1^3)(x_1y_2-x_2y_1)$.
\end{itemize}

Since
\begin{eqnarray*}
  \mathrm{d}(L_\mathbf{x})_\mathbf{0}=\left[\begin{matrix}    
      1 & 0  & 0& 0& 0& 0& 0\\
   0 & 1& 0& 0& 0& 0& 0\\
   -\frac{x_2}{2} & \frac{x_1}{2}& 1& 0& 0& 0& 0\\
   -\frac{x_1x_2}{12}-\frac{x_3}{2} &\frac{x_1^2}{12} & \frac{x_1}{2}& 1& 0& 0& 0\\
   -\frac{x_1x_3}{12}-\frac{x_4}{2} &  0& \frac{x_1^2}{12}& \frac{x_1}{2}& 1& 0& 0\\
 \frac{x_1^3x_2}{720}-\frac{x_1x_4}{12}-\frac{x_5}{2} &  -\frac{x_1^4}{720}& 0& \frac{x_1^2}{12}& \frac{x_1}{2}& 1& 0\\
 \frac{x_1^3x_3}{720}-\frac{x_1x_5}{12}-\frac{x_6}{2}& 0 &-\frac{x_1^4}{720} &0& \frac{x_1^2}{12}& \frac{x_1}{2}& 1
   \end{matrix}\right]\,,
   \end{eqnarray*}
the induced left-invariant vector fields \eqref{leftinvariant vf} are: 
\begin{itemize}
\item $X_1=\partial_{x_1} -\frac{x_2}{2}\partial_{x_3}-\big(\frac{x_3}{2}+\frac{x_1x_2}{12}\big)\partial_{x_4}-\big(\frac{x_4}{2}+\frac{x_1x_3}{12}\big)\partial_{x_5}+\big(\frac{x_1^3x_2}{720}-\frac{x_1x_4}{12}-\frac{x_5}{2}\big)\partial_{x_6}\saltox+\big(\frac{x_1^3x_3}{720}-\frac{x_1x_5}{12}-\frac{x_6}{2}\big)\partial_{x_7}\,;$
\item $X_2=\partial_{x_2}+\frac{x_1}{2}\partial_{x_3}+\frac{x_1^2}{12}\partial_{x_4}-\frac{x_1^4}{720}\partial_{x_6}\,;$
\item $X_3=\partial_{x_3}+\frac{x_1}{2}\partial_{x_4}+\frac{x_1^2}{12}\partial_{x_5}-\frac{x_1^4}{720}\partial_{x_7}\,;$
\item $X_4=\partial_{x_4}+\frac{x_1}{2}\partial_{x_5}+\frac{x_1^2}{12}\partial_{x_6}\,;$
\item $X_5=\partial_{x_5}+\frac{x_1}{2}\partial_{x_6}+\frac{x_1^2}{12}\partial_{x_7}$;
\item $X_6=\partial_{x_6}+\frac{x_1}{2}\partial_{x_7}$;
\item $X_7=\partial_{x_7}$,
\end{itemize}
and the respective left-invariant 1-forms \eqref{leftinvariant form} are: 
\begin{itemize}
\item $\theta_1=dx_1$;
\item $\theta_2=dx_2$;
\item $\theta_3=dx_3-\frac{x_1}{2}dx_2+\frac{x_2}{2}dx_1$;
\item $\theta_4=dx_4-\frac{x_1}{2}dx_3+\frac{x_1^2}{6}dx_2+\big(\frac{x_3}{2}-\frac{x_1x_2}{6}\big)dx_1$;
\item $\theta_5=dx_5-\frac{x_1}{2}dx_4+\frac{x_1^2}{6}dx_3-\frac{x_1^3}{24}dx_2+\big(\frac{x_1^2x_2}{24}-\frac{x_1x_3}{6}+\frac{x_4}{2}\big)dx_1$;
\item $\theta_6=dx_6-\frac{x_1}{2}dx_5+\frac{x_1^2}{6}dx_4-\frac{x_1^3}{24}dx_3+\frac{x_1^4}{120}dx_2+\big(\frac{x_5}{2}-\frac{x_1x_4}{6}+\frac{x_1^2x_3}{24}-\frac{x_1^3x_2}{120}\big)dx_1$;
\item $\theta_7=dx_7-\frac{x_1}{2}dx_6+\frac{x_1^2}{6}dx_5-\frac{x_1^3}{24}dx_4+\frac{x_1^4}{120}dx_3-\frac{x_1^5}{720}dx_2+\big(\frac{x_6}{2}-\frac{x_1x_5}{6}+\frac{x_1^2x_4}{24}-\frac{x_1^3x_3}{120}\saltot+\frac{x_1^4x_2}{720}\big)dx_1$.
\end{itemize}

Finally, we have
\begin{eqnarray*}
  \mathrm{d}(R_\mathbf{x})_\mathbf{0}=\left[\begin{matrix}    
      1 & 0  & 0& 0& 0& 0& 0\\
   0 & 1& 0& 0& 0& 0& 0\\
   \frac{x_2}{2} & -\frac{x_1}{2}& 1& 0& 0& 0& 0\\
   \frac{x_3}{2}-\frac{x_1x_2}{12} &\frac{x_1^2}{12} & -\frac{x_1}{2}& 1& 0& 0& 0\\
   \frac{x_4}{2}-\frac{x_1x_3}{12} &  0& \frac{x_1^2}{12}& -\frac{x_1}{2}& 1& 0& 0\\
 \frac{x_1^3x_2}{720}-\frac{x_1x_4}{12}+\frac{x_5}{2} &  -\frac{x_1^4}{720}& 0& \frac{x_1^2}{12}& -\frac{x_1}{2}& 1& 0\\
 \frac{x_1^3x_3}{720}-\frac{x_1x_5}{12}+\frac{x_6}{2}& 0 &-\frac{x_1^4}{720} &0& \frac{x_1^2}{12}& -\frac{x_1}{2}& 1
   \end{matrix}\right]\,.
   \end{eqnarray*}

  \newpage
  \section{Some free-nilpotent groups in low dimension}
  In this final section we analyze free-nilpotent groups of low dimension.
  We shall denote with $\mathbb F_{rs}$ the simply connected Lie group whose Lie algebra has rank $r$ and is free up to nilpotency step $s$. In the specific, we shall study 
  $\mathbb F_{23}$,
  $\mathbb F_{24}$,
  $\mathbb F_{33}$,
  $\mathbb F_{25}$, respectively.

 \subsection*{$\mathbb F_{23}$.}

The following is the free-nilpotent Lie algebra with 2 generators and nilpotency step 3. It has dimension 5. This Lie algebra is also denoted as $N_{5,2,3}$, see page \pageref{N523}.

The non-trivial brackets coming from the Hall basis are:
\begin{equation*}\label{F23}
     [X_1,X_2]=X_{3}\,,\,[X_1,X_3]=X_4\,,\,[X_2,X_3]=X_5\,.
 \end{equation*}
  
  A representation of the left-invariant vector fields $X_i^s$ defined on $\mathbb{R}^5$ with respect to exponential coordinates of the second kind is the following:
  \begin{itemize}
      \item $X_1^s=\partial_{x_1}\,;$
      \item $X_2^s=\partial_{x_2}+x_1\partial_{x_3}+\frac{x_1^2}{2}\partial_{x_4}+x_1x_2\partial_{x_5}$;
      \item $X_3^s=\partial_{x_3}+x_1\partial_{x_4}+x_2\partial_{x_5}$;
      \item $X_4^s=\partial_{x_4}$;
      \item $X_5^s=\partial_{x_5}$.
  \end{itemize}





One can relate the  exponential coordinates of first kind to the exponential coordinates of second kind.
If we denote by $\alpha_1,\ldots, \alpha_5$ the coordinates of second kind and by
$a_1,\ldots, a_5$ the coordinates of first kind, the change of coordinates are given as follows:

\begin{equation*}
\left\{\begin{array}{ccl}
a_1 &=&\alpha_1 \\
a_2 &=&\alpha_2\\
a_3 &=&\alpha_3-\dfrac{\alpha_1\alpha_2}{2}\\
a_4 &=&\alpha_4-\dfrac{\alpha_1\alpha_3}{2}+\dfrac{\alpha_1^2\alpha_2}{12}\\
a_5 &=&\alpha_5-\dfrac{\alpha_2\alpha_3}{2} - \dfrac{\alpha_1\alpha_2^2}{12}
\end{array}
\right. \;,\qquad
\left\{\begin{array}{ccl}
\alpha_1&=&a_1\\
\alpha_2    &=&a_2\\
\alpha_3    &=&a_3 +\dfrac{a_1 a_2}{2} \\
\alpha_4    &=&a_4 + \dfrac{a_1a_3}{2} +\dfrac{a_1^2a_2}{6}\\
\alpha_5    &=&a_5 + \dfrac{a_2 a_3}{2}+\dfrac{a_1 a_2^2 }{3} 
\end{array}
\right.\; .
\end{equation*}

More explicitly, we have
$$
\exp(\alpha_5 X_5)\exp(\alpha_4 X_4)\cdots\exp(\alpha_1 X_1)=\qquad\qquad\qquad\qquad\qquad\qquad\qquad\qquad\qquad\qquad$$
$$=
\exp\left( \alpha_1 X_1+\alpha_2 X_2 + (\alpha_3-\dfrac{\alpha_1\alpha_2}{2})X_3 +(\alpha_4-\dfrac{\alpha_1\alpha_3}{2}+\dfrac{\alpha_1^2\alpha_2}{12})X_4+
(\alpha_5-\dfrac{\alpha_2\alpha_3}{2} - \dfrac{\alpha_1\alpha_2^2}{12}) X_5\right),$$
and viceversa
$$
\exp\left( a_1 X_1+a_2 X_2 + a_3 X_3 +a_4 X_4+a_5 X_5\right) 
=\qquad\qquad\qquad\qquad\qquad\qquad\qquad\qquad$$
$$=
\exp((a_5+ \dfrac{a_2 a_3}{2}  +\dfrac{a_1 a_2^2}{3}  ) X_5 )\exp((a_4+ \dfrac{a_1 a_3}{2}  +\dfrac{a_1^2 a_2}{6}  )  X_4)
\cdot\exp((a_3 +\dfrac{1}{2} a_1 a_2)X_3) 
\cdot\exp(a_2 X_2) \cdot\exp(a_1 X_1) \,.$$


 \subsection*{$\mathbb F_{24}$.}

The following is the free-nilpotent Lie algebra with 2 generators and nilpotency step 4. It has dimension 8.

The non-trivial brackets coming from the Hall basis are:
\begin{equation*}
 \begin{aligned}
     &\qquad[X_1,X_2]=X_{3}\,,\,[X_1,X_3]=X_4\,,\,[X_2,X_3]=X_{5}\,,\\ &[X_1,X_4]=X_6\,,\,[X_1,X_5]= [X_2,X_4]=X_7\,,\,
     [X_2,X_5]=X_8\,.
     \end{aligned}
 \end{equation*}
 
  The composition law \eqref{group law in G} of $\mathbb{F}_{24}$ is given by:

\begin{itemize}
    \item $z_1=x_1+y_1$;
    \item $z_2=x_2+y_2$;
    \item $z_3=x_3+y_3+\frac{1}{2}(x_1y_2-x_2y_1)$;
    \item $z_4=x_4+y_4+\frac{1}{2}(x_1y_3-x_3y_1)+\frac{1}{12}(x_1-y_1)(x_1y_2-x_2y_1)$;
    \item $z_5=x_5+y_5+\frac{1}{2}(x_2y_3-x_3y_2)+\frac{1}{12}(x_2-y_2)(x_1y_2-x_2y_1)$;
    \item $z_6=x_6+y_6+\frac{1}{2}(x_1y_4-x_4y_1)+\frac{1}{12}(x_1-y_1)(x_1y_3-x_3y_1) -\frac{1}{24}x_1y_1(x_1y_2-x_2y_1)$;
    \item $z_7=x_7+y_7+\frac{1}{2}(x_1y_5-x_5y_1+x_2y_4-x_4y_2)+\frac{1}{12}(x_1-y_1)(x_2y_3-x_3y_2)\salto+\frac{1}{12}(x_2-y_2)(x_1y_3-x_3y_1)]-\frac{1}{24}(x_1y_2+x_2y_1)(x_1y_2-x_2y_1)$;
    \item $z_8=x_8+y_8+\frac{1}{2}(x_2y_5-x_5y_2)+\frac{1}{12}(x_2-y_2)(x_2y_3-x_3y_2)-\frac{1}{24}x_2y_2(x_1y_2-x_2y_1)$
\end{itemize}

Since 

\begin{eqnarray*}
  \mathrm{d}(L_\mathbf{x})_\mathbf{0}=\left[\begin{matrix} 
   1 & 0  & 0& 0& 0& 0& 0& 0\\
   0 & 1& 0& 0& 0& 0& 0 & 0\\
  -\frac{x_2}{2} & \frac{x_1}{2}& 1& 0& 0& 0& 0 & 0\\
   -\frac{x_3}{2}-\frac{x_1x_2}{12} & \frac{x_1^2}{12}& \frac{x_1}{2}& 1& 0& 0& 0& 0\\
   -\frac{x_2^2}{12} &  \frac{x_1x_2}{12}-\frac{x_3}{2}& \frac{x_2}{2}& 0& 1& 0& 0& 0\\
 -\frac{x_4}{2}-\frac{x_1x_3}{12} & 0 & \frac{x_1^2}{12}& \frac{x_1}{2}& 0& 1& 0 &0\\
 -\frac{x_5}{2}-\frac{x_2x_3}{12}&  -\frac{x_4}{2}-\frac{x_1x_3}{12}& \frac{x_1x_2}{6}& \frac{x_2}{2}& \frac{x_1}{2}& 0& 1 &0\\
 0 & -\frac{x_5}{2}-\frac{x_2x_3}{12} & \frac{x_2^2}{12} & 0 & \frac{x_2}{2}& 0& 0 &1
   \end{matrix}\right]\,,
   \end{eqnarray*}
the induced left-invariant vector fields \eqref{leftinvariant vf} are: 
\begin{itemize}
\item $X_1={\partial}_{x_1}-\frac{x_2}{2}{\partial}_{x_3}-\big(\frac{x_3}{2}+\frac{x_1x_2}{12}\big){\partial}_{x_4}-\frac{x_2^2}{12}\partial_{x_5}-\big(\frac{x_4}{2}+\frac{x_1x_3}{12}\big)\partial_{x_6}-\big(\frac{x_5}{2}+\frac{x_2x_3}{12}\big)\partial_{x_7}\,;$
\item $X_2={\partial}_{x_2}+\frac{x_1}{2}{\partial}_{x_3}+\frac{x_1^2}{12}{\partial}_{x_4}+\big(\frac{x_1x_2}{12}-\frac{x_3}{2}\big)\partial_{x_5}-\big(\frac{x_4}{2}+\frac{x_1x_3}{12}\big)\partial_{x_7}-\big(\frac{x_5}{2}+\frac{x_2x_3}{12}\big)\partial_{x_8}\,;$
\item $X_3={\partial}_{x_3}+\frac{x_1}{2}{\partial}_{x_4}+\frac{x_2}{2}{\partial}_{x_5}+\frac{x_1^2}{12}\partial_{x_6}+\frac{x_1x_2}{6}\partial_{x_7}+\frac{x_2^2}{12}\partial_{x_8}\,;$
\item $X_4={\partial}_{x_4}+\frac{x_1}{2}{\partial}_{x_6}+\frac{x_2}{2}{\partial}_{x_7}\,;$
\item $X_5={\partial}_{x_5}+\frac{x_1}{2}{\partial}_{x_7}+\frac{x_2}{2}{\partial}_{x_8}\,$;
\item $X_6={\partial}_{x_6}$;
\item $X_7={\partial}_{x_7}$;
\item $X_8={\partial}_{x_8}$,
\end{itemize}
and the respective left-invariant 1-forms \eqref{leftinvariant form} are: 
\begin{itemize}
\item $\theta_1=dx_1$;
\item $\theta_2=dx_2$;
\item $\theta_3=dx_3-\frac{x_1}{2}dx_2+\frac{x_2}{2}dx_1$;
\item $\theta_4=dx_4-\frac{x_1}{2}dx_3+\frac{x_1^2}{6}dx_2+\big(\frac{x_3}{2}-\frac{x_1x_2}{6}\big)dx_1$;
\item $\theta_5=dx_5-\frac{x_2}{2}dx_3+\big(\frac{x_3}{2}+\frac{x_1x_2}{6}\big)dx_2-\frac{x_2^2}{6}dx_1$;
\item $\theta_6=dx_6-\frac{x_1}{2}dx_4+\frac{x_1^2}{6}dx_3-\frac{x_1^3}{24}dx_2+\big(\frac{x_4}{2}-\frac{x_1x_3}{6}+\frac{x_1^2x_2}{24}\big)dx_1$;
\item $\theta_7=dx_7-\frac{x_1}{2}dx_5-\frac{x_2}{2}dx_4+\frac{x_1x_2}{3}dx_3+\big(\frac{x_4}{2}-\frac{x_1x_3}{6}-\frac{x_1^2x_2}{12}\big)dx_2+\big(\frac{x_5}{2}-\frac{x_2x_3}{6}+\frac{x_1x_2^2}{12}\big)dx_1$;
\item $\theta_8=dx_8-\frac{x_2}{2}dx_5+\frac{x_2^2}{6}dx_3+\big(\frac{x_5}{2}-\frac{x_2x_3}{6}-\frac{x_1x_2^2}{24}\big)dx_2+\frac{x_2^3}{24}dx_1$.
\end{itemize}

Finally, we have
\begin{eqnarray*}
  \mathrm{d}(R_\mathbf{x})_\mathbf{0}=\left[\begin{matrix} 
   1 & 0  & 0& 0& 0& 0& 0& 0\\
   0 & 1& 0& 0& 0& 0& 0 & 0\\
  \frac{x_2}{2} & -\frac{x_1}{2}& 1& 0& 0& 0& 0 & 0\\
   \frac{x_3}{2}-\frac{x_1x_2}{12} & \frac{x_1^2}{12}& -\frac{x_1}{2}& 1& 0& 0& 0& 0\\
   -\frac{x_2^2}{12} &  \frac{x_1x_2}{12}+\frac{x_3}{2}& -\frac{x_2}{2}& 0& 1& 0& 0& 0\\
 \frac{x_4}{2}-\frac{x_1x_3}{12} & 0 & \frac{x_1^2}{12}& -\frac{x_1}{2}& 0& 1& 0 &0\\
 \frac{x_5}{2}-\frac{x_2x_3}{12}&  \frac{x_4}{2}-\frac{x_1x_3}{12}& \frac{x_1x_2}{6}& -\frac{x_2}{2}& -\frac{x_1}{2}& 0& 1 &0\\
 0 & \frac{x_5}{2}-\frac{x_2x_3}{12} & \frac{x_2^2}{12} & 0 & -\frac{x_2}{2}& 0& 0 &1
   \end{matrix}\right]\,.
   \end{eqnarray*}

One can relate the  exponential coordinates of first kind to the exponential coordinates of second kind.
If we denote by $\alpha_1,\ldots, \alpha_8$ the coordinates of second kind and by
$a_1,\ldots, a_8$ the coordinates of first kind, the change of coordinates are given as follows:

\begin{equation*}
\left\{\begin{array}{ccl}
a_1 &=&\alpha_1 \\
a_2 &=&\alpha_2\\
a_3 &=&\alpha_3-\frac{\alpha_1\alpha_2}{2}\\
a_4 &=&\alpha_4-\frac{\alpha_1\alpha_3}{2}+\frac{\alpha_1^2\alpha_2}{12}\\
a_5 &=&\alpha_5-\frac{\alpha_2\alpha_3}{2} - \frac{\alpha_1\alpha_2^2}{12}\\
a_6 &=& \alpha_6-\frac{\alpha_1\alpha_4}{2}+\frac{\alpha_1^2\alpha_3}{12}\\
a_7 &=& \alpha_7-\frac{\alpha_1\alpha_5+\alpha_2\alpha_4}{2}+\frac{\alpha_1\alpha_2\alpha_3}{6}+\frac{\alpha_1^2\alpha_2^2}{24}\\
a_8&=&\alpha_8-\frac{\alpha_2\alpha_5}{2}+\frac{\alpha_2^2\alpha_3}{12}
\end{array}
\right. \qquad\qquad\qquad\qquad\qquad\qquad
\end{equation*}

and

\begin{equation*}
\left\{\begin{array}{ccl}
\alpha_1&=&a_1\\
\alpha_2    &=&a_2\\
\alpha_3    &=&a_3 +\frac{a_1 a_2}{2} \\
\alpha_4    &=&a_4 + \frac{a_1a_3}{2} +\frac{a_1^2a_2}{6}\\
\alpha_5    &=&a_5 + \frac{a_2 a_3}{2}+\frac{a_1 a_2^2 }{3} \\
\alpha_6 &=&a_6+\frac{a_1a_4}{2}+\frac{a_1^2a_3}{6}+\frac{a_1^3a_2}{24}\\
\alpha_7 &=&a_7+\frac{a_1a_5+a_2a_4}{2}+\frac{a_1a_2a_3}{3}+\frac{a_1^2a_2^2}{8}\\
\alpha_8 & =& a_8+\frac{a_2a_5}{2}+\frac{a_2^2a_3}{6}+\frac{a_1a_2^3}{8}\,. 
\end{array}
\right.\; .\qquad\qquad\qquad\qquad\qquad\qquad\qquad\quad
\end{equation*}

More explicitly, we have
$$
\exp(\alpha_8 X_8)\exp(\alpha_7 X_7)\cdots\exp(\alpha_1 X_1)=\qquad\qquad\qquad\qquad\qquad\qquad\qquad\qquad\qquad\qquad$$
$$=
\exp\left( \alpha_1 X_1+\alpha_2 X_2 + (\alpha_3-\frac{\alpha_1\alpha_2}{2})X_3 +(\alpha_4-\frac{\alpha_1\alpha_3}{2}+\frac{\alpha_1^2\alpha_2}{12})X_4+
\right.$$
$$+(\alpha_5-\frac{\alpha_2\alpha_3}{2} - \frac{\alpha_1\alpha_2^2}{12}) X_5+(\alpha_6-\frac{\alpha_1\alpha_4}{2}+\frac{\alpha_1^2\alpha_3}{12})X_6+\qquad\qquad$$
$$\left.\qquad\qquad+(\alpha_7-\frac{\alpha_1\alpha_5+\alpha_2\alpha_4}{2}+\frac{\alpha_1\alpha_2\alpha_3}{6}+\frac{\alpha_1^2\alpha_2^2}{24})X_7+(\alpha_8-\frac{\alpha_2\alpha_5}{2}+\frac{\alpha_2^2\alpha_3}{12})X_8\right) $$
and viceversa
$$
\exp\left( a_1 X_1+a_2 X_2 + \cdots+a_7X_7+a_8X_8\right)
=\qquad\qquad\qquad\qquad\qquad\qquad\qquad\qquad\quad$$
$$\quad\qquad=\exp\bigg((a_8+\frac{a_2a_5}{2}+\frac{a_2^2a_3}{6}+\frac{a_1a_2^3}{8})X_8\bigg)\cdot
\exp\bigg((a_7+\frac{a_1a_5+a_2a_4}{2}+\frac{a_1a_2a_3}{3}+\frac{a_1^2a_2^2}{8})X_7\bigg)\cdot $$
$$\quad\cdot\exp\bigg((a_6+\frac{a_1a_4}{2}+\frac{a_1^2a_3}{6}+\frac{a_1^3a_2}{24})X_6\bigg)\cdot\exp\bigg((a_5+ \frac{a_2 a_3}{2}  +\frac{a_1 a_2^2}{3}  ) X_5 \bigg)\cdot$$
$$\qquad\qquad\cdot\exp\bigg((a_4+ \frac{a_1 a_3}{2}  +\frac{a_1^2 a_2}{6}  )  X_4\bigg)\cdot\exp\bigg((a_3 +\frac{1}{2} a_1 a_2)X_3\bigg) 
\cdot\exp(a_2 X_2) \cdot\exp(a_1 X_1) \,.$$


 \subsection*{$\mathbb F_{33}$.}

The following is the free-nilpotent Lie algebra with 3 generators and nilpotency step 3. It has dimension 14.

The non-trivial brackets coming from the Hall basis are:
\begin{equation*}
 \begin{aligned}
     &\;\quad[X_1,X_2]=X_{4}\,,\,[X_1,X_3]=X_5\,,\,[X_2,X_3]=X_{6}\,,\,[X_1,X_4]=X_7\,,\\ &\quad[X_1,X_5]=X_8\,,\,[X_1,X_6]=X_9\,,\, [X_2,X_4]=X_{10}\,,\,
     [X_2,X_6]=X_{11}\,,\\
     &[X_3,X_4]=X_{12}\,,\,[X_3,X_5]=X_{13}\,,\,[X_3,X_6]=X_{14}\,,\,[X_2,X_5]=X_9+X_{12}\,.
     \end{aligned}
 \end{equation*}
 
   The composition law \eqref{group law in G} of $\mathbb{F}_{33}$ is given by:

\begin{itemize}
    \item $z_1=x_1+y_1$;
    \item $z_2=x_2+y_2$;
    \item $z_3=x_3+y_3$;
    \item $z_4=x_4+y_4+\frac{1}{2}(x_1y_2-x_2y_1)$;
    \item $z_5=x_5+y_5+\frac{1}{2}(x_1y_3-x_3y_1)$;
    \item $z_6=x_6+y_6+\frac{1}{2}(x_2y_3-x_3y_2)$;
    \item $z_7=x_7+y_7+\frac{1}{2}(x_1y_4-x_4y_1)+\frac{1}{12}(x_1-y_1)(x_1y_2-x_2y_1)$;
    \item $z_8=x_8+y_8+\frac{1}{2}(x_1y_5-x_5y_1)+\frac{1}{12}(x_1-y_1)(x_1y_3-x_3y_1)$;
    \item $z_9=x_9+y_9+\frac{1}{2}(x_1y_6-x_6y_1+x_2y_5-x_5y_2)+\frac{1}{12}(x_1-y_1)(x_2y_3-x_3y_2)\salto+\frac{1}{12}(x_2-y_2)(x_1y_3-x_3y_1)$;
    \item $z_{10}=x_{10}+y_{10}+\frac{1}{2}(x_2y_4-x_4y_2)+\frac{1}{12}(x_2-y_2)(x_1y_2-x_2y_1)$;
    \item $z_{11}=x_{11}+y_{11}+\frac{1}{2}(x_2y_6-x_6y_2)+\frac{1}{12}(x_2-y_2)(x_2y_3-x_3y_2)$;
    \item $z_{12}=x_{12}+y_{12}+\frac{1}{2}(x_2y_5-x_5y_2+x_3y_4-x_4y_3)+\frac{1}{12}(x_2-y_2)(x_1y_3-x_3y_1)\salto+\frac{1}{12}(x_2-y_2)(x_1y_2-x_2y_1)$;
    \item $z_{13}=x_{13}+y_{13}+\frac{1}{2}(x_3y_5-x_5y_3)+\frac{1}{12}(x_3-y_3)(x_1y_3-x_3y_1)$;
    \item $z_{14}=x_{14}+y_{14}+\frac{1}{2}(x_3y_6-x_6y_3)+\frac{1}{12}(x_3-y_3)(x_2y_3-x_3y_2)$.
\end{itemize}

Since 

\begin{eqnarray*}
  \mathrm{d}(L_\mathbf{x})_\mathbf{0}=\left[\begin{matrix}
   1 & 0  & 0& 0& 0& 0& 0& 0& 0& 0& 0& 0& 0& 0\\
   0 & 1& 0& 0& 0& 0& 0 & 0& 0& 0& 0& 0& 0& 0\\
  0 & 0& 1& 0& 0& 0& 0 & 0& 0& 0& 0& 0& 0& 0\\
   -\frac{x_2}{2} & \frac{x_1}{2}& 0& 1& 0& 0& 0& 0& 0& 0& 0& 0& 0& 0\\
   -\frac{x_3}{2} &  0& \frac{x_1}{2}& 0& 1& 0& 0& 0& 0& 0& 0& 0& 0& 0\\
 0 & -\frac{x_3}{2} & \frac{x_2}{2}& 0& 0& 1& 0 &0& 0& 0& 0& 0& 0& 0\\
 -\frac{x_4}{2}-\frac{x_1x_2}{12}&  \frac{x_1^2}{12}& 0& \frac{x_1}{2}& 0& 0& 1 &0& 0& 0& 0& 0& 0& 0\\
 -\frac{x_5}{2}-\frac{x_1x_3}{12} & 0 & \frac{x_1^2}{12} & 0 & \frac{x_1}{2}& 0& 0 &1& 0& 0& 0& 0& 0& 0\\
-\frac{x_6}{2}-\frac{x_2x_3}{12} & -\frac{x_5}{2}-\frac{x_1x_3}{12} & \frac{x_1x_2}{6} & 0 & \frac{x_2}{2} & \frac{x_1}{2} & 0 & 0& 1 & 0 & 0 & 0 & 0 & 0 \\
 -\frac{x_2^2}{12} & \frac{x_1x_2}{12}-\frac{x_4}{2} & 0 & \frac{x_2}{2} & 0 & 0 & 0 & 0 & 0 & 1 & 0 & 0 & 0 & 0 \\
 0 & -\frac{x_2x_3}{12}-\frac{x_6}{2} & \frac{x_2^2}{12} & 0 & 0 & \frac{x_2}{2} & 0 & 0 & 0 & 0 & 1& 0 & 0 & 0 \\
 -\frac{x_2x_3}{6} & \frac{x_1x_3}{12}-\frac{x_5}{2}& \frac{x_1x_2}{12}-\frac{x_4}{2} & \frac{x_3}{2} & \frac{x_2}{2} & 0 & 0 & 0 & 0 & 0 & 0 & 1 & 0 & 0\\
 -\frac{x_3^2}{12} & 0 & \frac{x_1x_3}{12}-\frac{x_5}{2} & 0 & \frac{x_3}{2} & 0 & 0 & 0 & 0 & 0 & 0 & 0 & 1 & 0\\
 0 & -\frac{x_3^2}{12} & \frac{x_2x_3}{12}-\frac{x_6}{2} & 0 &0 & \frac{x_3}{2} & 0 & 0 & 0 & 0 & 0 & 0 & 0 & 1
   \end{matrix}\right]\,,
   \end{eqnarray*}
the induced left-invariant vector fields \eqref{leftinvariant vf} are: 
\begin{itemize}
\item $X_1={\partial}_{x_1}-\frac{x_2}{2}{\partial}_{x_4}-\frac{x_3}{2}\partial_{x_5}-\big(\frac{x_4}{2}+\frac{x_1x_2}{12}\big){\partial}_{x_7}-\big(\frac{x_5}{2}+\frac{x_1x_3}{12}\big)\partial_{x_8}-\big(\frac{x_6}{2}+\frac{x_2x_3}{12}\big)\partial_{x_{9}}\saltox-\frac{x_2^2}{12}\partial_{x_{10}}-\frac{x_2x_3}{6}\partial_{x_{12}}-\frac{x_3^2}{12}\partial_{x_{13}}\,;$
\item $X_2={\partial}_{x_2}+\frac{x_1}{2}{\partial}_{x_4}-\frac{x_3}{2}\partial_{x_{6}}+\frac{x_1^2}{12}{\partial}_{x_7}-\big(\frac{x_1x_3}{12}+\frac{x_5}{2}\big)\partial_{x_9}+\big(\frac{x_1x_2}{12}-\frac{x_4}{2}\big)\partial_{x_{10}}-\big(\frac{x_6}{2}\saltox+\frac{x_2x_3}{12}\big)\partial_{x_{11}}+\big(\frac{x_1x_3}{12}-\frac{x_5}{2}\big)\partial_{x_{12}}-\frac{x_3^2}{12}\partial_{x_{14}}\,;$
\item $X_3={\partial}_{x_3}+\frac{x_1}{2}{\partial}_{x_5}+\frac{x_2}{2}{\partial}_{x_6}+\frac{x_1^2}{12}\partial_{x_8}+\frac{x_1x_2}{6}\partial_{x_9}+\frac{x_2^2}{12}\partial_{x_{11}}+\big(\frac{x_1x_2}{12}-\frac{x_4}{2}\big)\partial_{x_{12}}+\big(\frac{x_1x_3}{12}\saltox-\frac{x_5}{2}\big)\partial_{x_{13}}+\big(\frac{x_2x_3}{12}-\frac{x_6}{2}\big)\partial_{x_{14}}\,;$
\item $X_4={\partial}_{x_4}+\frac{x_1}{2}\partial_{x_7}+\frac{x_2}{2}\partial_{x_{10}}+\frac{x_3}{2}\partial_{x_{12}}\,;$
\item $X_5={\partial}_{x_5}+\frac{x_1}{2}\partial_{x_8}+\frac{x_2}{2}\partial_{x_9}+\frac{x_2}{2}\partial_{x_{12}}+\frac{x_3}{2}\partial_{x_{13}}\,$;
\item $X_6={\partial}_{x_6}+\frac{x_1}{2}\partial_{x_{9}}+\frac{x_2}{2}\partial_{x_{10}}$;
\item $X_7={\partial}_{x_7}+\frac{x_1}{2}\partial_{x_9}+\frac{x_2}{2}\partial_{x_{11}}+\frac{x_3}{2}\partial_{x_{14}}$;
\item $X_8={\partial}_{x_8}$;
\item $X_9=\partial_{x_{9}}$;
\item $X_{10}=\partial_{x_{10}}$;
\item $X_{11}=\partial_{x_{11}}$;
\item $X_{12}=\partial_{x_{12}}$;
\item $X_{13}=\partial_{x_{13}}$;
\item $X_{14}=\partial_{x_{14}}$,
\end{itemize}
and the respective left-invariant 1-forms \eqref{leftinvariant form} are: 
\begin{itemize}
\item $\theta_1=dx_1$;
\item $\theta_2=dx_2$;
\item $\theta_3=dx_3$;
\item $\theta_4=dx_4-\frac{x_1}{2}dx_2+\frac{x_2}{2}dx_1$;
\item $\theta_5=dx_5-\frac{x_1}{2}dx_3+\frac{x_3}{2}dx_1$;
\item $\theta_6=dx_6-\frac{x_2}{2}dx_3+\frac{x_3}{2}dx_2$;
\item $\theta_7=dx_7-\frac{x_1}{2}dx_4+\frac{x_1^2}{6}dx_2+\big(\frac{x_4}{2}-\frac{x_1x_2}{6}\big)dx_1$;
\item $\theta_8=dx_8-\frac{x_1}{2}dx_5+\frac{x_1^2}{6}dx_3+\big(\frac{x_5}{2}-\frac{x_1x_3}{6}\big)dx_1$;
\item $\theta_9=dx_9-\frac{x_1}{2}dx_6-\frac{x_2}{2}dx_5+\frac{x_1x_2}{3}dx_3+\big(\frac{x_5}{2}-\frac{x_1x_3}{6}\big)dx_2+\big(\frac{x_6}{2}-\frac{x_2x_3}{6}\big) dx_1$;
\item $\theta_{10}=dx_{10}-\frac{x_2}{2}dx_4+\big(\frac{x_4}{2}+\frac{x_1x_2}{6}\big)dx_2-\frac{x_2^2}{6}dx_1$;
\item $\theta_{11}=dx_{11}-\frac{x_2}{2}dx_6+\frac{x_2^2}{6}dx_3+\big(\frac{x_6}{2}-\frac{x_2x_3}{6}\big)dx_2$;
\item $\theta_{12}=dx_{12}-\frac{x_2}{2}dx_5-\frac{x_3}{2}dx_4+\big(\frac{x_1x_2}{6}+\frac{x_4}{2}\big)dx_3+\big(\frac{x_1x_3}{6}+\frac{x_5}{2}\big)dx_2-\frac{x_2x_3}{3}dx_1$;
\item $\theta_{13}=dx_{13}-\frac{x_3}{2}dx_5+\big(\frac{x_1x_3}{6}+\frac{x_5}{2}\big)dx_3-\frac{x_3^2}{6}dx_1$;
\item $\theta_{14}=dx_{14}-\frac{x_3}{2}dx_6+\big(\frac{x_6}{2}+\frac{x_2x_3}{6}\big)dx_3-\frac{x_3^2}{6}dx_2$.
\end{itemize}

Finally, we have
\begin{eqnarray*}
  \mathrm{d}(R_\mathbf{x})_\mathbf{0}=\left[\begin{matrix}
   1 & 0  & 0& 0& 0& 0& 0& 0& 0& 0& 0& 0& 0& 0\\
   0 & 1& 0& 0& 0& 0& 0 & 0& 0& 0& 0& 0& 0& 0\\
  0 & 0& 1& 0& 0& 0& 0 & 0& 0& 0& 0& 0& 0& 0\\
   \frac{x_2}{2} & -\frac{x_1}{2}& 0& 1& 0& 0& 0& 0& 0& 0& 0& 0& 0& 0\\
   \frac{x_3}{2} &  0& -\frac{x_1}{2}& 0& 1& 0& 0& 0& 0& 0& 0& 0& 0& 0\\
 0 & \frac{x_3}{2} & -\frac{x_2}{2}& 0& 0& 1& 0 &0& 0& 0& 0& 0& 0& 0\\
 \frac{x_4}{2}-\frac{x_1x_2}{12}&  \frac{x_1^2}{12}& 0& -\frac{x_1}{2}& 0& 0& 1 &0& 0& 0& 0& 0& 0& 0\\
 \frac{x_5}{2}-\frac{x_1x_3}{12} & 0 & \frac{x_1^2}{12} & 0 & -\frac{x_1}{2}& 0& 0 &1& 0& 0& 0& 0& 0& 0\\
\frac{x_6}{2}-\frac{x_2x_3}{12} & \frac{x_5}{2}-\frac{x_1x_3}{12} & \frac{x_1x_2}{6} & 0 & -\frac{x_2}{2} & -\frac{x_1}{2} & 0 & 0& 1 & 0 & 0 & 0 & 0 & 0 \\
 -\frac{x_2^2}{12} & \frac{x_1x_2}{12}+\frac{x_4}{2} & 0 & -\frac{x_2}{2} & 0 & 0 & 0 & 0 & 0 & 1 & 0 & 0 & 0 & 0 \\
 0 & \frac{x_6}{2}-\frac{x_2x_3}{12} & \frac{x_2^2}{12} & 0 & 0 & -\frac{x_2}{2} & 0 & 0 & 0 & 0 & 1& 0 & 0 & 0 \\
 -\frac{x_2x_3}{6} & \frac{x_1x_3}{12}+\frac{x_5}{2}& \frac{x_1x_2}{12}+\frac{x_4}{2} & -\frac{x_3}{2} & -\frac{x_2}{2} & 0 & 0 & 0 & 0 & 0 & 0 & 1 & 0 & 0\\
 -\frac{x_3^2}{12} & 0 & \frac{x_1x_3}{12}+\frac{x_5}{2} & 0 & -\frac{x_3}{2} & 0 & 0 & 0 & 0 & 0 & 0 & 0 & 1 & 0\\
 0 & -\frac{x_3^2}{12} & \frac{x_2x_3}{12}+\frac{x_6}{2} & 0 &0 & -\frac{x_3}{2} & 0 & 0 & 0 & 0 & 0 & 0 & 0 & 1
   \end{matrix}\right]\,.
   \end{eqnarray*}
   
  A representation of the left-invariant vector fields $X_i^s$ defined on $\mathbb{R}^{14}$ with respect to exponential coordinates of the second kind is the following:
  \begin{itemize}
      \item $X_1^s=\partial_{x_1}\,;$
      \item $X_2^s=\partial_{x_2}+x_1\partial_{x_4}+\frac{x_1^2}{2}\partial_{x_7}+x_1x_2\partial_{x_{10}}$;
      \item $X_3^s=\partial_{x_3}+x_1\partial_{x_5}+x_2\partial_{x_6}+\frac{x_1^2}{2}\partial_{x_8}+x_1x_2\partial_{x_9}+\frac{x_2^2}{2}\partial_{x_{11}}+({x_1x_2}-x_4)\partial_{x_{12}}+x_1x_3\partial_{x_{13}}+x_2x_3\partial_{x_{14}}$;
      \item $X_4^s=\partial_{x_4}+x_1\partial_{x_7}+x_2\partial_{x_{10}}$;
      \item $X_5^s=\partial_{x_5}+x_1\partial_{x_8}+x_2\partial_{x_9}+x_2\partial_{x_{12}}+x_3\partial_{x_{13}}$;
      \item $X_6^s=\partial_{x_6}+x_1\partial_{x_9}+x_2\partial_{x_{11}}+x_3\partial_{x_{14}}$;
      \item $X_7^s=\partial_{x_7}$;
      \item $X_8^s=\partial_{x_8}$;
      \item $X_9^s=\partial_{x_9}$;
      \item $X_{10}^s=\partial_{x_{10}}$;
      \item $X_{11}^s=\partial_{x_{11}}$;
      \item $X_{12}^s=\partial_{x_{12}}$;
      \item $X_{13}^s=\partial_{x_{13}}$;
      \item $X_{14}^s=\partial_{x_{14}}$.
  \end{itemize}

One can relate the  exponential coordinates of first kind to the exponential coordinates of second kind.
If we denote by $\alpha_1,\ldots, \alpha_{14}$ the coordinates of second kind and by
$a_1,\ldots, a_{14}$ the coordinates of first kind, the change of coordinates are given as follows:

\begin{equation*}
\left\{\begin{array}{ccl}
a_1 &=&\alpha_1 \\
a_2 &=&\alpha_2\\
a_3 &=&\alpha_3\\
a_4 &=&\alpha_4-\frac{\alpha_1\alpha_2}{2}\\
a_5 &=&\alpha_5-\frac{\alpha_1\alpha_3}{2} \\
a_6 &=& \alpha_6-\frac{\alpha_2\alpha_3}{2}\\
a_7 &=& \alpha_7-\frac{\alpha_1\alpha_4}{2}+\frac{\alpha_1^2\alpha_2}{12}\\
a_8&=&\alpha_8-\frac{\alpha_1\alpha_5}{2}+\frac{\alpha_1^2\alpha_{3}}{12}\\
a_9&=&\alpha_9-\frac{\alpha_1\alpha_6+\alpha_2\alpha_5}{2}+\frac{\alpha_1\alpha_2\alpha_3}{6}\\
a_{10}&=&\alpha_{10}-\frac{\alpha_2\alpha_4}{2}-\frac{\alpha_1\alpha_2^2}{12}\\
a_{11}&=&\alpha_{11}-\frac{\alpha_2\alpha_6}{2}+\frac{\alpha_2^2\alpha_3}{12}\\
a_{12}&=&\alpha_{12}+\frac{\alpha_3\alpha_4-\alpha_2\alpha_5}{2}-\frac{\alpha_1\alpha_2\alpha_3}{6}\\
a_{13}&=&\alpha_{13}-\frac{\alpha_3\alpha_5}{2}-\frac{\alpha_1\alpha_3^2}{12}\\
a_{14}&=&\alpha_{14}-\frac{\alpha_3\alpha_6}{2}-\frac{\alpha_2\alpha_3^2}{12}
\end{array}
\right. \qquad\qquad\qquad\qquad\qquad\qquad
\end{equation*}

and

\begin{equation*}
\left\{\begin{array}{ccl}
\alpha_1&=&a_1\\
\alpha_2    &=&a_2\\
\alpha_3    &=&a_3  \\
\alpha_4    &=&a_4 +\frac{a_1 a_2}{2}\\
\alpha_5    &=&a_5 +\frac{a_1 a_3 }{2} \\
\alpha_6 &=&a_6+ \frac{a_2 a_3}{2}\\
\alpha_7 &=&a_7+\frac{a_1a_4}{2}+\frac{a_1^2a_2}{6}\\
\alpha_8 & =& a_8+\frac{a_1a_5}{2}+\frac{a_1^2a_3}{6}\\
\alpha_9&=&a_9+\frac{a_1a_6+a_2a_5}{2}+\frac{a_1a_2a_3}{3}\\
\alpha_{10}&=&a_{10}+\frac{a_2a_4}{2}+\frac{a_1a_2^2}{3}\\
\alpha_{11}&=&a_{11}+\frac{a_2a_6}{2}+\frac{a_2^2a_3}{6}\\
\alpha_{12}&=&a_{12}+\frac{a_2a_5-a_3a_4}{2}+\frac{a_1a_2a_3}{6}\\
\alpha_{13}&=&a_{13}+\frac{a_3a_5}{2}+\frac{a_1a_3^2}{3}\\
\alpha_{14}&=&a_{14}+\frac{a_3a_6}{2}+\frac{a_2a_3^2}{3}\,. 
\end{array}
\right.\; .\qquad\qquad\qquad\qquad\qquad\qquad\qquad\quad
\end{equation*}

More explicitly, we have
$$
\exp(\alpha_{14} X_{14})\exp(\alpha_{13} X_{13})\cdots\exp(\alpha_1 X_1)=\qquad\qquad\qquad\qquad\qquad\qquad\qquad\qquad\qquad\qquad$$
$$=
\exp\left( \alpha_1 X_1+\alpha_2 X_2 + \alpha_3X_3 +(\alpha_4-\frac{\alpha_1\alpha_2}{2})X_4+(\alpha_5-\frac{\alpha_1\alpha_3}{2} ) X_5+\qquad
\right.$$
$$+(\alpha_6-\frac{\alpha_2\alpha_3}{2})X_6+(\alpha_7-\frac{\alpha_1\alpha_4}{2}+\frac{\alpha_1^2\alpha_2}{12})X_7+(\alpha_8-\frac{\alpha_1\alpha_5}{2}+\frac{\alpha_1^2\alpha_3}{12})X_8+ $$
$$+(\alpha_9-\frac{\alpha_1\alpha_6+\alpha_2\alpha_5}{2}+\frac{\alpha_1\alpha_2\alpha_3}{6})X_9+(\alpha_{10}-\frac{\alpha_2\alpha_4}{2}-\frac{\alpha_1\alpha_2^2}{12})X_{10}+\quad$$
$$+(\alpha_{11}-\frac{\alpha_2\alpha_6}{2}+\frac{\alpha_2^2\alpha_3}{12})X_{11}+(\alpha_{12}+\frac{\alpha_3\alpha_4-\alpha_2\alpha_5}{2}-\frac{\alpha_1\alpha_2\alpha_3}{6})X_{12}+ $$
$$\left.+(\alpha_{13}-\frac{\alpha_3\alpha_5}{2}-\frac{\alpha_1\alpha_3^2}{12})X_{13}+(\alpha_{14}-\frac{\alpha_3\alpha_6}{2}-\frac{\alpha_2\alpha_3^2}{12})X_{14}\right)\qquad  $$
and viceversa
$$
\exp\left( a_1 X_1+a_2 X_2 + \cdots+a_{13}X_{13}+a_{14}X_{14}\right)
=\qquad\qquad\qquad\qquad\qquad\qquad\qquad\qquad\quad$$
$$\qquad\qquad
=\exp\bigg((a_{14}+\frac{a_3a_6}{2}+\frac{a_2a_3^2}{3})X_{14}\bigg)\cdot\exp\bigg((a_{13}+\frac{a_3a_5}{2}+\frac{a_1a_3^2}{3} )X_{13}\bigg)\cdot\qquad\qquad\quad$$
$$\qquad\cdot \exp\bigg( (a_{12}+\frac{a_2a_5-a_3a_4}{2}+\frac{a_1a_2a_3}{6})X_{12}\bigg)\cdot \exp\bigg((a_{11}+\frac{a_2a_6}{2}+\frac{a_2^2a_3}{6}) X_{11}\bigg)\cdot\qquad $$
$$\cdot\exp\bigg((a_{10}+\frac{a_2a_4}{2}+\frac{a_1a_2^2}{3})X_{10}\bigg)\cdot\exp\bigg((a_9+\frac{a_1a_6+a_2a_5}{2}+\frac{a_1a_2a_3}{3})X_9\bigg)\cdot\qquad  $$
$$\cdot\exp\bigg((a_8+\frac{a_1a_5}{2}+\frac{a_1^2a_3}{6})X_8\bigg)\cdot
\exp\bigg((a_7+\frac{a_1a_4}{2}+\frac{a_1^2a_2}{6})X_7\bigg)\cdot \exp\bigg((a_6+ \frac{a_2 a_3}{2})X_6\bigg)\cdot$$
$$\cdot\exp\bigg((a_5 +\frac{a_1 a_3 }{2}  ) X_5 \bigg)\cdot\exp\bigg((a_4+ \frac{a_1 a_2}{2}   )  X_4\bigg)\cdot\exp(a_3 X_3) 
\cdot\exp(a_2 X_2) \cdot\exp(a_1 X_1) \,.$$


 \subsection*{$\mathbb F_{25}$.}

The following is the free-nilpotent Lie algebra with 2 generators and nilpotency step 5. It has dimension 14.

The non-trivial brackets coming from the Hall basis are:
\begin{equation*}
 \begin{aligned}
     &\;[X_1,X_2]=X_{3}\,,\,[X_1,X_3]=X_4\,,\,[X_2,X_3]=X_{5}\,,\,[X_1,X_4]=X_6\,,\\
     &
     [X_2,X_5]=X_8\,,\,[X_1,X_5]= [X_2,X_4]=X_7\,,\,[X_1,X_7]=X_{10}+X_{13}\,,\\
     &\qquad[X_1,X_6]=X_9\,,\,[X_1,X_8]=X_{11}+X_{14}\,,\,[X_2,X_6]=X_{10}\,,\\&[X_2,X_7]=X_{11}\,,\,[X_2,X_8]=X_{12}\,,\,[X_3,X_4]=X_{13}\,,\,[X_3,X_5]=X_{14}\,.
     \end{aligned}
 \end{equation*}

  The composition law \eqref{group law in G} of $\mathbb{F}_{25}$ is given by:

\begin{itemize}
    \item $z_1=x_1+y_1$;
    \item $z_2=x_2+y_2$;
    \item $z_3=x_3+y_3+\frac{1}{2}(x_1y_2-x_2y_1)$;
    \item $z_4=x_4+y_4+\frac{1}{2}(x_1y_3-x_3y_1)+\frac{1}{12}(x_1-y_1)(x_1y_2-x_2y_1)$;
    \item $z_5=x_5+y_5+\frac{1}{2}(x_2y_3-x_3y_2)+\frac{1}{12}(x_2-y_2)(x_1y_2-x_2y_1)$;
    \item $z_6=x_6+y_6+\frac{1}{2}(x_1y_4-x_4y_1)+\frac{1}{12}(x_1-y_1)(x_1y_3-x_3y_1) -\frac{1}{24}x_1y_1(x_1y_2-x_2y_1)$;
    \item $z_7=x_7+y_7+\frac{1}{2}(x_1y_5-x_5y_1+x_2y_4-x_4y_2)+\frac{1}{12}(x_1-y_1)(x_2y_3-x_3y_2)\salto+\frac{1}{12}(x_2-y_2)(x_1y_3-x_3y_1)-\frac{1}{24}(x_1y_2+x_2y_1)(x_1y_2-x_2y_1)$;
    \item $z_8=x_8+y_8+\frac{1}{2}(x_2y_5-x_5y_2)+\frac{1}{12}(x_2-y_2)(x_2y_3-x_3y_2)-\frac{1}{24}x_2y_2(x_1y_2-x_2y_1)$;
    \item $z_9=x_9+y_9+\frac{1}{2}(x_1y_6-x_6y_1)+\frac{1}{12}(x_1-y_1)(x_1y_4-x_4y_1)-\frac{1}{24}x_1y_1(x_1y_3-x_3y_1)\salto+\frac{1}{720}(y_1^3-x_1^3)(x_1y_2-x_2y_1)+\frac{1}{180}(x_1y_1^2-x_1^2y_1)(x_1y_2-x_2y_1)$;
    \item $z_{10}=x_{10}+y_{10}+\frac{1}{2}(x_1y_7-x_7y_1+x_2y_6-x_6y_2)+\frac{1}{12}(x_2-y_2)(x_1y_4-x_4y_1)\salto+\frac{1}{12}(x_1-y_1)(x_1y_5-x_5y_1+x_2y_4-x_4y_2)-\frac{1}{24}(x_1y_2+x_2y_1)(x_1y_3-x_3y_1)\salto-\frac{1}{24}x_1y_1(x_2y_3-x_3y_2)+\frac{1}{90}(x_1y_1y_2-x_1x_2y_1)(x_1y_2-x_2y_1)\salto+\frac{1}{240}(y_1^2y_2-x_1^2x_2)(x_1y_2-x_2y_1)+\frac{1}{180}(x_2y_1^2-x_1^2y_2)(x_1y_2-x_2y_1)$;
    \item $z_{11}=x_{11}+y_{11}+\frac{1}{2}(x_1y_8-x_8y_1+x_2y_7-x_7y_2)+\frac{1}{12}(x_1-y_1)(x_2y_5-x_5y_2)\salto+\frac{1}{12}(x_2-y_2)(x_1y_4-x_4y_2+x_1y_5-x_5y_1)-\frac{1}{24}(x_2y_1+x_1y_2)(x_2y_3-x_3y_2)\salto-\frac{1}{24}x_2y_2(x_1y_3-x_3y_1)+\frac{1}{240}(y_1y_2^2-x_1x_2^2)(x_1y_2-x_2y_1)\salto+\frac{1}{90}(x_2y_1y_2-x_1x_2y_2)(x_1y_2-x_2y_1)+\frac{1}{180}(x_1y_2^2-x_2^2y_1)(x_1y_2-x_2y_1)$;
    \item $z_{12}=x_{12}+y_{12}+\frac{1}{2}(x_2y_8-x_8y_2)+\frac{1}{12}(x_2-y_2)(x_2y_5-x_5y_2)-\frac{1}{24}x_2y_2(x_2y_3-x_3y_2)\salto+\frac{1}{720}(y_2^3-x_2^3)(x_1y_2-x_2y_1)+\frac{1}{180}(x_2y_2^2-x_2^2y_2)(x_1y_2-x_2y_1)$;
    \item $z_{13}=x_{13}+y_{13}+\frac{1}{2}(x_1y_7-x_7y_1+x_3y_4-x_4y_3)+\frac{1}{12}(x_3-y_3)(x_1y_3-x_3y_1)\salto+\frac{1}{12}(x_1-y_1)(x_1y_5-x_5y_1+x_2y_4-x_4y_2)+\frac{1}{12}(y_4-x_4)(x_1y_2-x_2y_1)]\salto-\frac{1}{24}[x_1y_3(x_1y_2-x_2y_1)+x_2y_1(x_1y_3-x_3y_1)+x_1y_1(x_2y_3-x_3y_2)]\salto+\frac{1}{360}(3x_1+y_1)(x_1y_2-x_2y_1)^2+\frac{1}{360}(y_1^2y_2-x_1^2x_2)(x_1y_2-x_2y_1)\salto-\frac{1}{90}x_1x_2y_1(x_1y_2-x_2y_1)+\frac{1}{180}(x_2y_1^2+x_1y_1y_2)(x_1y_2-x_2y_1)$;
    \item $z_{14}=x_{14}+y_{14}+\frac{1}{2}(x_1y_8-x_8y_1+x_3y_5-x_5y_3)+\frac{1}{12}(x_1-y_1)(x_2y_5-x_5y_2)\salto+\frac{1}{12}(x_3-y_3)(x_2y_3-x_3y_2)+\frac{1}{12}(y_5-x_5)(x_1y_2-x_2y_1)\salto-\frac{1}{24}x_2y_3(x_1y_2-x_2y_1)-\frac{1}{24}x_2y_1(x_2y_3-x_3y_2)+\frac{1}{360}(3x_2+y_2)(x_1y_2-x_2y_1)^2\salto+\frac{1}{720}(y_1y_2^2-x_1x_2^2)(x_1y_2-x_2y_1)
    +\frac{1}{180}(x_2y_1y_2-x_2^2y_1)(x_1y_2-x_2y_1)$.
\end{itemize}

Since 

\begin{eqnarray*}
  \mathrm{d}(L_\mathbf{x})_\mathbf{0}=\left[\begin{matrix}
   1 & 0  & 0& 0& 0& 0& 0& 0& 0& 0& 0& 0& 0& 0\\
   0 & 1& 0& 0& 0& 0& 0 & 0& 0& 0& 0& 0& 0& 0\\
  -\frac{x_2}{2} & \frac{x_1}{2}& 1& 0& 0& 0& 0 & 0& 0& 0& 0& 0& 0& 0\\
   -\frac{x_3}{2}-\frac{x_1x_2}{12} & \frac{x_1^2}{12}& \frac{x_1}{2}& 1& 0& 0& 0& 0& 0& 0& 0& 0& 0& 0\\
   -\frac{x_2^2}{12} &  \frac{x_1x_2}{12}-\frac{x_3}{2}& \frac{x_2}{2}& 0& 1& 0& 0& 0& 0& 0& 0& 0& 0& 0\\
 -\frac{x_4}{2}-\frac{x_1x_3}{12} & 0 & \frac{x_1^2}{12}& \frac{x_1}{2}& 0& 1& 0 &0& 0& 0& 0& 0& 0& 0\\
 -\frac{x_5}{2}-\frac{x_2x_3}{12}&  -\frac{x_4}{2}-\frac{x_1x_3}{12}& \frac{x_1x_2}{6}& \frac{x_2}{2}& \frac{x_1}{2}& 0& 1 &0& 0& 0& 0& 0& 0& 0\\
 0 & -\frac{x_5}{2}-\frac{x_2x_3}{12} & \frac{x_2^2}{12} & 0 & \frac{x_2}{2}& 0& 0 &1& 0& 0& 0& 0& 0& 0\\
 a_1 & -\frac{x_1^4}{720} & 0 & \frac{x_1^2}{12} & 0 & \frac{x_1}{2} & 0 & 0& 1 & 0 & 0 & 0 & 0 & 0 \\
 a_2 & a_3 & 0 & \frac{x_1x_2}{6} & \frac{x_1^2}{12} & \frac{x_2}{2} & \frac{x_1}{2} & 0 & 0 & 1 & 0 & 0 & 0 & 0 \\
 a_4 & a_5 & 0 & \frac{x_2^2}{12} & \frac{x_1x_2}{6} & 0 & \frac{x_2}{2} & \frac{x_1}{2} & 0 & 0 & 1& 0 & 0 & 0 \\
 \frac{x_2^4}{720} & a_6& 0 & 0 & \frac{x_2^2}{12} & 0 & 0 & \frac{x_2}{2} & 0 & 0 & 0 & 1 & 0 & 0\\
 a_7 & a_8 & a_9 & a_{10} & \frac{x_1^2}{12} & 0 & \frac{x_1}{2} & 0 & 0 & 0 & 0 & 0 & 1 & 0\\
 a_{11} & a_{12} & a_{13} & 0 & a_{14} & 0 & 0 &\frac{x_1}{2} & 0 & 0 & 0 & 0 & 0 & 1
   \end{matrix}\right]\,,
   \end{eqnarray*}
   where
   \begin{eqnarray*}
   a_1&=&\frac{x_1^3x_2}{720}-\frac{x_1x_4}{12}-\frac{x_6}{2}\,;\\
   a_2&=&\frac{x_1^2x_2^2}{240}-\frac{x_2x_4+x_1x_5}{12}-\frac{x_7}{2}\,;\\
   a_3&=&-\frac{x_1^3x_2}{240}-\frac{x_1x_4}{12}-\frac{x_6}{2}\,;\\
   a_4&=&\frac{x_1x_2^3}{240}-\frac{x_2x_5}{12}-\frac{x_8}{2}\,;\\
   a_5&=&-\frac{x_1^2x_2^2}{240}-\frac{x_2x_4+x_1x_5}{12}-\frac{x_7}{2}\,;\\
   a_6&=&-\frac{x_1x_2^3}{720}-\frac{x_2x_5}{12}-\frac{x_8}{2} \,;\\
   a_7&=&\frac{x_1^2x_2^2}{360}-\frac{x_3^2+x_1x_5-x_2x_4}{12}-\frac{x_7}{2}\,;\\
   a_8&=&-\frac{x_1^3x_2}{360}-\frac{x_1x_4}{6}\,;\\
   a_9&=&\frac{x_1x_3}{12}-\frac{x_4}{2}\,;\\
   a_{10}&=&\frac{x_1x_2}{12}+\frac{x_3}{2}\,;\\
   a_{11}&=&\frac{x_1x_2^3}{720}+\frac{x_2x_5}{12}-\frac{x_8}{2}\,;\\
   a_{12}&=&-\frac{x_1^2x_2^2}{720}-\frac{x_3^2+2x_1x_5}{12}\,;\\
   a_{13}&=&\frac{x_2x_3}{12}-\frac{x_5}{2}\,;\\
   a_{14}&=&\frac{x_1x_2}{12} +\frac{x_3}{2}\,,
   \end{eqnarray*}
   the induced left-invariant vector fields \eqref{leftinvariant vf} are: 
\begin{itemize}
\item $X_1={\partial}_{x_1}-\frac{x_2}{2}{\partial}_{x_3}-\big(\frac{x_3}{2}+\frac{x_1x_2}{12}\big){\partial}_{x_4}-\frac{x_2^2}{12}\partial_{x_5}-\big(\frac{x_4}{2}+\frac{x_1x_3}{12}\big)\partial_{x_6}-\big(\frac{x_5}{2}+\frac{x_2x_3}{12}\big)\partial_{x_7}\saltox+\big(\frac{x_1^3x_2}{720}-\frac{x_1x_4}{12}-\frac{x_6}{2}\big)\partial_{x_9}+\big(\frac{x_1^2x_2^2}{240}-\frac{x_2x_4+x_1x_5}{12}-\frac{x_7}{2}\big)\partial_{x_{10}}+\big(\frac{x_1x_2^3}{240}-\frac{x_2x_5}{12}\saltox-\frac{x_8}{2}\big)\partial_{x_{11}}+\frac{x_2^4}{720}\partial_{x_{12}}+\big(\frac{x_1^2x_2^2}{360}-\frac{x_3^2-x_2x_4+x_1x_5}{12}-\frac{x_7}{2}\big)\partial_{x_{13}}+\big(\frac{x_1x_2^3}{720}+\frac{x_2x_5}{12}-\frac{x_8}{2}\big)\partial_{x_{14}}\,;$
\item $X_2={\partial}_{x_2}+\frac{x_1}{2}{\partial}_{x_3}+\frac{x_1^2}{12}{\partial}_{x_4}+\big(\frac{x_1x_2}{12}-\frac{x_3}{2}\big)\partial_{x_5}-\big(\frac{x_4}{2}+\frac{x_1x_3}{12}\big)\partial_{x_7}-\big(\frac{x_5}{2}+\frac{x_2x_3}{12}\big)\partial_{x_8}\saltox-\frac{x_1^4}{720}\partial_{x_{9}}-\big(\frac{x_1^3x_2}{240}+\frac{x_1x_4}{12}+\frac{x_6}{2}\big)\partial_{x_{10}}-\big(\frac{x_1^2x_2^2}{240}+\frac{x_2x_4+x_1x_5}{12}+\frac{x_7}{2}\big)\partial_{x_{11}}-\big(\frac{x_1x_2^3}{720}\saltox+\frac{x_2x_5}{12}+\frac{x_8}{2}\big)\partial_{x_{12}}-\big(\frac{x_1^3x_2}{360}+\frac{x_1x_4}{6}\big)\partial_{x_{13}}-\big(\frac{x_1^2x_2^2}{720}+\frac{x_3^2+2x_1x_5}{12}\big)\partial_{x_{14}}\,;$
\item $X_3={\partial}_{x_3}+\frac{x_1}{2}{\partial}_{x_4}+\frac{x_2}{2}{\partial}_{x_5}+\frac{x_1^2}{12}\partial_{x_6}+\frac{x_1x_2}{6}\partial_{x_7}+\frac{x_2^2}{12}\partial_{x_8}+\big(\frac{x_1x_3}{12}-\frac{x_4}{2}\big)\partial_{x_{13}}+\big(\frac{x_1x_3}{12}\saltox-\frac{x_5}{2}\big)\partial_{x_{14}}\,;$
\item $X_4={\partial}_{x_4}+\frac{x_1}{2}{\partial}_{x_6}+\frac{x_2}{2}{\partial}_{x_7}+\frac{x_1^2}{12}\partial_{x_9}+\frac{x_1x_2}{6}\partial_{x_{10}}+\frac{x_2^2}{12}\partial_{x_{11}}+\big(\frac{x_1x_2}{12}+\frac{x_3}{2}\big)\partial_{x_{13}}\,;$
\item $X_5={\partial}_{x_5}+\frac{x_1}{2}{\partial}_{x_7}+\frac{x_2}{2}{\partial}_{x_8}+\frac{x_1^2}{12}\partial_{x_{10}}+\frac{x_1x_2}{6}\partial_{x_{11}}+\frac{x_2^2}{12}\partial_{x_{12}}+\frac{x_1^2}{12}\partial_{x_{13}}+\big(\frac{x_1x_2}{12}+\frac{x_3}{2}\big)\partial_{x_{14}}\,$;
\item $X_6={\partial}_{x_6}+\frac{x_1}{2}\partial_{x_{9}}+\frac{x_2}{2}\partial_{x_{10}}$;
\item $X_7={\partial}_{x_7}+\frac{x_1}{2}\partial_{x_{10}}+\frac{x_2}{2}\partial_{x_{11}}+\frac{x_1}{2}\partial_{x_{13}}$;
\item $X_8={\partial}_{x_8}+\frac{x_1}{2}\partial_{x_{11}}+\frac{x_2}{2}\partial_{x_{12}}+\frac{x_1}{2}\partial_{x_{14}}$;
\item $X_9=\partial_{x_{9}}$;
\item $X_{10}=\partial_{x_{10}}$;
\item $X_{11}=\partial_{x_{11}}$;
\item $X_{12}=\partial_{x_{12}}$;
\item $X_{13}=\partial_{x_{13}}$;
\item $X_{14}=\partial_{x_{14}}$,
\end{itemize}
and the respective left-invariant 1-forms \eqref{leftinvariant form} are: 
\begin{itemize}
\item $\theta_1=dx_1$;
\item $\theta_2=dx_2$;
\item $\theta_3=dx_3-\frac{x_1}{2}dx_2+\frac{x_2}{2}dx_1$;
\item $\theta_4=dx_4-\frac{x_1}{2}dx_3+\frac{x_1^2}{6}dx_2+\big(\frac{x_3}{2}-\frac{x_1x_2}{6}\big)dx_1$;
\item $\theta_5=dx_5-\frac{x_2}{2}dx_3+\big(\frac{x_3}{2}+\frac{x_1x_2}{6}\big)dx_2-\frac{x_2^2}{6}dx_1$;
\item $\theta_6=dx_6-\frac{x_1}{2}dx_4+\frac{x_1^2}{6}dx_3-\frac{x_1^3}{24}dx_2+\big(\frac{x_4}{2}-\frac{x_1x_3}{6}+\frac{x_1^2x_2}{24}\big)dx_1$;
\item $\theta_7=dx_7-\frac{x_1}{2}dx_5-\frac{x_2}{2}dx_4+\frac{x_1x_2}{3}dx_3+\big(\frac{x_4}{2}-\frac{x_1x_3}{6}-\frac{x_1^2x_2}{12}\big)dx_2+\big(\frac{x_5}{2}-\frac{x_2x_3}{6}+\frac{x_1x_2^2}{12}\big)dx_1$;
\item $\theta_8=dx_8-\frac{x_2}{2}dx_5+\frac{x_2^2}{6}dx_3+\big(\frac{x_5}{2}-\frac{x_2x_3}{6}-\frac{x_1x_2^2}{24}\big)dx_2+\frac{x_2^3}{24}dx_1$;
\item $\theta_9=dx_9-\frac{x_1}{2}dx_6+\frac{x_1^2}{6}dx_4-\frac{x_1^3}{24}dx_3+\frac{x_1^4}{120}dx_2+\big(\frac{x_6}{2}-\frac{x_1^3x_2}{120}+\frac{x_1^2x_3}{24}-\frac{x_1x_4}{6}\big) dx_1$;
\item $\theta_{10}=dx_{10}-\frac{x_1}{2}dx_7+\frac{x_1^2}{6}dx_5-\frac{x_2}{2}dx_6+\frac{x_1x_2}{3}dx_4-\frac{x_1^2x_2}{8}dx_3+\big(\frac{x_1^3x_2}{40}+\frac{x_1^2x_3}{24}-\frac{x_1x_4}{6}\saltot+\frac{x_6}{2}\big)dx_2+\big(\frac{x_1x_2x_3}{12}-\frac{x_1x_2^2}{40}-\frac{x_2x_4+x_1x_5}{6}-\frac{x_7}{2}\big)dx_1$;
\item $\theta_{11}=dx_{11}-\frac{x_1}{2}dx_8-\frac{x_2}{2}dx_7+\frac{x_1x_2}{3}dx_5+\frac{x_2^2}{6}dx_4-\frac{x_1x_2^2}{8}dx_3+\big(\frac{x_1^2x_2^2}{40}+\frac{x_1x_2x_3}{12}\saltot-\frac{x_2x_4+x_1x_5}{6}-\frac{x_7}{2}\big)dx_2+\big(\frac{x_2^2x_3}{24}-\frac{x_1x_2^3}{40}-\frac{x_2x_5}{6}-\frac{x_8}{2}\big)dx_1$;
\item $\theta_{12}=dx_{12}-\frac{x_2}{2}dx_8+\frac{x_2^2}{6}dx_5-\frac{x_2^3}{24}dx_3+\big(\frac{x_1x_2^3}{120}+\frac{x_2^2x_3}{24}-\frac{x_2x_5}{6}-\frac{x_8}{2}\big)dx_2+\frac{x_2^4}{120}dx_1$;
\item $\theta_{13}=dx_{13}-\frac{x_1}{2}dx_7+\frac{x_1^2}{6}dx_5+\big(\frac{x_1x_2}{6}-\frac{x_3}{2}\big)dx_4+\big(\frac{x_4}{2}+\frac{x_1x_3}{6}-\frac{x_1^2x_2}{12}\big)dx_3+\big(\frac{x_1^3x_2}{60}\saltot-\frac{x_1x_4}{3}\big)dx_2-\big(\frac{x_1^2x_2^2}{60}-\frac{x_1x_2x_3}{12}+\frac{x_3^2-x_2x_4+x_1x_5}{6}+\frac{x_7}{2}\big)dx_1$;
\item $\theta_{14}=dx_{14}-\frac{x_1}{2}dx_8+\big(\frac{x_1x_2}{6}-\frac{x_3}{2}\big)dx_5+\big(\frac{x_5}{2}+\frac{x_2x_3}{6}-\frac{x_1x_2^2}{24}\big)dx_3+\big(\frac{x_1^2x_2^2}{120}-\frac{x_3^2+2x_1x_5}{6}\big)dx_2\saltot-\big(\frac{x_1x_2^3}{120}-\frac{x_2^2x_3}{24}-\frac{x_2x_5}{6}+\frac{x_8}{2}\big)dx_1$.
\end{itemize}

Finally, we have
\begin{eqnarray*}
  \mathrm{d}(R_\mathbf{x})_\mathbf{0}=\left[\begin{matrix}
   1 & 0  & 0& 0& 0& 0& 0& 0& 0& 0& 0& 0& 0& 0\\
   0 & 1& 0& 0& 0& 0& 0 & 0& 0& 0& 0& 0& 0& 0\\
  \frac{x_2}{2} & -\frac{x_1}{2}& 1& 0& 0& 0& 0 & 0& 0& 0& 0& 0& 0& 0\\
   \frac{x_3}{2}-\frac{x_1x_2}{12} & \frac{x_1^2}{12}& -\frac{x_1}{2}& 1& 0& 0& 0& 0& 0& 0& 0& 0& 0& 0\\
   -\frac{x_2^2}{12} &  \frac{x_1x_2}{12}+\frac{x_3}{2}& \frac{x_2}{2}& 0& 1& 0& 0& 0& 0& 0& 0& 0& 0& 0\\
 \frac{x_4}{2}-\frac{x_1x_3}{12} & 0 & \frac{x_1^2}{12}& -\frac{x_1}{2}& 0& 1& 0 &0& 0& 0& 0& 0& 0& 0\\
 \frac{x_5}{2}-\frac{x_2x_3}{12}&  \frac{x_4}{2}-\frac{x_1x_3}{12}& \frac{x_1x_2}{6}& -\frac{x_2}{2}& -\frac{x_1}{2}& 0& 1 &0& 0& 0& 0& 0& 0& 0\\
 0 & \frac{x_5}{2}-\frac{x_2x_3}{12} & \frac{x_2^2}{12} & 0 & -\frac{x_2}{2}& 0& 0 &1& 0& 0& 0& 0& 0& 0\\
 a_1 & -\frac{x_1^4}{720} & 0 & \frac{x_1^2}{12} & 0 & -\frac{x_1}{2} & 0 & 0& 1 & 0 & 0 & 0 & 0 & 0 \\
 a_2 & a_3 & 0 & \frac{x_1x_2}{6} & \frac{x_1^2}{12} & -\frac{x_2}{2} & -\frac{x_1}{2} & 0 & 0 & 1 & 0 & 0 & 0 & 0 \\
 a_4 & a_5 & 0 & \frac{x_2^2}{12} & \frac{x_1x_2}{6} & 0 & -\frac{x_2}{2} & -\frac{x_1}{2} & 0 & 0 & 1& 0 & 0 & 0 \\
 \frac{x_2^4}{720} & a_6& 0 & 0 & \frac{x_2^2}{12} & 0 & 0 & -\frac{x_2}{2} & 0 & 0 & 0 & 1 & 0 & 0\\
 a_7 & a_8 & a_9 & a_{10} & \frac{x_1^2}{12} & 0 & -\frac{x_1}{2} & 0 & 0 & 0 & 0 & 0 & 1 & 0\\
 a_{11} & a_{12} & a_{13} & 0 & a_{14} & 0 & 0 &-\frac{x_1}{2} & 0 & 0 & 0 & 0 & 0 & 1
   \end{matrix}\right]\,,
   \end{eqnarray*}
   where
   \begin{eqnarray*}
   a_1&=&\frac{x_1^3x_2}{720}-\frac{x_1x_4}{12}+\frac{x_6}{2}\,;\\
   a_2&=&\frac{x_1^2x_2^2}{240}-\frac{x_2x_4+x_1x_5}{12}+\frac{x_7}{2}\,;\\
   a_3&=&-\frac{x_1^3x_2}{240}-\frac{x_1x_4}{12}+\frac{x_6}{2}\,;\\
   a_4&=&\frac{x_1x_2^3}{240}-\frac{x_2x_5}{12}+\frac{x_8}{2}\,;\\
   a_5&=&-\frac{x_1^2x_2^2}{240}-\frac{x_2x_4+x_1x_5}{12}+\frac{x_7}{2}\,;\\
   a_6&=&-\frac{x_1x_2^3}{720}-\frac{x_2x_5}{12}+\frac{x_8}{2} \,;\\
   a_7&=&\frac{x_1^2x_2^2}{360}-\frac{x_3^2+x_1x_5-x_2x_4}{12}+\frac{x_7}{2}\,; \\
   a_8&=&-\frac{x_1^3x_2}{360}-\frac{x_1x_4}{6}\,;\\
   a_9&=&\frac{x_1x_3}{12}+\frac{x_4}{2}\,;\\
   a_{10}&=&\frac{x_1x_2}{12}-\frac{x_3}{2}\,;\\
   a_{11}&=&\frac{x_1x_2^3}{720}+\frac{x_2x_5}{12}+\frac{x_8}{2}\,;\end{eqnarray*}
   \begin{eqnarray*}
   a_{12}&=&-\frac{x_1^2x_2^2}{720}-\frac{x_3^2+2x_1x_5}{12}\,;\\
   a_{13}&=&\frac{x_2x_3}{12}+\frac{x_5}{2}\,;\\
   a_{14}&=&\frac{x_1x_2}{12} -\frac{x_3}{2}\,,
   \end{eqnarray*}

  A representation of the left-invariant vector fields $X_i^s$ defined on $\mathbb{R}^{14}$ with respect to exponential coordinates of the second kind is the following:
  \begin{itemize}
      \item $X_1^s=\partial_{x_1}\,;$
      \item $X_2^s=\partial_{x_2}+x_1\partial_{x_3}+\frac{x_1^2}{2}\partial_{x_4}+x_1x_2\partial_{x_5}+\frac{x_1^3}{6}\partial_{x_6}+\frac{x_1^2x_2}{2}\partial_{x_7}+\frac{x_1x_2^2}{2}\partial_{x_8}+\frac{x_1^4}{24}\partial_{x_9}+\frac{x_1^3x_2}{6}\partial_{x_{10}}+\frac{x_1^2x_2^2}{4}\partial_{x_{11}}+\frac{x_1x_2^3}{6}\partial_{x_{12}}+\frac{x_1^2x_3}{2}\partial_{x_{13}}+x_1x_2x_3\partial_{x_{14}}$;
      \item $X_3^s=\partial_{x_3}+x_1\partial_{x_4}+x_2\partial_{x_5}+\frac{x_1^2}{2}\partial_{x_6}+x_1x_2\partial_{x_7}+\frac{x_2^2}{2}\partial_{x_8}+\frac{x_1^3}{6}\partial_{x_9}+\frac{x_1^2x_2}{2}\partial_{x_{10}}+\frac{x_1x_2^2}{2}\partial_{x_{11}}+\frac{x_2^3}{6}\partial_{x_{12}}+x_1x_3\partial_{x_{13}}+x_2x_3\partial_{x_{14}}$;
      \item $X_4^s=\partial_{x_4}+x_1\partial_{x_6}+x_2\partial_{x_7}+\frac{x_1^2}{2}\partial_{x_{9}}+x_1x_2\partial_{x_{10}}+\frac{x_2^2}{2}\partial_{x_{11}}+x_3\partial_{x_{13}}$;
      \item $X_5^s=\partial_{x_5}+x_1\partial_{x_7}+x_2\partial_{x_8}+\frac{x_1^2}{2}\partial_{x_{10}}+x_1x_2\partial_{x_{11}}+\frac{x_2^2}{2}\partial_{x_{12}}+\frac{x_1^2}{2}\partial_{x_{13}}+x_3\partial_{x_{14}}$;
      \item $X_6^s=\partial_{x_6}+x_1\partial_{x_9}+x_2\partial_{x_{10}}$;
      \item $X_7^s=\partial_{x_7}+x_1\partial_{x_{10}}+x_2\partial_{x_{11}}+x_1\partial_{}x_{13}$;
      \item $X_8^s=\partial_{x_8}+x_1\partial_{x_{11}}+x_2\partial_{x_{12}}+x_1\partial_{x_{14}}$;
      \item $X_9^s=\partial_{x_9}$;
      \item $X_{10}^s=\partial_{x_{10}}$;
      \item $X_{11}^s=\partial_{x_{11}}$;
      \item $X_{12}^s=\partial_{x_{12}}$;
      \item $X_{13}^s=\partial_{x_{13}}$;
      \item $X_{14}^s=\partial_{x_{14}}$.
  \end{itemize}

One can relate the  exponential coordinates of first kind to the exponential coordinates of second kind.
If we denote by $\alpha_1,\ldots, \alpha_{14}$ the coordinates of second kind and by
$a_1,\ldots, a_{14}$ the coordinates of first kind, the change of coordinates are given as follows:

\begin{equation*}
\left\{\begin{array}{ccl}
a_1 &=&\alpha_1 \\
a_2 &=&\alpha_2\\
a_3 &=&\alpha_3-\frac{\alpha_1\alpha_2}{2}\\
a_4 &=&\alpha_4-\frac{\alpha_1\alpha_3}{2}+\frac{\alpha_1^2\alpha_2}{12}\\
a_5 &=&\alpha_5-\frac{\alpha_2\alpha_3}{2} - \frac{\alpha_1\alpha_2^2}{12}\\
a_6 &=& \alpha_6-\frac{\alpha_1\alpha_4}{2}+\frac{\alpha_1^2\alpha_3}{12}\\
a_7 &=& \alpha_7-\frac{\alpha_1\alpha_5+\alpha_2\alpha_4}{2}+\frac{\alpha_1\alpha_2\alpha_3}{6}+\frac{\alpha_1^2\alpha_2^2}{24}\\
a_8&=&\alpha_8-\frac{\alpha_2\alpha_5}{2}+\frac{\alpha_2^2\alpha_3}{12}\\
a_9&=&\alpha_9-\frac{\alpha_1\alpha_6}{2}+\frac{\alpha_1^2\alpha_4}{12}-\frac{\alpha_1^4\alpha_2}{720}\\
a_{10}&=&\alpha_{10}-\frac{\alpha_1\alpha_7+\alpha_2\alpha_6}{2}+\frac{\alpha_1^2\alpha_5}{12}+\frac{\alpha_1\alpha_2\alpha_4}{6}-\frac{\alpha_1^3\alpha_2^2}{180}\\
a_{11}&=&\alpha_{11}-\frac{\alpha_1\alpha_8+\alpha_2\alpha_7}{2}+\frac{\alpha_1\alpha_2\alpha_5}{6}+\frac{\alpha_2^2\alpha_4}{12}+\frac{\alpha_1^2\alpha_2^3}{180}\\
a_{12}&=&\alpha_{12}-\frac{\alpha_2\alpha_8}{2}+\frac{\alpha_2^2\alpha_5}{12}+\frac{\alpha_1\alpha_2^4}{720}\\
a_{13}&=&\alpha_{13}-\frac{\alpha_1\alpha_7+\alpha_3\alpha_4}{2}+\frac{\alpha_1^2\alpha_5-\alpha_1\alpha_3^2}{12}+\frac{\alpha_1\alpha_2\alpha_4}{3}-\frac{\alpha_1^3\alpha_2^2}{120}\\
a_{14}&=&\alpha_{14}-\frac{\alpha_1\alpha_8+\alpha_3\alpha_5}{2}+\frac{\alpha_1\alpha_2\alpha_5}{3}-\frac{\alpha_2\alpha_3^2+\alpha_1\alpha_2^2\alpha_3}{12}-\frac{\alpha_1^2\alpha_2^3}{360}
\end{array}
\right. \qquad\qquad\qquad\qquad\qquad
\end{equation*}

and

\begin{equation*}
\left\{\begin{array}{ccl}
\alpha_1&=&a_1\\
\alpha_2    &=&a_2\\
\alpha_3    &=&a_3 +\frac{a_1 a_2}{2} \\
\alpha_4    &=&a_4 + \frac{a_1a_3}{2} +\frac{a_1^2a_2}{6}\\
\alpha_5    &=&a_5 + \frac{a_2 a_3}{2}+\frac{a_1 a_2^2 }{3} \\
\alpha_6 &=&a_6+\frac{a_1a_4}{2}+\frac{a_1^2a_3}{6}+\frac{a_1^3a_2}{24}\\
\alpha_7 &=&a_7+\frac{a_1a_5+a_2a_4}{2}+\frac{a_1a_2a_3}{3}+\frac{a_1^2a_2^2}{8}\\
\alpha_8 & =& a_8+\frac{a_2a_5}{2}+\frac{a_2^2a_3}{6}+\frac{a_1a_2^3}{8}\\
\alpha_9&=&a_9+\frac{a_1a_6}{2}+\frac{a_1^2a_4}{6}+\frac{a_1^3a_3}{24}+\frac{a_1^4a_2}{120}\\
\alpha_{10}&=&a_{10}+\frac{a_1a_7+a_2a_6}{2}+\frac{a_1^2a_5}{6}+\frac{a_1a_2a_4}{3}+\frac{a_1^2a_2a_3}{8}+\frac{a_1^3a_2^2}{30}\\
\alpha_{11}&=&a_{11}+\frac{a_1a_8+a_2a_7}{2}+\frac{a_1a_2a_5}{3}+\frac{a_2^2a_4}{6}+\frac{a_1a_2^2a_3}{8}+\frac{a_1^2a_2^3}{20}\\
\alpha_{12}&=&a_{12}+\frac{a_2a_8}{2}+\frac{a_2^2a_5}{6}+\frac{a_2^3a_3}{24}+\frac{a_1a_2^4}{30}\\
\alpha_{13}&=&a_{13}+\frac{a_1a_7+a_3a_4}{2}+\frac{a_1^2a_5+a_1a_2a_4}{6}+\frac{a_1a_3^2}{3}+\frac{a_1^2a_2a_3}{4}+\frac{a_1^3a_2^2}{20}\\
\alpha_{14}&=&a_{14}+\frac{a_1a_8+a_3a_5}{2}+\frac{a_1a_2a_5}{6}+\frac{a_2a_3^2}{3}+\frac{3a_1a_2^2a_3}{8}+\frac{a_1^2a_2^3}{10}
\end{array}
\right.\; .\qquad\qquad\qquad\qquad\qquad\qquad\quad
\end{equation*}

More explicitly, we have
$$
\exp(\alpha_{14} X_{14})\exp(\alpha_{13} X_{13})\cdots\exp(\alpha_1 X_1)=\qquad\qquad\qquad\qquad\qquad\qquad\qquad\qquad\qquad\qquad$$
$$=
\exp\left( \alpha_1 X_1+\alpha_2 X_2 + (\alpha_3-\frac{\alpha_1\alpha_2}{2})X_3 +(\alpha_4-\frac{\alpha_1\alpha_3}{2}+\frac{\alpha_1^2\alpha_2}{12})X_4+
\right.$$
$$+(\alpha_5-\frac{\alpha_2\alpha_3}{2} - \frac{\alpha_1\alpha_2^2}{12}) X_5+(\alpha_6-\frac{\alpha_1\alpha_4}{2}+\frac{\alpha_1^2\alpha_3}{12})X_6+\qquad\qquad$$
$$\qquad\qquad+(\alpha_7-\frac{\alpha_1\alpha_5+\alpha_2\alpha_4}{2}+\frac{\alpha_1\alpha_2\alpha_3}{6}+\frac{\alpha_1^2\alpha_2^2}{24})X_7+(\alpha_8-\frac{\alpha_2\alpha_5}{2}+\frac{\alpha_2^2\alpha_3}{12})X_8 +$$
$$ +(\alpha_9-\frac{\alpha_1\alpha_6}{2}+\frac{\alpha_1^2\alpha_4}{12}-\frac{\alpha_1^4\alpha_2}{720})X_9 +(\alpha_{10}-\frac{\alpha_1\alpha_7+\alpha_2\alpha_6}{2}+\frac{\alpha_1^2\alpha_5+2\alpha_1\alpha_2\alpha_4}{12}-\frac{\alpha_1^3\alpha_2^2}{180})X_{10}+$$
$$ +(\alpha_{11}-\frac{\alpha_1\alpha_8+\alpha_2\alpha_7}{2}+\frac{2\alpha_1\alpha_2\alpha_5+\alpha_2^2\alpha_4}{12}+\frac{\alpha_1^2\alpha_2^3}{180})X_{11} +(\alpha_{12}-\frac{\alpha_2\alpha_8}{2}+\frac{\alpha_2^2\alpha_5}{12}+\frac{\alpha_1\alpha_2^4}{720})X_{12}+ $$
$$+(\alpha_{13}-\frac{\alpha_1\alpha_7+\alpha_3\alpha_4}{2}+\frac{\alpha_1^2\alpha_5-\alpha_1\alpha_3^2}{12}+\frac{\alpha_1\alpha_2\alpha_4}{3}-\frac{\alpha_1^3\alpha_2^2}{120})X_{13}+  $$
$$\left.+(\alpha_{14}-\frac{\alpha_1\alpha_8+\alpha_3\alpha_5}{2}+\frac{\alpha_1\alpha_2\alpha_5}{3}-\frac{\alpha_2\alpha_3^2+\alpha_1\alpha_2^2\alpha_3}{12}-\frac{\alpha_1^2\alpha_2^3}{360})X_{14}\right) $$
and viceversa
$$
\exp\left( a_1 X_1+a_2 X_2 + \cdots+a_{13}X_{13}+a_{14}X_{14}\right)
=\qquad\qquad\qquad\qquad\qquad\qquad\qquad\qquad\quad$$
$$\qquad=\exp\bigg((a_{14}+\frac{a_1a_8+a_3a_5}{2}+\frac{a_1a_2a_5}{6}+\frac{a_2a_3^2}{3}+\frac{3a_1a_2^2a_3}{8}+\frac{a_1^2a_2^3}{10})X_{14}\bigg)\cdot $$
$$\cdot\exp\bigg((a_{13}+\frac{a_1a_7+a_3a_4}{2}+\frac{a_1^2a_5+a_1a_2a_4}{6}+\frac{a_1a_3^2}{3}+\frac{a_1^2a_2a_3}{4}+\frac{a_1^3a_2^2}{20})X_{13}\bigg)\cdot  $$
$$ \cdot\exp\bigg((a_{12}+\frac{a_2a_8}{2}+\frac{a_2^2a_5}{6}+\frac{a_2^3a_3}{24}+\frac{a_1a_2^4}{30})X_{12}\bigg)\cdot $$
$$ \cdot\exp\bigg((a_{11}+\frac{a_1a_8+a_2a_7}{2}+\frac{a_1a_2a_5}{3}+\frac{a_2^2a_4}{6}+\frac{a_1a_2^2a_3}{8}+\frac{a_1^2a_2^3}{20})X_{11}\bigg)\cdot $$
$$ \cdot\exp\bigg((a_{10}+\frac{a_1a_7+a_2a_6}{2}+\frac{a_1^2a_5}{6}+\frac{a_1a_2a_4}{3}+\frac{a_1^2a_2a_3}{8}+\frac{a_1^3a_2^2}{30})X_{10}\bigg)\cdot $$
$$\cdot\exp\bigg((a_9+\frac{a_1a_6}{2}+\frac{a_1^2a_4}{6}+\frac{a_1^3a_3}{24}+\frac{a_1^4a_2}{120})X_{9}\bigg)\cdot  $$
$$\cdot\exp\bigg((a_8+\frac{a_2a_5}{2}+\frac{a_2^2a_3}{6}+\frac{a_1a_2^3}{8})X_8\bigg)\cdot
\exp\bigg((a_7+\frac{a_1a_5+a_2a_4}{2}+\frac{a_1a_2a_3}{3}+\frac{a_1^2a_2^2}{8})X_7\bigg)\cdot $$
$$\quad\cdot\exp\bigg((a_6+\frac{a_1a_4}{2}+\frac{a_1^2a_3}{6}+\frac{a_1^3a_2}{24})X_6\bigg)\cdot\exp\bigg((a_5+ \frac{a_2 a_3}{2}  +\frac{a_1 a_2^2}{3}  ) X_5 \bigg)\cdot$$
$$\qquad\qquad\cdot\exp\bigg((a_4+ \frac{a_1 a_3}{2}  +\frac{a_1^2 a_2}{6}  )  X_4\bigg)\cdot\exp\bigg((a_3 +\frac{1}{2} a_1 a_2)X_3\bigg) 
\cdot\exp(a_2 X_2) \cdot\exp(a_1 X_1) \,.$$


  \newpage

\bibliography{general_bibliography_tripaldi}
\bibliographystyle{amsalpha}
\end{document}